\documentclass[a4paper, 12pt]{article}
\usepackage[T2A]{fontenc}
\usepackage[english,russian]{babel}
\usepackage{amsmath}
\usepackage{amsfonts,amssymb,mathrsfs,amscd}
\usepackage{bm}
\renewcommand{\Im}{\operatorname{Im}}
\newcommand{\Ker}{\operatorname{Ker}}
\renewcommand{\dim}{\operatorname{dim}}
\newcommand{\codim}{\operatorname{codim}}
\newcommand{\Kappa}{\operatorname{\mathcal K}}
\newcommand{\sgn}{\operatorname{sign}}

\newcommand{\inter}{\operatorname{int}}
\newcommand{\e}{\mathfrak e}

\newcommand{\mx}{\operatorname{\mathfrak M}}
\newcommand{\p}{\operatorname{\mathfrak p}}
\renewcommand{\P}{\operatorname{\mathfrak P}}
\newcommand{\Q}{\operatorname{\mathfrak Q}}
\renewcommand{\span}{\operatorname{span}}
\newcommand{\cmx}{\operatorname{\mathfrak C}}
\newcommand{\mn}{\operatorname{\mathfrak m}}
\newcommand{\cmn}{\operatorname{\mathfrak c}}
\newcommand{\m}{\mathfrak m}
\newcommand{\n}{\mathfrak n}
\newcommand{\M}{\mathfrak M}
\newcommand{\s}{\bm s}

\newcommand{\card}{\operatorname{card}}

\newcommand{\supp}{\operatorname{supp}}
\newcommand{\supvrai}{\operatornamewithlimits{sup\,vrai}}
\newcommand{\N}{\mathbb N}
\newcommand{\Z}{\mathbb Z}
\newcommand{\R}{\mathbb R}
\newcommand{\Nu}{\mathcal N}

\newcommand{\D}{\mathcal D}

\renewcommand{\L}{\mathcal L}
\newcommand{\mes}{\operatorname{mes}}
\allowdisplaybreaks

\begin{document}

\author{ С. Н. Кудрявцев }
\title{Восстановление производных по значениям функций из классов Никольского -- Бесова \\
смешанной гладкости в областях определённого вида}
\date{}
\maketitle
\begin{abstract}
В работе рассмотрены пространства Никольского и Бесова с нормами, в
определении которых вместо смешанных модулей непрерывности известных порядков
определённых смешанных производных функций используются "$L_p$-усреднённые"
смешанные модули непрерывности функций соответствующих порядков.
Установлены оценки сверху и снизу величины наилучшей точности
восстановления производных по значениям функций в заданном числе точек
для классов таких функций в ограниченных областях определённого вида. Эти оценки
не слабее, а в некоторых случаях являются более сильными, чем
соответствующие оценки, полученные автором ранее в рассматриваемой задаче
для указанных классов функций на кубе $ I^d. $ Помимо этого существенно расширен
класс пространств Никольского -- Бесова смешанной гладкости, для которых
получены упомянутые оценки в рассматриваемой задаче.
\end{abstract}

Ключевые слова: точность, восстановление, производная,
значения функций, классы Никольского -- Бесова смешанной гладкости
\bigskip

\centerline{Введение}
\bigskip

В работе при соответствующих условиях на $ \lambda \in \Z_+^d, d \in \N, $
установлены оценки сверху и снизу величины наилучшей
точности восстановления в пространстве $ L_q(D), 1 \le q \le \infty, $
производной $ \D^\lambda f $ по значениям в $ n $ точках функций $ f $
из классов Никольского $ (\mathcal S_p^\alpha \mathcal H)^\prime(D) $ и Бесова
$ (\mathcal S_{p,\theta}^\alpha \mathcal B)^\prime(D),
\alpha \in \R_+^d, 1 \le p < \infty, 1 \le \theta < \infty, $ заданных в
ограниченных областях $ D \subset \R^d, $ определённого вида.

Настоящая работа продолжает исследования, проводившиеся автором в
[1] -- [7] для упомянутой задачи в отношении различных классов функций
конечной гладкости. Отметим, что средства, использовавшиеся в [1] -- [7]
для получения как верхних, так и нижних оценок указанной выше величины для
рассматривавшихся в этих работах классов функций, неприменимы в
отношении классов, изучаемых в настоящей статье. Помимо этого, схемы,
применявшиеся в [1] -- [7] при получении метрических соотношений, используемых
для вывода верхних оценок точности приближений производных функций из
исследуемых в [1] -- [7] классов, непригодны для изучаемых в настоящей
работе классов. Предложенные ниже средства и схемы вывода оценок
позволили найти при определённых условиях слабую асимптотику
отмеченной выше величины для классов Никольского и Бесова смешанной гладкости
функций, заданных в ограниченных областях определённого вида.
При этом в некоторых случаях получены более сильные чем установленные ранее
оценки сверху и снизу рассматриваемой здесь величины для таких классов функций
на кубе $ i^d. $

В заключение в качестве иллюстрации приведём полученный в \S 2 результат,
касающийся величины наилучшей точности
восстановления в $ L_q(D) $ производных $ \D^\lambda f $ по значениям в $ n $
точках функций $ f $ из класса Бесова
$ (\mathcal S_{p, \theta}^\alpha \mathcal B)^\prime(D) $ в ограниченной
области $ D, $ имеющей определённую структуру, определяемой равенством
$$
\sigma_n(\D^\lambda,(\mathcal S_{p,\theta}^\alpha \mathcal B)^\prime(D),L_q(D))
= \inf_{A,\phi} \sup_{f \in (\mathcal S_{p,\theta}^\alpha \mathcal B)^\prime(D)}
\|\D^\lambda f - A \circ \phi(f)\|_{L_q(D)},
$$
где $ \inf $ берётся по всем отображениям $ A: \R^n \mapsto
L_q(D) $ и отображениям $ \phi: C(D) \mapsto \R^n $ вида $
\phi(f) = (f(x^1),\ldots,f(x^n)), x^j \in D, j =1,\ldots.n. $
Если $ d \in \N, D \subset \R^d $ -- ограниченная область определённого
строения, $ \alpha \in \R_+^d, 1 \le p < \infty, 1 \le \theta, q \le \infty, \lambda \in
\Z_+^d, $ установлено (см. теорему 2.2.7), что при $ \alpha_j -p^{-1} >0,
\alpha_j -\lambda_j -(p^{-1} -q^{-1})_+ >0, j =1,\ldots,d, $ справедливо
неравенство
$$
\sigma_n(\D^\lambda,(\mathcal S_{p,\theta}^\alpha \mathcal B)^\prime(D),L_q(D))
\le c_1 n^{-\bm m} (\log n)^{(\bm m +1 -1/\max(p,\theta))(\bm c -1)},
$$
где
\begin{multline*}
\bm m = \min\{\alpha_j -\lambda_j -(p^{-1} -q^{-1})_+:
j =1,\ldots,d\}, \\
\bm c = \card\{j =1,\ldots,d: \alpha_j -\lambda_j -(p^{-1} -q^{-1})_+ = \bm m\},
\end{multline*}

$ c_1 $ -- положительная константа, не зависящая от $ n. $
Отметим, что в аналогичной оценке из [4] слагаемое $ -1/\max(p,\theta) $
в показателе степени $ \log n $ отсутствует. При этом полученные
ниже оценки при $ \bm c = 1 $ дают порядок изучаемой величины. Таким образом
существенно расширена совокупность классов функций Никольского -- Бесова
смешанной гладкости, для которых установлена слабая асимптотика поведения
в зависимости от $n$ величины наилучшей точности восстановления в
рассматриваемой задаче.

Работа состоит из введения и трёх параграфов, в первом из которых приведены
предварительные сведения, во втором установлеа верхняя оценка, а в третьем --
нижняя оценка изучаемой величины.
\bigskip

\centerline{\S 1. Предварительные сведения и вспомогательные
утверждения}
\bigskip

1.1. В этом пункте вводятся обозначения, относящиеся к
функциональным классам и пространствам, рассматриваемым в
настоящей работе, а также приводятся некоторые факты, необходимые
в дальнейшем.

Для $ d \in \N $ через $ \Z_+^d $ обозначим множество
$$
\Z_+^d = \{\lambda = (\lambda_1, \ldots, \lambda_d) \in \Z^d:
\lambda_j \ge0, j=1, \ldots, d\}.
$$
При $ d \in \N $ для $ l \in \Z_+^d $ обозначим также через $ \Z_+^d(l) $
множество
\begin{equation*} \tag{1.1.1}
\Z_+^d(l) = \{ \lambda  \in \Z_+^d: \lambda_j \le l_j, j=1, \ldots, d\}.
\end{equation*}

Для $  d \in  \N, l \in \Z_+^d $ через $ \mathcal P^{d,l} $ будем
обозначать пространство вещественных  полиномов, состоящее из всех
функций $ f: \R^d \mapsto \R $ вида
$$
f(x) = \sum_{\lambda \in \Z_+^d(l)} a_{\lambda} \cdot x^{\lambda}, x \in \R^d,
$$
где  $   a_{\lambda}   \in   \R,  x^{\lambda} = x_1^{\lambda_1}
\ldots x_d^{\lambda_d}, \lambda \in \Z_+^d(l). $

При $ d \in \N $ в $ \R^d $ фиксируем норму
\begin{equation*}
\|x\| = \max_{j =1, \ldots,d} |x_j|.
\end{equation*}

Для множества $ A $ из топологического пространства $ T $ через $ \overline A $
обозначается замыкание множества $ A, $ а через $ \inter A $  -- его внутренность.

Через $ \chi_A $ будем обозначать характеристическую функцию множества
$ A \subset \R^d. $

Для топологического пространства $ T $ через $ C(T) $
обозначим пространство непрерывных вещественных
функций, заданных на $ T. $ Если $ T $ --
компактное топологическое пространство, то в $ C(T) $ фиксируем норму,
определяемую равенством
$$
\|f\|_{C(T)} = \max_{x \in T} |f(x)|, f \in C(T).
$$

При $ d \in \N $ через $ C_0(\R^d) $ обозначим банахово пространство
вещественных непрерывных функций $ f, $ заданных на $ \R^d, $ для каждой из
которых $ \lim_{ \| x \| \to \infty} f(x) =0, $ с нормой
\begin{equation*}
\| f \|_{C_0(\R^d)} = \max_{x \in \R^d} | f(x) |.
\end{equation*}

Для измеримого по Лебегу множества $ D \subset \R^d $ при $ 1 \le p \le \infty $
через $ L_p(D), $ как обычно, обозначается
пространство  всех  вещественных измеримых на $ D $ функций $f,$
для которых определена норма
$$
\|f\|_{L_p(D)} = \begin{cases} (\int_D |f(x)|^p dx)^{1/p}, 1 \le p < \infty; \\
\supvrai_{x \in D}|f(x)|, p = \infty.
\end{cases}
$$

Отметим здесь неоднократно используемое в дальнейшем неравенство
\begin{equation*} \tag{1.1.2}
| \sum_{j =1}^n x_j|^a \le \sum_{j =1}^n |x_j|^a, x_j \in \R, j =1,\ldots,n,
n \in \N, 0 \le a \le 1.
\end{equation*}

Для $ x,y \in \R^d $ положим $ xy = x \cdot y = (x_1 y_1, \ldots, x_d y_d), $
а для $ x \in \R^d $ и $ A \subset \R^d $ определим
$$
x A = x \cdot A = \{xy: y \in A\}.
$$

При $ d \in \N $ для $ x,y \in \R^d $ будем обозначать
$$
(x,y) = \sum_{j =1}^d x_j y_j.
$$

Для $ x \in \R^d: x_j \ne 0, $ при $ j=1,\ldots,d,$ положим
$ x^{-1} = (x_1^{-1},\ldots,x_d^{-1}). $

При $ d \in \N $ для $ x,y \in \R^d $ будем писать $ x \le y (x < y), $ если
для каждого $ j=1,\ldots,d $ выполняется неравенство $ x_j \le y_j (x_j < y_j). $

При $ d \in \N $ для $ x \in \R^d $ положим
$$
x_+ = ((x_1)_+, \ldots, (x_d)_+),
$$
где $ t_+ = \frac{1} {2} (t +|t|), t \in \R. $

Обозначим через $ \R_+^d $ множество $ x \in \R^d, $ для которых
$ x_j >0 $ при $ j=1,\ldots,d,$ и для $ a \in \R_+^d, x \in \R^d $
положим $ a^x = a_1^{x_1} \ldots a_d^{x_d}.$

При $ d \in \N $ определим множества
$$
I^d = \{x \in \R^d: 0 < x_j < 1,j=1,\ldots,d\},
$$
$$
\overline I^d = \{x \in \R^d: 0 \le x_j \le 1,j=1,\ldots,d\},
$$
$$
B^d = \{x \in \R^d: -1 \le x_j \le 1,j=1,\ldots,d\}.
$$

Через $ \e $ будем обозначать вектор в $ \R^d, $ задаваемый
равенством $ \e = (1,\ldots,1). $

При $ d \in \N $ для $ \lambda \in \Z_+^d $ через $ \D^\lambda $
будем обозначать оператор дифференцирования $ \D^\lambda =
\frac{\D^{|\lambda|}} {\D x_1^{\lambda_1} \ldots \D x_d^{\lambda_d}}, $
где $ |\lambda| = \sum_{j=1}^d \lambda_j. $

Теперь приведём некоторые факты, относящиеся к полиномам, которыми
мы будем пользоваться ниже.

В [3] содержится такое утверждение.

    Лемма 1.1.1

Пусть $ d \in \N, l \in \Z_+^d, \lambda \in \Z_+^d, 1 \le p,q \le \infty,
\rho, \sigma \in \R_+^d. $ Тогда существует константа
$ c_1(d,l,\lambda,\rho, \sigma) >0 $ такая, что для любых
измеримых по Лебегу множеств $ D,Q \subset \R^d, $  для которых
можно найти $ \delta \in \R_+^d $ и $ x^0 \in \R^d $ такие, что
$ D \subset (x^0 +\rho \delta B^d) $ и $ (x^0 +\sigma \delta I^d) \subset Q, $
для любого полинома $ f \in \mathcal P^{d,l} $
выполняется неравенство
    \begin{equation*} \tag{1.1.3}
\| \D^\lambda f\|_{L_q(D)} \le c_1 \delta^{-\lambda -p^{-1} \e +q^{-1} \e}
\|f\|_{L_p(q)}.
   \end{equation*}

Далее, напомним, что для открытого множеста $ D \subset \R^d $ и вектора
$ h \in \R^d $ через $ D_h $ обозначается множество
$$
D_h = \{x \in D: x +th \in D \ \forall t \in \overline I\}.
$$

Для функции $ f, $ заданной на открытом множестве $ D \subset \R^d, $ и
вектора $ h \in \R^d $ определим на $ D_h $ её разность $ \Delta_h f $
с шагом $ h, $ полагая
$$
(\Delta_h f)(x) = f(x +h) -f(x), x \in D_h,
$$
а для $ l \in \N: l \ge 2, $ на $ D_{lh} $ определим $l$-ую
разность $ \Delta_h^l f $ функции $ f $ с шагом $ h $ равенством
$$
(\Delta_h^l f)(x) = (\Delta_h (\Delta_h^{l-1} f))(x), x \in
D_{lh},
$$
положим также $ \Delta_h^0 f = f. $

Как известно, справедливо равенство
$$
(\Delta_h^l f)(\cdot) = \sum_{k=0}^l C_l^k (-1)^{l-k} f(\cdot +kh),
C_l^k =\frac{l!} {k! (l-k)!}.
$$

При $ d \in \N $ для $ j=1,\ldots,d$ через $ e_j $ будем
обозначать вектор $ e_j = (0,\ldots,0,1_j,0,\ldots,0).$

Как показано в [8], справедлива следующая лемма.

Лемма 1.1.2

   Пусть $ d \in \N, l \in \Z_+^d. $ Тогда
для любого $ \delta \in \R_+^d $ и $ x^0 \in \R^d $ для $ Q = x^0 +\delta I^d $
существует единственный линейный оператор
$ P_{\delta, x^0}^{d,l}: L_1(Q) \mapsto \mathcal P^{d,l}, $
обладающий следующими свойствами:

1) для $ f \in \mathcal P^{d,l} $ имеет место равенство
\begin{equation*} \tag{1.1.4}
P_{\delta, x^0}^{d,l}(f \mid_Q) = f,
  \end{equation*}

2)
\begin{equation*}
\Ker P_{\delta,x^0}^{d,l} = \biggl\{\,f \in L_1(Q):
\int \limits_{Q} f(x) g(x) \,dx =0\ \forall g \in \mathcal P^{d,l}\,\biggr\},
\end{equation*}

причём существуют константы $ c_2(d,l) >0 $ и $ c_3(d,l) >0 $
такие, что

   3) при $ 1 \le p \le \infty $ для $ f \in L_p(Q) $ справедливо неравенство
  \begin{equation*} \tag{1.1.5}
\|P_{\delta, x^0}^{d,l} f \|_{L_p(Q)} \le c_2 \|f\|_{L_p(Q)},
  \end{equation*}

4) при $ 1 \le p < \infty $ для $ f \in L_p(Q) $ выполняется неравенство
   \begin{equation*} \tag{1.1.6}
  \| f -P_{\delta, x^0}^{d,l} f \|_{L_p(Q)} \le c_3 \sum_{j=1}^d
\delta_j^{-1/p} \biggl(\int_{\delta_j B^1} \int_{Q_{(l_j +1) \xi e_j}}
|\Delta_{\xi e_j}^{l_j +1} f(x)|^p dx d\xi\biggr)^{1/p}.
\end{equation*}

Теперь определим пространства и классы функций, рассматриваемые в настоящей
работе (ср. с [9], [10]). Но прежде введём некоторые обозначения.

При $ d \in \N $ для $ x \in \R^d $ через $\s(x) $ обозначим
множество $ \s(x) = \{j =1,\ldots,d: x_j \ne 0\}, $ а для
множества $ J \subset \{1,\ldots,d\} $ через $ \chi_J $ обозначим
вектор из $ \R^d $ с компонентами
$$
(\chi_J)_j = \begin{cases} 1, & \text{ для } j \in J; \\
0, & \text{ для } j \in (\{1,\ldots,d\} \setminus J).
\end{cases}
$$

При $ d \in \N $ для $ x \in \R^d $ и $ J = \{j_1,\ldots,j_k\}
\subset \N: 1 \le j_1 < j_2 < \ldots < j_k \le d, $ через $ x^J $
обозначим вектор $ x^J = (x_{j_1},\ldots,x_{j_k}) \in \R^k, $ а
для множества $ A \subset \R^d $ положим $ A^J = \{x^J: x \in A\}. $

Для открытого множества $ D \subset \R^d $ и векторов $ h \in \R^d $ и $ l \in
\Z_+^d $ через $ D_h^l $ обозначим множество
\begin{multline*}
D_h^l = (\ldots (D_{l_d h_d e_d})_{l_{d-1} h_{d-1} e_{d-1}}
\ldots)_{l_1 h_1 e_1} = \{ x \in D: x +tlh \in D \ \forall t \in
\overline I^d\} = \\ \{ x \in D: (x +\sum_{j \in \s(l)} t_j l_j h_j e_j) \in D \
\forall t^{\s(l)} \in (\overline I^d)^{\s(l)} \}.
\end{multline*}

Пусть $ d \in \N, D $ -- открытое множество в $ \R^d $ и $ 1 \le p \le \infty. $
Тогда для $ f \in L_p(D), h \in \R^d $ и $ l \in \Z_+^d $ определим в $ D_h^l $
смешанную разность функции $ f $ порядка $ l, $ соответствующую вектору $ h, $
равенством
\begin{multline*}
(\Delta_h^l f)(x) = \biggl(\biggl(\prod_{j=1}^d \Delta_{h_j e_j}^{l_j}\biggr) f\biggr)(x)
= \biggl(\biggl(\prod_{j \in \s(l)} \Delta_{h_j e_j}^{l_j}\biggr) f\biggr)(x) = \\
\sum_{k \in \Z_+^d(l)} (-\e)^{l-k} C_l^k f(x+kh), x \in D_h^l,
\end{multline*}
где $ C_l^k = \prod_{j=1}^d C_{l_j}^{k_j}. $

Имея в виду, что для $ f \in L_p(D), l \in \Z_+^d $ и векторов
$ h,h^\prime \in \R^d: h^{\s(l)} = (h^\prime)^{\s(l)}, $ соблюдается
соотношение
$$
\| \Delta_h^l f\|_{L_p(D_h^l)} = \| \Delta_{h^\prime}^l
f\|_{L_p(D_{h^\prime}^l)}, 1 \le p \le \infty,
$$
определим при $ 1 \le p \le \infty $ для функции $ f \in L_p(D) $ смешанный
модуль непрерывности в $ L_p(D) $ порядка $ l \in \Z_+^d $ равенством
$$
\Omega^l (f,t^{\s(l)})_{L_p(D)} = \supvrai_{ \{ h \in \R^d:
h^{\s(l)} \in t^{\s(l)} (B^d)^{\s(l)} \}} \| \Delta_h^l f\|_{L_p(D_h^l)},
t^{\s(l)} \in (\R_+^d)^{\s(l)}.
$$
Кроме того, при тех же условиях введём в рассмотрение для функции $ f $
"усреднённый" смешанный модуль непрерывности в $ L_p(D) $ порядка $ l, $
полагая
\begin{multline*}
\Omega^{\prime l} (f, t^{\s(l)})_{L_p(D)} = \begin{cases}
\biggl((2 t^{\s(l)})^{-\e^{\s(l)}}
\int_{ t^{\s(l)} (B^d)^{\s(l)}} \| \Delta_\xi^l f\|_{L_p(D_\xi^l)}^p
d \xi^{\s(l)}\biggr)^{1 /p} = \\
\biggl((2 t^{\s(l)})^{-\e^{\s(l)}}
\int_{ (t B^d)^{\s(l)}} \int_{D_\xi^{l \chi_{\s(l)}}}
| \Delta_\xi^{l \chi_{\s(l)}} f(x)|^p dx
d \xi^{\s(l)}\biggr)^{1 /p}, p \ne \infty; \\
\Omega^l (f,t^{\s(l)})_{L_p(D)}, p = \infty,
\end{cases} \\
t^{\s(l)} \in (\R_+^d)^{\s(l)}.
\end{multline*}

Из приведенных определений видно, что
\begin{multline*} \tag{1.1.7}
\Omega^{\prime l} (f, t^{\s(l)})_{L_p(D)} \le
\Omega^l (f, t^{\s(l)})_{L_p(D)}, t^{\s(l)} \in (\R_+^d)^{\s(l)}, \\
f \in L_p(D), 1 \le p \le \infty, l \in \Z_+^d, \\
D \text{ -- произвольное открытое множество в } \R^d.
\end{multline*}

Пусть теперь $ \alpha \in \R_+^d, 1 \le p \le \infty $ и $ D $ --
область в $ \R^d. $ Тогда зададим вектор $ l = l(\alpha) \in \N^d, $ полагая
$ l_j = \min \{m \in \N: \alpha_j < m \}, j =1,\ldots,d, $ и обозначим через
$ (S_p^\alpha H)^\prime(D) ((\mathcal S_p^\alpha \mathcal H)^\prime(D)) $
множество всех функций $ f \in L_p(D), $ обладающих тем свойством, что для
любого непустого множества $ J \subset \{1,\ldots,d\} $ выполняется неравенство
$$
\sup_{t^J \in (\R_+^d)^J}
(t^J)^{-\alpha^J} \Omega^{\prime l \chi_J}(f,t^J)_{L_p(D)}
= \sup_{t^J \in (\R_+^d)^J} (\prod_{j \in J} t_j^{-\alpha_j})
\Omega^{\prime l \chi_J}(f, t^{\s(l\chi_J)})_{L_p(D)} < \infty (\le 1).
$$

Пусть $ \alpha,p,D $ и $ l = l(\alpha) $ -- те же, что и выше, и
$ \theta \in \R: 1 \le \theta < \infty. $ Тогда обозначим через
$ (S_{p,\theta}^\alpha B)^\prime(D) ((\mathcal S_{p,\theta}^\alpha
\mathcal B)^\prime(D)) $ множество всех функций $ f \in L_p(D), $ которые для
любого непустого множества $ J \subset \{1,\ldots,d\} $
удовлетворяют условию
\begin{multline*}
\biggl(\int_{(\R_+^d)^J} (t^J)^{-\e^J -\theta \alpha^J}
(\Omega^{\prime l \chi_J}(f, t^J)_{L_p(D)})^\theta dt^J\biggr)^{1/\theta} =\\
\biggl(\int_{(\R_+^d)^J} (\prod_{j \in J} t_j^{-1 -\theta \alpha_j})
(\Omega^{\prime l \chi_J}(f, t^{\s(l \chi_J)})_{L_p(D)})^\theta \prod_{j
\in J} dt_j\biggr)^{1/\theta} < \infty (\le 1).
\end{multline*}

При $ \theta = \infty $ положим $ (S_{p,\infty}^\alpha B)^\prime(D) =
(S_p^\alpha H)^\prime(D), (\mathcal S_{p,\infty}^\alpha \mathcal B)^\prime(D) =
(\mathcal S_p^\alpha \mathcal H)^\prime(D). $

Учитывая, что для $ f \in (S_{p,\theta}^\alpha B)^{\prime}(D),
t^J \in (\R_+^d)^J (J \subset \{1, \ldots, d\}: J \ne \emptyset) $
выполняется неравенство
\begin{multline*} \tag{1.1.8}
(t^J)^{-\alpha^J} \Omega^{\prime l \chi_J}(f, t^J)_{L_p(D)} = \\
\biggl(\int_{t^J +t^J (I^d)^J} (t^J)^{-\e^J -\theta \alpha^J} (\Omega^{\prime l \chi_J}(f,
t^J)_{L_p(D)})^\theta d\tau^J\biggr)^{1/\theta} = \\
\biggl(\int_{t^J +t^J (I^d)^J} (t^J)^{-\e^J -\theta \alpha^J} ((2 t^J)^{-\e^J}
\int_{ (t B^d)^J} \| \Delta_\xi^{l \chi_J} f\|_{L_p(D_\xi^{l \chi_J})}^p
d \xi^J)^{\theta /p} d\tau^J\biggr)^{1/\theta} \le \\
\biggl(\int_{t^J +t^J (I^d)^J} (\prod_{j \in J} 2^{1 +\theta \alpha_j})
(\tau^J)^{-\e^J -\theta \alpha^J} ((\tau^J)^{-\e^J}
\int_{ (\tau B^d)^J} \| \Delta_\xi^{l \chi_J} f\|_{L_p(D_\xi^{l \chi_J})}^p
d \xi^J)^{\theta /p} d\tau^J\biggr)^{1/\theta} = \\
\biggl(\int_{t^J +t^J (I^d)^J} (\prod_{j \in J} 2^{1 +\theta \alpha_j})
(\tau^J)^{-\e^J -\theta \alpha^J} (\prod_{j \in J} 2^{\theta /p})
(\Omega^{\prime l \chi_J}(f, \tau^J)_{L_p(D)})^\theta
d\tau^J\biggr)^{1/\theta} \le \\
(\prod_{j \in J} 2^{\alpha_j +1/\theta +1 /p}) \biggl(\int_{(\R_+^d)^J}
(\tau^J)^{-\e^J -\theta \alpha^J} (\Omega^{\prime l \chi_J}(f,
\tau^J)_{L_p(D)})^\theta d\tau^J\biggr)^{1/\theta},
\end{multline*}
заключаем, что
\begin{equation*}
(\mathcal S_{p, \theta}^\alpha \mathcal B)^\prime(D) \subset c_4(\alpha)
(\mathcal S_p^\alpha \mathcal H)^\prime(D),
\end{equation*}
где $ c_4(\alpha) = \prod_{j=1}^d 2^{2+\alpha_j}. $

Пусть $ \alpha,p,D $ и $ l = l(\alpha) $ -- те же, что и выше, и
$ \theta \in \R: 1 \le \theta < \infty. $ Тогда обозначим через
$ (S_{p,\theta}^\alpha B)^0(D) ((\mathcal S_{p,\theta}^\alpha \mathcal B)^0(D)) $
множество всех функций $ f \in L_p(D), $ которые для любого непустого множества
$ J \subset \{1,\ldots,d\} $ подчинены условию
\begin{multline*}
\biggl(\int_{(\R_+^d)^J} (t^J)^{-\e^J -\theta \alpha^J}
(\Omega^{l \chi_J}(f, t^J)_{L_p(D)})^\theta dt^J\biggr)^{1/\theta} = \\
\biggl(\int_{(\R_+^d)^J} (\prod_{j \in J} t_j^{-1 -\theta \alpha_j})
(\Omega^{l \chi_J}(f, t^{\s(l \chi_J)})_{L_p(D)})^\theta \prod_{j
\in J} dt_j\biggr)^{1/\theta} < \infty (\le 1),
\end{multline*}
а через $ (S_p^\alpha H)^0(D) ((\mathcal S_p^\alpha \mathcal H)^0(D)) $ --
множество всех функций $ f \in L_p(D), $ обладающих тем свойством, что для
любого непустого множества $ J \subset \{1,\ldots,d\} $ выполняется неравенство
$$
\sup_{t^J \in (\R_+^d)^J}
(t^J)^{-\alpha^J} \Omega^{l \chi_J}(f,t^J)_{L_p(D)}
= \sup_{t^J \in (\R_+^d)^J} (\prod_{j \in J} t_j^{-\alpha_j})
\Omega^{l \chi_J}(f, t^{\s(l\chi_J)})_{L_p(D)} < \infty (\le 1).
$$
При $ \theta = \infty $ положим $ (S_{p,\infty}^\alpha B)^0(D) =
(S_p^\alpha H)^0(D), (\mathcal S_{p,\infty}^\alpha \mathcal B)^0(D) =
(\mathcal S_p^\alpha \mathcal H)^0(D). $
Из (1.1.7) вытекает, что
\begin{multline*} \tag{1.1.9}
\| f\|_{(S_{p, \theta}^\alpha B)^\prime(D)} \le
\| f\|_{(S_{p, \theta}^\alpha B)^{0}(D)}, \\
\ f \in (S_{p, \theta}^\alpha B)^{0}(D), \
\alpha \in \R_+^d, 1 \le p \le \infty, 1 \le \theta \le \infty, \\
D \text{ -- произвольная область в } \R^d.
\end{multline*}

Обозначим через $ C^\infty(D) $ пространство бесконечно
дифференцируемых функций в области $ D \subset \R^d, $ а через
$ C_0^\infty(D) $ -- пространство функций $ f \in C^\infty(\R^d), $
каждая из которых имеет компактный носитель $ \supp f \subset D. $

В заключение этого пункта введём ещё несколько обозначений.

Для банахова пространства $ X $ (над $ \R$) обозначим
$ B(X) = \{x \in X: \|x\|_X \le 1\}. $

Для банаховых пространств $ X,Y $ через $ \mathcal B(X,Y) $
обозначим банахово пространство, состоящее из непрерывных линейных
операторов $ T: X \mapsto Y, $ с нормой
$$
\|T\|_{\mathcal B(X,Y)} = \sup_{x \in B(X)} \|Tx\|_Y.
$$
Отметим, что если $ X=Y,$ то $ \mathcal B(X,Y) $ является банаховой алгеброй.
\bigskip

1.2. В этом пункте содержатся сведения о кратных рядах, которыми будем
пользоваться в дальнейшем.

Для $ d \in \N, y \in \R^d $ положим
$$
\mn(y) = \min_{j=1,\ldots,d} y_j
$$
и для банахова пространства $ X, $ вектора $ x \in X $ и семейства
$ \{x_\kappa \in X, \kappa \in \Z_+^d\} $ будем писать $ x =
\lim_{ \mn(\kappa) \to \infty} x_\kappa, $ если для любого $ \epsilon >0 $
существует $ n_0 \in \N $ такое, что для любого $ \kappa \in \Z_+^d, $ для
которого $ \mn(\kappa) > n_0, $
справедливо неравенство $ \|x -x_\kappa\|_X  < \epsilon. $

Пусть $ X $ -- банахово пространство, $ d \in \N $ и
$ \{ x_\kappa \in X, \kappa \in \Z_+^d\} $ -- семейство векторов.
Тогда под суммой ряда $ \sum_{\kappa \in \Z_+^d} x_\kappa $ будем
понимать вектор $ x \in X, $ для которого выполняется равенство $
 x = \lim_{\mn(k) \to \infty} \sum_{\kappa \in \Z_+^d(k)} x_\kappa. $

При $ d \in \N $ через $ \Upsilon^d $ обозначим множество
$$
\Upsilon^d = \{ \epsilon \in \Z^d: \epsilon_j \in \{0,1\},
j=1,\ldots,d\}.
$$

Имеет место

Лемма 1.2.1

Пусть $ X $ -- банахово пространство, а вектор $ x \in X $ и
семейство $ \{x_\kappa \in X, \kappa \in \Z_+^d\} $ таковы, что
$ x = \lim_{ \mn(\kappa) \to \infty} x_\kappa, $ Тогда для семейства
$ \{ \mathcal X_\kappa \in X, \kappa \in \Z_+^d \}, $ элементы которого
определяются равенством
$$
\mathcal X_\kappa = \sum_{\epsilon \in \Upsilon^d: \s(\epsilon)
\subset \s(\kappa)} (-\e)^\epsilon x_{\kappa -\epsilon}, \kappa \in \Z_+^d,
$$
справедливо равенство
$$
x = \sum_{\kappa \in \Z_+^d} \mathcal X_\kappa.
$$

Лемма является следствием того, что при $ k \in \Z_+^d $
выполняется равенство
\begin{equation*}
\sum|{\kappa \in \Z_+^d(k)}
\mathcal X_\kappa = x_k \text{ (см. [8])}.
\end{equation*}

Замечание.

Легко заметить, что для любого семейства чисел
$ \{x_\kappa \in \R: x_\kappa \ge 0, \kappa \in \Z_+^d\}, $ если ряд
$ \sum_{\kappa \in \Z_+^d} x_\kappa $ сходится, т.е. существует предел
$ \lim_{\mn(k) \to \infty} \sum_{\kappa \in \Z_+^d(k)} x_\kappa, $
что эквивалентно соотношению
$ \sup_{k \in \Z_+^d} \sum_{\kappa \in \Z_+^d(k)} x_\kappa < \infty, $

то для любой последовательности подмножеств
$ \{Z_n \subset \Z_+^d, n \in \Z_+\}, $
 таких, что $ \card Z_n < \infty,
Z_n \subset Z_{n+1},  n \in \Z_+, $
и $ \cup_{ n \in \Z_+} Z_n = \Z_+^d, $
справедливо равенство
$ \sum_{\kappa \in \Z_+^d} x_\kappa =
\lim_{ n \to \infty} \sum_{\kappa \in Z_n} x_\kappa =
\sup_{k \in \Z_+^d} \sum_{\kappa \in \Z_+^d(k)} x_\kappa. $
Отсюда несложно понять, что если для семейства векторов
$ \{x_\kappa \in X, \kappa \in \Z_+^d\}  $ банахова пространства $ X $
ряд $ \sum_{\kappa \in \Z_+^d} \| x_\kappa \|_X $ сходится,
то для любой последовательности подмножеств
$ \{Z_n \subset \Z_+^d, n \in \Z_+\}, $
 таких, что $ \card Z_n < \infty,
Z_n \subset Z_{n+1},  n \in \Z_+, $
и $ \cup_{ n \in \Z_+} Z_n = \Z_+^d, $
в $ X $ соблюдается равенство
$ \sum_{\kappa \in \Z_+^d} x_\kappa =
\lim_{ n \to \infty} \sum_{\kappa \in Z_n} x_\kappa. $

При $ d \in \N $ для $ x \in \R^d $ обозначим
\begin{eqnarray*}
\mx(x)  &=& \max_{j=1,\ldots,d} x_j, \\
\cmx(x) &=& \card \{j \in \{1,\ldots,d\}: x_j = \mx(x)\}, \\
\cmn(x) &=& \card \{j \in \{1,\ldots,d\}: x_j = \mn(x)\}.
\end{eqnarray*}

Лемма 1.2.2

Пусть $ d \in \N, \beta \in \R_+^d, \alpha \in \R^d $ и
$ \mx(\beta^{-1} \alpha) >0. $ Тогда существуют константы
$ c_1(d,\alpha,\beta)  >0 $ и $ c_2(d,\alpha,\beta) >0 $ такие, что
для $ r \in \N $ соблюдается неравенство
\begin{equation*} \tag{1.2.1}
c_1 2^{ \mx(\beta^{-1} \alpha) r} r^{ \cmx(\beta^{-1} \alpha) -1}
\le \sum_{ \kappa \in \Z_+^d: (\kappa, \beta) \le r} 2^{(\kappa, \alpha)} \le
c_2 2^{ \mx(\beta^{-1} \alpha) r} r^{ \cmx(\beta^{-1} \alpha) -1}.
\end{equation*}

Лемма 1.2.3

Пусть $ d \in \N, \alpha, \beta \in \R_+^d. $ Тогда существуют
константы $ c_3(d,\alpha,\beta) >0 $ и $ c_4(d,\alpha,\beta) >0 $
такие, что при $ r \in \N $ справедливо неравенство
\begin{equation*} \tag{1.2.2}
c_3 2^{-\mn(\beta^{-1} \alpha) r} r^{\cmn(\beta^{-1} \alpha) -1}
\le \sum_{\kappa \in \Z_+^d: (\kappa, \beta) > r} 2^{-(\kappa, \alpha)} \le
c_4 2^{-\mn(\beta^{-1} \alpha) r} r^{\cmn(\beta^{-1} \alpha) -1}.
\end{equation*}

Доказательство лемм 1.2.2 и 1.2.3 приведено в [11].
\bigskip

1.3. В этом пункте приведём некоторые вспомогательные утверждения,
которые используются в следующем пункте и далее.

Следующее утверждение установлено в [12].

Лемма 1.3.1

Пусть $ d \in \N, 1 \le p < \infty. $ Тогда

1) при  $ j=1,\ldots,d$ для любого непрерывного линейного оператора
$ T: L_p(\R) \mapsto L_p(\R) $ существует единственный непрерывный
линейный оператор $ \mathcal T^j: L_p(\R^d) \mapsto L_p(\R^d), $ для которого
для любой функции $ f \in L_p(\R^d) $ почти для всех $ (x_1,\ldots,x_{j-1},
x_{j+1},\ldots,x_d) \in \R^{d-1} $ в $ L_p(\R) $ имеет место равенство
\begin{multline*} \tag{1.3.1}
(\mathcal T^j f)(x_1,\ldots, x_{j-1},\cdot,x_{j+1},\ldots,x_d) =\\
(T(f(x_1,\ldots,x_{j-1},\cdot,x_{j+1},
\ldots,x_d)))(\cdot),
\end{multline*}

2) при этом для каждого $ j=1,\ldots,d $ отображение $ V_j^{L_p}, $ которое
каждому оператору $ T \in \mathcal B(L_p(\R), L_p(\R)) $ ставит в соответствие
оператор $ V_j^{L_p}(T) = \mathcal T^j \in \mathcal B(L_p(\R^d), L_p(\R^d)), $
удовлетворяющий (1.3.1), является непрерывным гомоморфизмом банаховой алгебры
$ \mathcal B(L_p(\R), L_p(\R)) $ в банахову алгебру $ \mathcal B(L_p(\R^d),
L_p(\R^d)), $

3) причём для любых операторов $ S,T \in \mathcal B(L_p(\R), L_p(\R)) $
при любых $ i,j =1,\ldots,d: i \ne j, $ выполняется равенство
\begin{equation*} \tag{1.3.2}
(V_i^{L_p}(S) V_j^{L_p}(T)) f = (V_j^{L_p}(T) V_i^{L_p}(S)) f, f \in L_p(\R^d).
\end{equation*}

Лемма 1.3.2

Пусть $ d \in \N. $ Тогда

1) при $ j=1,\ldots,d$ для любого непрерывного линейного оператора
$ T: C_0(\R) \mapsto C_0(\R) $ существует единственный непрерывный
линейный оператор
$ \mathcal T^j: C_0(\R^d) \mapsto C_0(\R^d), $ обладающий тем свойством, что
для любой функции $ f \in C_0(\R^d) $ для всех $ (x_1,\ldots,x_{j-1},x_{j+1},
\ldots,x_d) \in \R^{d-1} $ при любом $ x_j \in \R $
имеет место равенство
\begin{multline*} \tag{1.3.3}
(\mathcal T^j f)(x_1,\ldots, x_{j-1},x_j,x_{j+1},\ldots,x_d) =\\
(T(f(x_1,\ldots,x_{j-1},\cdot,x_{j+1},\ldots,x_d)))(x_j),
\end{multline*}

2) при этом для каждого $ j=1,\ldots,d $ отображение $ V_j^C, $ которое каждому
оператору $ T \in \mathcal B(C_0(\R), C_0(\R)) $ ставит в
соответствие оператор $ V_j^C(T) = \mathcal T^j \in \mathcal B(C_0(\R^d),
C_0(\R^d)), $
удовлетворяющий (1.3.3), является непрерывным гомоморфизмом банаховой
алгебры $ \mathcal B(C_0(\R), C_0(\R)) $ в банахову алгебру
$ \mathcal B(C_0(\R^d), C_0(\R^d)), $

3) причём для любых операторов $ S,T \in \mathcal B(C_0(\R), C_0(\R)) $
при любых $ i,j =1,\ldots,d: i \ne j, $ справедливо равенство
\begin{equation*} \tag{1.3.4}
(V_i^C(S) V_j^C(T)) f = (V_j^C(T) V_i^C(S)) f, f \in C_0(\R^d).
\end{equation*}

Доказательство.

Пусть $ T \in \mathcal B(C_0(\R), C_0(\R)), j =1,\ldots,d. $
Каждой функции $ f \in C_0(\R^d) $ сопоставим функцию $ \mathcal
T^j f, $ значение которой в точке $(x_1,\ldots,x_{j-1},x_j,x_{j+1},\ldots,x_d) \in \R^d $
определяется равенством (1.3.3), и проверим, что функция $
\mathcal T^j f \in C_0(\R^d). $ Для этого заметим, что всякая
функция $ f \in C_0(\R^d) $ равномерно непрерывна на $ \R^d. $ В
самом деле, для произвольного $ \epsilon > 0 $ выберем $ R > 0  $
так, чтобы при $ \|x\| \ge R $ выполнялось неравенство $ |f(x)| <
\epsilon.$ Принимая во внимание, что $ f $ непрерывна на компакте
$ \{ x \in \R^d: \|x\| \le R\}, $ заключаем, что $ f $ -- равномерно
непрерывна на этом компакте. Исходя из этого, возьмём $ \delta >
0, $ для которого при $ x, y \in \R^d: \|x\| \le R, \|y\| \le R,
\|x -y\| < \delta, $ соблюдается неравенство $ | f(x) -f(y)| <
\epsilon. $ Тогда для $ x, y \in \R^d: \|x\| \ge R, \|y\| \ge R,
\| x -y\| < \delta, $ имеет место неравенство $ | f(x) -f(y)| \le
| f(x)| +| f(y)| < 2 \epsilon. $ А для $ x, y \in \R^d: \|x\| \le
R, \|y\| \ge R, \|x -y\| < \delta, $ находя $ z = x +t (y -x), t
\in \overline I, $ для которого $ \|z\| = R, $ имеем $ \| x -z\| =
t \|y -x\| < \delta, \| y -z\| = (1 -t) \|y -x\| < \delta, $ что в
силу сказанного выше влечёт неравенство $ | f(x) -f(y)| \le | f(x)
-f(z)| +| f(z) -f(y)| < \epsilon +2 \epsilon = 3 \epsilon. $ Тем
самым, равномерная непрерывность $ f $ на $ \R^d $ установлена.

Пользуясь этим, проверим непрерывность $ \mathcal T^j f $ на $ \R^d $ для $ f
\in C_0(\R^d. $
Пусть $ x \in \R^d $ и $ \epsilon >0. $ Учитывая условия леммы, выберем
$ \delta >0 $ так, чтобы при $ | x_j -y_j| < \delta $ выполнялось неравенство
$$
| (T f(x_1,\ldots, x_{j-1}, \cdot, x_{j+1}, \ldots,x_d))(x_j) -
(T f(x_1,\ldots, x_{j-1}, \cdot, x_{j+1}, \ldots,x_d))(y_j)| < \epsilon,
$$
а при $ \max_{\{i =1,\ldots,d: i \ne j\}} | x_i -y_i| < \delta $
для $ z \in \R $ соблюдалось неравенство
$$
| f(x_1,\ldots, x_{j-1},z, x_{j+1}, \ldots,x_d) -
f(y_1,\ldots, y_{j-1},z,y_{j+1},\ldots,y_d) | < \epsilon.
$$
Тогда отсюда для $ y \in \R^d $ таких, что $ \max_{i=1,\ldots,d} | x_i -y_i|
= \| x -y \| < \delta, $ имеем
\begin{multline*}
| (\mathcal T^j f)(x) -(\mathcal T^j f)(y)| \le \\
| (T f(x_1,\ldots,x_{j-1},\cdot, x_{j+1}, \ldots,x_d))(x_j) -
(T f(x_1,\ldots, x_{j-1}, \cdot, x_{j+1}, \ldots,x_d))(y_j)| +\\
| (T f(x_1,\ldots, x_{j-1}, \cdot, x_{j+1}, \ldots,x_d))(y_j) -
(T f(y_1,\ldots, y_{j-1}, \cdot, y_{j+1}, \ldots,y_d))(y_j)| < \\ \epsilon +
\| T(f(x_1,\ldots, x_{j-1}, \cdot, x_{j+1}, \ldots,x_d) -
f(y_1,\ldots, y_{j-1}, \cdot, y_{j+1}, \ldots,y_d))\|_{C_0(\R)} \le\\
\epsilon +\|T\|_{\mathcal B(C_0(\R), C_0(\R))}
\|f(x_1,\ldots, x_{j-1}, \cdot, x_{j+1}, \ldots,x_d) -\\
f(y_1,\ldots, y_{j-1}, \cdot, y_{j+1}, \ldots,y_d)\|_{C_0(\R)} \le
\epsilon +\|T\|_{\mathcal B(C_0(\R), C_0(\R))} \epsilon.
\end{multline*}
Тем самым непрерывность функции $ \mathcal T^j f $ в $ \R^d $
установлена.

Теперь проверим, что для $ f \in C_0(\R^d) $ предел
$ \lim_{ \|x\| \to \infty} (\mathcal T^j f)(x) =0. $ Прежде всего для
произвольного $ \epsilon > 0 $ выберем $ R_0 > 0 $ так, чтобы при
$ \|x\| > R_0 $ соблюдалось неравенство $ | f(x) | < \epsilon. $ Тогда для
$ x \in \R^d: \max_{ i =1,\ldots,d: i \ne j} | x_i| > R_0, $ выполняется
неравенство
\begin{multline*}
| (\mathcal T^j f)(x)| =
| (T(f(x_1,\ldots, x_{j-1}, \cdot, x_{j+1}, \ldots,x_d)))(x_j)| \le \\
\| T(f(x_1,\ldots, x_{j-1}, \cdot, x_{j+1}, \ldots,x_d)) \|_{C_0(\R)} \le\\
\| T\|_{\mathcal B(C_0(\R), C_0(\R))}
\| f(x_1,\ldots, x_{j-1}, \cdot, x_{j +1}, \ldots,x_d)\|_{C_0(\R)} \le \\
\| T\|_{\mathcal B(C_0(\R), C_0(\R))} \epsilon.
\end{multline*}

Далее, для каждого $ y^\prime = (y_1,\ldots,y_{j -1}, y_{j +1},\ldots,y_d) \in
\R^{d -1}: \max_{i =1,\ldots,d: i \ne j}
| y_i| \le R_0, $ фиксируем $ R_{y^\prime}, $ для которого при
$ | x_j| > R_{y^\prime} $ справедливо нервенство
$ | (T(f(y_1,\ldots,y_{j -1},\cdot,y_{j +1},\ldots,y_d)))(x_j) | <
\epsilon. $ В силу равномерной непрерывности на $ \R^d $ функции
$ f $ выберем $ \delta > 0, $ для которого при $ \| \xi -\eta\| \le
\delta $ выполняется неравенство $ | f(\xi) -f(\eta)| < \epsilon, $
и построим конечную $ \delta$-сеть $ \{y^{\prime s} \in \R^{d -1}:
\| y^{\prime s} \| \le R_0, s =1,\ldots,S\} $ для компактного
множества $ \{y^\prime \in \R^{d -1}: \| y^\prime \| \le R_0\}. $
Тогда для $ x \in \R^d: \max_{i =1,\ldots,d: i \ne j} | x_i| =
\| x^\prime \| \le R_0, | x_j| > \max_{s =1,\ldots,S} R_{y^{\prime s}}, $ беря
$ y^{\prime s}, s \in \{1,\ldots,S\}, $ для которого $ \| x^\prime -y^{\prime s} \|
= \max_{i =1,\ldots,d: i \ne j} | x_i -y_i^s| \le \delta, $ на основании
сказанного выше имеем
\begin{multline*}
| (\mathcal T^j f)(x_1,\ldots, x_{j-1},x_j,x_{j+1},\ldots,x_d)| =\\
| (T(f(x_1,\ldots,x_{j-1},\cdot,x_{j+1},\ldots,x_d)))(x_j)| \le \\
| (T(f(x_1,\ldots,x_{j-1},\cdot,x_{j+1},\ldots,x_d)))(x_j) -
(T(f(y_1^s,\ldots,y_{j -1}^s,\cdot,y_{j +1}^s,\ldots,y_d^s)))(x_j)| + \\
| (T(f(y_1^s,\ldots,y_{j -1}^s,\cdot,y_{j +1}^s,\ldots,y_d^s)))(x_j)| < \\
\epsilon +| (T(f(x_1,\ldots,x_{j-1},\cdot,x_{j+1},\ldots,x_d) -
f(y_1^s,\ldots,y_{j -1}^s,\cdot,y_{j +1}^s,\ldots,y_d^s)))(x_j)| \le \\
\epsilon +\| T(f(x_1,\ldots,x_{j-1},\cdot,x_{j+1},\ldots,x_d) -
f(y_1^s,\ldots,y_{j-1}^s,\cdot,y_{j +1}^s,\ldots,y_d^s))\|_{C_0(\R)} \le \\
\epsilon +\|T\|_{\mathcal B(C_0(\R), C_0(\R))}
\| f(x_1,\ldots,x_{j-1},\cdot,x_{j+1},\ldots,x_d) -\\
f(y_1^s,\ldots,y_{j-1}^s,\cdot,y_{j +1}^s,\ldots,y_d^s)\|_{C_0(\R)} \le
\epsilon +\|T\|_{\mathcal B(C_0(\R), C_0(\R))} \epsilon.
\end{multline*}
Таким образом, определённое выше соответствие $ \mathcal T^j, $ является
отображением $ C_0(\R^d) $ в себя.

Линейность отображения $ \mathcal T^j: C_0(\R^d) \mapsto C_0(\R^d) $
проверяется непосредственно выкладкой с использованием (1.3.3),
опираясь на линейность отображения $ T. $ Проверим его
непрерывность.

Пусть $ f \in C_0(\R^d). $ Тогда
\begin{multline*}
\| \mathcal T^j f \|_{C_0(\R^d)} = \sup_{x \in \R^d}
| (\mathcal T^j f)(x)| = \\
\sup_{(x_1,\ldots, x_{j-1}, x_{j+1}, \ldots,x_d)
\in \R^{d-1}} \sup_{x_j \in \R}
| (\mathcal T^j f)(x_1,\ldots, x_{j-1}, x_j, x_{j+1}, \ldots,x_d)| = \\
\sup_{(x_1,\ldots, x_{j-1}, x_{j+1}, \ldots,x_d) \in \R^{d-1}} \sup_{x_j \in \R}
| (T f(x_1,\ldots, x_{j-1}, \cdot, x_{j+1}, \ldots,x_d))(x_j)| \le \\
\sup_{(x_1,\ldots, x_{j-1}, x_{j+1}, \ldots,x_d) \in \R^{d-1}}
\|T\|_{\mathcal B(C_0(\R), C_0(\R))} \sup_{x_j \in \R}
| f(x_1,\ldots, x_{j-1}, x_j, x_{j+1}, \ldots,x_d)| = \\
\|T\|_{\mathcal B(C_0(\R), C_0(\R))}
\sup_{(x_1,\ldots, x_{j-1}, x_{j+1}, \ldots,x_d) \in \R^{d -1}} \sup_{x_j \in \R}
| f(x_1,\ldots, x_{j-1}, x_j, x_{j+1}, \ldots,x_d)| = \\
\|T\|_{\mathcal B(C_0(\R), C_0(\R))} \sup_{x \in \R^d} | f(x)|.
\end{multline*}
Откуда заключаем, что $ \mathcal T^j \in \mathcal B(C_0(\R^d), C_0(\R^d)) $ и
$ \|\mathcal T^j\|_{\mathcal B(C_0(\R^d), C_0(\R^d))} \le
\| T \|_{\mathcal B(C_0(\R), C_0(\R))}. $

Единственность оператора $ \mathcal T^j, $ удовлетворяющего (1.3.3), очевидна.

Соблюдение свойств $ V_j^C(\alpha S +\beta T) = \alpha V_j^C(S)
+\beta V_j^C(T) $ и $ V_j^C(ST) = V_j^C(S) V_j^C(T) $ устанавливается
с помощью (1.3.3).

Для доказательства утверждения п. 3) достаточно показать, что равенство
(1.3.4) имеет место для любой функции $ f $ вида $ f = \phi g, $ где
$ g \in \mathcal P^{d,l}, l \in \Z_+^d, $ а $ \phi(x_1,\ldots,x_d) =
\prod_{j =1}^d \phi_j(x_j), $ причём для каждого $ j =1,\ldots,d $ функция
$ \phi_j $ непрерывна на $ \R, $ и удовлетворяет условиям:
\begin{multline*}
\phi_J(x_j) =1, \text{ при } | x_j| < R_0;
\phi_j(x_j) =0, \text{ при } | x_j| > R_1; \\
0 \le \phi_j(x_j) \le 1, \text{ при } R_0 \le | x_j| \le R_1, 0 < R_0 < R_1 < \infty.
\end{multline*}
А для получения (1.3.4) в общем случае достаточно заметить, что множество
функций указанного вида плотно в $ C_0(\R^d). \square $

Замечание.

Если при $ d \in \N, 1 \le p, q < \infty, $ оператор
$ T \in \mathcal B(L_p(\R), L_p(\R)) \cap \mathcal B(L_q(\R), L_q(\R))
( T \in \mathcal B(L_p(\R), L_p(\R)) \cap \mathcal B(C_0(\R), C_0(\R))), $ то
при $ j =1,\ldots,d $ для $ f \in L_p(\R^d) \cap L_q(\R^d)
(f \in L_p(\R^d) \cap C_0(\R^d)) $
справедливо равенство $ (V_j^{L_p} T) f = (V_j^{L_q} T) f
((V_j^{L_p} T) f = (V_j^C T) f). $ Поэтому символы $ L_p, L_q, C $ в
качестве индексов у $ V_j $ можно опускать.

Для $ \phi \in L_\infty(\R) $ через $ M_\phi $ обозначим линейный оператор
в пространстве локально суммируемых функций на $ \R, $ определяемый
равенством $ M_\phi g = \phi g, $ где $ g $ -- локально суммируемая на $ \R $
функция.

Так же, как соответствующее утверждение в [13], устанавливается лемма 1.3.3.

Лемма 1.3.3

Пусть $ d \in \N, l \in \N^d, 1 \le p < \infty. $ Тогда существует константа
$ c_1(d,l) >0 $ такая, что для любых $ x^0, X^0 \in \R^d, \delta, \Delta \in
\R_+^d $ таких, что $ Q = (x^0 +\delta I^d) \subset (X^0 +\Delta I^d), $ для
любой функции $ f \in L_p(\R^d), $ для любого $ \xi \in \R^d, $ для любых
множеств $ J, J^\prime \subset \{1,\ldots,d\} $ имеет место неравенство
\begin{multline*} \tag{1.3.5}
\biggl\| \Delta_\xi^{l \chi_J} ((\prod_{j \in J^\prime}
V_j(E -M_{\chi^j_{\Delta_j, X^0_j}}
P_{\delta_j, x_j^0}^{1,l_j -1})) f)\biggr\|_{L_p(Q_\xi^{l \chi_J})} \le \\
c_1 \biggl(\prod_{j \in J^\prime \setminus J} \delta_j^{-1/p} \biggr)
\biggl(\int_{(\delta B^d)^{J^\prime \setminus J}} \int_{ Q_\xi^{l
\chi_{J \cup J^\prime}}} |(\Delta_\xi^{l \chi_{J \cup J^\prime}} f)(x)|^p
dx d\xi^{J^\prime \setminus J} \biggr)^{1/p},
\end{multline*}
где $ E $ -- тождественный оператор, а для каждого $ j \in J^\prime $ функция
$ \chi^j_{\Delta_j, X^0_j} $ непрерывна на $ \R, $ имеет компактный носитель и
$ \chi^j_{\Delta_j, X^0_j}(x) =1 $ для $ x \in X^0_j +\Delta_j I. $

Доказательство.

Сначала убедимся в справедливости (1.3.5) в случае, когда
множества $ J $ и $ J^\prime $ не пересекаются. А именно, покажем,
что если выполнены условия леммы и $ J \cap J^\prime = \emptyset, $ то
\begin{multline*} \tag{1.3.6}
\|\Delta_\xi^{l \chi_J} ((\prod_{j \in J^\prime} V_j(E -M_{\chi^j_{\Delta_j, X^0_j}}
P_{\delta_j, x_j^0}^{1, l_j -1})) f)\|_{L_p(Q_\xi^{l \chi_J})} \le\\
c_2^k \biggl(\prod_{j \in J^\prime} \delta_j^{-1/p}\biggr)
\biggl(\int_{(\delta B^d)^{J^\prime}} \int_{ Q_\xi^{l \chi_{J \cup
J^\prime}}} | (\Delta_\xi^{l \chi_{J \cup J^\prime}} f)(x)|^p dx
d\xi^{J^\prime}\biggr)^{1/p},
\end{multline*}
где $ k = \card J^\prime, $ а $ c_2(d,l) = \max_{j=1,\ldots,d}
c_3(1, l_j-1) (см. (1.1.6)). $

Доказательство (1.3.6) проведём по индукции относительно $ k. $ При
$ k=1, $ а, следовательно, $ J^\prime = \{j\} (j \in \{1,\ldots,d\}
\setminus J), $ используя теорему Фубини, (1.3.1), (1.1.6), в условиях леммы имеем
\begin{multline*}
\|\Delta_\xi^{l \chi_J} (V_j(E -M_{\chi^j_{\Delta_j, X^0_j}}
P_{\delta_j, x_j^0}^{1, l_j -1}) f)\|_{L_p (Q_\xi^{l \chi_J})}^p  = \\
\int\limits_{\substack{ \prod_{i=1,\ldots,d: i \ne j} \\ (x_i^0
+\delta_i I)_{l_i (\chi_J)_i \xi_i}}} \int\limits_{(x_j^0 +\delta_j I)}
|((E -M_{\chi^j_{\Delta_j, X^0_j}} P_{\delta_j, x_j^0}^{1, l_j -1})
(\Delta_\xi^{l \chi_J} f)(x_1, \ldots,
x_{j-1}, \cdot, x_{j+1}, \ldots,x_d))(x_j)|^p \\
\times dx_j dx_1 \ldots dx_{j-1} dx_{j+1} \ldots dx_d = \\
\int\limits_{\substack{ \prod_{i=1,\ldots,d: i \ne j} \\ (x_i^0
+\delta_i I)_{l_i (\chi_J)_i \xi_i}}} \int\limits_{(x_j^0 +\delta_j I)}
|((E -P_{\delta_j, x_j^0}^{1, l_j -1}) (\Delta_\xi^{l \chi_J} f)
(x_1, \ldots, x_{j-1}, \cdot, x_{j+1}, \ldots,x_d))(x_j)|^p \\
\times dx_j dx_1 \ldots dx_{j-1} dx_{j+1} \ldots dx_d \le \\
c_2^p \delta_j^{-1} \int_{\delta_j B^1} \int_{ Q_\xi^{l \chi_{J \cup
J^\prime}}} |(\Delta_\xi^{l \chi_{J \cup J^\prime}} f)(x)|^p dx
d\xi_j.
\end{multline*}
Отсюда следует (1.3.6) в случае, когда $ J^\prime = \{j\}, $ а $ j \in
\{1,\ldots,d\} \setminus J. $

Предположим, что при $ k \in \N $ неравенство (1.3.6) имеет место
при соблюдении условий леммы 1.3.3 для $ J,J^\prime \subset \{1.\ldots,d\}:
J \cap J^\prime = \emptyset $ и $ \card J^\prime \le k. $ Покажем, что тогда оно
справедливо в условиях леммы 1.3.3 в ситуации, когда $ J,J^\prime \subset
\{1,\ldots,d\}: J \cap J^\prime = \emptyset, $ и $ \card J^\prime = k+1. $

Для этого, фиксируя для таких множеств $ J,J^\prime $ элемент $ j \in J^\prime $
и применяя (1.3.6) сначала с множествами $ J $ и $ \{j\} $ (см. (1.3.2)), а
затем с множествами $ J \cup \{j\} $ и $ J^\prime \setminus \{j\}, $ получаем
(1.3.6) в рассматриваемой ситуации, что и завершает доказательство (1.3.6) в
общей ситуации.

Для вывода (1.3.5) заметим, что при $ d \in \N, l \in \N^d, 1 \le p < \infty,
\delta, \Delta \in \R_+^d, x^0, X^0 \in \R^d $ для $ f \in L_p(\R^d) $ при
$ j=1,\ldots,d, \xi_j \in \R $ для $ \chi^j_{\Delta_j, X^0_j}, $
удовлетворяющих условиям леммы 1.3.3,  ввиду п. 2) леммы 1.3.1 в
$ (X^0 +\Delta I^d)_{l_j \xi_j e_j} $ справедливо равенство
\begin{multline*} \tag{1.3.7}
\Delta_{\xi_j e_j}^{l_j} (V_j(E -M_{\chi^j_{\Delta_j, X^0_j}}
P_{\delta_j, x_j^0}^{1,l_j-1}) f) =\\
\Delta_{\xi_j e_j}^{l_j} f -
\Delta_{\xi_j e_j}^{l_j} (V_j(M_{\chi^j_{\Delta_j, X^0_j}}
P_{\delta_j, x_j^0}^{1,l_j-1}) f).
\end{multline*}
При этом в силу (1.3.1) почти для всех
$ (x_1,\ldots,x_{j-1},x_{j+1},\ldots,x_d) \in \R^{d-1} $ для почти
всех $ x_j \in (X_j^0 +\Delta_j I)_{l_j \xi_j} $ выполняется равенство
\begin{multline*}
(\Delta_{\xi_j e_j}^{l_j} (V_j (M_{\chi^j_{\Delta_j, X^0_j}}
P_{\delta_j, x_j^0}^{1,l_j-1}) f))(x_1,\ldots,x_{j-1},x_j,x_{j+1},\ldots,x_d)
= \\
(\Delta_{\xi_j}^{l_j} ((M_{\chi^j_{\Delta_j, X^0_j}}
P_{\delta_j,x_j^0}^{1,l_j-1}) f(x_1,\ldots,x_{j-1},
\cdot,x_{j+1},\ldots,x_d)))(x_j) = \\
(\Delta_{\xi_j}^{l_j} (P_{\delta_j,x_j^0}^{1,l_j-1}
f(x_1,\ldots,x_{j-1},\cdot,x_{j+1},\ldots,x_d)))(x_j) =0.
\end{multline*}
Сопоставляя сказанное с (1.3.7), приходим к выводу, что почти для
всех $ x \in (X^0 +\Delta I^d)_{l_j \xi_j e_j} $ имеет место равенство
$ (\Delta_{\xi_j e_j}^{l_j} (V_j(E -M_{\chi^j_{\Delta_j, X^0_j}}
P_{\delta_j, x_j^0}^{1,l_j-1}) f))(x) = (\Delta_{\xi_j e_j}^{l_j} f)(x). $

Для доказательства (1.3.5) покажем, что при $ d \in \N, l \in \N^d,
1 \le p < \infty, x^0, X^0 \in \R^d, \delta, \Delta \in \R_+^d $ для $ f \in
L_p(\R^d), \xi \in \R^d, \mathcal J \subset \{1,\ldots,d\} $ и любого набора
функций $ \chi^j_{\Delta_j, X^0_j}, j \in \mathcal J, $ со свойствами,
указанными в лемме, почти для всех $ x \in (X^0 +\Delta I^d)_\xi^{l \chi_{\mathcal J}} $
выполняется равенство
\begin{multline*} \tag{1.3.8}
\biggl((\prod_{j \in \mathcal J} \Delta_{\xi_j e_j}^{l_j}) ((\prod_{j \in
\mathcal J} V_j(E -M_{\chi^j_{\Delta_j, X^0_j}}
P_{\delta_j, x_j^0}^{1,l_j-1})) f)\biggr)(x) = \\
\biggl((\prod_{j \in \mathcal J} \Delta_{\xi_j e_j}^{l_j}) f\biggr)(x).
\end{multline*}

Равенство (1.3.8) установим по индукции относительно $ \card \mathcal J. $

Как показано выше, в случае, когда $ \card \mathcal J =1, $
равенство (1.3.8) справедливо. Предположим, что при $ k \in \N $
оно верно в случае, когда $ \card \mathcal J \le k. $ Проверим,
что тогда оно соблюдается и в ситуации, когда $ \card \mathcal J = k+1. $

Фиксируя $ j \in \mathcal J $ и применяя (1.3.8) сначала с
множеством $ \{j\}, $ а затем с множеством $ \mathcal J \setminus
\{j\}, $ с учётом (1.3.2) получаем (1.3.8) в случае, когда $ \card
\mathcal J = k+1, $ что и завершает вывод (1.3.8) в общем случае.

Теперь при $ d \in \N, l \in \N^d, 1 \le p < \infty, x^0, X^0 \in \R^d,
\delta, \Delta \in \R_+^d $ для $ f \in L_p(\R^d), \xi \in \R^d $ и любых
множеств $ J,J^\prime \subset \{1,\ldots,d\}, $ функций
$ \chi^j_{\Delta_j, X^0_j}, j \in J^\prime, $ описанных в лемме 1.3.3,
используя (1.3.2) и (1.3.8), имеем
\begin{multline*} \tag{1.3.9}
\biggl(\Delta_\xi^{l \chi_J} ((\prod_{j \in J^\prime} V_j(E -
M_{\chi^j_{\Delta_j, X^0_j}} P_{\delta_j, x_j^0}^{1,l_j -1})) f)\biggr)(x) = \\
\biggl(\Delta_\xi^{l \chi_J} ((\prod_{j \in J^\prime \setminus J}
V_j(E -M_{\chi^j_{\Delta_j, X^0_j}} P_{\delta_j, x_j^0}^{1,l_j-1})) f)\biggr)(x), x \in
(X^0 +\Delta I^d)_\xi^{l \chi_J}.
\end{multline*}

Соединяя (1.3.9) с (1.3.6), получаем (1.3.5). $ \square $
\bigskip

1.4. В этом пункте вводятся в рассмотрение пространства кусочно-
полиномиальных функций и операторы в них, которые используются для построения
средств приближения функций из изучаемых нами пространств. Но сначала приведём
некоторые вспомогательные сведения.

Введём в рассмотрение систему разбиений единицы на открытых множествах, с
помощью которой строятся средства приближения функций из интересующих нас
пространств. Для этого обозначим через $ \psi^{1,0} $ характеристическую функцию
интервала $ I, $ т.е. функцию, определяемую равенством
$$
\psi^{1,0}(x) = \begin{cases} 1, & \text{ для } x \in I; \\
0, & \text{ для } x \in \R \setminus I.
\end{cases}
$$
При $ m \in \N $ положим
$$
\psi^{1,m}(x) = \int_I \psi^{1, m-1}(x -y) dy \ (\text{см., например, } [14]),
$$
а для $ d \in \N, m \in \Z_+^d $ определим
$$
\psi^{d,m}(x) = \prod_{j=1}^d \psi^{1,m_j}(x_j), x =
(x_1,\ldots,x_d) \in \R^d.
$$

Для $ d \in \N, m,n \in \Z^d: m \le n, $ обозначим
\begin{equation*} \tag{1.4.1}
\Nu_{m,n}^d = \{ \nu \in \Z^d: m \le \nu \le n \} = \prod_{j=1}^d
\Nu_{m_j,n_j}^1.
\end{equation*}

Опираясь на определения, используя индукцию, нетрудно проверить
следующие свойства функций $ \psi^{d,m}, d \in \N, m \in \Z_+^d. $

1) При $ d \in \N, m \in \Z_+^d $
$$
\sgn \psi^{d,m}(x) = \begin{cases} 1, \text{ для } x \in ((m +\e) I^d); \\
0, \text{ для } x \in \R^d \setminus ((m +\e) I^d),
\end{cases}
$$

2) при $ d \in \N, m \in \Z_+^d $ для каждого $ \lambda \in \Z_+^d(m) $
(обобщённая) производная $ \D^\lambda \psi^{d,m} \in L_\infty(\R^d), $

3) при $ d \in \N, m \in \Z_+^d $ почти для всех $ x \in \R^d $
справедливо равенство
$$
\sum_{\nu \in \Z^d} \psi^{d,m}(x -\nu) =1,
$$

4) при $ m \in \N $ для всех $ x \in \R $ (при $ m =0 $ почти для всех $ x \in \R $)
имеет место равенство
\begin{equation*} \tag{1.4.2}
\psi^{1,m}(x) = \sum_{\mu \in \Nu_{0, m+1}^1} a_{\mu}^m
\psi^{1,m}(2x -\mu),
\end{equation*}
где $ a_\mu^m = 2^{-m} C_{m+1}^\mu. $
Используя разложение Ньютона для $ (1+1)^{m+1} $ и $ (-1+1)^{m+1}, $ легко
проверить, что при $ m \in \Z_+ $ выполняются равенства
\begin{equation*} \tag{1.4.3}
\sum_{\mu \in \Nu_{0,m +1}^1 \cap (2 \Z)} a_\mu^m =1,
\sum_{\mu \in \Nu_{0,m +1}^1 \cap (2 \Z +1)} a_\mu^m =1.
\end{equation*}

При $ d \in \N $ для $ t \in \R^d $ через $ 2^t $ будем обозначать
вектор $ 2^t = (2^{t_1}, \ldots, 2^{t_d}). $
Для $ d \in \N, m,\kappa \in \Z_+^d, \nu \in \Z^d $ обозначим
$$
g_{\kappa, \nu}^{d,m}(x) = \psi^{d,m}(2^\kappa x -\nu) =
\prod_{j =1}^d \psi^{1,m_j}( 2^{\kappa_j} x_j -\nu_j), x \in \R^d.
$$
Из первого среди приведенных выше свойств функций $ \psi^{d,m} $
следует, что при $ d \in \N, m,\kappa \in \Z_+^d, \nu \in \Z^d $ носитель
\begin{equation*} \tag{1.4.4}
\supp g_{\kappa,\nu}^{d,m} =
2^{-\kappa} \nu +2^{-\kappa} (m +\e) \overline I^d.
\end{equation*}
При $ d \in \N, \kappa \in \Z_+^d, \nu \in \Z^d $ обозначим
\begin{equation*} \tag{1.4.5}
Q_{\kappa, \nu}^d = 2^{-\kappa} \nu +2^{-\kappa} I^d,
\overline Q_{\kappa, \nu}^d = 2^{-\kappa} \nu +2^{-\kappa} \overline I^d.
\end{equation*}

Отметим некоторые полезные для нас свойства носителей функций
$ g_{\kappa,\nu}^{d,m}. $

При $ d \in \N, m,\kappa \in \Z_+^d $ для каждого $ \nu^\prime \in \Z^d $
имеет место равенство
\begin{equation*} \tag{1.4.6}
\{ \nu \in \Z^d: Q_{\kappa, \nu^\prime}^d \cap
\supp g_{\kappa, \nu}^{d,m} \ne \emptyset\} = \nu^\prime +\Nu_{-m,0}^d.
\end{equation*}

Из свойства 3) функций $ \psi^{d,m} $ вытекает, что при $ d \in \N,
m, \kappa \in \Z_+^d $ для любого открытого множества $ U \subset \R^d $
почти для всех $ x \in U $ соблюдается равенство
\begin{equation*}
\sum_{ \nu \in \Z^d: \supp g_{\kappa, \nu}^{d,m} \cap U \ne \emptyset}
g_{\kappa, \nu}^{d,m}(x) =1.
\end{equation*}

Имея в виду свойство 2) функций $ \psi^{d,m}, $ отметим, что при
$ d \in \N, m,\kappa \in \Z_+^d, \nu \in \Z^d, \lambda \in \Z_+^d(m)$
выполняется равенство
\begin{multline*} \tag{1.4.7}
\| \D^\lambda g_{\kappa, \nu}^{d,m} \|_{L_\infty (\R^d)} =
2^{(\kappa, \lambda)} \| \D^\lambda \psi^{d,m} \|_{L_\infty(\R^d)} =
c_1(d,m,\lambda) 2^{(\kappa, \lambda)}.
\end{multline*}

Введём в рассмотрение следующие пространства кусочно-полиномиальных
функций.
При $ d \in \N, l \in \Z_+^d, m \in \N^d, \kappa \in \Z_+^d $ и открытого
множества $ U \subset \R^d, $ полагая
\begin{equation*} \tag{1.4.8}
N_\kappa^{d,m,U} = \{\nu \in \Z^d: \supp g_{\kappa, \nu}^{d,m}
\cap U \ne \emptyset\},
\end{equation*}
через $ \mathcal P_\kappa^{d,l,m,U} $ обозначим линейное пространство,
состоящее из функций $ f: \R^d \mapsto \R, $ для каждой из которых существует
набор полиномов $ \{f_\nu \in \mathcal P^{d,l}, \nu \in N_\kappa^{d,m,U}\} $
такой, что для $ x \in \R^d $ выполняется равенство
\begin{equation*} \tag{1.4.9}
f(x) = \sum_{\nu \in N_\kappa^{d,m,U}} f_\nu(x) g_{\kappa,\nu}^{d,m}(x).
\end{equation*}

Нетрудно проверить, что при $ d \in \N, l \in \Z_+^d, m \in \N^d,
\kappa \in \Z_+^d $ и ограниченного открытого множества $ U \subset \R^d $
отображение, которое каждому набору полиномов $ \{f_\nu \in \mathcal P^{d, l},
\nu \in N_\kappa^{d,m,U} \} $ ставит в соответствие функцию $ f, $ задаваемую
равенством (1.4.9), является изоморфизмом прямого произведения
$ \card N_\kappa^{d,m,U} $ экземпляров пространства $ \mathcal P^{d,l} $
на пространство $ \mathcal P_\kappa^{d,l,m,U}. $

Опираясь на (1.4.2) устанавливается следующая лемма (см. [6]).

Лемма 1.4.1

Пусть $ d \in \N, l \in \Z_+^d, m \in \N^d, \kappa \in \Z_+^d, U $ -- открытое
ограниченное множество в $ \R^d. $ Тогда при $ j =1,\ldots,d $ линейный
оператор $ H_\kappa^{j,d,l,m,U}:
\mathcal P_\kappa^{d,l,m,U} \mapsto \mathcal P_{\kappa +e_j}^{d,l,m,U}, $
значение которого на функции $ f \in \mathcal P_\kappa^{d,l,m,U}, $ задаваемой
равенством (1.4.9), определяется соотношением
\begin{multline*} \tag{1.4.10}
(H_\kappa^{j,d,l,m,U} f)(x) = \\
\sum_{\nu \in N_{\kappa +e_j}^{d,m,U}}
\biggl(\sum_{\substack{\nu^\prime \in N_\kappa^{d,m,U}, \mu_j \in \Nu_{0, m_j +1}^1: \\
 2 \nu^\prime_j +\mu_j = \nu_j, \nu^\prime_i = \nu_i, i = 1,\ldots,d, i \ne j }}
a_{\mu_j}^{m_j} f_{\nu^\prime}(x)\biggr) g_{\kappa +e_j,\nu}^{d,m}(x), x \in \R^d,
\end{multline*}
обладает тем свойством, что для $ f \in \mathcal P_\kappa^{d,l,m,U} $ выполняется равенство
\begin{equation*}
(H_\kappa^{j,d,l,m,U} f) \mid_{U} = f \mid_{U}.
\end{equation*}

При $ m \in \N^d, \epsilon \in \Upsilon^d, \nu \in \Z^d $ обозначим
через $ \M_{\epsilon}^m(\nu) $ множество наборов чисел
\begin{multline*} \tag{1.4.11}
\M_{\epsilon}^m(\nu) = \{ \m^{\epsilon} = \{ \m_j \in \Nu_{0, m_j +1}^1,
j \in \s(\epsilon)\}: \\
(\nu_j -\m_j) /2 \in \Z \ \forall j \in \s(\epsilon)\} = \\
\prod_{j \in \s(\epsilon)} \{ \m_j \in \Nu_{0, m_j +1}^1: (\nu_j -\m_j) /2 \in \Z\} =
\prod_{j \in \s(\epsilon)} \M_1^{m_j}(\nu_j),
\end{multline*}
и каждой паре $ \nu \in \Z^d, \m^{\epsilon} \in \M_{\epsilon}^m(\nu) $ сопоставим
$ \n_{\epsilon}(\nu,\m^{\epsilon}) \in \Z^d, $ полагая
\begin{equation*} \tag{1.4.12}
(\n_{\epsilon}(\nu,\m^{\epsilon}))_j = \begin{cases} (\nu_j -\m_j) /2,
j \in \s(\epsilon); \\
\nu_j, j \in \{1,\ldots,d\} \setminus \s(\epsilon).
\end{cases}
\end{equation*}

В дальнейшем будет полезен следующий факт. При $ d \in \N, m \in \N^d $ для
$ \nu \in \Z^d, \epsilon, \epsilon^\prime \in \Upsilon^d: \s(\epsilon) \cap
\s(\epsilon^\prime) = \emptyset, $ (а, значит, $ \epsilon +\epsilon^\prime \in
\Upsilon^d, \s(\epsilon +\epsilon^\prime) = \s(\epsilon) \cup
\s(\epsilon^\prime) $) и любых $ \m^{\epsilon +\epsilon^\prime}
\in \M_{\epsilon +\epsilon^\prime}^m(\nu), \m^{\epsilon} \in
\M_{\epsilon}^m(\nu), \m^{\epsilon^\prime} \in
\M_{\epsilon^\prime}^m(\nu): (\m^{\epsilon +\epsilon^\prime})_j =
(\m^{\epsilon})_j, j \in \s(\epsilon), (\m^{\epsilon
+\epsilon^\prime})_j = (\m^{\epsilon^\prime})_j, j \in
\s(\epsilon^\prime), $ соблюдается равенство
\begin{equation*} \tag{1.4.13}
\n_{\epsilon +\epsilon^\prime}(\nu, \m^{\epsilon +\epsilon^\prime}) =
\n_{\epsilon^\prime}(\n_{\epsilon}(\nu, \m^{\epsilon}), \m^{\epsilon^\prime}).
\end{equation*}

Замечание.

При $ d \in \N, m \in \N^d $ для $ \kappa \in \Z_+^d,
\epsilon \in \Upsilon^d: \s(\epsilon) \subset \s(\kappa), $ для открытого
множества $ U \subset \R^d $ и $ \nu \in N_{\kappa}^{d,m,U}, \m^{\epsilon} \in
\M_{\epsilon}^m(\nu) $ имеет место включение
\begin{equation*} \tag{1.4.14}
\n_{\epsilon}(\nu,\m^{\epsilon}) \in N_{\kappa -\epsilon}^{d,m,U}.
\end{equation*}

В самом деле, при соблюдении условий замечания, выбирая $ x \in U \cap
\supp g_{\kappa,\nu}^{d,m}, $ ввиду (1.4.4) получаем, что
\begin{equation*}
2^{-\kappa_i} \nu_i \le x_i \le 2^{-\kappa_i} \nu_i +2^{-\kappa_i}
(m_i +1) \text{ при } i = 1,\ldots,d,
\end{equation*}
откуда с учётом (1.4.12) имеем
\begin{multline*}
2^{-(\kappa -\epsilon)_i} (\n_{\epsilon}(\nu,\m^{\epsilon}))_i =
2^{-\kappa_i} \nu_i \le x_i \le 2^{-\kappa_i} \nu_i +
2^{-\kappa_i} (m_i +1) = \\
2^{-(\kappa -\epsilon)_i} (\n_{\epsilon}(\nu,\m^{\epsilon}))_i +
2^{-(\kappa -\epsilon)_i} (m_i +1) \text{ при } i = 1,\ldots,d: i \notin \s(\epsilon); \\
2^{-(\kappa -\epsilon)_j} (\n_{\epsilon}(\nu,\m^{\epsilon}))_j =
2^{-\kappa_j +1} (\nu_j -\m_j) /2 =
2^{-\kappa_j} \nu_j -2^{-\kappa_j} \m_j \le \\
2^{-\kappa_j} \nu_j \le x_j \le 2^{-\kappa_j} \nu_j +2^{-\kappa_j} (m_j +1) = \\
2^{-\kappa_j +1} (\nu_j -\m_j) /2 +
2^{-\kappa_j +1} (\m_j /2) +2^{-\kappa_j +1} (m_j +1) /2 = \\
2^{-\kappa_j +1} (\nu_j -\m_j) /2 +2^{-\kappa_j +1} (\m_j +m_j +1) /2 \le
2^{-\kappa_j +1} (\nu_j -\m_j) /2 +2^{-\kappa_j +1} (m_j +1) = \\
2^{-(\kappa -\epsilon)_j} (\n_{\epsilon}(\nu,\m^{\epsilon}))_j +
2^{-(\kappa -\epsilon)_j} (m_j +1) \text{ при } j \in \s(\epsilon),
\end{multline*}
т.е. $ x \in U \cap \supp g_{\kappa -\epsilon, \n_{\epsilon}(\nu,\m^{\epsilon})}^{d,m}, $
а, следовательно, $ \n_{\epsilon}(\nu,\m^{\epsilon}) \in N_{\kappa -\epsilon}^{d,m,U}. $

Частным случаем леммы 1.2.2 из [6] является лемма 1.4.2, для формулировки
которой при $ d \in \N $ для $ j \in \{1,\ldots,d\} $ обозначим через
$ \eta^j: \R^d \times \R^d \mapsto \R^d $ отображение, определяемое соотношением
$$
(\eta^j(\xi,x))_i = \begin{cases} \xi_i, i =1, \ldots, j; \\
x_i, i = j +1, \ldots, d,
\end{cases} \xi,x \in \R^d.
$$

Отметим также, что доказательство леммы 1.4.2 опирается на лемму 1.4.1.

Лемма 1.4.2

Пусть $ d \in \N, l \in \Z_+^d, m \in \N^d, \kappa \in \Z_+^d, \epsilon \in
\Upsilon^d: \s(\epsilon) \subset \s(\kappa), U $ -- ограниченное открытое
множество в $ \R^d. $ Тогда линейный оператор
$ H_{\kappa, \kappa -\epsilon}^{d,l,m,U}:
\mathcal P_{\kappa -\epsilon}^{d,l,m,U} \mapsto \mathcal P_\kappa^{d,l,m,U}, $
значение которого для $ f \in \mathcal P_{\kappa -\epsilon}^{d,l,m,U} $
определяется равенством
\begin{equation*} \tag{1.4.15}
H_{\kappa, \kappa -\epsilon}^{d,l,m,U} f = \begin{cases} f, \text{ при }
\epsilon =0; \\
(\prod_{j \in \s(\epsilon)} H_{\eta^j(\kappa -\epsilon, \kappa)}^{j,d,l,m,U}) f,
\text{ при } \epsilon \ne 0, (\text{ см. } (1.4.10)),
\end{cases}
\end{equation*}
обладает следующими свойствами:

1) для $ f \in \mathcal P_{\kappa -\epsilon}^{d,l,m,U} $ выполняется равенство
\begin{equation*} \tag{1.4.16}
(H_{\kappa, \kappa -\epsilon}^{d,l,m,U} f) \mid_{U} = f \mid_{U};
\end{equation*}

2) для $ f \in \mathcal P_{\kappa -\epsilon}^{d,l,m,U} $ вида
\begin{equation*} \tag{1.4.17}
f = \sum_{\nu^\prime \in N_{\kappa -\epsilon}^{d,m,U}}
f_{\kappa -\epsilon, \nu^\prime} g_{\kappa -\epsilon, \nu^\prime}^{d,m},
\{f_{\kappa -\epsilon, \nu^\prime} \in \mathcal P^{d,l},
\nu^\prime \in N_{\kappa -\epsilon}^{d,m,U}\},
\end{equation*}
имеет место представление
\begin{equation*} \tag{1.4.18}
H_{\kappa, \kappa -\epsilon}^{d,l,m,U} f = \sum_{\nu \in N_{\kappa}^{d,m,U}}
f_{\kappa,\nu} g_{\kappa,\nu}^{d,m},
\end{equation*}
где
\begin{equation*} \tag{1.4.19}
f_{\kappa,\nu} = \sum_{\m^{\epsilon} \in \M_{\epsilon}^m(\nu)} A_{\m^{\epsilon}}^m
f_{\kappa -\epsilon, \n_{\epsilon}(\nu,\m^{\epsilon})},
\end{equation*}
а
\begin{equation*} \tag{1.4.20}
A_{\m^{\epsilon}}^m = \prod_{i \in \s(\epsilon)}
a_{\m_i}^{m_i}, (\text{ см. } (1.4.2)), \m^{\epsilon} \in \M_{\epsilon}^m(\nu),
\nu \in N_{\kappa}^{d,m,U}.
\end{equation*}

Отметим ещё частный случай леммы 1.2.3 из [6] (см. (1.4.3), (1.4.20)).

Лемма 1.4.3

При $ d \in \N, \nu \in \Z^d, \epsilon \in \Upsilon^d, m \in \N^d $ имеет место
равенство
\begin{equation*} \tag{1.4.21}
\sum_{\m^{\epsilon} \in \M_{\epsilon}^m(\nu)} A_{\m^{\epsilon}}^m =1.
\end{equation*}
\bigskip

\centerline{\S 2. Оценка сверху наилучшей точности восстановления
в $ L_q(D) $ производных $ \D^\lambda f$}
\centerline{по значениям в $n$ точках функций $ f $ из
$ (\mathcal S_{p,\theta}^\alpha \mathcal B)^\prime(D) $}
\bigskip

2.1. В этом пункте будут построены средства приближения функций из
рассматриваемых нами пространств, на которые опирается вывод основных
результатов работы.

При $ d \in \N, l, \kappa \in \Z_+^d, \nu \in \Z^d $ определим линейный
оператор $ S_{\kappa,\nu}^{d,l}: L_1(Q_{\kappa,\nu}^d) \mapsto \mathcal P^{d,l}, $
полагая $ S_{\kappa, \nu}^{d,l} = P_{\delta, x^0}^{d,l} $ при
$ \delta = 2^{-\kappa}, x^0 = 2^{-\kappa} \nu $ (см. лемму 1.1.2 и (1.4.5)).
Отметим, что в ситуации, когда $ Q_{\kappa,\nu}^d \subset D, $ где $ D $ --
область в $ \R^d, $ для $ f \in L_1(D) $ вместо $ S_{\kappa, \nu}^{d,l}(f \mid_{Q_{\kappa, \nu}^d}) $
будем писать $ S_{\kappa, \nu}^{d,l} f.$

При $ d \in \N $ для ограниченной области $ D \subset \R^d, $ её открытого
подмножества $ U \subset D $ и $ \kappa \in \Z_+^d $ таких, что множество
\begin{equation*} \tag{2.1.1}
\{\nu^\prime \in \Z^d: Q_{\kappa,\nu^\prime}^d \subset D\} \ne \emptyset, \
(\text{см. } (1.4.5)),
\end{equation*}
а $ m \in \N^d, $ фиксируем некоторое отображение
\begin{equation*} \tag{2.1.2}
\nu_\kappa = \nu_\kappa^{d,m,D,U}: N_\kappa^{d,m,U} \ni \nu \mapsto
\nu_\kappa^{d,m,D,U}(\nu) \in \{\nu^\prime \in \Z^d:
Q_{\kappa,\nu^\prime}^d \subset D\} \text{ (см. } (1.4.8)),
\end{equation*}
и при $ l \in \Z_+^d $ определим линейный непрерывный оператор
$ E_\kappa^{d,l,m,D,U,\nu_\kappa}: L_1(D) \mapsto \mathcal P_\kappa^{d,l,m,U} \cap
L_\infty(\R^d) $ (см. п. 1.4.) равенством
\begin{equation*} \tag{2.1.3}
E_\kappa^{d,l,m,D,U,\nu_\kappa} f = \sum_{\nu \in N_\kappa^{d,m,U}}
(S_{\kappa, \nu_\kappa^{d,m,D,U}(\nu)}^{d,l} f ) g_{\kappa, \nu}^{d,m}, \
f \in L_1(D).
\end{equation*}

Замечание.

Учитывая, что при $ d \in \N, \kappa \in \Z_+^d, \nu \in \Nu_{0, 2^\kappa -\e}^d $
клетка $ (2^{-\kappa} \nu +2^{-\kappa} I^d) \subset I^d, $ а, значит,
\begin{multline*}
Q_{\kappa^0,\nu^0}^d = (2^{-\kappa^0}  \nu^0 +2^{-\kappa^0} I^d) \supset
(2^{-\kappa^0} \nu^0 +2^{-\kappa^0} (2^{-\kappa} \nu +2^{-\kappa} I^d)) = \\
(2^{-\kappa^0 -\kappa} (2^\kappa \nu^0 +\nu) +2^{-\kappa^0 -\kappa} I^d) =
Q_{\kappa^0 +\kappa, 2^\kappa \nu^0 +\nu}^d, \kappa^0 \in \Z_+^d,
\nu^0 \in \Z^d,
\end{multline*}
замечаем, что если при $ d \in \N $ для области $ D \subset \R^d $ и $ \kappa^0
\in \Z_+^d $ соблюдается (2.1.1) при $ \kappa^0 $ вместо $ \kappa, $
то для $ \kappa \in \Z_+^d $ имеет место (2.1.1) при $ \kappa^0 +\kappa $ вместо
$ \kappa. $

Как покозано в [13], имеет место следующее утверждение.

Предложение 2.1.1
.

Пусть $ d \in \N, l \in \N^d, m \in \N^d, 1 \le p < \infty $
а ограниченная область $ D \subset \R^d $ и её открытое подмножество
$ U \subset D $ таковы, что существуют $ \kappa^0 = \kappa^0(d,m,D,U) \in \Z_+^d,
\gamma^0 = \gamma^0(d,m,D,U) \in \R_+^d, $ для которых при любом $ \kappa \in \Z_+^d $
существует отображение $ \nu_{\kappa^0 +\kappa} =
\nu_{\kappa^0 +\kappa}^{d,m,D,U}: N_{\kappa^0 +\kappa}^{d,m,U} \mapsto \Z^d, $
обладающее тем свойством, что для каждого
$ \nu \in N_{\kappa^0 +\kappa}^{d,m,U} $ справедливо включение
\begin{equation*} \tag{2.1.4}
Q_{\kappa^0 +\kappa,\nu_{\kappa^0 +\kappa}^{d,m,D,U}(\nu)}^d \subset D \cap
(2^{-\kappa^0 -\kappa} \nu +\gamma^0 2^{-\kappa^0 -\kappa} B^d).
\end{equation*}
Тогда для любой функции $ f \in L_p(D) $ в $ L_p(U) $ имеет место равенство
\begin{equation*} \tag{2.1.5}
f \mid_U = \lim_{\mn(\kappa) \to \infty} (E_{\kappa^0 +\kappa}^{d,l -\e,m,D,U,
\nu_{\kappa^0 +\kappa}} f) \mid_U \text{ (см. } (2.1.3)).
\end{equation*}

Ещё будет нужна следующая теорема.

Теорема 2.1.2

Пусть $ d \in \N, \alpha \in \R_+^d, 1 \le p < \infty, 1 \le q \le \infty $ и
$ \lambda \in \Z_+^d $ удовлетворяют условию
\begin{equation*} \tag{2.1.6}
\alpha -\lambda -(p^{-1} -q^{-1})_+ \e >0.
\end{equation*}
Тогда существуют константы $ c_{1}(d,\alpha) >0, c_2(d,\alpha,p,q,\lambda) >0 $
и $ c_{3}(d,\alpha,p,q,\lambda) >0 $ такие, что при любых $ \delta \in
\R_+^d, x^0 \in \R^d $ для $ Q = x^0 +\delta I^d $ для $ f \in
(S_p^\alpha H)^\prime(Q) $ при $ l = l(\alpha) $ выполняются неравенства
\begin{multline*} \tag{2.1.7}
\| \D^\lambda f -\D^\lambda P_{\delta, x^0}^{d, l -\e} f\|_{L_q(Q)} \le
c_{2} \delta^{-\lambda -p^{-1} \e +q^{-1} \e}
\sum_{ J \subset \Nu_{1,d}^1: J \ne \emptyset} (\prod_{j \in J}
\delta_j^{\lambda_j +(p^{-1} -q^{-1})_+}) \\ \times
\int_{ (c_{1} \delta I^d)^J} (\prod_{j \in J} t_j^{-\lambda_j -p^{-1} -(p^{-1}
-q^{-1})_+ -1}) \biggl(\int_{(t B^d)^J} \int_{ Q_\xi^{l \chi_J}}
|\Delta_\xi^{l \chi_J} f(x)|^p dx d\xi^J\biggr)^{1/p} dt^J,
\end{multline*}

\begin{multline*} \tag{2.1.8}
\| \D^\lambda f \|_{L_q(Q)} \le c_{3} \delta^{-\lambda -p^{-1} \e +q^{-1} \e}
\biggl(\|f\|_{L_p(Q)} +\sum_{ J \subset \Nu_{1,d}^1: J \ne \emptyset}
(\prod_{j \in J} \delta_j^{\lambda_j +(p^{-1} -q^{-1})_+}) \\ \times
\int_{ (c_{1} \delta I^d)^J} (\prod_{j \in J}
t_j^{-\lambda_j -p^{-1} -(p^{-1} -q^{-1})_+ -1})
\biggl(\int_{(t B^d)^J} \int_{ Q_\xi^{l \chi_J}}
|\Delta_\xi^{l \chi_J} f(x)|^p dx d\xi^J\biggr)^{1/p} dt^J\biggr).
\end{multline*}

Доказательство.

Прежде всего отметим, что в силу (2.3.16) из [13] имеет место включение
\begin{equation*}
(S_p^\alpha H)^\prime(Q) \subset (S_p^\alpha H)^0(Q).
\end{equation*}
Используя определённые в [4] объекты, применяя (1.1.3), (1.1.4), (1.1.5),
а также (1.1.3) и (1.1.2) из [4], для $ f \in (S_p^\alpha H)^\prime(Q) $ имеем
\begin{multline*}
\| \D^\lambda f -\D^\lambda P_{\delta, x^0}^{d, l -\e} f\|_{L_q(Q)} \le
\| \D^\lambda f -\D^\lambda P_{\delta, x^0}^{d, l -\e,\lambda} f\|_{L_q(Q)} + \\
\| \D^\lambda P_{\delta, x^0}^{d, l -\e,\lambda} f -
\D^\lambda P_{\delta, x^0}^{d, l -\e,0} f\|_{L_q(Q)} +
\| \D^\lambda P_{\delta, x^0}^{d, l -\e,0} f -
\D^\lambda P_{\delta, x^0}^{d, l -\e} f\|_{L_q(Q)} \le \\
\| \D^\lambda f -\D^\lambda P_{\delta, x^0}^{d, l -\e,\lambda} f\|_{L_q(Q)} + \\
c_4 \delta^{-\lambda -p^{-1} \e +q^{-1} \e}
\| P_{\delta, x^0}^{d, l -\e,\lambda} f -
P_{\delta, x^0}^{d, l -\e,0} f\|_{L_p(Q)} + \\
c_4 \delta^{-\lambda -p^{-1} \e +q^{-1} \e}
\| P_{\delta, x^0}^{d, l -\e,0} f -
P_{\delta, x^0}^{d, l -\e} f\|_{L_p(Q)}  = \\
\| \D^\lambda f -\D^\lambda P_{\delta, x^0}^{d, l -\e,\lambda} f\|_{L_q(Q)} +
c_4 \delta^{-\lambda -p^{-1} \e +q^{-1} \e}
\| P_{\delta, x^0}^{d, l -\e,\lambda} (f -
P_{\delta, x^0}^{d, l -\e,0} f)\|_{L_p(Q)} + \\
c_4 \delta^{-\lambda -p^{-1} \e +q^{-1} \e}
\| P_{\delta, x^0}^{d, l -\e} ((P_{\delta, x^0}^{d, l -\e,0} f) -f)\|_{L_p(Q)} \le \\
\| \D^\lambda f -\D^\lambda P_{\delta, x^0}^{d, l -\e,\lambda} f\|_{L_q(Q)} +
c_5 \delta^{-\lambda -p^{-1} \e +q^{-1} \e}
\| f -P_{\delta, x^0}^{d, l -\e,0} f\|_{L_p(Q)}.
\end{multline*}
Ввиду отмеченного выше включения, применяя (1.3.52) из [4], получаем
\begin{multline*}
\| \D^\lambda f -\D^\lambda P_{\delta, x^0}^{d, l -\e,\lambda} f\|_{L_q(Q)} \le
c_{6} \delta^{-\lambda -p^{-1} \e +q^{-1} \e}
\sum_{ J \subset \Nu_{1,d}^1: J \ne \emptyset} (\prod_{j \in J}
\delta_j^{\lambda_j +(p^{-1} -q^{-1})_+}) \\ \times
\int_{ (c_{1} \delta I^d)^J} (\prod_{j \in J} t_j^{-\lambda_j -p^{-1} -(p^{-1}
-q^{-1})_+ -1}) \biggl(\int_{(t B^d)^J} \int_{ Q_\xi^{l \chi_J}}
|\Delta_\xi^{l \chi_J} f(x)|^p dx d\xi^J\biggr)^{1/p} dt^J,
\end{multline*}
\begin{multline*}
\| f -P_{\delta, x^0}^{d, l -\e,0} f\|_{L_p(Q)} \le \\
c_{7} \sum_{ J \subset \Nu_{1,d}^1: J \ne \emptyset} \int_{ (c_{1} \delta I^d)^J}
(\prod_{j \in J} t_j^{-p^{-1} -1}) \biggl(\int_{(t B^d)^J}
\int_{ Q_\xi^{l \chi_J}} |\Delta_\xi^{l \chi_J} f(x)|^p dx d\xi^J\biggr)^{1/p} dt^J \le \\
c_7 \sum_{ J \subset \Nu_{1,d}^1: J \ne \emptyset} \int_{ (c_{1} \delta I^d)^J}
\biggl(\prod_{j \in J} (c_1 \delta_j)^{\lambda_j +(p^{-1} -q^{-1})_+}
t_j^{-\lambda_j -p^{-1} -(p^{-1} -q^{-1})_+ -1}\biggr) \\
\times \biggl(\int_{(t B^d)^J}
\int_{ Q_\xi^{l \chi_J}} |\Delta_\xi^{l \chi_J} f(x)|^p dx d\xi^J\biggr)^{1/p} dt^J \le \\
c_8 \sum_{ J \subset \Nu_{1,d}^1: J \ne \emptyset}
(\prod_{j \in J} \delta_j^{\lambda_j +(p^{-1} -q^{-1})_+})
\int_{ (c_{1} \delta I^d)^J}
(\prod_{j \in J} t_j^{-\lambda_j -p^{-1} -(p^{-1} -q^{-1})_+ -1}) \\
\biggl(\int_{(t B^d)^J} \int_{ Q_\xi^{l \chi_J}}
|\Delta_\xi^{l \chi_J} f(x)|^p dx d\xi^J\biggr)^{1/p} dt^J.
\end{multline*}
Объединяя полученные оценки, приходим к (2.1.7). Неравенство (2.1.8) является
следствием (1.3.53) из [4]. $ \square $

Предложение 2.1.3

Пусть $ d \in \N, \alpha \in \R_+^d, 1 \le p < \infty $ удовлетворяют условию
\begin{equation*} \tag{2.1.9}
\alpha -p^{-1} \e >0.
\end{equation*}
Тогда для любой области $ D \subset \R^d $ для любой функции $ f \in
(S_p^\alpha H)^\prime(D) $ существует функция $ F \in C(D), $ для которой
почти для всех $ x \in D $ соблюдается равенство $ F(x) = f(x). $

Доказательство.

Для произвольной области $ D \subset \R^d $ и $ f \in (S_p^\alpha H)^\prime(D), $
представляя $ D $ в виде счётного объединения кубов
$ D = \cup_{i \in \N} Q^i, $ где $ Q^i = x^i +\delta^i I^d,
x^i \in \R^d, \delta^i > 0, i \in \N, $ для каждой функции
$ f \mid_{Q^i}(x^i +\delta^i \cdot) \in (S_p^\alpha h)^\prime(I^d), i \in \N, $
на основании предложения 1.3.10 из [4] (см. также (2.3.16) из [13]) возьмём
функцию $ F^i \in C(Q^i) $ такую, что почти для всех $ x \in Q^i $ верно
равенство $ F^i(x) = f \mid_{Q^i}(x) = f(x), i \in \N. $ Принимая во внимание,
что при $ i, j \in \N $ почти для всех $ x \in Q^i \cap Q^j $ имеет место
равенство $ F^i(x) = f(x) = F^j(x), $ а, следовательно, $ F^i(x) = F^j(x) $ для
всех $ x \in Q^i \cap Q^j, $ определим функцию $ F: D \mapsto \R $ равенством
$ F(x) = F^i(x) $ при $ x \in Q^i, i \in \N. $ Понятно, что построенная функия
$ F \in C(D) $ и почти для всех $ x \in D $ выпоняется равенство $ F(x) = f(x).  \square $

Теперь напомним некоторые сведения, касающиеся интерполяционных полиномов.

При $ l \in \Z_+ $ в интервале $ I $ фиксируем систему из $ (l+1) $
различных точек $ \{ \xi_{1,0}^{1,l,\lambda} \in I, \lambda = 0,\ldots, l \} $
и построим систему полиномов $ \{\pi_{1,0}^{1,l,\lambda} \in
\mathcal P^{1,l}, \lambda =0, \ldots,l \}, $ обладающую тем свойством, что при
$ \lambda, \mu =0,\ldots, l $ выполняются равенства
$$
\pi_{1,0}^{1,l,\lambda}(\xi_{1,0}^{1,l,\mu}) =
\begin{cases} 1, \text{ при } \lambda = \mu; \\
0, \text{ при } \lambda \ne \mu.
\end{cases}
$$

При $ l \in \Z_+, \delta \in \R_+, x^0 \in \R $ определим систему
точек $ \{ \xi_{\delta, x^0}^{1,l,\lambda} \in x^0 +\delta I,
\lambda =0,\ldots,l\} $ и систему полиномов
$ \{ \pi_{\delta,x^0}^{1,l,\lambda} \in \mathcal P^{1,l}, \lambda =0,\ldots,l\}, $
полагая
$$
\xi_{\delta, x^0}^{1,l,\lambda} =
x^0 +\delta \xi_{1,0}^{1,l,\lambda}, \\
\pi_{\delta, x^0}^{1,l,\lambda}(x)  =
\pi_{1,0}^{1,l,\lambda}(\delta^{-1} (x -x^0)).
$$
Ясно, что при $ \lambda, \mu =0, \ldots,l $ соблюдаются равенства
\begin{equation*} \tag{2.1.10}
\pi_{\delta, x^0}^{1,l,\lambda}(\xi_{\delta, x^0}^{1,l,\mu}) =
\pi_{1,0}^{1,l,\lambda}(\xi_{1,0}^{1,l,\mu}) = \begin{cases}
1, \text{ при } \lambda = \mu; \\
0, \text{ при } \lambda \ne \mu.
\end{cases}
\end{equation*}

При $ d \in \N, l \in \Z_+^d, \delta \in \R_+^d, x^0 \in \R^d $
построим систему точек $ \{ \xi_{\delta, x^0}^{d,l,\lambda} \in
x^0 +\delta I^d, \lambda \in \Z_+^d(l) \} $ и систему полиномов
$ \{ \pi_{\delta, x^0}^{d,l,\lambda} \in \mathcal P^{d,l}, \lambda \in \Z_+^d(l) \}, $
задавая
$$
( \xi_{\delta, x^0}^{d,l,\lambda})_j =
\xi_{\delta_j, x_j^0}^{1,l_j,\lambda_j}, j =1,\ldots,d, \\
\pi_{\delta, x^0}^{d,l,\lambda}(x) =
\prod_{j =1}^d
\pi_{\delta_j, x_j^0}^{1,l_j,\lambda_j}(x_j).
$$

Из этих определений и (2.1.10) следует, что при $ \lambda, \mu \in \Z_+^d(l) $
справедливы равенства
\begin{equation*} \tag{2.1.11}
\pi_{\delta, x^0}^{d,l,\lambda}(\xi_{\delta, x^0}^{d,l,\mu}) =
\prod_{j =1}^d \pi_{\delta_j,x_j^0}^{1,l_j,\lambda_j}
(\xi_{\delta_j, x_j^0}^{1,l_j,\mu_j}) =
\begin{cases}
1, \text{ при } \lambda = \mu; \\
0, \text{ при } \lambda \ne \mu.
\end{cases}
\end{equation*}

Понятно, что
\begin{multline*}
\pi_{\delta, x^0}^{d,l,\lambda}(x) = \prod_{j =1}^d
\pi_{1, 0}^{1,l_j,\lambda_j}(\delta_j^{-1} (x_j -x_j^0)) =
\pi_{\e, 0}^{d,l,\lambda}(\delta^{-1} (x -x^0)), \\
\xi_{\delta, x^0}^{d,l,\lambda} = x^0 +\delta \xi_{\e,0}^{d,l,\lambda}, \lambda
\in \Z_+^d(l).
\end{multline*}

При $ d \in \N, l \in \Z_+^d, \delta \in \R_+^d, x^0 \in \R^d $ обозначим через
$ A_{\delta, x^0}^{d,l}: \R^{ (l+\e)^{\e} } \mapsto \mathcal P^{d,l} $ линейный
оператор, который каждому набору чисел $ t = \{ t_\lambda \in \R, \lambda \in
\Z_+^d(l) \} $ ставит в соответствие полином
$$
A_{\delta, x^0}^{d,l} t = \sum_{ \lambda \in \Z_+^d(l)} t_\lambda
\pi_{\delta, x^0}^{d,l,\lambda} \in \mathcal P^{d,l},
$$
а через $ \phi_{\delta, x^0}^{d,l} $ обозначим
линейное отображение пространства $ C(x^0 +\delta I^d) $
в $ \R^{(l +\e)^{\e}}, $ которое каждой функции $ f \in C(x^0 +\delta I^d) $
сопоставляет набор её значений $ \{ f(\xi_{\delta, x^0}^{d,l,\lambda}),
\lambda \in \Z_+^d(l) \}, $ и определим линейный оператор
$ \mathcal P_{\delta,x^0}^{d,l}: C(x^0 +\delta I^d) \mapsto \mathcal P^{d,l} $
равенством
$$
\mathcal P_{\delta, x^0}^{d,l} =
A_{\delta, x^0}^{d,l} \circ \phi_{\delta, x^0}^{d,l} = \sum_{\lambda \in
\Z_+^d(l)} f(\xi_{\delta, x^0}^{d,l,\lambda}) \pi_{\delta, x^0}^{d,l,\lambda}.
$$

Как видно из опредений и (2.1.11), при $ \lambda \in \Z_+^d(l) $ имеет место
равенство
\begin{equation*} \tag{2.1.12}
(\mathcal P_{\delta, x^0}^{d,l} f) (\xi_{\delta,x^0}^{d,l,\lambda}) =
f(\xi_{\delta, x^0}^{d,l,\lambda}).
\end{equation*}

Заметим, что вследствие (2.1.11) система полиномов
$ \{\pi_{\delta, x^0}^{d,l,\lambda} \in \mathcal P^{d,l}, \lambda \in
\Z_+^d(l) \} $ -- линейно независима, и число элементов этой системы
$ \card \Z_+^d(l) = \dim \mathcal P^{d,l}, $ и, следовательно, эта система
является базисом в $ \mathcal P^{d,l}.$
Поэтому для любого полинома $ f \in \mathcal P^{d,l} $ существует набор чисел
$ \{ t_\lambda \in \R, \lambda \in\Z_+^d(l) \} $ такой, что
$$
f = \sum_{ \lambda \in \Z_+^d(l)} t_\lambda \pi_{\delta,x^0}^{d,l,\lambda}.
$$
Отсюда ввиду (2.1.11) получаем
$$
t_\lambda = f(\xi_{\delta, x^0}^{d,l,\lambda}), \lambda \in \Z_+^d(l),
$$
т.е.
\begin{equation*} \tag{2.1.13}
f = \mathcal P_{\delta, x^0}^{d,l} f, f \in \mathcal P^{d,l}.
\end{equation*}

Лемма 2.1.4

Пусть $ d \in \N, l \in \Z_+^d, \delta \in \R_+^d, x^0 \in \R^d $
и для каждого $ j =1,\ldots,d $ функция $ \chi^j $ непрерывна на $ \R, $ имеет
компактный носитель, а $ \chi(x) = \prod_{j =1}^d \chi^j(x_j). $
Тогда для любой функции $ f \in C_0(\R^d) $ для всех $ x \in \R^d $
справедливо равенство
\begin{equation*} \tag{2.1.14}
\chi(x) (\mathcal P_{\delta, x^0}^{d,l}(f \mid_{x^0 +\delta I^d}))(x)
= \biggl((\prod_{j =1}^d V_j(M_{\chi^j}
\mathcal P_{\delta_j, x_j^0}^{1,l_j})) f\biggr)(x).
\end{equation*}

Доказательство.

Введём следующее обозначение. При $ d \in \N $ для $ J \subset \Nu_{1,d}^1 $
обозначим через $ \eta_J: \R^d \times \R^d \mapsto \R^d $ отображение, у
которого
$$
(\eta_J(\xi,x))_j =
\begin{cases}
\xi_j, j \in J; \\
x_j, j \in \Nu_{1,d}^1 \setminus J,
\end{cases}
\xi,x \in \R^d.
$$

В условиях леммы для любой функции $ f \in C_0(\R^d) $ и любого
множества $ J \subset \Nu_{1,d}^1 $ имеет место равенство
\begin{multline*} \tag{2.1.15}
((\prod_{j \in J} V_j(M_{\chi^j}
\mathcal P_{\delta_j, x_j^0}^{1, l_j})) f)(x) =\\
(\prod_{j \in J} \chi^j(x_j))
\sum_{ \lambda \in \Z_+^d(l \chi_J)}
f(\eta_J(\xi_{\delta,x^0}^{d,l,\lambda}, x))
\prod_{j \in J} \pi_{\delta_j,x_j^0}^{1,l_j,\lambda_j}(x_j),
x \in \R^d.
\end{multline*}

Равенство (2.1.15) установим по индукции относительно $ \card J. $

В случае, когда $ \card J =1, $ и, следовательно, $ J = \{j\}, $
где $ j \in \Nu_{1,d}^1, $ согласно (1.3.3) для $ f \in C_0(\R^d), x \in
\R^d $ имеем
\begin{multline*}
(V_j(M_{\chi^j} \mathcal P_{\delta_j, x_j^0}^{1, l_j}) f)(x) =
((M_{\chi^j} \mathcal P_{\delta_j, x_j^0}^{1, l_j})
f(x_1, \ldots, x_{j-1}, \cdot, x_{j+1}, \ldots, x_d))(x_j) = \\
\chi^j(x_j) \sum_{\lambda_j =0}^{l_j}
f(x_1, \ldots, x_{j-1}, \xi_{\delta_j,x_j^0}^{1,l_j,\lambda_j}, x_{j+1}, \ldots, x_d)
\pi_{\delta_j,x_j^0}^{1,l_j,\lambda_j}(x_j),
\end{multline*}
что совпадает с (2.1.15) при $ J = \{j\}. $

Предположим теперь, что равенство (2.1.15) справедливо для любого
множества $ J \subset \Nu_{1,d}^1, $ у которого $ \card J \le k, (1 \le k \le d-1), $
и покажем, что тогда оно верно для любого множества $ J \subset
\Nu_{1,d}^1, $ у которого $ \card J = k+1. $
В этой ситуации, представляя $ J $ в виде $ J = J^\prime \cup \{j\}, $
где $ j \not \in J^\prime, $ на основании предположения индукции
с учётом (1.3.4) для $ x \in \R^d $ получаем
\begin{multline*}
\biggl((\prod_{i \in J} V_i(M_{\chi^i}
\mathcal P_{\delta_i, x_i^0}^{1, l_i})) f\biggr)(x) = \\
\biggl(V_j(M_{\chi^j} \mathcal P_{\delta_j, x_j^0}^{1, l_j})
((\prod_{i \in J^\prime} V_i(M_{\chi^i}
\mathcal P_{\delta_i, x_i^0}^{1, l_i})) f)\biggr)(x) = \\
\biggl(V_j(M_{\chi^j} \mathcal P_{\delta_j, x_j^0}^{1, l_j})
((\prod_{i \in J^\prime} \chi^i(y_i))
\sum_{ \lambda^\prime \in \Z_+^d(l \chi_{J^\prime})}
f(\eta_{J^\prime}(\xi_{\delta,x^0}^{d,l,\lambda^\prime}, y))
\prod_{i \in J^\prime}
\pi_{\delta_i,x_i^0}^{1,l_i,\lambda_i^\prime}(y_i))\biggr)(x) = \\
\sum_{ \lambda^\prime \in \Z_+^d(l \chi_{J^\prime})}
(V_j(M_{\chi^j} \mathcal P_{\delta_j, x_j^0}^{1, l_j})
(f(\eta_{J^\prime}(\xi_{\delta,x^0}^{d,l,\lambda^\prime}, y))
\prod_{i \in J^\prime} (\chi^i(y_i)
\pi_{\delta_i,x_i^0}^{1,l_i,\lambda_i^\prime}(y_i))))(x) = \\
\sum_{ \lambda^\prime \in \Z_+^d(l \chi_{J^\prime})}
\chi^j(x_j) \sum_{\lambda_j =0}^{l_j}
f(\eta_{J^\prime}(\xi_{\delta,x^0}^{d,l,\lambda^\prime},
x_1, \ldots, x_{j-1}, \xi_{\delta_j,x_j^0}^{1,l_j,\lambda_j}, x_{j +1}, \ldots, x_d))\\
\times (\prod_{i \in J^\prime} (\chi^i(x_i)
\pi_{\delta_i,x_i^0}^{1,l_i,\lambda_i^\prime}(x_i)))
\pi_{\delta_j,x_j^0}^{1,l_j,\lambda_j}(x_j) = \\
(\prod_{i \in J} \chi^i(x_i))
\sum_{ \lambda \in \Z_+^d(l \chi_J)}
f(\eta_J(\xi_{\delta,x^0}^{d,l,\lambda}, x))
\prod_{i \in J} \pi_{\delta_i,x_i^0}^{1,l_i,\lambda_i}(x_i),
\end{multline*}
что и завершает вывод (2.1.15).

Полагая в (2.1.15) $ J = \Nu_{1,d}^1, $ заключаем, что в условиях леммы
при $ x \in \R^d $ выполняется равенство
\begin{multline*}
((\prod_{j =1}^d V_j(M_{\chi^j}
\mathcal P_{\delta_j, x_j^0}^{1, l_j})) f)(x) =\\
(\prod_{j =1}^d \chi^j(x_j))
\sum_{ \lambda \in \Z_+^d(l)}
f(\xi_{\delta,x^0}^{d,l,\lambda})
\prod_{j =1}^d  \pi_{\delta_j,x_j^0}^{1,l_j,\lambda_j}(x_j) =
\chi(x) (\mathcal P_{\delta,x^0}^{d,l} f)(x),
\end{multline*}
что совпадает с (2.1.14). $ \square $

Для доказательства леммы 2.1.6 понадобится следующая лемма (см. [4]).

Лемма 2.1.5

Пусть $ d \in \N, l \in \Z_+^d. $ Тогда  существует константа $ c_9(d,l) >0 $
такая, что для любого $ \delta \in \R_+^d $ и любой точки $ x^0 \in \R^d $
для $ Q = x^0 +\delta I^d $ для любой функции $ f \in C(Q) \cap L_\infty(Q) $
справедлива оценка
\begin{equation*} \tag{2.1.16}
\| \mathcal P_{\delta, x^0}^{d,l} f\|_{L_\infty(Q)} \le
c_9 \| f\|_{L_\infty(Q)}.
\end{equation*}

Доказательство.

В условиях леммы ввиду определений имеем
\begin{multline*}
\| \mathcal P_{\delta, x^0}^{d,l} f\|_{L_\infty(Q)} =
\| \sum_{\lambda \in \Z_+^d(l)} f(\xi_{\delta, x^0}^{d,l,\lambda})
\pi_{\delta, x^0}^{d,l,\lambda} \|_{L_\infty(Q)} \le \\
\sum_{\lambda \in \Z_+^d(l)} | f(\xi_{\delta, x^0}^{d,l,\lambda})| \cdot
\| \pi_{\delta, x^0}^{d,l,\lambda} \|_{L_\infty(Q)} \le
\biggl(\sum_{\lambda \in \Z_+^d(l)} \| \pi_{\delta, x^0}^{d,l,\lambda} \|_{L_\infty(Q)}\biggr)
\|f\|_{L_\infty(Q)}.
\end{multline*}
Отсюда, пользуясь тем, что
\begin{multline*}
\| \pi_{\delta, x^0}^{d,l,\lambda} \|_{L_\infty(Q)} =
\sup_{x \in Q} | \pi_{\delta, x^0}^{d,l,\lambda}(x)| =
\sup_{x \in (x^0 +\delta I^d)} | \pi_{\delta, x^0}^{d,l,\lambda}(x)| = \\
\sup_{x \in I^d} | \pi_{\e, 0}^{d,l,\lambda}(\delta^{-1}(x^0 +\delta x -x^0))| =
\sup_{x \in I^d} | \pi_{\e, 0}^{d,l,\lambda}(x)|, \lambda \in
\Z_+^d(l),
\end{multline*}
приходим к (2.1.16). $ \square $

Лемма 2.1.6

Пусть $ d \in \N, l \in \Z_+^d $ и $ \rho, \sigma \in \R_+^d. $
Тогда существует константа $ c_{10}(d,l,\rho,\sigma) >0 $ такая, что для любых
$ x^0, \bm x^0, X^0 \in \R^d $ и $ \delta, \bm \delta, \Delta \in \R_+^d $ таких, что
$ (x^0 +\sigma \delta I^d) \subset (\bm x^0 +\bm \delta I^d) \subset
(\bm x^0 +\rho \delta I^d) \cap (X^0 +\Delta I^d), $
для любого множества $ J \subset \{1,\ldots,d\}, $ и любого набора функций
$ \chi^j_{\Delta_j, X^0_j}, j \in J, $ каждая из которых непрерывна на $ \R, $
имеет компактный носитель и $ \chi^j_{\Delta_j, X^0_j}(x) =1 $ для $ x \in
(X^0_j +\Delta_j I), $ для $ f \in C_0(\R^d) $ имеет место неравенство
\begin{equation*} \tag{2.1.17}
\biggl\| (\prod_{j \in J} V_j(E -M_{\chi^j_{\Delta_j, X^0_j}}
\mathcal P_{\sigma_j \delta_j,x_j^0}^{1,l_j})) f \biggr\|_{L_\infty(D)} \le
c_{10} \|f\|_{L_\infty(D)},
\end{equation*}
где $ D = (\bm x^0 +\bm \delta I^d). $

Доказательство.

Ввиду леммы 1.3.2 понятно, что (2.1.17) достаточно доказать в случае, когда
$ J = \{j\}, j \in \{1,\ldots,d\}. $ Переходя к выводу (2.1.17) в этом
случае, заметим, что в условиях леммы имеет место включение
\begin{equation*} \tag{2.1.18}
(\bm x^0_j +\bm \delta_j I) \subset (\bm x^0_j +\rho_j \delta_j I) \subset
(x_j^0 +\rho_j \delta_j B^1), j \in \Nu_{1,d}^1.
\end{equation*}

Далее, ввиду п. 2) леммы 1.3.2 имеем
\begin{equation*} \tag{2.1.19}
\| V_j(E -M_{\chi^j_{\Delta_j, X^0_j}} \mathcal P_{\sigma_j \delta_j,
x_j^0}^{1,l_j}) f\|_{L_\infty(D)} \le
\|f\|_{L_\infty(D)} +\| V_j(M_{\chi^j_{\Delta_j, X^0_j}}
\mathcal P_{\sigma_j \delta_j,x_j^0}^{1,l_j}) f\|_{L_\infty(D)}.
\end{equation*}

Оценивая второе слагаемое в правой части (2.1.19), с помощью соотношений
(1.3.3), (1.1.3) (с учётом (2.1.18)), (2.1.16) находим, что
\begin{multline*} \tag{2.1.20}
\| V_j(M_{\chi^j_{\Delta_j, X^0_j}}
\mathcal P_{\sigma_j \delta_j, x_j^0}^{1,l_j}) f\|_{L_\infty(D)} = \\
\sup_{ x \in D} | ((M_{\chi^j_{\Delta_j, X^0_j}}
\mathcal P_{\sigma_j \delta_j, x_j^0}^{1,l_j}) f(x_1,\ldots,x_{j -1},
\cdot,x_{j +1},\ldots,x_d))(x_j)| = \\
\sup_{\substack{ (x_1,\ldots,x_{j -1},x_{j +1},\ldots,x_d) \in \\ \prod_{i =1,\ldots,d:
i \ne j} (\bm x_i^0 +\bm \delta_i I)}}
\sup_{ x_j \in (\bm x^0_j +\bm \delta_j I)} | \chi^j_{\Delta_j, X^0_j}(x_j)
(\mathcal P_{\sigma_j \delta_j,x_j^0}^{1,l_j}
f(x_1,\ldots,x_{j-1},\cdot,x_{j+1},\ldots,x_d))(x_j)| = \\
\sup_{\substack{ (x_1,\ldots,x_{j -1},x_{j +1},\ldots,x_d) \in \\ \prod_{i =1,\ldots,d:
i \ne j} (\bm x_i^0 +\bm \delta_i I)}}
\sup_{ x_j \in (\bm x^0_j +\bm \delta_j I)}
| (\mathcal P_{\sigma_j \delta_j,x_j^0}^{1,l_j}
f(x_1,\ldots,x_{j-1},\cdot,x_{j+1},\ldots,x_d))(x_j)| \le \\
\sup_{ (x_1,\ldots,x_{j -1},x_{j +1},\ldots,x_d) \in \prod_{i =1,\ldots,d:
i \ne j} (\bm x_i^0 +\bm \delta_i I)}
c_1(1,l_j,0,\rho_j,\sigma_j) \\
\times \sup_{ x_j \in (x_j^0 +\sigma_j \delta_j I)}
|(\mathcal P_{\sigma_j \delta_j,x_j^0}^{1,l_j}
f(x_1,\ldots,x_{j-1},\cdot,x_{j+1},\ldots,x_d))(x_j)| \le \\
c_{11} \sup_{ (x_1,\ldots,x_{j -1},x_{j +1},\ldots,x_d) \in \prod_{i =1,\ldots,d:
i \ne j} (\bm x_i^0 +\bm \delta_i I)} \\
\sup_{ x_j \in (x_j^0 +\sigma_j \delta_j I)}
| f(x_1,\ldots,x_{j-1},x_j,x_{j+1},\ldots,x_d)| \le \\
c_{11} \sup_{ (x_1,\ldots,x_{j -1},x_{j +1},\ldots,x_d) \in \prod_{i =1,\ldots,d:
i \ne j} (\bm x_i^0 +\bm \delta_i I)} \\
\sup_{ x_j \in (\bm x_j^0 +\bm \delta_j I)}
| f(x_1,\ldots,x_{j-1},x_j,x_{j+1},\ldots,x_d)| = \\
c_{11} \sup_{x \in D} |f(x)| = c_{11} \|f\|_{L_\infty(D)}.
\end{multline*}
Соединяя (2.1.19) с (2.1.20), приходим к (2.1.17) в ситуации, когда $ J = \{j\},
j \in \{1,\ldots,d\}. \square $

Определим операторы, с которыми мы будем иметь дело в дальнейшем, и отметим их
свойства полезные для нас.

Для $ d \in \N, l, \kappa \in \Z_+^d, \nu \in \Z^d $ определим линейный
оператор $ \mathtt R_{\kappa,\nu}^{d,l}: C(Q_{\kappa,\nu}^d) \mapsto \mathcal P^{d,l} $
равенством
$$
\mathtt R_{\kappa, \nu}^{d,l} f =
\mathcal P_{\delta, x^0}^{d,l} f, f \in C(x^0 +\delta I^d),
$$
при $ \delta = 2^{-\kappa}, x^0 = 2^{-\kappa} \nu. $

Как видно из (2.1.14), для $ d \in \N, l, \kappa \in \Z_+^d, \nu \in \Z^d $
для любого набора функций $ \chi^j, j = 1,\ldots,d, \chi, $ удовлетворяющих
условиям леммы 2.1.4, для $ f \in C_0(\R^d) $ и $ x \in \R^d $
верно равенство
\begin{equation*} \tag{2.1.21}
\chi(x) (\mathtt R_{\kappa, \nu}^{d,l} (f \mid_{Q_{\kappa,\nu}^d}))(x) =
((\prod_{j =1}^d V_j(M_{\chi^j} \mathtt R_{\kappa_j, \nu_j}^{1,l_j})) f)(x).
\end{equation*}

При $ d \in \N $ для ограниченной области $ D \subset \R^d, $ её открытого
подмножества $ U \subset D $ и $ \kappa \in \Z_+^d $ таких, что выполняется
(2.1.1), для $ m \in \N^d $ и некоторого фиксированного отображения $ \nu_\kappa
= \nu_\kappa^{d,m,D,U} $ вида (2.1.2), при $ l \in \Z_+^d $ определим линейный
оператор
$$
R_\kappa^{d,l,m,D,U,\nu_\kappa}: C(D) \mapsto \mathcal P_\kappa^{d,l,m,U}
$$
равенством
\begin{equation*} \tag{2.1.22}
R_\kappa^{d,l,m,D,U,\nu_\kappa} f = \sum_{ \nu \in N_\kappa^{d,m,U}}
(\mathtt R_{\kappa, \nu_\kappa^{d,m,D,U}(\nu)}^{d,l} f) g_{\kappa, \nu}^{d,m}.
\end{equation*}

Предложение 2.1.7

Пусть $ d \in \N, m \in \N^d, $ а ограниченная область $ D \subset \R^d $
и её открытое подмножество $ U \subset D $ удовлетворяют условиям предложения
2.1.1. Пусть ещё для $ \alpha \in \R_+^d, 1 \le p < \infty $ соблюдается
неравенство (2.1.9). Тогда для любой функции $ f \in (S_p^\alpha H)^\prime(D) $
при $ l = l(\alpha) $ в $ L_p(U) $ имеет место равенство
\begin{equation*} \tag{2.1.23}
f \mid_U = \lim_{\mn(\kappa) \to \infty} (R_{\kappa^0 +\kappa}^{d,l -\e,m,D,U,\nu_{\kappa^0 +\kappa}} f) \mid_U.
\end{equation*}

Доказательство.

Сначала отметим свойства некоторых вспомогательных множеств и других
объектов, которые понадобятся для доказательства предложения.

Для $ d \in \N, \kappa \in \Z_+^d, m \in \N^d $ и открытого множества
$ U \subset \R^d $ обозначим множество
$$
G_\kappa^{d,m,U} = \cup_{\nu \in N_\kappa^{d,m,U}} \supp g_{\kappa,\nu}^{d,m},
$$
для которого справедливо представление
\begin{multline*} \tag{2.1.24}
G_\kappa^{d,m,U} = (\cup_{n \in \Z^d: Q_{\kappa,n}^d \cap G_\kappa^{d,m,U} \ne
\emptyset} Q_{\kappa,n}^d) \cup (G_\kappa^{d,m,U} \cap A_\kappa^d), \\
\text{ где } \mes A_\kappa^d =0, A_\kappa^d \cap Q_{\kappa,n}^d = \emptyset,
Q_{\kappa,n}^d \cap Q_{\kappa,n^\prime}^d = \emptyset, n,n^\prime \in
\Z^d: n \ne n^\prime.
\end{multline*}
Понятно, что для $ n \in \Z^d $ такого, что пересечение $ Q_{\kappa,n}^d \cap
G_\kappa^{d,m,U} \ne \emptyset, $ множество
$ \{\nu \in N_\kappa^{d,m,U}: Q_{\kappa,n}^d \cap
\supp g_{\kappa,\nu}^{d,m} \ne \emptyset\} \ne \emptyset, $ при этом, используя
(1.4.4), (1.4.5), легко проверить, что для $ \nu \in N_\kappa^{d,m,U}:
Q_{\kappa,n}^d \cap \supp g_{\kappa,\nu}^{d,m} \ne \emptyset, $ верно
соотношение
\begin{equation*} \tag{2.1.25}
Q_{\kappa,n}^d \subset \overline Q_{\kappa,n}^d \subset
\supp g_{\kappa,\nu}^{d,m} \subset G_\kappa^{d,m,U}.
\end{equation*}

Отметим ещё, что в условиях предложения при $ \kappa \in \Z_+^d $
для $ n \in \Z^d: Q_{\kappa^0 +\kappa,n}^d \cap G_{\kappa^0 +\kappa}^{d,m,U} \ne \emptyset,
\nu \in N_{\kappa^0 +\kappa}^{d,m,U}: \supp g_{\kappa^0 +\kappa, \nu}^{d,m} \cap Q_{\kappa^0 +\kappa,n}^d
\ne \emptyset, $
вследствие (2.1.4) и замкнутости $ B^d, $ а также благодаря (2.1.25) и (1.4.4)
(при $ \kappa^0 +\kappa $ вместо $ \kappa $ ) справедливо включение
\begin{equation*} \tag{2.1.26}
Q_{\kappa^0 +\kappa,\nu_{\kappa^0 +\kappa}^{d,m,D,U}(\nu)}^d \subset
\overline Q_{\kappa^0 +\kappa,\nu_{\kappa^0 +\kappa}^{d,m,D,U}(\nu)}^d
\subset (2^{-\kappa^0 -\kappa} n +\gamma^1 2^{-\kappa^0 -\kappa} B^d),
\end{equation*}
\begin{equation*} \tag{2.1.27}
Q_{\kappa^0 +\kappa,n}^d \subset \overline Q_{\kappa^0 +\kappa,n}^d
\subset (2^{-\kappa^0 -\kappa} \nu_{\kappa^0 +\kappa}^{d,m,D,U}(\nu) +
(\gamma^1 +\e) 2^{-\kappa^0 -\kappa} B^d),
\end{equation*}
с $ \gamma^1 = \gamma^1(d,m,D,U) > \e. $

В условиях предложения при $ \kappa \in \Z_+^d, n \in \Z^d:
Q_{\kappa^0 +\kappa,n}^d \cap G_{\kappa^0 +\kappa}^{d,m,U} \ne \emptyset, $
зададим
\begin{equation*}
x_{\kappa^0 +\kappa,n}^{d,m,D,U} = 2^{-\kappa^0 -\kappa} n -\gamma^1 2^{-\kappa^0 -\kappa}; \\
\delta_{\kappa^0 +\kappa,n}^{d,m,D,U} = 2 \gamma^1 2^{-\kappa^0 -\kappa}
\end{equation*}
и определим клетку $ D_{\kappa^0 +\kappa,n}^{d,m,D,U} $ равенством
\begin{equation*} \tag{2.1.28}
D_{\kappa^0 +\kappa,n}^{d,m,D,U} = x_{\kappa^0 +\kappa,n}^{d,m,D,U} +
\delta_{\kappa^0 +\kappa,n}^{d,m,D,U} I^d = \inter (2^{-\kappa^0 -\kappa} n +\gamma^1 2^{-\kappa^0 -\kappa} B^d).
\end{equation*}

Из приведенных определений с учётом того, что $ \gamma^1 > \e, $ видно, что
соблюдается включение
\begin{equation*} \tag{2.1.29}
Q_{\kappa^0 +\kappa, n}^d \subset D_{\kappa^0 +\kappa,n}^{d,m,D,U},
n \in \Z^d: Q_{\kappa^0 +\kappa,n}^d \cap G_{\kappa^0 +\kappa}^{d,m,U} \ne \emptyset,
\kappa \in \Z_+^d.
\end{equation*}

Учитывая (2.1.29), (2.1.28), нетрудно видеть, что в условиях предложения
существует константа $ c_{12}(d,m,D,U) >0 $ такая, что при $ \kappa \in \Z_+^d $
для каждого $ x \in \R^d $ число
\begin{equation*} \tag{2.1.30}
\card \{ n \in \Z^d: Q_{\kappa^0 +\kappa,n}^d \cap G_{\kappa^0 +\kappa}^{d,m,U} \ne \emptyset,
x \in D_{\kappa^0 +\kappa,n}^{d,m,D,U} \} \le c_{12}.
\end{equation*}

Из (2.1.26) и (2.1.28) следует, что при $ \kappa \in \Z_+^d$ для $ n \in \Z^d:
Q_{\kappa^0 +\kappa,n}^d \cap G_{\kappa^0 +\kappa}^{d,m,U} \ne \emptyset, \nu \in N_{\kappa^0 +\kappa}^{d,m,U}:
\supp g_{\kappa^0 +\kappa, \nu}^{d,m} \cap Q_{\kappa^0 +\kappa,n}^d \ne
\emptyset, $ имеет место включение
\begin{equation*} \tag{2.1.31}
Q_{\kappa^0 +\kappa,\nu_{\kappa^0 +\kappa}^{d,m,D,U}(\nu)}^d \subset
D_{\kappa^0 +\kappa,n}^{d,m,D,U}.
\end{equation*}

Из (1.4.6) с учётом (1.4.1) вытекает, что при $ \kappa \in \Z_+^d $ для
$ n \in \Z^d: Q_{\kappa^0 +\kappa,n}^d \cap G_{\kappa^0 +\kappa}^{d,m,U} \ne
\emptyset, $ верно неравенство
\begin{equation*} \tag{2.1.32}
\card \{ \nu \in N_{\kappa^0 +\kappa}^{d,m,U}: \supp g_{\kappa^0 +\kappa, \nu}^{d,m} \cap
Q_{\kappa^0 +\kappa,n}^d \ne \emptyset \} \le c_{13}(d,m).
\end{equation*}

Пусть теперь в условиях предложения $ f \in (S_p^\alpha H)^\prime(D), l = l(\alpha),
\kappa \in \Z_+^d. $ Тогда с учётом предложения 2.1.3 имеем
\begin{multline*} \tag{2.1.33}
\| f \mid_U -(R_{\kappa^0 +\kappa}^{d,l -\e,m,D,U,\nu_{\kappa^0 +\kappa}} f) \mid_U\|_{L_p(U)} \le 
\| f \mid_U -(E_{\kappa^0 +\kappa}^{d,l -\e,m,D,U,\nu_{\kappa^0 +\kappa}} f) \mid_U \|_{L_p(U)} +\\
\| (E_{\kappa^0 +\kappa}^{d,l -\e,m,D,U,\nu_{\kappa^0 +\kappa}} f) \mid_U
-(R_{\kappa^0 +\kappa}^{d,l -\e,m,D,U,\nu_{\kappa^0 +\kappa}} f) \mid_U \|_{L_p(U)}.
\end{multline*}

В силу (2.1.5) для получения (2.1.23) остаётся показать, что второе слагаемое
в правой части (2.1.33) сходится к нулю при $ \mn(\kappa) \to \infty. $
Для этого ввиду (2.1.3), (2.1.22) имеем
\begin{multline*} \tag{2.1.34}
\| (E_{\kappa^0 +\kappa}^{d,l -\e,m,D,U,\nu_{\kappa^0 +\kappa}} f) \mid_U
-(R_{\kappa^0 +\kappa}^{d,l -\e,m,D,U,\nu_{\kappa^0 +\kappa}} f) \mid_U \|_{L_p(U)}^p \le \\
\| E_{\kappa^0 +\kappa}^{d,l -\e,m,D,U,\nu_{\kappa^0 +\kappa}} f
-R_{\kappa^0 +\kappa}^{d,l -\e,m,D,U,\nu_{\kappa^0 +\kappa}} f\|_{L_p(\R^d)}^p = \\
\biggl\| \sum_{\nu \in N_{\kappa^0 +\kappa}^{d,m,U}}
(S_{\kappa^0 +\kappa, \nu_{\kappa^0 +\kappa}^{d,m,D,U}(\nu)}^{d,l -\e} f)
g_{\kappa^0 +\kappa, \nu}^{d,m} -
\sum_{\nu \in N_{\kappa^0 +\kappa}^{d,m,U}}
(\mathtt R_{\kappa^0 +\kappa, \nu_{\kappa^0 +\kappa}^{d,m,D,U}(\nu)}^{d,l -\e} f)
g_{\kappa^0 +\kappa, \nu}^{d,m} \biggr\|_{L_p(\R^d)}^p = \\
\biggl\| \sum_{\nu \in N_{\kappa^0 +\kappa}^{d,m,U}}
((S_{\kappa^0 +\kappa, \nu_{\kappa^0 +\kappa}^{d,m,D,U}(\nu)}^{d,l -\e} f) -
(\mathtt R_{\kappa^0 +\kappa, \nu_{\kappa^0 +\kappa}^{d,m,D,U}(\nu)}^{d,l -\e} f))
g_{\kappa^0 +\kappa, \nu}^{d,m} \biggr\|_{L_p(\R^d)}^p = \\
\int_{\R^d} \biggl| \sum_{\nu \in N_{\kappa^0 +\kappa}^{d,m,U}}
((S_{\kappa^0 +\kappa, \nu_{\kappa^0 +\kappa}^{d,m,D,U}(\nu)}^{d,l -\e} f)(x) -
(\mathtt R_{\kappa^0 +\kappa, \nu_{\kappa^0 +\kappa}^{d,m,D,U}(\nu)}^{d,l -\e} f)(x))
g_{\kappa^0 +\kappa, \nu}^{d,m}(x) \biggr|^p dx = \\
\int_{G_{\kappa^0 +\kappa}^{d,m,U}}
\biggl| \sum_{\nu \in N_{\kappa^0 +\kappa}^{d,m,U}}
((S_{\kappa^0 +\kappa, \nu_{\kappa^0 +\kappa}^{d,m,D,U}(\nu)}^{d,l -\e} f)(x) \textendash\
(\mathtt R_{\kappa^0 +\kappa, \nu_{\kappa^0 +\kappa}^{d,m,D,U}(\nu)}^{d,l -\e} f)(x))
g_{\kappa^0 +\kappa, \nu}^{d,m}(x) \biggr|^p dx.
\end{multline*}

В силу (2.1.24) выводим
\begin{multline*} \tag{2.1.35}
\int_{G_{\kappa^0 +\kappa}^{d,m,U}}
\biggl| \sum_{\nu \in N_{\kappa^0 +\kappa}^{d,m,U}}
((S_{\kappa^0 +\kappa, \nu_{\kappa^0 +\kappa}^{d,m,D,U}(\nu)}^{d,l -\e} f)(x) -
(\mathtt R_{\kappa^0 +\kappa, \nu_{\kappa^0 +\kappa}^{d,m,D,U}(\nu)}^{d,l -\e} f)(x))
g_{\kappa^0 +\kappa, \nu}^{d,m}(x) \biggr|^p dx = \\
\sum_{\substack{n \in \Z^d: Q_{\kappa^0 +\kappa,n}^d \\ \cap
G_{\kappa^0 +\kappa}^{d,m,U} \ne \emptyset}}
\int\limits_{Q_{\kappa^0 +\kappa,n}^d} \biggl| \sum_{\nu \in N_{\kappa^0 +\kappa}^{d,m,U}}
((S_{\kappa^0 +\kappa, \nu_{\kappa^0 +\kappa}^{d,m,D,U}(\nu)}^{d,l -\e} f)(x) -
(\mathtt R_{\kappa^0 +\kappa, \nu_{\kappa^0 +\kappa}^{d,m,D,U}(\nu)}^{d,l -\e} f)(x))
g_{\kappa^0 +\kappa, \nu}^{d,m}(x) \biggr|^p dx = \\
\sum_{\substack{n \in \Z^d: Q_{\kappa^0 +\kappa,n}^d \\
\cap G_{\kappa^0 +\kappa}^{d,m,U} \ne \emptyset}}
\biggl\| \sum_{\substack{\nu \in N_{\kappa^0 +\kappa}^{d,m,U}:\\
\supp g_{\kappa^0 +\kappa, \nu}^{d,m} \cap Q_{\kappa^0 +\kappa,n}^d \ne \emptyset}}
(S_{\kappa^0 +\kappa, \nu_{\kappa^0 +\kappa}^{d,m,D,U}(\nu)}^{d,l -\e} f -
\mathtt R_{\kappa^0 +\kappa, \nu_{\kappa^0 +\kappa}^{d,m,D,U}(\nu)}^{d,l -\e} f)
g_{\kappa^0 +\kappa, \nu}^{d,m} \biggr\|_{L_p(Q_{\kappa^0 +\kappa,n}^d)}^p \le \\
\sum_{\substack{n \in \Z^d: Q_{\kappa^0 +\kappa,n}^d \\
\cap G_{\kappa^0 +\kappa}^{d,m,U} \ne \emptyset}}
\biggl( \sum_{\substack{\nu \in N_{\kappa^0 +\kappa}^{d,m,U}:\\
\supp g_{\kappa^0 +\kappa, \nu}^{d,m} \cap Q_{\kappa^0 +\kappa,n}^d \ne \emptyset}}
\| (S_{\kappa^0 +\kappa, \nu_{\kappa^0 +\kappa}^{d,m,D,U}(\nu)}^{d,l -\e} f -
\mathtt R_{\kappa^0 +\kappa, \nu_{\kappa^0 +\kappa}^{d,m,D,U}(\nu)}^{d,l -\e} f)
g_{\kappa^0 +\kappa, \nu}^{d,m} \|_{L_p(Q_{\kappa^0 +\kappa,n}^d)} \biggr)^p.
\end{multline*}

Далее, для $ n \in \Z^d: Q_{\kappa^0 +\kappa,n}^d \cap G_{\kappa^0 +\kappa}^{d,m,U} \ne \emptyset,
\nu \in N_{\kappa^0 +\kappa}^{d,m,U}: \supp g_{\kappa^0 +\kappa, \nu}^{d,m}
\cap Q_{\kappa^0 +\kappa,n}^d \ne \emptyset, $ учитывая (1.4.7), а затем
благодаря (2.1.27) используя (1.1.3), получаем
\begin{multline*} \tag{2.1.36}
\| (S_{\kappa^0 +\kappa, \nu_{\kappa^0 +\kappa}^{d,m,D,U}(\nu)}^{d,l -\e} f -
\mathtt R_{\kappa^0 +\kappa, \nu_{\kappa^0 +\kappa}^{d,m,D,U}(\nu)}^{d,l -\e} f)
g_{\kappa^0 +\kappa, \nu}^{d,m} \|_{L_p(Q_{\kappa^0 +\kappa,n}^d)} \le \\
\| S_{\kappa^0 +\kappa, \nu_{\kappa^0 +\kappa}^{d,m,D,U}(\nu)}^{d,l -\e} f -
\mathtt R_{\kappa^0 +\kappa, \nu_{\kappa^0 +\kappa}^{d,m,D,U}(\nu)}^{d,l -\e} f
\|_{L_p(Q_{\kappa^0 +\kappa,n}^d)} \le \\
c_{14} (2^{-\kappa^0 -\kappa})^{p^{-1} \e}
\| S_{\kappa^0 +\kappa, \nu_{\kappa^0 +\kappa}^{d,m,D,U}(\nu)}^{d,l -\e} f -
\mathtt R_{\kappa^0 +\kappa, \nu_{\kappa^0 +\kappa}^{d,m,D,U}(\nu)}^{d,l -\e} f
\|_{L_\infty(Q_{\kappa^0 +\kappa,\nu_{\kappa^0 +\kappa}^{d,m,D,U}(\nu)}^d)} = \\
c_{15} 2^{-(\kappa, p^{-1} \e)}
\| S_{\kappa^0 +\kappa, \nu_{\kappa^0 +\kappa}^{d,m,D,U}(\nu)}^{d,l -\e} f -
\mathtt R_{\kappa^0 +\kappa, \nu_{\kappa^0 +\kappa}^{d,m,D,U}(\nu)}^{d,l -\e} f
\|_{L_\infty(Q_{\kappa^0 +\kappa,\nu_{\kappa^0 +\kappa}^{d,m,D,U}(\nu)}^d)}.
\end{multline*}

Применяя (2.1.13) и (2.1.16), на основании (2.1.7) при $ \lambda =0,
q = \infty $ (см. также (2.1.4) и предложение 2.1.3) заключаем, что для
$ \nu \in N_{\kappa^0 +\kappa}^{d,m,U} $ выполняется неравенство
\begin{multline*} \tag{2.1.37}
\| S_{\kappa^0 +\kappa, \nu_{\kappa^0 +\kappa}^{d,m,D,U}(\nu)}^{d,l -\e} f -
\mathtt R_{\kappa^0 +\kappa, \nu_{\kappa^0 +\kappa}^{d,m,D,U}(\nu)}^{d,l -\e} f
\|_{L_\infty(Q_{\kappa^0 +\kappa,\nu_{\kappa^0 +\kappa}^{d,m,D,U}(\nu)}^d)} = \\
\| \mathtt R_{\kappa^0 +\kappa, \nu_{\kappa^0 +\kappa}^{d,m,D,U}(\nu)}^{d,l -\e}
(f -S_{\kappa^0 +\kappa, \nu_{\kappa^0 +\kappa}^{d,m,D,U}(\nu)}^{d,l -\e} f)
\|_{L_\infty(Q_{\kappa^0 +\kappa,\nu_{\kappa^0 +\kappa}^{d,m,D,U}(\nu)}^d)} \le \\
c_9 \| f -S_{\kappa^0 +\kappa, \nu_{\kappa^0 +\kappa}^{d,m,D,U}(\nu)}^{d,l -\e} f
\|_{L_\infty(Q_{\kappa^0 +\kappa,\nu_{\kappa^0 +\kappa}^{d,m,D,U}(\nu)}^d)} \le \\
c_9 c_2 (2^{-\kappa^0 -\kappa})^{-p^{-1} \e }
\sum_{ J \subset \Nu_{1,d}^1: J \ne \emptyset}
\biggl(\prod_{j \in J} (2^{-\kappa^0_j -\kappa_j})^{p^{-1}}\biggr)
\int_{ (c_1 2^{-\kappa^0 -\kappa} I^d)^J}
\biggl(\prod_{j \in J} t_j^{-2 p^{-1} -1}\biggr)\times \\
\biggl(\int_{(t B^d)^J}
\int_{ (Q_{\kappa^0 +\kappa,\nu_{\kappa^0 +\kappa}^{d,m,D,U}(\nu)}^d)_\xi^{l \chi_J}}
|\Delta_\xi^{l \chi_J} f(x)|^p dx d\xi^J\biggr)^{1/p} dt^J \le \\
c_{16} 2^{(\kappa, p^{-1} \e) }
\sum_{ J \subset \Nu_{1,d}^1: J \ne \emptyset}
\biggl(\prod_{j \in J} 2^{-\kappa_j p^{-1}}\biggr)
\int_{ (c_{17} 2^{-\kappa} I^d)^J}
\biggl(\prod_{j \in J} t_j^{-2 p^{-1} -1}\biggr)\times \\
\biggl(\int_{(t B^d)^J}
\int_{ (Q_{\kappa^0 +\kappa,\nu_{\kappa^0 +\kappa}^{d,m,D,U}(\nu)}^d)_\xi^{l \chi_J}}
|\Delta_\xi^{l \chi_J} f(x)|^p dx d\xi^J\biggr)^{1/p} dt^J.
\end{multline*}

Соединяя (2.1.36), (2.1.37) и принимая во внимание (2.1.31), (2.1.4),
для $ n \in \Z^d: Q_{\kappa^0 +\kappa,n}^d \cap
G_{\kappa^0 +\kappa}^{d,m,U} \ne \emptyset,
\nu \in N_{\kappa^0 +\kappa}^{d,m,U}: \supp g_{\kappa^0 +\kappa, \nu}^{d,m}
\cap Q_{\kappa^0 +\kappa,n}^d \ne \emptyset, $ находим, что
\begin{multline*}
\| (S_{\kappa^0 +\kappa, \nu_{\kappa^0 +\kappa}^{d,m,D,U}(\nu)}^{d,l -\e} f -
\mathtt R_{\kappa^0 +\kappa, \nu_{\kappa^0 +\kappa}^{d,m,D,U}(\nu)}^{d,l -\e} f)
g_{\kappa^0 +\kappa, \nu}^{d,m} \|_{L_p(Q_{\kappa^0 +\kappa,n}^d)} \le \\
c_{15} 2^{-(\kappa, p^{-1} \e)}
c_{16} 2^{(\kappa, p^{-1} \e) }
\sum_{ J \subset \Nu_{1,d}^1: J \ne \emptyset}
\biggl(\prod_{j \in J} 2^{-\kappa_j p^{-1}}\biggr)
\int_{ (c_{17} 2^{-\kappa} I^d)^J}
\biggl(\prod_{j \in J} t_j^{-2 p^{-1} -1}\biggr)\times \\
\biggl(\int_{(t B^d)^J}
\int_{ (Q_{\kappa^0 +\kappa,\nu_{\kappa^0 +\kappa}^{d,m,D,U}(\nu)}^d)_\xi^{l \chi_J}}
|\Delta_\xi^{l \chi_J} f(x)|^p dx d\xi^J\biggr)^{1/p} dt^J = \\
c_{18} \sum_{ J \subset \Nu_{1,d}^1: J \ne \emptyset}
\biggl(\prod_{j \in J} 2^{-\kappa_j p^{-1}}\biggr)
\int_{ (c_{17} 2^{-\kappa} I^d)^J} \biggl(\prod_{j \in J} t_j^{-2 p^{-1} -1}\biggr)\times \\
\biggl(\int_{(t B^d)^J}
\int_{ (Q_{\kappa^0 +\kappa,\nu_{\kappa^0 +\kappa}^{d,m,D,U}(\nu)}^d)_\xi^{l \chi_J}}
|\Delta_\xi^{l \chi_J} f(x)|^p dx d\xi^J\biggr)^{1/p} dt^J \le \\
c_{18} \sum_{ J \subset \Nu_{1,d}^1: J \ne \emptyset}
\biggl(\prod_{j \in J} 2^{-\kappa_j p^{-1}}\biggr)
\int_{ (c_{17} 2^{-\kappa} I^d)^J} \biggl(\prod_{j \in J} t_j^{-2 p^{-1} -1}\biggr)\times \\
\biggr(\int_{(t B^d)^J}
\int_{ (D \cap D_{\kappa^0 +\kappa,n}^{d,m,D,U})_\xi^{l \chi_J}}
|\Delta_\xi^{l \chi_J} f(x)|^p dx d\xi^J\biggr)^{1/p} dt^J \le \\
c_{18} \sum_{ J \subset \Nu_{1,d}^1: J \ne \emptyset}
\biggl(\prod_{j \in J} 2^{-\kappa_j p^{-1}}\biggr)
\int_{ (c_{17} 2^{-\kappa} I^d)^J} \biggl(\prod_{j \in J} t_j^{-2 p^{-1} -1}\biggr)\times \\
\biggl(\int_{(t B^d)^J} \int_{D_\xi^{l \chi_J} \cap D_{\kappa^0 +\kappa,n}^{d,m,D,U}}
|\Delta_\xi^{l \chi_J} f(x)|^p dx d\xi^J\biggr)^{1/p} dt^J.
\end{multline*}

Подставляя эту оценку в (2.1.35) и применяя (2.1.32) и неравенство
Гёльдера, выводим
\begin{multline*} \tag{2.1.38}
\int_{G_{\kappa^0 +\kappa}^{d,m,U}}
\biggl| \sum_{\nu \in N_{\kappa^0 +\kappa}^{d,m,U}}
((S_{\kappa^0 +\kappa, \nu_{\kappa^0 +\kappa}^{d,m,D,U}(\nu)}^{d,l -\e} f)(x) -
(\mathtt R_{\kappa^0 +\kappa, \nu_{\kappa^0 +\kappa}^{d,m,D,U}(\nu)}^{d,l -\e} f)(x))
g_{\kappa^0 +\kappa, \nu}^{d,m}(x) \biggr|^p dx \le \\
\sum_{\substack{n \in \Z^d: \\ Q_{\kappa^0 +\kappa,n}^d \cap
G_{\kappa^0 +\kappa}^{d,m,U} \ne \emptyset}}
\biggl( \sum_{\substack{\nu \in N_{\kappa^0 +\kappa}^{d,m,U}:\\
\supp g_{\kappa^0 +\kappa, \nu}^{d,m} \cap Q_{\kappa^0 +\kappa,n}^d \ne \emptyset}}
c_{18} \sum_{ J \subset \Nu_{1,d}^1: J \ne \emptyset}
(\prod_{j \in J} 2^{-\kappa_j p^{-1}})\times\\
\int_{ (c_{17} 2^{-\kappa} I^d)^J} (\prod_{j \in J} t_j^{-2 p^{-1} -1})
(\int_{(t B^d)^J} \int_{D_\xi^{l \chi_J} \cap D_{\kappa^0 +\kappa,n}^{d,m,D,U}}
|\Delta_\xi^{l \chi_J} f(x)|^p dx d\xi^J)^{1/p} dt^J \biggr)^p \le \\
\sum_{\substack{n \in \Z^d: \\ Q_{\kappa^0 +\kappa,n}^d \cap
G_{\kappa^0 +\kappa}^{d,m,U} \ne \emptyset}}
\biggl( c_{13} c_{18} \sum_{ J \subset \Nu_{1,d}^1: J \ne \emptyset}
(\prod_{j \in J} 2^{-\kappa_j p^{-1}})
\int_{ (c_{17} 2^{-\kappa} I^d)^J} (\prod_{j \in J} t_j^{-2 p^{-1} -1})\times \\
(\int_{(t B^d)^J} \int_{D_\xi^{l \chi_J} \cap D_{\kappa^0 +\kappa,n}^{d,m,D,U}}
|\Delta_\xi^{l \chi_J} f(x)|^p dx d\xi^J)^{1/p} dt^J \biggr)^p \le \\
c_{19} \sum_{\substack{n \in \Z^d: \\ Q_{\kappa^0 +\kappa,n}^d \cap
G_{\kappa^0 +\kappa}^{d,m,U} \ne \emptyset}}
\sum_{ J \subset \Nu_{1,d}^1: J \ne \emptyset}
(\prod_{j \in J} 2^{-\kappa_j})
\biggl(\int_{ (c_{17} 2^{-\kappa} I^d)^J}
(\prod_{j \in J} t_j^{-2 p^{-1} -1})\times\\
(\int_{(t B^d)^J} \int_{D_\xi^{l \chi_J} \cap D_{\kappa^0 +\kappa,n}^{d,m,D,U}}
|\Delta_\xi^{l \chi_J} f(x)|^p dx d\xi^J)^{1/p} dt^J \biggr)^p = \\
c_{19} \sum_{ J \subset \Nu_{1,d}^1: J \ne \emptyset}
(\prod_{j \in J} 2^{-\kappa_j})
\sum_{n \in \Z^d: Q_{\kappa^0 +\kappa,n}^d \cap
G_{\kappa^0 +\kappa}^{d,m,U} \ne \emptyset}
\biggl(\int_{ (c_{17} 2^{-\kappa} I^d)^J}
(\prod_{j \in J} t_j^{-2 p^{-1} -1})\times\\
(\int_{(t B^d)^J} \int_{D_\xi^{l \chi_J} \cap D_{\kappa^0 +\kappa,n}^{d,m,D,U}}
|\Delta_\xi^{l \chi_J} f(x)|^p dx d\xi^J)^{1/p} dt^J \biggr)^p.
\end{multline*}

Теперь фиксируем $ \epsilon \in \R_+^d $ так, что $ \alpha -p^{-1} \e -\epsilon
> 0. $ Тогда, используя неравенство Гёльдера, для $ J \subset \Nu_{1,d}^1:
J \ne \emptyset, n \in \Z^d: Q_{\kappa^0 +\kappa,n}^d \cap
G_{\kappa^0 +\kappa}^{d,m,U} \ne \emptyset, $ имеем
\begin{multline*} \tag{2.1.39}
\biggl(\int_{ (c_{17} 2^{-\kappa} I^d)^J}
(\prod_{j \in J} t_j^{-2 p^{-1} -1})
(\int_{(t B^d)^J} \int_{D_\xi^{l \chi_J} \cap D_{\kappa^0 +\kappa,n}^{d,m,D,U}}
|\Delta_\xi^{l \chi_J} f(x)|^p dx d\xi^J)^{1/p} dt^J \biggr)^p = \\
\biggl(\int_{ (c_{17} 2^{-\kappa} I^d)^J}
(\prod_{j \in J} (t_j^{\epsilon_j -1 /p^\prime} t_j^{-\epsilon_j -2 p^{-1} -1 /p}))\times\\
(\int_{(t B^d)^J} \int_{D_\xi^{l \chi_J} \cap D_{\kappa^0 +\kappa,n}^{d,m,D,U}}
|\Delta_\xi^{l \chi_J} f(x)|^p dx d\xi^J)^{1/p} dt^J \biggr)^p \le \\
\biggl(\int_{ (c_{17} 2^{-\kappa} I^d)^J}
\prod_{j \in J} t_j^{p^\prime \epsilon_j -1} dt^J \biggr)^{p /p^\prime}
\int_{ (c_{17} 2^{-\kappa} I^d)^J} (\prod_{j \in J} t_j^{-p \epsilon_j -3})\times\\
\int_{(t B^d)^J} \int_{ D_\xi^{l \chi_J} \cap D_{\kappa^0 +\kappa,n}^{d,m,D,U}}
|\Delta_\xi^{l \chi_J} f(x)|^p dx d\xi^J dt^J \le \\
c_{20} 2^{-(\kappa^J, \epsilon^J) p}
\int_{ (c_{17} 2^{-\kappa} I^d)^J} (\prod_{j \in J} t_j^{-p \epsilon_j -3})
\int_{(t B^d)^J} \int_{ D_\xi^{l \chi_J} \cap D_{\kappa^0 +\kappa,n}^{d,m,D,U}}
|\Delta_\xi^{l \chi_J} f(x)|^p dx d\xi^J dt^J.
\end{multline*}

Используя (2.1.39) для оценки правой части (2.1.38), для $ J \subset
\Nu_{1,d}^1: J \ne \emptyset, $ с учётом (2.1.30) получаем
\begin{multline*} \tag{2.1.40}
\sum_{n \in \Z^d: Q_{\kappa^0 +\kappa,n}^d \cap
G_{\kappa^0 +\kappa}^{d,m,U} \ne \emptyset}
\biggl(\int_{ (c_{17} 2^{-\kappa} I^d)^J}
(\prod_{j \in J} t_j^{-2 p^{-1} -1})\times\\
(\int_{(t B^d)^J} \int_{D_\xi^{l \chi_J} \cap D_{\kappa^0 +\kappa,n}^{d,m,D,U}}
|\Delta_\xi^{l \chi_J} f(x)|^p dx d\xi^J)^{1/p} dt^J \biggr)^p \le \\
\sum_{n \in \Z^d: Q_{\kappa^0 +\kappa,n}^d \cap
G_{\kappa^0 +\kappa}^{d,m,U} \ne \emptyset}
c_{20} 2^{-(\kappa^J, \epsilon^J) p}
\int_{ (c_{17} 2^{-\kappa} I^d)^J} (\prod_{j \in J} t_j^{-p \epsilon_j -3})\times\\
\int_{(t B^d)^J} \int_{ D_\xi^{l \chi_J} \cap D_{\kappa^0 +\kappa,n}^{d,m,D,U}}
|\Delta_\xi^{l \chi_J} f(x)|^p dx d\xi^J dt^J = \\
c_{20} 2^{-(\kappa^J, \epsilon^J) p}
\sum_{n \in \Z^d: Q_{\kappa^0 +\kappa,n}^d \cap
G_{\kappa^0 +\kappa}^{d,m,U} \ne \emptyset}
\int_{ (c_{17} 2^{-\kappa} I^d)^J} (\prod_{j \in J} t_j^{-p \epsilon_j -3})\times\\
\int_{(t B^d)^J} \int_{ D_\xi^{l \chi_J} \cap D_{\kappa^0 +\kappa,n}^{d,m,D,U}}
|\Delta_\xi^{l \chi_J} f(x)|^p dx d\xi^J dt^J = \\
c_{20} 2^{-(\kappa^J, \epsilon^J) p}
\int_{ (c_{17} 2^{-\kappa} I^d)^J} (\prod_{j \in J} t_j^{-p \epsilon_j -3})\times\\
\sum_{n \in \Z^d: Q_{\kappa^0 +\kappa,n}^d \cap
G_{\kappa^0 +\kappa}^{d,m,U} \ne \emptyset}
\int_{(t B^d)^J} \int_{ D_\xi^{l \chi_J} \cap D_{\kappa^0 +\kappa,n}^{d,m,D,U}}
|\Delta_\xi^{l \chi_J} f(x)|^p dx d\xi^J dt^J = \\
c_{20} 2^{-(\kappa^J, \epsilon^J) p}
\int_{ (c_{17} 2^{-\kappa} I^d)^J} (\prod_{j \in J} t_j^{-p \epsilon_j -3})\times\\
\int_{(t B^d)^J} \sum_{n \in \Z^d: Q_{\kappa^0 +\kappa,n}^d \cap
G_{\kappa^0 +\kappa}^{d,m,U} \ne \emptyset}
\int_{ D_\xi^{l \chi_J}} \chi_{D_{\kappa^0 +\kappa,n}^{d,m,D,U}}(x)
|\Delta_\xi^{l \chi_J} f(x)|^p dx d\xi^J dt^J = \\
c_{20} 2^{-(\kappa^J, \epsilon^J) p}
\int_{ (c_{17} 2^{-\kappa} I^d)^J} (\prod_{j \in J} t_j^{-p \epsilon_j -3})\times\\
\int_{(t B^d)^J} \int_{ D_\xi^{l \chi_J}}
(\sum_{n \in \Z^d: Q_{\kappa^0 +\kappa,n}^d \cap
G_{\kappa^0 +\kappa}^{d,m,U} \ne \emptyset}
\chi_{D_{\kappa^0 +\kappa,n}^{d,m,D,U}}(x))
|\Delta_\xi^{l \chi_J} f(x)|^p dx d\xi^J dt^J \le \\
c_{20} 2^{-(\kappa^J, \epsilon^J) p}
\int_{ (c_{17} 2^{-\kappa} I^d)^J} (\prod_{j \in J} t_j^{-p \epsilon_j -3})\times\\
\int_{(t B^d)^J} \int_{ D_\xi^{l \chi_J}}
c_{12} |\Delta_\xi^{l \chi_J} f(x)|^p dx d\xi^J dt^J \le \\
c_{21} 2^{-(\kappa^J, \epsilon^J) p}
\int_{ (c_{17} 2^{-\kappa} I^d)^J} (\prod_{j \in J} t_j^{-p \epsilon_j -2})
(\Omega^{\prime l \chi_J}(f, t^J)_{L_p(D)})^p dt^J \le \\
c_{22}(d,\alpha,p,m,D,U,f) 2^{-(\kappa^J, \epsilon^J) p}
\int_{ (c_{17} 2^{-\kappa} I^d)^J} (\prod_{j \in J} t_j^{-p \epsilon_j -2})
(t^J)^{p \alpha^J} dt^J = \\
c_{22} 2^{-(\kappa^J, \epsilon^J) p}
\int_{ (c_{17} 2^{-\kappa} I^d)^J} (\prod_{j \in J}
t_j^{p(\alpha_j -\epsilon_j -p^{-1}) -1}) dt^J \le \\
c_{23} 2^{-(\kappa^J, \epsilon^J) p}
2^{-p(\kappa^J, \alpha^J -\epsilon^J -p^{-1} \e^J)} =
c_{23} 2^{-p(\kappa^J, \alpha^J -p^{-1} \e^J)}.
\end{multline*}

Подстановка (2.1.40) в (2.1.38) и применение (1.1.2) при $ a = 1 /p $ даёт
оценку
\begin{multline*} \tag{2.1.41}
\int_{G_{\kappa^0 +\kappa}^{d,m,U}}
\biggl| \sum_{\nu \in N_{\kappa^0 +\kappa}^{d,m,U}}
((S_{\kappa^0 +\kappa, \nu_{\kappa^0 +\kappa}^{d,m,D,U}(\nu)}^{d,l -\e} f)(x) -
(\mathtt R_{\kappa^0 +\kappa, \nu_{\kappa^0 +\kappa}^{d,m,D,U}(\nu)}^{d,l -\e} f)(x))
g_{\kappa^0 +\kappa, \nu}^{d,m}(x) \biggr|^p dx \le \\
c_{19} \sum_{ J \subset \Nu_{1,d}^1: J \ne \emptyset}
(\prod_{j \in J} 2^{-\kappa_j}) c_{23} 2^{-p(\kappa^J, \alpha^J -p^{-1} \e^J)} = \\
c_{24} \sum_{ J \subset \Nu_{1,d}^1: J \ne \emptyset}
2^{-p(\kappa^J, \alpha^J)} \le
c_{24} (\sum_{ J \subset \Nu_{1,d}^1: J \ne \emptyset}
2^{-(\kappa^J, \alpha^J)})^p.
\end{multline*}
Соединяя (2.1.34) с (2.1.41), приходим к соотношению
\begin{multline*}
\| (E_{\kappa^0 +\kappa}^{d,l -\e,m,D,U,\nu_{\kappa^0 +\kappa}} f) \mid_U
-(R_{\kappa^0 +\kappa}^{d,l -\e,m,D,U,\nu_{\kappa^0 +\kappa}} f) \mid_U \|_{L_p(U)} \le \\
c_{25} \sum_{ J \subset \Nu_{1,d}^1: J \ne \emptyset} 2^{-(\kappa^J, \alpha^J)}
\to 0 \text{ при } \mn(\kappa) \to \infty. \square
\end{multline*}

Для дальнейшего изложения понадобятся ещё некоторые факты.
Как известно, имеет место следующее утверждение.

Лемма 2.1.8

Пусть $ d \in \N, \lambda \in \Z_+^d, D $ --- область в $ \R^d $ и
функция $ f \in C^\infty(D), $ а $ g \in L_1(D), $ причём
для каждого $ \mu \in \Z_+^d(\lambda) $ (см. (1.1.1)) обобщённая производная
$ \D^\mu g \in L_1(D). $ Тогда в пространстве обобщённых функций в
области $ D $ имеет место равенство
\begin{equation*} \tag{2.1.42}
\D^\lambda (fg) = \sum_{ \mu \in \Z_+^d(\lambda)} C_\lambda^\mu
\D^{\lambda -\mu} f \D^\mu g.
\end{equation*}

Для формулировки следующего утверждения введём обозначение.

Пусть $ d \in \N, D -- $ ограниченная область в $ \R^d $ и $ \kappa^0 \in
\Z_+^d $ таковы, что выполняется (2.1.1) с $ \kappa^0 $ вместо $ \kappa, $ а
$ U \subset D $ -- открытое подмножество $ D. $ Имея в виду замечание после
(2.1.3), при $ m \in \N^d $ рассмотрим некоторое семейство отображений
$ \Nu = \{ \nu_{\kappa^0 +\kappa}^{d,m,D,U}, \kappa \in \Z_+^d \} $
вида (2.1.2) с $ \kappa^0 +\kappa $ вместо $ \kappa $ и при $ \kappa, l \in
\Z_+^d, $ исходя из (2.1.22), определим линейный оператор
$ \mathcal R_{\kappa^0, \kappa}^{d,l,m,D,U,\Nu}: C(D) \mapsto
\mathcal P_{\kappa^0 +\kappa}^{d,l,m,U}, $ полагая
\begin{equation*} \tag{2.1.43}
\mathcal R_{\kappa^0, \kappa}^{d,l,m,D,U,\Nu} = \sum_{\epsilon \in \Upsilon^d:
\s(\epsilon) \subset \s(\kappa)} (-\e)^\epsilon
H_{\kappa^0 +\kappa, \kappa^0 +\kappa -\epsilon}^{d,l,m,U}
R_{\kappa^0 +\kappa -\epsilon}^{d,l,m,D,U,\nu_{\kappa^0 +\kappa -\epsilon}}.
\end{equation*}

При этом, принимая во внимание (2.1.43), (2.1.22), (1.4.17), (1.4.18), (1.4.19),
имеем
\begin{multline*} \tag{2.1.44}
\mathcal R_{\kappa^0, \kappa}^{d,l,m,D,U,\Nu} f =\\
\sum_{\epsilon \in \Upsilon^d: \s(\epsilon) \subset \s(\kappa)} (-\e)^\epsilon
\sum_{\nu \in N_{\kappa^0 +\kappa}^{d,m,U}}
\biggl(\sum_{\m^{\epsilon} \in \M_{\epsilon}^m(\nu)} A_{\m^{\epsilon}}^m
\mathtt R_{\kappa^0 +\kappa -\epsilon, \nu_{\kappa^0 +\kappa -\epsilon}^{d,m,D,U}
(\n_{\epsilon}(\nu,\m^{\epsilon}))}^{d,l} f\biggr)
g_{\kappa^0 +\kappa,\nu}^{d,m} = \\
\sum_{\epsilon \in \Upsilon^d: \s(\epsilon) \subset \s(\kappa)}
\sum_{\nu \in N_{\kappa^0 +\kappa}^{d,m,U}} \biggl((-\e)^\epsilon
\sum_{\m^{\epsilon} \in \M_{\epsilon}^m(\nu)} A_{\m^{\epsilon}}^m
\mathtt R_{\kappa^0 +\kappa -\epsilon, \nu_{\kappa^0 +\kappa -\epsilon}^{d,m,D,U}
(\n_{\epsilon}(\nu,\m^{\epsilon}))}^{d,l} f \biggr)
g_{\kappa^0 +\kappa,\nu}^{d,m} = \\
\sum_{\nu \in N_{\kappa^0 +\kappa}^{d,m,U}}
\sum_{\epsilon \in \Upsilon^d: \s(\epsilon) \subset \s(\kappa)}
\biggl((-\e)^\epsilon \sum_{\m^{\epsilon} \in \M_{\epsilon}^m(\nu)} A_{\m^{\epsilon}}^m
\mathtt R_{\kappa^0 +\kappa -\epsilon, \nu_{\kappa^0 +\kappa -\epsilon}^{d,m,D,U}
(\n_{\epsilon}(\nu,\m^{\epsilon}))}^{d,l} f\biggr)
g_{\kappa^0 +\kappa,\nu}^{d,m} = \\
\sum_{\nu \in N_{\kappa^0 +\kappa}^{d,m,U}}
\biggl(\sum_{\epsilon \in \Upsilon^d: \s(\epsilon) \subset \s(\kappa)} (-\e)^\epsilon
\sum_{\m^{\epsilon} \in \M_{\epsilon}^m(\nu)} A_{\m^{\epsilon}}^m
\mathtt R_{\kappa^0 +\kappa -\epsilon, \nu_{\kappa^0 +\kappa -\epsilon}^{d,m,D,U}
(\n_{\epsilon}(\nu,\m^{\epsilon}))}^{d,l} f \biggr)
g_{\kappa^0 +\kappa,\nu}^{d,m} = \\
\sum_{ \nu \in N_{\kappa^0 +\kappa}^{d,m,U}}
(V_{\kappa^0, \kappa,\nu}^{d,l,m,D,U,\Nu} f) g_{\kappa^0 +\kappa, \nu}^{d,m}, f \in C(D),
\end{multline*}
где $ V_{\kappa^0, \kappa,\nu}^{d,l,m,D,U,\Nu}: C(D) \mapsto
\mathcal P^{d,l} $ -- линейный оператор, определяемый при $ d \in \N, D \subset
\R^d, U \subset D, l \in \Z_+^d, m \in \N^d, \kappa \in \Z_+^d, \nu \in
N_{\kappa^0 +\kappa}^{d,m,U}, \Nu = \{\nu_{\kappa^0 +\kappa}, \kappa \in \Z_+^d \}, $
удовлетворяющих указанным выше условиям, равенством
\begin{equation*} \tag{2.1.45}
V_{\kappa^0, \kappa,\nu}^{d,l,m,D,U,\Nu} f =
\sum_{\epsilon \in \Upsilon^d: \s(\epsilon) \subset \s(\kappa)} (-\e)^\epsilon
\sum_{\m^{\epsilon} \in \M_{\epsilon}^m(\nu)} A_{\m^{\epsilon}}^m
\mathtt R_{\kappa^0 +\kappa -\epsilon, \nu_{\kappa^0 +\kappa -\epsilon}^{d,m,D,U}
(\n_{\epsilon}(\nu,\m^{\epsilon}))}^{d,l} f, f \in C(D).
\end{equation*}

При $ d \in \N $ обозначим через
$$
\Sigma^d = \{ \sigma \in \Z^d: \sigma_j \in \{-1,1\}, j =1,\ldots,d\}.
$$

Приведём уточнённый вариант теоремы, установленной в [13], которым
мы воспользуемся в следующем пункте.

Теорема 2.1.9

При $ d \in \N $ пусть $ D $ -- ограниченная область в $ \R^d, $
обладающая тем свойством, что существует система открытых подмножеств
$ \{ U_i \subset D, i =1,\ldots,\mathcal I\}, $ для которых существует число
$ \bm \delta > 0 $ такое, что набор открытых множеств
$$
V_i = \{ x \in U_i: \inf_{ y \in D \setminus U_i} \| x -y\| >
\bm \delta\}, i =1,\ldots,\mathcal I,
$$
образует покрытие области $ D, $ и  таких, что при $ i
=1,\ldots,\mathcal I $ существуют $ \delta_i
\in \R_+^d $ и $ \sigma^i \in \Sigma^d, $ для которых
имеет место включение $ (U_i +\sigma^i \delta^i I^d) \subset D. $
Пусть ещё $ \alpha \in \R_+^d, 1 \le p < \infty, 1 \le \theta \le \infty. $
Тогда существует непрерывное линейное отображение
$ \mathcal E^{d,\alpha,p,\theta,D}: (S_{p,\theta}^\alpha B)^\prime(D)
\mapsto L_p(\R^d), $ обладающее следующими свойствами:

1) для $ f \in (S_{p,\theta}^\alpha B)^\prime(D) $ выполняется равенство
\begin{equation*} \tag{2.1.46}
f = (\mathcal E^{d,\alpha,p,\theta,D} f) \mid_{D},
\end{equation*}

2) существует константа $ c_{26}(d,\alpha,p,\theta,D) >0 $ такая, что
при $ l = l(\alpha), \lambda \in \Z_+^d: \lambda < \alpha, $ для
$ f \in (S_{p,\theta}^\alpha B)^\prime(D) $ и любого множества
$ J \subset \{1,\ldots,d\} $ верно неравенство
\begin{multline*} \tag{2.1.47}
\biggl(\int_{(\R_+^d)^J} (t^J)^{-\e^J -\theta (\alpha^J -\lambda^J)}
(\Omega^{(l -\lambda) \chi_J}(\D^\lambda (\mathcal E^{d,\alpha,p,\theta,D} f),
t^J)_{L_p(\R^d)})^\theta dt^J\biggr)^{1/\theta} \le \\
c_{26} \| f \|_{(S_{p,\theta}^\alpha B)^\prime(D)}, \ \theta \ne \infty; \\
\sup_{t^J \in (\R_+^d)^J} (t^J)^{-(\alpha^J -\lambda^J)}
\Omega^{(l -\lambda) \chi_J}(\D^\lambda (\mathcal E^{d,\alpha,p,\theta,D} f),
t^J)_{L_p(\R^d)} \le \\
c_{26} \| f \|_{(S_{p,\theta}^\alpha B)^\prime(D)}, \ \theta = \infty,
\end{multline*}

3) при соблюдении (2.1.9) для $ f \in (S_{p,\theta}^\alpha B)^\prime(D) $
имеет место включение
\begin{equation*} \tag{2.1.48}
(\mathcal E^{d,\alpha,p,\theta,D} f) \in C_0(\R^d).
\end{equation*}

Доказательство.

В условиях теоремы 2.1.9 оператор $ \mathcal E^{d,\alpha,p,\theta,D}, $
построенный при доказательстве теоремы 2.3.1 из [13],
обладает свойствами, описанными в пп. 1), 2). Для проверки соблюдения
включения (2.1.48) при условии (2.1.9) достаточно принять во внимание (2.1.47),
(1.1.7), (1.1.8), предложение 2.1.3 и тот факт, что в случае ограниченности
области $ D $ для $ f \in (S_{p,\theta}^\alpha B)^\prime(D) $ функция
$ \mathcal E^{d,\alpha,p,\theta,D} f $ имеет компактный носитель (см.
доказательство теорем 2.3.1, 2.2.7 и предложения 2.2.4 из [13]).
$ \square $

Замечание.

Если область $ D \subset \R^d $ удовлетворяет условиям теоремы 2.1.9, то
при любом $ \sigma \in \Sigma^d $ область $ \sigma D $ удовлетворяет
условиям теоремы 2.1.9.
\bigskip

2.2. В этом пункте будет установлена оценка сверху изучаемой величины. Для
этого понадобится лемма 2.2.2, которая использует лемму 2.2.1.

Лемма 2.2.1

При $ d \in \N $ пусть область $ D \subset \R^d $ и её открытое подмножество
$ U \subset D $ таковы, что существует $ \delta \in \R_+^d, $ для которого имеет
место включение $ (U +\delta I^d) \subset D. $ Тогда при любом $ m \in \N^d $ 
существуют константы $ \Kappa^0 = \Kappa^0(d,m,D,U) \in \Z_+^d, \Gamma^0 =
\Gamma^0(d,m,D,U) \in \R_+^d, $ для которых существуют семейства отображений
$$
\Nu = \Nu^{d,m,D,U} = \{ \nu_{\Kappa^0 +\kappa}^{d,m,D,U}: N_{\Kappa^0 +\kappa}^{d,m,U}
\mapsto \Z^d, \kappa \in \Z_+^d\}, \\
\{n_{\Kappa^0 +\kappa}^{d,m,D,U}: N_{\Kappa^0 +\kappa}^{d,m,U} \mapsto \Z^d,
\kappa \in \Z_+^d\}, 
$$
обладающие следующими свойствами:

1) при $ \kappa \in \Z_+^d $ для каждого $ \nu \in N_{\Kappa^0 +\kappa}^{d,m,U} $
справедливо включение
\begin{equation*} \tag{2.2.1}
(Q_{\Kappa^0 +\kappa,\nu_{\Kappa^0 +\kappa}^{d,m,D,U}(\nu)}^d \cup
Q_{\Kappa^0 +\kappa,n_{\Kappa^0 +\kappa}^{d,m,D,U}(\nu)}^d) \subset D \cap
(2^{-\Kappa^0 -\kappa} \nu +\Gamma^0 2^{-\Kappa^0 -\kappa} B^d);
\end{equation*}

2) при $ \kappa \in \Z_+^d, \nu \in N_{\Kappa^0 +\kappa}^{d,m,U},
\epsilon \in \Upsilon^d: \s(\epsilon) \subset \s(\kappa), \m^{\epsilon} \in
\M_{\epsilon}^m(\nu) $ для
$$
\mathcal D_{\Kappa^0 +\kappa,\nu,\epsilon,\m^{\epsilon}}^{d,m,D,U} =
\bm x_{\Kappa^0 +\kappa,\nu,\epsilon,\m^{\epsilon}}^{d,m,D,U} +
\bm \delta_{\Kappa^0 +\kappa,\nu,\epsilon,\m^{\epsilon}}^{d,m,D,U} I^d,
$$
где точка $ \bm x_{\Kappa^0 +\kappa,\nu,\epsilon,\m^{\epsilon}}^{d,m,D,U} \in \R^d $
и вектор $ \bm \delta_{\Kappa^0 +\kappa,\nu,\epsilon,\m^{\epsilon}}^{d,m,D,U} \in \R_+^d $
определяются равенствами
\begin{multline*}
(\bm x_{\Kappa^0 +\kappa,\nu,\epsilon,\m^{\epsilon}}^{d,m,D,U})_j =\\
\min(2^{-\Kappa^0_j -\kappa_j} (n_{\Kappa^0 +\kappa}^{d,m,D,U}(\nu))_j,
2^{-\Kappa^0_j -\kappa_j +\epsilon_j} (\nu_{\Kappa^0 +\kappa -\epsilon}^{d,m,D,U}
(\n_{\epsilon}(\nu,\m^{\epsilon})))_j),  j \in \Nu_{1,d}^1; \\
(\bm \delta_{\Kappa^0 +\kappa,\nu,\epsilon,\m^{\epsilon}}^{d,m,D,U})_j =
\max(2^{-\Kappa^0_j -\kappa_j} (n_{\Kappa^0 +\kappa}^{d,m,D,U}(\nu))_j
+2^{-\Kappa^0_j -\kappa_j}, \\
2^{-\Kappa^0_j -\kappa_j +\epsilon_j} (\nu_{\Kappa^0 +\kappa -\epsilon}^{d,m,D,U}
(\n_{\epsilon}(\nu,\m^{\epsilon})))_j +2^{-\Kappa^0_j -\kappa_j +\epsilon_j})
-(\bm x_{\Kappa^0 +\kappa,\nu,\epsilon,\m^{\epsilon}}^{d,m,D,U})_j,
j \in \Nu_{1,d}^1,
\end{multline*}
и
\begin{multline*}
\mathcal D_{\Kappa^0 +\kappa,\nu}^{d,m,D,U} = \{\bm x \in \R^d: 
\min_{\bm \epsilon \in \Upsilon^d: \s(\bm \epsilon) \subset \s(\kappa), \m^{\bm \epsilon} 
\in \M_{\bm \epsilon}^m(\nu)}
(\bm x_{\Kappa^0 +\kappa,\nu,\bm \epsilon,\m^{\bm \epsilon}}^{d,m,D,U})_j < \bm x_j < \\
\max_{\bm \epsilon \in \Upsilon^d: \s(\bm \epsilon) \subset 
\s(\kappa), \m^{\bm \epsilon} \in \M_{\bm \epsilon}^m(\nu)} 
(\bm x_{\Kappa^0 +\kappa,\nu,\bm \epsilon,\m^{\bm \epsilon}}^{d,m,D,U})_j +
(\bm \delta_{\Kappa^0 +\kappa,\nu,\bm \epsilon,\m^{\bm \epsilon}}^{d,m,D,U})_j, j =1,\ldots,d\},
\end{multline*}
соблюдается включение
\begin{equation*} \tag{2.2.2}
\mathcal D_{\Kappa^0 +\kappa,\nu,\epsilon,\m^{\epsilon}}^{d,m,D,U} \subset 
\mathcal D_{\Kappa^0 +\kappa,\nu}^{d,m,D,U} \subset D;
\end{equation*}

3) при $ \kappa \in \Z_+^d, \nu \in N_{\Kappa^0 +\kappa}^{d,m,U} $ для любых
$ \epsilon \in \Upsilon^d: \s(\epsilon) \subset \s(\kappa), $ и $ \m^{\epsilon}
\in \M_{\epsilon}^m(\nu) $ при $ j \in \Nu_{1,d}^1 \setminus \s(\epsilon) $
выполняется равенство
\begin{equation*} \tag{2.2.3}
(\nu_{\Kappa^0 +\kappa -\epsilon}^{d,m,D,U}
(\n_{\epsilon}(\nu,\m^{\epsilon})))_j =
(\nu_{\Kappa^0 +\kappa}^{d,m,D,U}(\nu))_j.
\end{equation*}

Доказательство.

В условиях леммы при $ m \in \N^d $ и $ \Kappa^0, \kappa \in \Z_+^d $ для $ \nu
\in N_{\Kappa^0 +\kappa}^{d,m,U}, $ учитывая открытость $ U $ и (1.3.4),
выберем $ x \in U \cap (2^{-\Kappa^0 -\kappa} \nu +2^{-\Kappa^0 -\kappa} (m +\e) I^d). $
Тогда $ x \in U $ и при $ j = 1,\ldots,d $ выполняются соотношения
\begin{equation*} \tag{2.2.4}
2^{-\Kappa^0_j -\kappa_j} \nu_j < x_j < 2^{-\Kappa^0_j -\kappa_j} \nu_j +
2^{-\Kappa^0_j -\kappa_j} (m_j +1),
\end{equation*}
а для $ y \in (2^{-\Kappa^0 -\kappa} (\nu +m +\e) +2^{-\Kappa^0 -\kappa} I^d) $
при $ j = 1,\ldots,d $ соблюдаются неравенства
$$
2^{-\Kappa^0_j -\kappa_j} (\nu_j +m_j +1) < y_j < 2^{-\Kappa^0_j -\kappa_j}
(\nu_j +m_j +1) +2^{-\Kappa^0_j -\kappa_j},
$$
и, значит,
\begin{multline*} \tag{2.2.5}
x_j < 2^{-\Kappa^0_j -\kappa_j} (\nu_j +m_j +1) < y_j <
2^{-\Kappa^0_j -\kappa_j} \nu_j +2^{-\Kappa^0_j -\kappa_j} (m_j +2) < \\
x_j+2^{-\Kappa^0_j -\kappa_j} (m_j +2),
\end{multline*}
т.е.
\begin{equation*} \tag{2.2.6}
Q_{\Kappa^0 +\kappa, \nu +m +\e}^d \subset x +2^{-\Kappa^0 -\kappa} (m +2 \e) I^d
\subset (U +2^{-\Kappa^0 -\kappa} (m +2 \e) I^d).
\end{equation*}

Теперь в условиях леммы при $ m \in \N^d $ фиксируем $ \Kappa^0 =
\Kappa^0(d,m,D,U) \in \Z_+^d $ такое, что $ 2^{-\Kappa^0} (m +2 \e) < \delta, $
и определим при $ \kappa \in \Z_+^d $ отображения
$ \nu_{\Kappa^0 +\kappa}^{d,m,D,U}: N_{\Kappa^0 +\kappa}^{d,m,U} \mapsto \Z^d,
n_{\Kappa^0 +\kappa}^{d,m,D,U}: N_{\Kappa^0 +\kappa}^{d,m,U} \mapsto \Z^d, $
полагая для $ \nu \in N_{\Kappa^0 +\kappa}^{d,m,U} $ значение
$$
\nu_{\Kappa^0 +\kappa}^{d,m,D,U}(\nu) = n_{\Kappa^0 +\kappa}^{d,m,D,U}(\nu) =
\nu +m +\e.
$$

Проверим, что для рассматриваемых объектов соблюдаются условия пп. 1) -- 3). 
Для этого при $ \kappa \in \Z_+^d $ для $ \nu \in N_{\Kappa^0 +\kappa}^{d,m,U} $ 
в силу (2.2.6) имеем
\begin{multline*}
Q_{\Kappa^0 +\kappa,\nu_{\Kappa^0 +\kappa}^{d,m,D,U}(\nu)}^d =
Q_{\Kappa^0 +\kappa,n_{\Kappa^0 +\kappa}^{d,m,D,U}(\nu)}^d =
Q_{\Kappa^0 +\kappa, \nu +m +\e}^d \subset\\
 (U +2^{-\Kappa^0 -\kappa} (m +2 \e) I^d)
\subset (U +2^{-\Kappa^0} (m +2 \e) I^d) \subset (U +\delta I^d) \subset D
\end{multline*}
и
\begin{multline*}
Q_{\Kappa^0 +\kappa,\nu_{\Kappa^0 +\kappa}^{d,m,D,U}(\nu)}^d =
Q_{\Kappa^0 +\kappa,n_{\Kappa^0 +\kappa}^{d,m,D,U}(\nu)}^d =
Q_{\Kappa^0 +\kappa, \nu +m +\e}^d =\\
(2^{-\Kappa^0 -\kappa} (\nu +m +\e) +2^{-\Kappa^0 -\kappa} I^d) \subset
(2^{-\Kappa^0 -\kappa} \nu +2^{-\Kappa^0 -\kappa} (m +\e +I^d)) \subset\\
(2^{-\Kappa^0 -\kappa} \nu +2^{-\Kappa^0 -\kappa} (m +2 \e) \overline I^d) \subset
(2^{-\Kappa^0 -\kappa} \nu +2^{-\Kappa^0 -\kappa} (m +2 \e) B^d),
\end{multline*}
т.е. выполняется (2.2.1) с $ \Gamma^0 = m +2 \e. $

Далее, при $ \kappa \in \Z_+^d, \nu \in N_{\Kappa^0 +\kappa}^{d,m,U}, $ выбирая 
$ x \in (U \cap (2^{-\Kappa^0 -\kappa} \nu +2^{-\Kappa^0 -\kappa} (m +\e) I^d)), $
ввиду (2.2.5), (2.2.4) для $ j =1.\ldots,d $ имеем
\begin{multline*} \tag{2.2.7}
x_j < 2^{-\Kappa^0_j -\kappa_j} (\nu_{\Kappa^0 +\kappa}^{d,m,D,U}(\nu))_j <
2^{-\Kappa^0_j -\kappa_j} (\nu_{\Kappa^0 +\kappa}^{d,m,D,U}(\nu))_j
+2^{-\Kappa^0_j -\kappa_j} =\\
2^{-\Kappa^0_j -\kappa_j} \nu_j + 2^{-\Kappa^0_j -\kappa_j} (m_j +2) <
x_j +2^{-\Kappa^0_j -\kappa_j} (m_j +2) < \\
x_j +2^{-\Kappa^0_j} (m_j +2) < x_j +\delta_j.
\end{multline*}
Принимая во внимание, что в описанных условиях при $ \epsilon \in \Upsilon^d: 
\s(\epsilon) \subset \s(\kappa), \m^{\epsilon} \in \M_{\epsilon}^m(\nu)$, имеет 
место включение
$$
x \in (U \cap (2^{-\Kappa^0 -\kappa +\epsilon} \n_{\epsilon}(\nu,\m^{\epsilon})
+2^{-\Kappa^0 -\kappa +\epsilon} (m +\e) I^d)), \text{ (см. вывод (1.4.14))},
$$
заключаем, что для $ j = 1,\ldots,d $ соблюдается (2.2.7) при
$ \Kappa^0 +\kappa -\epsilon $ вместо $ \Kappa^0 +\kappa $ и
$ \n_{\epsilon}(\nu,\m^{\epsilon}) \in N_{\Kappa^0 +\kappa -\epsilon}^{d,m,U} $
вместо $ \nu, $ т.е.
\begin{multline*} \tag{2.2.8}
x_j < 2^{-\Kappa^0_j -\kappa_j +\epsilon_j}
(\nu_{\Kappa^0 +\kappa -\epsilon}^{d,m,D,U}(\n_{\epsilon}(\nu,\m^{\epsilon})))_j <\\
2^{-\Kappa^0_j -\kappa_j +\epsilon_j}
(\nu_{\Kappa^0 +\kappa -\epsilon}^{d,m,D,U}(\n_{\epsilon}(\nu,\m^{\epsilon})))_j +
2^{-\Kappa^0_j -\kappa_j +\epsilon_j} =\\
2^{-\Kappa^0_j -\kappa_j +\epsilon_j} (\n_{\epsilon}(\nu,\m^{\epsilon}))_j +
2^{-\Kappa^0_j -\kappa_j +\epsilon_j} (m_j +2) <\\
x_j +2^{-\Kappa^0_j -\kappa_j +\epsilon_j} (m_j +2) \le
x_j +2^{-\Kappa^0_j} (m_j +2) < x_j +\delta_j.
\end{multline*}
Из (2.2.7) и (2.2.8) с учётом определений получаем, что при
$ \epsilon \in \Upsilon^d: \s(\epsilon) \subset \s(\kappa), \m^{\epsilon} \in 
\M_{\epsilon}^m(\nu), j =1,\ldots,d $ справедливо соотношение
\begin{multline*}
x_j < \min(2^{-\Kappa^0_j -\kappa_j} (n_{\Kappa^0 +\kappa}^{d,m,D,U}(\nu))_j,
2^{-\Kappa^0_j -\kappa_j +\epsilon_j} (\nu_{\Kappa^0 +\kappa -\epsilon}^{d,m,D,U}
(\n_{\epsilon}(\nu,\m^{\epsilon})))_j) = \\
(\bm x_{\Kappa^0 +\kappa,\nu,\epsilon,\m^{\epsilon}}^{d,m,D,U})_j <
(\bm x_{\Kappa^0 +\kappa,\nu,\epsilon,\m^{\epsilon}}^{d,m,D,U})_j +
(\bm \delta_{\Kappa^0 +\kappa,\nu,\epsilon,\m^{\epsilon}}^{d,m,D,U})_j =\\
\max(2^{-\Kappa^0_j -\kappa_j} (n_{\Kappa^0 +\kappa}^{d,m,D,U}(\nu))_j
+2^{-\Kappa^0_j -\kappa_j}, \\
2^{-\Kappa^0_j -\kappa_j +\epsilon_j} (\nu_{\Kappa^0 +\kappa -\epsilon}^{d,m,D,U}
(\n_{\epsilon}(\nu,\m^{\epsilon})))_j +2^{-\Kappa^0_j -\kappa_j +\epsilon_j}) <
x_j +\delta_j,
\end{multline*}
и
\begin{multline*}
x_j < \min_{\epsilon \in \Upsilon^d: \s(\epsilon) \subset \s(\kappa), 
\m^{\epsilon} \in \M_{\epsilon}^m(\nu)}
(\bm x_{\Kappa^0 +\kappa,\nu,\epsilon,\m^{\epsilon}}^{d,m,D,U})_j <\\
\max_{\epsilon \in \Upsilon^d: \s(\epsilon) \subset \s(\kappa),  
\m^{\epsilon} \in \M_{\epsilon}^m(\nu)} 
(\bm x_{\Kappa^0 +\kappa,\nu,\epsilon,\m^{\epsilon}}^{d,m,D,U})_j +
(\bm \delta_{\Kappa^0 +\kappa,\nu,\epsilon,\m^{\epsilon}}^{d,m,D,U})_j 
< x_j +\delta_j,
\end{multline*}
а, значит,
\begin{multline*}
\mathcal D_{\Kappa^0 +\kappa,\nu,\epsilon,\m^{\epsilon}}^{d,m,D,U} =
\bm x_{\Kappa^0 +\kappa,\nu,\epsilon,\m^{\epsilon}}^{d,m,D,U} +
\bm \delta_{\Kappa^0 +\kappa,\nu,\epsilon,\m^{\epsilon}}^{d,m,D,U} I^d
\subset \mathcal D_{\Kappa^0 +\kappa,\nu}^{d,m,D,U} \subset x +\delta I^d 
\subset (U +\delta I^d) \subset D,
\end{multline*}
что приводит к (2.2.2).

Наконец, при $ \kappa \in \Z_+^d, \nu \in N_{\Kappa^0 +\kappa}^{d,m,U} $ для
любых $ \epsilon \in \Upsilon^d: \s(\epsilon) \subset \s(\kappa), $ и
$ \m^{\epsilon} \in \M_{\epsilon}^m(\nu) $ при $ j \in \Nu_{1,d}^1 \setminus \s(\epsilon) $
в силу определений и (1.4.12) имеем
\begin{multline*}
(\nu_{\Kappa^0 +\kappa -\epsilon}^{d,m,D,U}
(\n_{\epsilon}(\nu,\m^{\epsilon})))_j = (\n_{\epsilon}(\nu,\m^{\epsilon}) +m +\e)_j =
(\n_{\epsilon}(\nu,\m^{\epsilon}))_j  +m_j +1 = \\
\nu_j +m_j +1 = (\nu +m +\e)_j =
(\nu_{\Kappa^0 +\kappa}^{d,m,D,U}(\nu))_j,
\end{multline*}
что совпадает с (2.2.3). $ \square $

Замечание.

В условиях леммы 2.2.1 при $ m \in \N^d $ существует константа 
$ \Gamma^1(d,m,D,U) \in \R_+^d $ такая, что при $ \kappa \in \Z_+^d, \nu \in 
N_{\Kappa^0 +\kappa}^{d,m,U}, \epsilon \in \Upsilon^d: \s(\epsilon) \subset 
\s(\kappa), \m^{\epsilon} \in \M_{\epsilon}^m(\nu) $ имеют место следующие 
соотношения:

\begin{equation*} \tag{2.2.9}
2^{-\Kappa^0 -\kappa} \le \bm \delta_{\Kappa^0 +\kappa,\nu,\epsilon,\m^{\epsilon}}^{d,m,D,U}
\le \Gamma^1 2^{-\Kappa^0 -\kappa},
\end{equation*}
\begin{equation*} \tag{2.2.10}
(Q_{\Kappa^0 +\kappa, n_{\Kappa^0 +\kappa}^{d,m,D,U}(\nu)}^d \cup
Q_{\Kappa^0 +\kappa -\epsilon, \nu_{\Kappa^0 +\kappa -\epsilon}^{d,m,D,U}
(\n_{\epsilon}(\nu,\m^{\epsilon}))}^d) \subset
\mathcal D_{\Kappa^0 +\kappa,\nu,\epsilon,\m^{\epsilon}}^{d,m,D,U} \text{ (см. } [13]).
\end{equation*}

Лемма 2.2.2

Пусть $ d \in \N, $ а ограниченная область $ D \subset \R^d $ и её открытое 
подмножество $ U \subset D $ удовлетворяют условиям леммы 2.2.1. 
Пусть ещё $ \alpha \in \R_+^d, 1 \le p < \infty, 1 \le q \le \infty, \lambda \in 
\Z_+^d $ и выполняется неравенство (2.1.9), а также $ m \in \N^d, \lambda 
\in \Z_+^d(m). $ Тогда существуют константы $ c_{1}(d,\alpha,p,q,\lambda,m,D,U) >0,
c_{2}(d,m,D,U) >0, c_3(d,\alpha,m,D,U) >0 $ и $ \bm \epsilon \in \R_+^d $
такие, что для любой функции $ f \in (S_p^\alpha H)^\prime(D) $ при
$ l = l(\alpha), \kappa \in \Z_+^d \setminus \{0\} $ соблюдается неравенство
\begin{multline*} \tag{2.2.11}
\| \D^\lambda \mathcal R_{\Kappa^0,\kappa}^{d,l -\e,m,D,U,\Nu} f
\|_{L_q(\R^d)} \le c_{1} 2^{(\kappa, \lambda +(p^{-1} -q^{-1})_+ \e)}
\biggl(\Omega^{\prime l \chi_{\s(\kappa)}}(f, (c_{2} 2^{-\kappa})^{\s(\kappa)})_{L_p(D)} \\
+\sum_{J \subset \Nu_{1,d}^1: J \ne \emptyset}
\biggl(\int_{ (c_{3} I^d)^J} (\prod_{j \in J} u_j^{-p(\alpha_j -\bm \epsilon_j) -1})
(\Omega^{\prime l \chi_{J \cup \s(\kappa)}}(f, (u \chi_{J \setminus \s(\kappa)} \\
+2^{-\kappa} u \chi_{\s(\kappa) \cap J} +c_{2} 2^{-\kappa}
\chi_{\s(\kappa) \setminus J})^{J \cup \s(\kappa)})_{L_p(D)})^p du^J \biggr)^{1/p} \biggr)\\
\text{(см. (2.1.43) с $ \Kappa^0, \Nu $ из леммы 2.2.1).}
\end{multline*}

Доказательство.

В условиях леммы пусть $ f \in (S_p^\alpha H)^\prime(D), l = l(\alpha),
\kappa \in \Z_+^d \setminus \{0\}, p \le q. $
Тогда, принимая во внимание предложение 2.1.3, а также (2.1.44), (2.1.42), имеем
\begin{multline*} \tag{2.2.12}
\| \D^\lambda \mathcal R_{\Kappa^0,\kappa}^{d,l -\e,m,D,U,\Nu} f \|_{L_q(\R^d)} = \\
\biggl\| \D^\lambda \biggl(\sum_{ \nu \in N_{\Kappa^0 +\kappa}^{d,m,U}}
(V_{\Kappa^0, \kappa,\nu}^{d,l -\e,m,D,U,\Nu} f) g_{\Kappa^0 +\kappa, \nu}^{d,m}\biggr) \biggr\|_{L_q(\R^d)} = \\
\biggl\| \sum_{ \nu \in N_{\Kappa^0 +\kappa}^{d,m,U}}
\D^\lambda \biggl((V_{\Kappa^0, \kappa,\nu}^{d,l -\e,m,D,U,\Nu} f)
g_{\Kappa^0 +\kappa, \nu}^{d,m} \biggr) \biggr\|_{L_q(\R^d)} = \\
\biggl\| \sum_{ \nu \in N_{\Kappa^0 +\kappa}^{d,m,U}}
\sum_{ \mu \in \Z_+^d(\lambda)} C_\lambda^\mu
\D^\mu (V_{\Kappa^0, \kappa,\nu}^{d,l -\e,m,D,U,\Nu} f)
\D^{\lambda -\mu} g_{\Kappa^0 +\kappa, \nu}^{d,m} \biggr\|_{L_q(\R^d)} = \\
\biggl\| \sum_{ \mu \in \Z_+^d(\lambda)} C_\lambda^\mu
\sum_{ \nu \in N_{\Kappa^0 +\kappa}^{d,m,U}}
\D^\mu (V_{\Kappa^0, \kappa,\nu}^{d,l -\e,m,D,U,\Nu} f)
\D^{\lambda -\mu} g_{\Kappa^0 +\kappa, \nu}^{d,m} \biggr\|_{L_q(\R^d)} \le \\
\sum_{ \mu \in \Z_+^d(\lambda)} C_\lambda^\mu
\biggl\| \sum_{ \nu \in N_{\Kappa^0 +\kappa}^{d,m,U}}
\D^\mu (V_{\Kappa^0, \kappa,\nu}^{d,l -\e,m,D,U,\Nu} f)
\D^{\lambda -\mu} g_{\Kappa^0 +\kappa, \nu}^{d,m} \biggr\|_{L_q(\R^d)}.
\end{multline*}

Оценивая правую часть (2.2.12), при $ \mu \in \Z_+^d(\lambda) $ с учётом
(2.1.24) получаем (ср. с (2.1.35))
\begin{multline*} \tag{2.2.13}
\biggl\| \sum_{ \nu \in N_{\Kappa^0 +\kappa}^{d,m,U}}
\D^\mu (V_{\Kappa^0, \kappa,\nu}^{d,l -\e,m,D,U,\Nu} f)
\D^{\lambda -\mu} g_{\Kappa^0 +\kappa, \nu}^{d,m} \biggr\|_{L_q(\R^d)}^q = \\
\int_{G_{\Kappa^0 +\kappa}^{d,m,U}}
\biggl| \sum_{ \nu \in N_{\Kappa^0 +\kappa}^{d,m,U}}
\D^\mu (V_{\Kappa^0, \kappa,\nu}^{d,l -\e,m,D,U,\Nu} f)
\D^{\lambda -\mu} g_{\Kappa^0 +\kappa, \nu}^{d,m} \biggr|^q dx = \\
\sum_{\substack{n \in \Z^d:\\ Q_{\Kappa^0 +\kappa,n}^d \cap
G_{\Kappa^0 +\kappa}^{d,m,U} \ne \emptyset}}
\int\limits_{Q_{\Kappa^0 +\kappa,n}^d} \biggl| \sum_{\substack{\nu \in
N_{\Kappa^0 +\kappa}^{d,m,U}:\\ Q_{\Kappa^0 +\kappa,n}^d \cap
\supp g_{\Kappa^0 +\kappa, \nu}^{d,m} \ne \emptyset}}
\D^\mu (V_{\Kappa^0, \kappa,\nu}^{d,l -\e,m,D,U,\Nu} f)
\D^{\lambda -\mu} g_{\Kappa^0 +\kappa, \nu}^{d,m} \biggr|^q dx \le \\
\sum_{\substack{n \in \Z^d:\\ Q_{\Kappa^0 +\kappa,n}^d \cap G_{\Kappa^0 +\kappa}^{d,m,U} \ne \emptyset}}
\biggl(\sum_{\substack{\nu \in N_{\Kappa^0 +\kappa}^{d,m,U}: \\ Q_{\Kappa^0 +\kappa,n}^d \cap
\supp g_{\Kappa^0 +\kappa, \nu}^{d,m} \ne \emptyset}}
\biggl\| \D^\mu (V_{\Kappa^0, \kappa,\nu}^{d,l -\e,m,D,U,\Nu} f)
\D^{\lambda -\mu} g_{\Kappa^0 +\kappa, \nu}^{d,m} \biggr\|_{L_q(Q_{\Kappa^0 +\kappa,n}^d)} \biggr)^q.
\end{multline*}

Для проведения оценки правой части (2.2.13) приведём некоторые полезные для
нас факты. Прежде всего, в условиях леммы, пользуясь (1.4.5), (2.2.1),
видим, что при $ \kappa \in \Z_+^d $ для $ \nu \in
N_{\Kappa^0 +\kappa}^{d,m,U} $ имеет место включение
\begin{equation*} \tag{2.2.14}
2^{-\Kappa^0 -\kappa} \nu -2^{-\Kappa^0 -\kappa}
n_{\Kappa^0 +\kappa}^{d,m,D,U}(\nu) \in \Gamma^0 2^{-\Kappa^0 -\kappa} B^d.
\end{equation*}

Кроме того, при $ \kappa \in \Z_+^d $ для $ n \in \Z^d:
Q_{\Kappa^0 +\kappa,n}^d \cap G_{\Kappa^0 +\kappa}^{d,m,U} \ne \emptyset, $
и $ \nu \in N_{\Kappa^0 +\kappa}^{d,m,U}: Q_{\Kappa^0 +\kappa,n}^d
\cap \supp g_{\Kappa^0 +\kappa, \nu}^{d,m} \ne \emptyset, $ учитывая
(2.1.25), (1.4.4), заключаем, что
\begin{equation*}
Q_{\Kappa^0 +\kappa, n}^d \subset 2^{-\Kappa^0 -\kappa} \nu +
2^{-\Kappa^0 -\kappa} (m +\e) \overline I^d \subset 2^{-\Kappa^0 -\kappa} \nu +
2^{-\Kappa^0 -\kappa} (m +\e) B^d
\end{equation*}
или
\begin{equation*}
(Q_{\Kappa^0 +\kappa, n}^d -2^{-\Kappa^0 -\kappa} \nu) \subset
2^{-\Kappa^0 -\kappa} (m +\e) B^d.
\end{equation*}
Из последнего включения и (2.2.14) вытекает, что при $ \kappa \in \Z_+^d $
для $ n \in \Z^d: Q_{\Kappa^0 +\kappa,n}^d \cap
G_{\Kappa^0 +\kappa}^{d,m,U} \ne \emptyset, $
и $ \nu \in N_{\Kappa^0 +\kappa}^{d,m,U}: Q_{\Kappa^0 +\kappa,n}^d
\cap \supp g_{\Kappa^0 +\kappa, \nu}^{d,m} \ne \emptyset, $ справедливо
соотношение
\begin{multline*}
Q_{\Kappa^0 +\kappa, n}^d -2^{-\Kappa^0 -\kappa} n_{\Kappa^0 +\kappa}^{d,m,D,U}(\nu) =\\
Q_{\Kappa^0 +\kappa, n}^d -2^{-\Kappa^0 -\kappa} \nu +2^{-\Kappa^0 -\kappa} \nu
-2^{-\Kappa^0 -\kappa} n_{\Kappa^0 +\kappa}^{d,m,D,U}(\nu) \subset \\
2^{-\Kappa^0 -\kappa} (m +\e) B^d +\Gamma^0 2^{-\Kappa^0 -\kappa} B^d \subset \\
(m +\e +\Gamma^0) 2^{-\Kappa^0 -\kappa} B^d =
\Gamma^2(d,m,D,U) 2^{-\Kappa^0 -\kappa} B^d
\end{multline*}
с $ \Gamma^2 = m +\e +\Gamma^0 > \e, $
или
\begin{equation*} \tag{2.2.15}
Q_{\Kappa^0 +\kappa, n}^d \subset
2^{-\Kappa^0 -\kappa} n_{\Kappa^0 +\kappa}^{d,m,D,U}(\nu) +
\Gamma^2 2^{-\Kappa^0 -\kappa} B^d,
\end{equation*}
а также ввиду замкнутости $ B^d $ --
\begin{equation*}
2^{-\Kappa^0 -\kappa} n +2^{-\Kappa^0 -\kappa} \overline I^d =
\overline Q_{\Kappa^0 +\kappa, n}^d \subset
2^{-\Kappa^0 -\kappa} n_{\Kappa^0 +\kappa}^{d,m,D,U}(\nu) +
\Gamma^2 2^{-\Kappa^0 -\kappa} B^d,
\end{equation*}
в частности,
\begin{equation*}
2^{-\Kappa^0 -\kappa} n \in
2^{-\Kappa^0 -\kappa} n_{\Kappa^0 +\kappa}^{d,m,D,U}(\nu) +
\Gamma^2 2^{-\Kappa^0 -\kappa} B^d,
\end{equation*}
или
\begin{equation*}
2^{-\Kappa^0 -\kappa} n -
2^{-\Kappa^0 -\kappa} n_{\Kappa^0 +\kappa}^{d,m,D,U}(\nu) \in
\Gamma^2 2^{-\Kappa^0 -\kappa} B^d,
\end{equation*}
а, значит,
\begin{equation*} \tag{2.2.16}
2^{-\Kappa^0 -\kappa} n_{\Kappa^0 +\kappa}^{d,m,D,U}(\nu) -
2^{-\Kappa^0 -\kappa} n \in \Gamma^2 2^{-\Kappa^0 -\kappa} B^d.
\end{equation*}

Далее, учитывая (2.2.9), понятно, что  при $ \kappa \in \Z_+^d $
для $ \nu \in N_{\Kappa^0 +\kappa}^{d,m,U}, \epsilon \in \Upsilon^d:
\s(\epsilon) \subset \s(\kappa), \m^{\epsilon} \in \M_{\epsilon}^m(\nu) $
соблюдается включение
\begin{equation*} \tag{2.2.17}
\overline {\mathcal D}_{\Kappa^0 +\kappa,\nu,\epsilon,\m^{\epsilon}}^{d,m,D,U} -
\bm x_{\Kappa^0 +\kappa,\nu,\epsilon,\m^{\epsilon}}^{d,m,D,U} \subset
\Gamma^1 2^{-\Kappa^0 -\kappa} B^d,
\end{equation*}
кроме того, принимая во внимание (1.4.5), (2.2.10), получаем, что
\begin{equation*}
2^{-\Kappa^0 -\kappa} n_{\Kappa^0 +\kappa}^{d,m,D,U}(\nu) \in
\overline Q_{\Kappa^0 +\kappa, n_{\Kappa^0 +\kappa}^{d,m,D,U}(\nu)}^d \subset
\overline {\mathcal D}_{\Kappa^0 +\kappa,\nu,\epsilon,\m^{\epsilon}}^{d,m,D,U},
\end{equation*}
что в силу (2.2.17) приводит к включению
\begin{equation*}
(2^{-\Kappa^0 -\kappa} n_{\Kappa^0 +\kappa}^{d,m,D,U}(\nu) -
\bm x_{\Kappa^0 +\kappa,\nu,\epsilon,\m^{\epsilon}}^{d,m,D,U}) \in
\Gamma^1 2^{-\Kappa^0 -\kappa} B^d,
\end{equation*}
а, значит,
\begin{equation*} \tag{2.2.18
}
(\bm x_{\Kappa^0 +\kappa,\nu,\epsilon,\m^{\epsilon}}^{d,m,D,U} -
2^{-\Kappa^0 -\kappa} n_{\Kappa^0 +\kappa}^{d,m,D,U}(\nu)) \in
\Gamma^1 2^{-\Kappa^0 -\kappa} B^d,
\end{equation*}
наконец, из (2.2.17), (2.2.18) вытекает, что
\begin{multline*} \tag{2.2.19}
\mathcal D_{\Kappa^0 +\kappa,\nu,\epsilon,\m^{\epsilon}}^{d,m,D,U} -
2^{-\Kappa^0 -\kappa} n_{\Kappa^0 +\kappa}^{d,m,D,U}(\nu) \subset
(\overline {\mathcal D}_{\Kappa^0 +\kappa,\nu,\epsilon,\m^{\epsilon}}^{d,m,D,U} -
\bm x_{\Kappa^0 +\kappa,\nu,\epsilon,\m^{\epsilon}}^{d,m,D,U}) + \\
(\bm x_{\Kappa^0 +\kappa,\nu,\epsilon,\m^{\epsilon}}^{d,m,D,U} -
2^{-\Kappa^0 -\kappa} n_{\Kappa^0 +\kappa}^{d,m,D,U}(\nu)) \subset \\
\Gamma^1 2^{-\Kappa^0 -\kappa} B^d +\Gamma^1 2^{-\Kappa^0 -\kappa} B^d
\subset 2 \Gamma^1 2^{-\Kappa^0 -\kappa} B^d.
\end{multline*}
Таким образом, при $ \kappa \in \Z_+^d $ для $ n \in \Z^d: Q_{\Kappa^0 +\kappa,n}^d
\cap G_{\Kappa^0 +\kappa}^{d,m,U} \ne \emptyset, $
и $ \nu \in N_{\Kappa^0 +\kappa}^{d,m,U}: Q_{\Kappa^0 +\kappa,n}^d
\cap \supp g_{\Kappa^0 +\kappa, \nu}^{d,m} \ne \emptyset, \epsilon \in
\Upsilon^d: \s(\epsilon) \subset \s(\kappa), \m^{\epsilon} \in
\M_{\epsilon}^m(\nu) $ согласно (2.2.19), (2.2.16) справедливо включение
\begin{multline*}
\mathcal D_{\Kappa^0 +\kappa,\nu,\epsilon,\m^{\epsilon}}^{d,m,D,U} -
2^{-\Kappa^0 -\kappa} n =
(\mathcal D_{\Kappa^0 +\kappa,\nu,\epsilon,\m^{\epsilon}}^{d,m,D,U} -
2^{-\Kappa^0 -\kappa} n_{\Kappa^0 +\kappa}^{d,m,D,U}(\nu)) + \\
(2^{-\Kappa^0 -\kappa} n_{\Kappa^0 +\kappa}^{d,m,D,U}(\nu) -
2^{-\Kappa^0 -\kappa} n) \subset \\
2 \Gamma^1 2^{-\Kappa^0 -\kappa} B^d +\Gamma^2 2^{-\Kappa^0 -\kappa} B^d \subset
\Gamma^3(d,m,D,U) 2^{-\Kappa^0 -\kappa} B^d
\end{multline*}
с $ \Gamma^3 = 2 \Gamma^1 +\Gamma^2 > \e, $ которое в соединении с (2.2.2) даёт
\begin{equation*}
\mathcal D_{\Kappa^0 +\kappa,\nu,\epsilon,\m^{\epsilon}}^{d,m,D,U} \subset D
\cap (2^{-\Kappa^0 -\kappa} n +\Gamma^3 2^{-\Kappa^0 -\kappa} B^d)
\end{equation*}
или
\begin{equation*} \tag{2.2.20}
\mathcal D_{\Kappa^0 +\kappa,\nu,\epsilon,\m^{\epsilon}}^{d,m,D,U} \subset
D \cap D_{\Kappa^0 +\kappa,n}^{\prime d,m,D,U},
\end{equation*}

где
$$
D_{\Kappa^0 +\kappa,n}^{\prime d,m,D,U} = 2^{-\Kappa^0 -\kappa} n +\Gamma^3
2^{-\Kappa^0 -\kappa} B^d,
\kappa \in \Z_+^d, n \in \Z^d: Q_{\Kappa^0 +\kappa,n}^d \cap
G_{\Kappa^0 +\kappa}^{d,m,U} \ne \emptyset.
$$

Из приведенных определений с учётом того, что $ \Gamma^3 > \e, $ видно, что
при $ \kappa \in \Z_+^d, n \in \Z^d: Q_{\Kappa^0 +\kappa,n}^d \cap
G_{\Kappa^0 +\kappa}^{d,m,U} \ne \emptyset, $
справедливо включение
\begin{equation*} \tag{2.2.21}
Q_{\Kappa^0 +\kappa, n}^d \subset D_{\Kappa^0 +\kappa, n}^{\prime d,m,D,U}.
\end{equation*}

Используя (2.2.21), легко проверить, что существует константа
$ c_4(d,m,D,U) >0 $ такая, что при $ \kappa \in \Z_+^d $ для каждого
$ x \in \R^d $ число
\begin{equation*} \tag{2.2.22}
\card \{ n \in \Z^d: Q_{\Kappa^0 +\kappa,n}^d \cap G_{\Kappa^0 +\kappa}^{d,m,U} \ne \emptyset,
x \in D_{\Kappa^0 +\kappa,n}^{\prime d,m,D,U} \} \le c_4.
\end{equation*}

Перейдём к оценке правой части (2.2.13). При $ n \in \Z^d:
Q_{\Kappa^0 +\kappa,n}^d \cap G_{\Kappa^0 +\kappa}^{d,m,U} \ne \emptyset,
\nu \in N_{\Kappa^0 +\kappa}^{d,m,U}:  Q_{\Kappa^0 +\kappa,n}^d \cap
\supp g_{\Kappa^0 +\kappa, \nu}^{d,m} \ne \emptyset, \mu \in \Z_+^d(\lambda), $
используя (1.4.7), а также, ввиду (2.1.45), (2.2.15) применяя (1.1.3), выводим
\begin{multline*} \tag{2.2.23}
\biggl\| \D^\mu (V_{\Kappa^0, \kappa,\nu}^{d,l -\e,m,D,U,\Nu} f)
\D^{\lambda -\mu} g_{\Kappa^0 +\kappa, \nu}^{d,m} \biggr\|_{L_q(Q_{\Kappa^0 +\kappa,n}^d)} \le \\
\| \D^{\lambda -\mu} g_{\Kappa^0 +\kappa, \nu}^{d,m} \|_{L_\infty(\R^d)} \cdot
\biggl\| \D^\mu (V_{\Kappa^0, \kappa,\nu}^{d,l -\e,m,D,U,\Nu} f) \biggr\|_{L_q(Q_{\Kappa^0 +\kappa,n}^d)} \le \\
c_5 2^{(\kappa, \lambda -\mu)}
\biggl\| \D^\mu (V_{\Kappa^0, \kappa,\nu}^{d,l -\e,m,D,U,\Nu} f) \biggr\|_{L_q(Q_{\Kappa^0 +\kappa,n}^d)} \le \\
c_5 2^{(\kappa, \lambda -\mu)} c_6 (2^{-\kappa})^{-\mu +q^{-1} \e}
\biggl\| V_{\Kappa^0, \kappa,\nu}^{d,l -\e,m,D,U,\Nu} f
\biggr\|_{L_\infty(Q_{\Kappa^0 +\kappa, n_{\Kappa^0 +\kappa}^{d,m,D,U}(\nu)}^d)} = \\
c_7 2^{(\kappa, \lambda -q^{-1} \e)}
\biggl\| V_{\Kappa^0, \kappa,\nu}^{d,l -\e,m,D,U,\Nu} f
\biggr\|_{L_\infty(Q_{\Kappa^0 +\kappa, n_{\Kappa^0 +\kappa}^{d,m,D,U}(\nu)}^d)}.
\end{multline*}

В условиях леммы преобразуем выражение (2.1.45) (с $ l -\e $ вместо $ l,
\Kappa^0 $ вместо $ \kappa^0 $ ) к виду, подходящему для получения интересующих
нас оценок. Для этого при $ d \in \N, m \in \N^d $ для ограниченной области
$ D \subset \R^d $ построим клетку $ Q^{d,m,D} = X^0 +\Delta I^d, X^0 =
X^{0,d,m,D} \in \R^d, \Delta = \Delta^{d,m,D} \in \R_+^d, $ такую, что
при любом $ \kappa \in \Z_+^d $ выполняется включение
\begin{equation*} \tag{2.2.24}
D \subset G_\kappa^{d,m,D} \subset Q^{d,m,D}.
\end{equation*}
Тогда в условиях леммы, фиксируя набор функций $ \chi^j = \chi^j_{\Delta_j, X^0_j}, j =1,\ldots,d, $
удовлетворяющих условиям леммы 1.3.3 (при $ J^\prime = \Nu_{1,d}^1 $ ), и
полагая $ \chi(x) = \prod_{j =1}^d \chi^j(x_j), x \in \R^d, $ а также при
$ \kappa \in \Z_+^d, \nu \in N_{\Kappa^0 +\kappa}^{d,m,U} $ обозначая 
$ F = \mathcal E^{d,\alpha,p,\infty,\mathcal D} (f \mid_{\mathcal D}) $ при
$ \mathcal D = \mathcal D_{\Kappa^0 +\kappa,\nu}^{d,m,D,U}, f \in (S_p^\alpha H)^\prime(D) $
(см. теорему 2.1.9 и (2.2.2)), принимая во внимание (2.2.10), (2.2.2), предложение 2.1.3 и
(2.1.46), а также с учётом (2.1.48) применяя (2.1.21), заключаем, что 
для $ f \in (S_p^\alpha H)^\prime(D) $ для всех $ x \in \R^d $ имеет место 
равенство

\begin{multline*}
\chi(x) (\mathtt R_{\Kappa^0 +\kappa -\epsilon,
\nu_{\Kappa^0 +\kappa -\epsilon}^{d,m,D,U}
(\n_{\epsilon}(\nu,\m^{\epsilon}))}^{d,l -\e}
(f \mid_{(Q_{\Kappa^0 +\kappa -\epsilon,\nu_{\Kappa^0 +\kappa -\epsilon}^{d,m,D,U}
(\n_{\epsilon}(\nu,\m^{\epsilon}))}^d}))(x) = \\
\chi(x) (\mathtt R_{\Kappa^0 +\kappa -\epsilon,
\nu_{\Kappa^0 +\kappa -\epsilon}^{d,m,D,U}
(\n_{\epsilon}(\nu,\m^{\epsilon}))}^{d,l -\e}
(F \mid_{(Q_{\Kappa^0 +\kappa -\epsilon,\nu_{\Kappa^0 +\kappa -\epsilon}^{d,m,D,U}
(\n_{\epsilon}(\nu,\m^{\epsilon}))}^d}))(x) = \\
\biggl((\prod_{j =1}^d V_j(M_{\chi^j_{\Delta_j, X^0_j}}
\mathtt R_{\Kappa^0_j +\kappa_j -\epsilon_j,
(\nu_{\Kappa^0 +\kappa -\epsilon}^{d,m,D,U}
(\n_{\epsilon}(\nu,\m^{\epsilon})))_j}^{1,l_j -1})) F \biggr)(x), \\
\kappa \in \Z_+^d, \nu \in N_{\Kappa^0 +\kappa}^{d,m,U}, \epsilon \in 
\Upsilon^d: \s(\epsilon) \subset \s(\kappa), \m^{\epsilon} \in \M_{\epsilon}^m(\nu).
\end{multline*}
Поэтому благодаря п. 2) леммы 1.3.2 и (1.3.4), в условиях леммы для
$ f \in (S_p^\alpha H)^\prime(D) $ при $ \kappa \in \Z_+^d,
\nu \in N_{\Kappa^0 +\kappa}^{d,m,U}, \epsilon \in \Upsilon^d: \s(\epsilon)
\subset \s(\kappa), \m^{\epsilon} \in \M_{\epsilon}^m(\nu) $ для всех
$ x \in \R^d $ выполняется равенство
\begin{multline*} \tag{2.2.25}
\chi(x) (\mathtt R_{\Kappa^0 +\kappa -\epsilon,
\nu_{\Kappa^0 +\kappa -\epsilon}^{d,m,D,U}
(\n_{\epsilon}(\nu,\m^{\epsilon}))}^{d, l -\e} f)(x) = \\
\biggl((\prod_{j =1}^d V_j(M_{\chi^j_{\Delta_j, X^0_j}}
\mathtt R_{\Kappa^0_j +\kappa_j -\epsilon_j,
(\nu_{\Kappa^0 +\kappa -\epsilon}^{d,m,D,U}
(\n_{\epsilon}(\nu,\m^{\epsilon})))_j}^{1,l_j -1})) F \biggr)(x) = \\
\biggl((\prod_{j =1}^d (E -V_j(E -M_{\chi^j_{\Delta_j, X^0_j}}
\mathtt R_{\Kappa^0_j +\kappa_j -\epsilon_j,
(\nu_{\Kappa^0 +\kappa -\epsilon}^{d,m,D,U}
(\n_{\epsilon}(\nu,\m^{\epsilon})))_j}^{1,l_j -1}))) F \biggr)(x) = \\
\sum_{ \gamma \in \Upsilon^d} (-\e)^\gamma \biggl((\prod_{j \in \s(\gamma)}
V_j(E -M_{\chi^j_{\Delta_j, X^0_j}} \mathtt R_{\Kappa^0_j +\kappa_j -\epsilon_j,
(\nu_{\Kappa^0 +\kappa -\epsilon}^{d,m,D,U}
(\n_{\epsilon}(\nu,\m^{\epsilon})))_j}^{1,l_j -1})) F \biggr)(x).
\end{multline*}

Исходя из (2.1.45), пользуясь (2.2.25), находим, что в условиях леммы при
$ \kappa \in \Z_+^d, \nu \in N_{\Kappa^0 +\kappa}^{d,m,U} $
для $ f \in (S_p^\alpha H)^\prime(D) $ имеет место равенство
\begin{multline*} \tag{2.2.26}
(V_{\Kappa^0, \kappa,\nu}^{d,l -\e,m,D,U,\Nu} f) \mid_{Q^{d,m,D}} =\\
\biggl(\sum_{\epsilon \in \Upsilon^d: \s(\epsilon) \subset \s(\kappa)} (-\e)^\epsilon
\sum_{\m^{\epsilon} \in \M_{\epsilon}^m(\nu)} A_{\m^{\epsilon}}^m
\mathtt R_{\Kappa^0 +\kappa -\epsilon, \nu_{\Kappa^0 +\kappa -\epsilon}^{d,m,D,U}
(\n_{\epsilon}(\nu,\m^{\epsilon}))}^{d,l -\e} f \biggr) \mid_{Q^{d,m,D}} = \\
\sum_{\epsilon \in \Upsilon^d: \s(\epsilon) \subset \s(\kappa)} (-\e)^\epsilon
\sum_{\m^{\epsilon} \in \M_{\epsilon}^m(\nu)} A_{\m^{\epsilon}}^m
(\chi \cdot (\mathtt R_{\Kappa^0 +\kappa -\epsilon, \nu_{\Kappa^0 +\kappa -\epsilon}^{d,m,D,U}
(\n_{\epsilon}(\nu,\m^{\epsilon}))}^{d,l -\e} f)) \mid_{Q^{d,m,D}} = \\
\sum_{\epsilon \in \Upsilon^d: \s(\epsilon) \subset \s(\kappa)} (-\e)^\epsilon
\sum_{\m^{\epsilon} \in \M_{\epsilon}^m(\nu)} A_{\m^{\epsilon}}^m \\
\biggl(\sum_{ \gamma \in \Upsilon^d} (-\e)^\gamma ((\prod_{j \in \s(\gamma)}
V_j(E -M_{\chi^j_{\Delta_j, X^0_j}} \mathtt R_{\Kappa^0_j +\kappa_j -\epsilon_j,
(\nu_{\Kappa^0 +\kappa -\epsilon}^{d,m,D,U}
(\n_{\epsilon}(\nu,\m^{\epsilon})))_j}^{1,l_j -1})) F) \biggr) \mid_{Q^{d,m,D}} = \\
\biggl(\sum_{\epsilon \in \Upsilon^d: \s(\epsilon) \subset \s(\kappa)} (-\e)^\epsilon
\sum_{\m^{\epsilon} \in \M_{\epsilon}^m(\nu)} A_{\m^{\epsilon}}^m
\sum_{ \gamma \in \Upsilon^d} (-\e)^\gamma \\
\times \biggl((\prod_{j \in \s(\gamma)}
V_j(E -M_{\chi^j_{\Delta_j, X^0_j}} \mathtt R_{\Kappa^0_j +\kappa_j -\epsilon_j,
(\nu_{\Kappa^0 +\kappa -\epsilon}^{d,m,D,U}
(\n_{\epsilon}(\nu,\m^{\epsilon})))_j}^{1,l_j -1})) F \biggr) \biggr) \mid_{Q^{d,m,D}} = \\
\biggl(\sum_{ \gamma \in \Upsilon^d} (-\e)^\gamma
\sum_{\epsilon \in \Upsilon^d: \s(\epsilon) \subset \s(\kappa)} (-\e)^\epsilon
\sum_{\m^{\epsilon} \in \M_{\epsilon}^m(\nu)} A_{\m^{\epsilon}}^m \\
\biggl((\prod_{j \in \s(\gamma)}
V_j(E -M_{\chi^j_{\Delta_j, X^0_j}} \mathtt R_{\Kappa^0_j +\kappa_j -\epsilon_j,
(\nu_{\Kappa^0 +\kappa -\epsilon}^{d, m,D,U}
(\n_{\epsilon}(\nu,\m^{\epsilon})))_j}^{1,l_j -1})) F \biggr) \biggr) \mid_{Q^{d,m,D}} = \\
\biggl(\sum_{ \gamma \in \Upsilon^d: \s(\kappa) \setminus \s(\gamma) = \emptyset} (-\e)^\gamma
\sum_{\epsilon \in \Upsilon^d: \s(\epsilon) \subset \s(\kappa)} (-\e)^\epsilon
\sum_{\m^{\epsilon} \in \M_{\epsilon}^m(\nu)} A_{\m^{\epsilon}}^m \\
\biggl((\prod_{j \in \s(\gamma)}
V_j(E -M_{\chi^j_{\Delta_j, X^0_j}} \mathtt R_{\Kappa^0_j +\kappa_j -\epsilon_j,
(\nu_{\Kappa^0 +\kappa -\epsilon}^{d, m,D,U}
(\n_{\epsilon}(\nu,\m^{\epsilon})))_j}^{1,l_j -1})) F \biggr) + \\
\sum_{ \gamma \in \Upsilon^d: \s(\kappa) \setminus \s(\gamma) \ne \emptyset} (-\e)^\gamma
\sum_{\epsilon \in \Upsilon^d: \s(\epsilon) \subset \s(\kappa)} (-\e)^\epsilon
\sum_{\m^{\epsilon} \in \M_{\epsilon}^m(\nu)} A_{\m^{\epsilon}}^m \\
\biggl((\prod_{j \in \s(\gamma)}
V_j(E -M_{\chi^j_{\Delta_j, X^0_j}} \mathtt R_{\Kappa^0_j +\kappa_j -\epsilon_j,
(\nu_{\Kappa^0 +\kappa -\epsilon}^{d, m,D,U}
(\n_{\epsilon}(\nu,\m^{\epsilon})))_j}^{1,l_j -1})) F \biggr) \biggr) \mid_{Q^{d,m,D}}.
\end{multline*}

При $ \kappa \in \Z_+^d, \nu \in N_{\Kappa^0 +\kappa}^{d,m,U} $
для $ \gamma \in \Upsilon^d: \s(\kappa) \setminus \s(\gamma) \ne \emptyset, $
ввиду (1.4.11), (1.4.20) соблюдается равенство
\begin{multline*} \tag{2.2.27}
\sum_{\epsilon \in \Upsilon^d: \s(\epsilon) \subset \s(\kappa)} (-\e)^\epsilon
\sum_{\m^{\epsilon} \in \M_{\epsilon}^m(\nu)} A_{\m^{\epsilon}}^m \\
\biggl(\prod_{j \in \s(\gamma)}
V_j(E -M_{\chi^j_{\Delta_j, X^0_j}} \mathtt R_{\Kappa^0_j +\kappa_j -\epsilon_j,
(\nu_{\Kappa^0 +\kappa -\epsilon}^{d, m,D,U}
(\n_{\epsilon}(\nu,\m^{\epsilon})))_j}^{1,l_j -1}) \biggr) = \\
\sum_{\epsilon, \epsilon^\prime \in \Upsilon^d: \s(\epsilon) \subset
(\s(\kappa) \cap \s(\gamma)), \s(\epsilon^\prime) \subset
(\s(\kappa) \setminus \s(\gamma))} (-\e)^{\epsilon +\epsilon^\prime}
\sum_{\m^{\epsilon +\epsilon^\prime} \in \M_{\epsilon +\epsilon^\prime}^m(\nu)} A_{\m^{\epsilon +\epsilon^\prime}}^m \\
\biggl(\prod_{j \in \s(\gamma)}
V_j(E -M_{\chi^j_{\Delta_j, X^0_j}} \mathtt R_{\Kappa^0_j +\kappa_j -\epsilon_j -\epsilon_j^\prime,
(\nu_{\Kappa^0 +\kappa -\epsilon -\epsilon^\prime}^{d, m,D,U}
(\n_{\epsilon +\epsilon^\prime}(\nu,\m^{\epsilon +\epsilon^\prime})))_j}^{1,l_j -1}) \biggr) = \\
\sum_{\substack{\epsilon \in \Upsilon^d:\\ \s(\epsilon) \subset
(\s(\kappa) \cap \s(\gamma))}} (-\e)^\epsilon
\sum_{\substack{\epsilon^\prime \in \Upsilon^d:\\
\s(\epsilon^\prime) \subset (\s(\kappa) \setminus \s(\gamma))}}
(-\e)^{\epsilon^\prime} \sum_{\substack{ \m^{\epsilon +\epsilon^\prime} =
(\m^{\epsilon}, \m^{\epsilon^\prime}):\\ \m^{\epsilon} \in \M_{\epsilon}^m(\nu),
\m^{\epsilon^\prime} \in \M_{\epsilon^\prime}^m(\nu)}} (\prod_{j \in
\s(\epsilon +\epsilon^\prime)} a_{\m_j}^{m_j}) \\ \times
\biggl(\prod_{j \in \s(\gamma)}
V_j(E -M_{\chi^j_{\Delta_j, X^0_j}} \mathtt R_{\Kappa^0_j +\kappa_j -\epsilon_j -\epsilon_j^\prime,
(\nu_{\Kappa^0 +\kappa -\epsilon -\epsilon^\prime}^{d, m,D,U}
(\n_{\epsilon +\epsilon^\prime}(\nu,\m^{\epsilon +\epsilon^\prime})))_j}^{1,l_j -1}) \biggr) = \\
\sum_{\substack{\epsilon \in \Upsilon^d:\\ \s(\epsilon) \subset
(\s(\kappa) \cap \s(\gamma))}} (-\e)^\epsilon
\sum_{\substack{\epsilon^\prime \in \Upsilon^d:\\
\s(\epsilon^\prime) \subset (\s(\kappa) \setminus \s(\gamma))}}
(-\e)^{\epsilon^\prime} \sum_{\substack{\m^{\epsilon} \in \M_{\epsilon}^m(\nu),
\m^{\epsilon^\prime} \in \M_{\epsilon^\prime}^m(\nu),\\ \m^{\epsilon +\epsilon^\prime} =
(\m^{\epsilon}, \m^{\epsilon^\prime})}} (\prod_{j \in \s(\epsilon)} a_{\m_j}^{m_j})
(\prod_{j \in \s(\epsilon^\prime)} a_{\m_j}^{m_j}) \\ \times
\biggl(\prod_{j \in \s(\gamma)}
V_j(E -M_{\chi^j_{\Delta_j, X^0_j}} \mathtt R_{\Kappa^0_j +\kappa_j -\epsilon_j -\epsilon_j^\prime,
(\nu_{\Kappa^0 +\kappa -\epsilon -\epsilon^\prime}^{d, m,D,U}
(\n_{\epsilon +\epsilon^\prime}(\nu,\m^{\epsilon +\epsilon^\prime})))_j}^{1,l_j -1}) \biggr).
\end{multline*}

Далее, заметим, что при $ \kappa \in \Z_+^d $ для $ \nu \in
N_{\Kappa^0 +\kappa}^{d,m,U}, \epsilon, \epsilon^\prime \in \Upsilon^d:
\s(\epsilon) \cap \s(\epsilon^\prime) = \emptyset,
\s(\epsilon +\epsilon^\prime) \subset \s(\kappa), $ и любых
$ \m^{\epsilon +\epsilon^\prime} \in \M_{\epsilon +\epsilon^\prime}^m(\nu),
\m^{\epsilon} \in \M_{\epsilon}^m(\nu),
\m^{\epsilon^\prime} \in \M_{\epsilon^\prime}^m(\nu):
(\m^{\epsilon +\epsilon^\prime})_j = (\m^{\epsilon})_j, j \in \s(\epsilon),
(\m^{\epsilon +\epsilon^\prime})_j = (\m^{\epsilon^\prime})_j, j \in
\s(\epsilon^\prime), $ согласно (1.4.13), (1.4.14) и (2.2.3)
при $ j \in \Nu_{1,d}^1 \setminus \s(\epsilon^\prime) $ имеет место равенство
\begin{equation*}
(\nu_{\Kappa^0 +\kappa -\epsilon -\epsilon^\prime}^{d,m,D,U}
(\n_{\epsilon +\epsilon^\prime}(\nu,\m^{\epsilon +\epsilon^\prime})))_j =
(\nu_{\Kappa^0 +\kappa -\epsilon -\epsilon^\prime}^{d,m,D,U}
(\n_{\epsilon^\prime}(\n_{\epsilon}(\nu, \m^{\epsilon}), \m^{\epsilon^\prime})))_j = \\
(\nu_{\Kappa^0 +\kappa -\epsilon}^{d,m,D,U}(\n_{\epsilon}(\nu,\m^{\epsilon})))_j.
\end{equation*}
Принимая во внимание это обстоятельство,
при $ \kappa \in \Z_+^d, \nu \in
N_{\Kappa^0 +\kappa}^{d,m,U}, \gamma \in \Upsilon^d: \s(\kappa) \setminus \s(\gamma) \ne \emptyset,
\epsilon \in \Upsilon^d: \s(\epsilon) \subset (\s(\kappa) \cap \s(\gamma)),
\epsilon^\prime \in \Upsilon^d: \s(\epsilon^\prime) \subset (\s(\kappa)
\setminus \s(\gamma)), \m^{\epsilon} \in \M_{\epsilon}^m(\nu),
\m^{\epsilon^\prime} \in \M_{\epsilon^\prime}^m(\nu), \m^{\epsilon +\epsilon^\prime}: $
\begin{equation*}
(\m^{\epsilon +\epsilon^\prime})_j = \begin{cases} (\m^{\epsilon})_j,
\text{ при } j \in \s(\epsilon); \\
(\m^{\epsilon^\prime})_j, \text{ при } j \in \s(\epsilon^\prime),
\end{cases}
\end{equation*}
с учётом (1.4.20), а также включения $ \s(\gamma) \subset
\Nu_{1,d}^1 \setminus \s(\epsilon^\prime) $ выводим
\begin{multline*} \tag{2.2.28}
(\prod_{j \in \s(\epsilon)} a_{\m_j}^{m_j})
(\prod_{j \in \s(\epsilon^\prime)} a_{\m_j}^{m_j}) \\ \times
(\prod_{j \in \s(\gamma)}
V_j(E -M_{\chi^j_{\Delta_j, X^0_j}} \mathtt R_{\Kappa^0_j +\kappa_j -\epsilon_j -\epsilon_j^\prime,
(\nu_{\Kappa^0 +\kappa -\epsilon -\epsilon^\prime}^{d,m,D,U}
(\n_{\epsilon +\epsilon^\prime}(\nu,\m^{\epsilon +\epsilon^\prime})))_j}^{1,l_j -1})) = \\
A_{\m^{\epsilon}}^m A_{\m^{\epsilon^\prime}}^m
(\prod_{j \in \s(\gamma)}
V_j(E -M_{\chi^j_{\Delta_j, X^0_j}} \mathtt R_{\Kappa^0_j +\kappa_j -\epsilon_j,
(\nu_{\Kappa^0 +\kappa -\epsilon}^{d,m,D,U}
(\n_{\epsilon}(\nu,\m^{\epsilon})))_j}^{1,l_j -1})).
\end{multline*}

Кроме того, используя (2.2.28), (1.4.21), при $ \kappa \in \Z_+^d,
\nu \in N_{\Kappa^0 +\kappa}^{d,m,U}, \gamma \in \Upsilon^d:
\s(\kappa) \setminus \s(\gamma) \ne \emptyset, \epsilon \in \Upsilon^d:
\s(\epsilon) \subset (\s(\kappa) \cap \s(\gamma)), \epsilon^\prime \in \Upsilon^d:
\s(\epsilon^\prime) \subset (\s(\kappa) \setminus \s(\gamma)), $ находим, что
\begin{multline*} \tag{2.2.29}
\sum_{ \m^{\epsilon} \in \M_{\epsilon}^m(\nu), \m^{\epsilon^\prime}
\in \M_{\epsilon^\prime}^m(\nu), \m^{\epsilon +\epsilon^\prime} =
(\m^{\epsilon}, \m^{\epsilon^\prime})}
(\prod_{j \in \s(\epsilon)} a_{\m_j}^{m_j})
(\prod_{j \in \s(\epsilon^\prime)} a_{\m_j}^{m_j}) \\ \times
(\prod_{j \in \s(\gamma)}
V_j(E -M_{\chi^j_{\Delta_j, X^0_j}} \mathtt R_{\Kappa^0_j +\kappa_j -\epsilon_j -\epsilon_j^\prime,
(\nu_{\Kappa^0 +\kappa -\epsilon -\epsilon^\prime}^{d, m,D,U}
(\n_{\epsilon +\epsilon^\prime}(\nu,\m^{\epsilon +\epsilon^\prime})))_j}^{1,l_j -1})) = \\
\sum_{ (\m^{\epsilon}, \m^{\epsilon^\prime}): \m^{\epsilon} \in \M_{\epsilon}^m(\nu),
\m^{\epsilon^\prime} \in \M_{\epsilon^\prime}^m(\nu)}
A_{\m^{\epsilon}}^m A_{\m^{\epsilon^\prime}}^m
(\prod_{j \in \s(\gamma)}
V_j(E -M_{\chi^j_{\Delta_j, X^0_j}} \mathtt R_{\Kappa^0_j +\kappa_j -\epsilon_j,
(\nu_{\Kappa^0 +\kappa -\epsilon}^{d, m,D,U}
(\n_{\epsilon}(\nu,\m^{\epsilon})))_j}^{1,l_j -1})) = \\
\sum_{\m^{\epsilon^\prime} \in \M_{\epsilon^\prime}^m(\nu)} A_{\m^{\epsilon^\prime}}^m
(\sum_{ \m^{\epsilon} \in \M_{\epsilon}^m(\nu)} A_{\m^{\epsilon}}^m
(\prod_{j \in \s(\gamma)}
V_j(E -M_{\chi^j_{\Delta_j, X^0_j}} \mathtt R_{\Kappa^0_j +\kappa_j -\epsilon_j,
(\nu_{\Kappa^0 +\kappa -\epsilon}^{d, m,D,U}
(\n_{\epsilon}(\nu,\m^{\epsilon})))_j}^{1,l_j -1}))) = \\
(\sum_{\m^{\epsilon^\prime} \in \M_{\epsilon^\prime}^m(\nu)} A_{\m^{\epsilon^\prime}}^m)
(\sum_{ \m^{\epsilon} \in \M_{\epsilon}^m(\nu)} A_{\m^{\epsilon}}^m
(\prod_{j \in \s(\gamma)}
V_j(E -M_{\chi^j_{\Delta_j, X^0_j}} \mathtt R_{\Kappa^0_j +\kappa_j -\epsilon_j,
(\nu_{\Kappa^0 +\kappa -\epsilon}^{d, m,D,U}
(\n_{\epsilon}(\nu,\m^{\epsilon})))_j}^{1,l_j -1}))) = \\
\sum_{ \m^{\epsilon} \in \M_{\epsilon}^m(\nu)} A_{\m^{\epsilon}}^m
(\prod_{j \in \s(\gamma)}
V_j(E -M_{\chi^j_{\Delta_j, X^0_j}} \mathtt R_{\Kappa^0_j +\kappa_j -\epsilon_j,
(\nu_{\Kappa^0 +\kappa -\epsilon}^{d, m,D,U}
(\n_{\epsilon}(\nu,\m^{\epsilon})))_j}^{1,l_j -1})).
\end{multline*}

Подставляя (2.2.29) в (2.2.27), приходим к равенству
\begin{multline*} \tag{2.2.30}
\sum_{\epsilon \in \Upsilon^d: \s(\epsilon) \subset \s(\kappa)} (-\e)^\epsilon
\sum_{\m^{\epsilon} \in \M_{\epsilon}^m(\nu)} A_{\m^{\epsilon}}^m \\
(\prod_{j \in \s(\gamma)}
V_j(E -M_{\chi^j_{\Delta_j, X^0_j}} \mathtt R_{\Kappa^0_j +\kappa_j -\epsilon_j,
(\nu_{\Kappa^0 +\kappa -\epsilon}^{d, m,D,U}
(\n_{\epsilon}(\nu,\m^{\epsilon})))_j}^{1,l_j -1})) = \\
\sum_{\substack{\epsilon \in \Upsilon^d:\\ \s(\epsilon) \subset
(\s(\kappa) \cap \s(\gamma))}} (-\e)^\epsilon
\sum_{\substack{\epsilon^\prime \in \Upsilon^d:\\
\s(\epsilon^\prime) \subset (\s(\kappa) \setminus \s(\gamma))}}
(-\e)^{\epsilon^\prime} (\sum_{ \m^{\epsilon} \in \M_{\epsilon}^m(\nu)} A_{\m^{\epsilon}}^m \\
(\prod_{j \in \s(\gamma)}
V_j(E -M_{\chi^j_{\Delta_j, X^0_j}} \mathtt R_{\Kappa^0_j +\kappa_j -\epsilon_j,
(\nu_{\Kappa^0 +\kappa -\epsilon}^{d, m,D,U}
(\n_{\epsilon}(\nu,\m^{\epsilon})))_j}^{1,l_j -1}))) = \\
\sum_{\substack{\epsilon^\prime \in \Upsilon^d:\\
\s(\epsilon^\prime) \subset (\s(\kappa) \setminus \s(\gamma))}} (-\e)^{\epsilon^\prime}
(\sum_{\substack{\epsilon \in \Upsilon^d:\\ \s(\epsilon) \subset
(\s(\kappa) \cap \s(\gamma))}} (-\e)^\epsilon
\sum_{ \m^{\epsilon} \in \M_{\epsilon}^m(\nu)} A_{\m^{\epsilon}}^m \\
(\prod_{j \in \s(\gamma)}
V_j(E -M_{\chi^j_{\Delta_j, X^0_j}} \mathtt R_{\Kappa^0_j +\kappa_j -\epsilon_j,
(\nu_{\Kappa^0 +\kappa -\epsilon}^{d, m,D,U}
(\n_{\epsilon}(\nu,\m^{\epsilon})))_j}^{1,l_j -1}))) = \\
(\sum_{\substack{\epsilon^\prime \in \Upsilon^d: \\
\s(\epsilon^\prime) \subset (\s(\kappa) \setminus \s(\gamma))}} (-\e)^{\epsilon^\prime})
(\sum_{\substack{\epsilon \in \Upsilon^d: \\ \s(\epsilon) \subset
(\s(\kappa) \cap \s(\gamma))}} (-\e)^\epsilon
\sum_{ \m^{\epsilon} \in \M_{\epsilon}^m(\nu)} A_{\m^{\epsilon}}^m \\
(\prod_{j \in \s(\gamma)}
V_j(E -M_{\chi^j_{\Delta_j, X^0_j}} \mathtt R_{\Kappa^0_j +\kappa_j -\epsilon_j,
(\nu_{\Kappa^0 +\kappa -\epsilon}^{d, m,D,U}
(\n_{\epsilon}(\nu,\m^{\epsilon})))_j}^{1,l_j -1}))) = 0,
\end{multline*}
при $ \kappa \in \Z_+^d, \nu \in N_{\Kappa^0 +\kappa}^{d,m,U},
\gamma \in \Upsilon^d: \s(\kappa) \setminus \s(\gamma) \ne \emptyset, $ ибо
в этом случае
$$
\sum_{\epsilon^\prime \in \Upsilon^d: \s(\epsilon^\prime) \subset
(\s(\kappa) \setminus \s(\gamma))} (-\e)^{\epsilon^\prime} =0.
$$

Соединяя (2.2.30) с (2.2.26), заключаем, что в условиях леммы
при $ \kappa \in \Z_+^d, \nu \in N_{\Kappa^0 +\kappa}^{d,m,U} $ для $ f \in
(S_p^\alpha H)^\prime(D) $ имеет место равенство
\begin{multline*} \tag{2.2.31}
(V_{\Kappa^0, \kappa,\nu}^{d,l -\e,m,D,U,\Nu} f) \mid_{Q^{d,m,D}} = \\
\biggl(\sum_{ \gamma \in \Upsilon^d: \s(\kappa) \setminus \s(\gamma) = \emptyset} (-\e)^\gamma
\sum_{\epsilon \in \Upsilon^d: \s(\epsilon) \subset \s(\kappa)} (-\e)^\epsilon
\sum_{\m^{\epsilon} \in \M_{\epsilon}^m(\nu)} A_{\m^{\epsilon}}^m \\
((\prod_{j \in \s(\gamma)}
V_j(E -M_{\chi^j_{\Delta_j, X^0_j}} \mathtt R_{\Kappa^0_j +\kappa_j -\epsilon_j,
(\nu_{\Kappa^0 +\kappa -\epsilon}^{d, m,D,U}
(\n_{\epsilon}(\nu,\m^{\epsilon})))_j}^{1,l_j -1})) F)\biggr) \mid_{Q^{d,m,D}} = \\
\biggl(\sum_{ \gamma \in \Upsilon^d: \s(\kappa) \subset \s(\gamma)} (-\e)^\gamma
\sum_{\epsilon \in \Upsilon^d: \s(\epsilon) \subset \s(\kappa)} (-\e)^\epsilon
\sum_{\m^{\epsilon} \in \M_{\epsilon}^m(\nu)} A_{\m^{\epsilon}}^m \\
((\prod_{j \in \s(\gamma)}
V_j(E -M_{\chi^j_{\Delta_j, X^0_j}} \mathtt R_{\Kappa^0_j +\kappa_j -\epsilon_j,
(\nu_{\Kappa^0 +\kappa -\epsilon}^{d, m,D,U}
(\n_{\epsilon}(\nu,\m^{\epsilon})))_j}^{1,l_j -1})) F)\biggr) \mid_{Q^{d,m,D}}.
\end{multline*}

Возвращаясь к оценке правой части (2.2.23), при $ \kappa \in \Z_+^d,
\nu \in N_{\Kappa^0 +\kappa}^{d,m,U} $ для $ f \in (S_p^\alpha H)^\prime(D), $
опираясь на (2.2.31), с учётом (2.2.10), (2.2.2) и (2.2.24), получаем
\begin{multline*} \tag{2.2.32}
\biggl\| V_{\Kappa^0, \kappa,\nu}^{d,l -\e,m,D,U,\Nu} f
\biggr\|_{L_\infty(Q_{\Kappa^0 +\kappa, n_{\Kappa^0 +\kappa}^{d,m,D,U}(\nu)}^d)} = \\
\biggl\| \sum_{ \gamma \in \Upsilon^d: \s(\kappa) \subset \s(\gamma)} (-\e)^\gamma
\sum_{\epsilon \in \Upsilon^d: \s(\epsilon) \subset \s(\kappa)} (-\e)^\epsilon
\sum_{\m^{\epsilon} \in \M_{\epsilon}^m(\nu)} A_{\m^{\epsilon}}^m \\
((\prod_{j \in \s(\gamma)}
V_j(E -M_{\chi^j_{\Delta_j, X^0_j}} \mathtt R_{\Kappa^0_j +\kappa_j -\epsilon_j,
(\nu_{\Kappa^0 +\kappa -\epsilon}^{d, m,D,U}
(\n_{\epsilon}(\nu,\m^{\epsilon})))_j}^{1,l_j -1})) F)
\biggr\|_{L_\infty(Q_{\Kappa^0 +\kappa, n_{\Kappa^0 +\kappa}^{d,m,D,U}(\nu)}^d)} \le \\
\sum_{ \gamma \in \Upsilon^d: \s(\kappa) \subset \s(\gamma)}
\sum_{\epsilon \in \Upsilon^d: \s(\epsilon) \subset \s(\kappa)}
\sum_{\m^{\epsilon} \in \M_{\epsilon}^m(\nu)} A_{\m^{\epsilon}}^m \\
\biggl\| (\prod_{j \in \s(\gamma)}
V_j(E -M_{\chi^j_{\Delta_j, X^0_j}} \mathtt R_{\Kappa^0_j +\kappa_j -\epsilon_j,
(\nu_{\Kappa^0 +\kappa -\epsilon}^{d, m,D,U}
(\n_{\epsilon}(\nu,\m^{\epsilon})))_j}^{1,l_j -1})) F
\biggr\|_{L_\infty(Q_{\Kappa^0 +\kappa, n_{\Kappa^0 +\kappa}^{d,m,D,U}(\nu)}^d)} \le \\
\sum_{ \gamma \in \Upsilon^d: \s(\kappa) \subset \s(\gamma)}
\sum_{\epsilon \in \Upsilon^d: \s(\epsilon) \subset \s(\kappa)}
\sum_{\m^{\epsilon} \in \M_{\epsilon}^m(\nu)} A_{\m^{\epsilon}}^m \\
\biggl\| (\prod_{j \in \s(\gamma)}
V_j(E -M_{\chi^j_{\Delta_j, X^0_j}} \mathtt R_{\Kappa^0_j +\kappa_j -\epsilon_j,
(\nu_{\Kappa^0 +\kappa -\epsilon}^{d, m,D,U}
(\n_{\epsilon}(\nu,\m^{\epsilon})))_j}^{1,l_j -1})) F
\biggr\|_{L_\infty(\mathcal D_{\Kappa^0 +\kappa,\nu,\epsilon,\m^{\epsilon}}^{d,m,D,U})}.
\end{multline*}
Оценивая правую часть (2.2.32), при $ \kappa \in \Z_+^d, \nu \in
N_{\Kappa^0 +\kappa}^{d,m,U}, \gamma \in \Upsilon^d: \s(\kappa) \subset \s(\gamma),
\epsilon \in \Upsilon^d: \s(\epsilon) \subset \s(\kappa),
\m^{\epsilon} \in \M_{\epsilon}^m(\nu),  $ ввиду (2.1.48) используя (1.3.4), а 
затем с учётом (2.2.10), (2.2.2), (2.2.24), (2.2.9) применяя (2.1.17), 
заключаем, что
\begin{multline*} \tag{2.2.33}
\biggl\| (\prod_{j \in \s(\gamma)}
V_j(E -M_{\chi^j_{\Delta_j, X^0_j}} \mathtt R_{\Kappa^0_j +\kappa_j -\epsilon_j,
(\nu_{\Kappa^0 +\kappa -\epsilon}^{d, m,D,U}
(\n_{\epsilon}(\nu,\m^{\epsilon})))_j}^{1,l_j -1})) F
\biggr\|_{L_\infty(\mathcal D_{\Kappa^0 +\kappa,\nu,\epsilon,\m^{\epsilon}}^{d,m,D,U})} = \\
\biggl\| (\prod_{j \in \s(\gamma) \setminus \s(\kappa)}
V_j(E -M_{\chi^j_{\Delta_j, X^0_j}} \mathtt R_{\Kappa^0_j +\kappa_j -\epsilon_j,
(\nu_{\Kappa^0 +\kappa -\epsilon}^{d, m,D,U}
(\n_{\epsilon}(\nu,\m^{\epsilon})))_j}^{1,l_j -1})) \\
((\prod_{j \in \s(\kappa)}
V_j(E -M_{\chi^j_{\Delta_j, X^0_j}} \mathtt R_{\Kappa^0_j +\kappa_j -\epsilon_j,
(\nu_{\Kappa^0 +\kappa -\epsilon}^{d, m,D,U}
(\n_{\epsilon}(\nu,\m^{\epsilon})))_j}^{1,l_j -1})) F)
\biggr\|_{L_\infty(\mathcal D_{\Kappa^0 +\kappa,\nu,\epsilon,\m^{\epsilon}}^{d,m,D,U})} \le \\
c_8 \biggl\| (\prod_{j \in \s(\kappa)}
V_j(E -M_{\chi^j_{\Delta_j, X^0_j}} \mathtt R_{\Kappa^0_j +\kappa_j -\epsilon_j,
(\nu_{\Kappa^0 +\kappa -\epsilon}^{d, m,D,U}
(\n_{\epsilon}(\nu,\m^{\epsilon})))_j}^{1,l_j -1})) F
\biggr\|_{L_\infty(\mathcal D_{\Kappa^0 +\kappa,\nu,\epsilon,\m^{\epsilon}}^{d,m,D,U})}.
\end{multline*}

Для проведения оценки правой части (2.2.33) при $ \kappa \in \Z_+^d, \nu
\in N_{\Kappa^0 +\kappa}^{d,m,U}, \epsilon \in \Upsilon^d:
\s(\epsilon) \subset \s(\kappa), \m^{\epsilon} \in \M_{\epsilon}^m(\nu),
j =1,\ldots,d $ обозначим через
$ \mathcal S_{\Kappa^0 +\kappa,\nu,\epsilon,\m^{\epsilon}}^{j,d,l -\e,m,D,U}:
L_1((\bm x_{\Kappa^0 +\kappa,\nu,\epsilon,\m^{\epsilon}}^{d,m,D,U})_j +
(\bm \delta_{\Kappa^0 +\kappa,\nu,\epsilon,\m^{\epsilon}}^{d,m,D,U})_j I) \mapsto
\mathcal P^{1, l_j -1}, $ линейный оператор, определяемый равенством
$$
\mathcal S_{\Kappa^0 +\kappa,\nu,\epsilon,\m^{\epsilon}}^{j,d,l -\e,m,D,U} =
P_{\delta, x^0}^{1, l_j -1}
$$
при $ \delta = (\bm \delta_{\Kappa^0 +\kappa,\nu,\epsilon,\m^{\epsilon}}^{d,m,D,U})_j,
x^0 = (\bm x_{\Kappa^0 +\kappa,\nu,\epsilon,\m^{\epsilon}}^{d,m,D,U})_j $
(см. лемму 1.1.2).

Далее, пользуясь тем, что при $ j =1,\ldots,d $ ввиду (2.2.10), (2.2.2),
(2.2.24), (2.1.13) в $ \mathcal B(C_0(\R), C_0(\R)) $ соблюдено равенство
\begin{multline*}
(E -M_{\chi^j_{\Delta_j, X^0_j}} \mathtt R_{\Kappa^0_j +\kappa_j -\epsilon_j,
(\nu_{\Kappa^0 +\kappa -\epsilon}^{d, m,D,U}
(\n_{\epsilon}(\nu,\m^{\epsilon})))_j}^{1,l_j -1})
(E -M_{\chi^j_{\Delta_j, X^0_j}}
\mathcal S_{\Kappa^0 +\kappa,\nu,\epsilon,\m^{\epsilon}}^{j,d,l -\e,m,D,U}) = \\
E -M_{\chi^j_{\Delta_j, X^0_j}} \mathtt R_{\Kappa^0_j +\kappa_j -\epsilon_j,
(\nu_{\Kappa^0 +\kappa -\epsilon}^{d, m,D,U}
(\n_{\epsilon}(\nu,\m^{\epsilon})))_j}^{1,l_j -1}
-M_{\chi^j_{\Delta_j, X^0_j}}
\mathcal S_{\Kappa^0 +\kappa,\nu,\epsilon,\m^{\epsilon}}^{j,d,l -\e,m,D,U} +\\
M_{\chi^j_{\Delta_j, X^0_j}} \mathtt R_{\Kappa^0_j +\kappa_j -\epsilon_j,
(\nu_{\Kappa^0 +\kappa -\epsilon}^{d, m,D,U}
(\n_{\epsilon}(\nu,\m^{\epsilon})))_j}^{1,l_j -1}
(M_{\chi^j_{\Delta_j, X^0_j}}
\mathcal S_{\Kappa^0 +\kappa,\nu,\epsilon,\m^{\epsilon}}^{j,d,l -\e,m,D,U}) = \\
E -M_{\chi^j_{\Delta_j, X^0_j}} \mathtt R_{\Kappa^0_j +\kappa_j -\epsilon_j,
(\nu_{\Kappa^0 +\kappa -\epsilon}^{d, m,D,U}
(\n_{\epsilon}(\nu,\m^{\epsilon})))_j}^{1,l_j -1}
-M_{\chi^j_{\Delta_j, X^0_j}}
\mathcal S_{\Kappa^0 +\kappa,\nu,\epsilon,\m^{\epsilon}}^{j,d,l -\e,m,D,U} +\\
M_{\chi^j_{\Delta_j, X^0_j}} \mathtt R_{\Kappa^0_j +\kappa_j -\epsilon_j,
(\nu_{\Kappa^0 +\kappa -\epsilon}^{d, m,D,U}
(\n_{\epsilon}(\nu,\m^{\epsilon})))_j}^{1,l_j -1}
(\mathcal S_{\Kappa^0 +\kappa,\nu,\epsilon,\m^{\epsilon}}^{j,d,l -\e,m,D,U}) = \\
E -M_{\chi^j_{\Delta_j, X^0_j}} \mathtt R_{\Kappa^0_j +\kappa_j -\epsilon_j,
(\nu_{\Kappa^0 +\kappa -\epsilon}^{d, m,D,U}
(\n_{\epsilon}(\nu,\m^{\epsilon})))_j}^{1,l_j -1}
-M_{\chi^j_{\Delta_j, X^0_j}}
\mathcal S_{\Kappa^0 +\kappa,\nu,\epsilon,\m^{\epsilon}}^{j,d,l -\e,m,D,U} +\\
M_{\chi^j_{\Delta_j, X^0_j}}
\mathcal S_{\Kappa^0 +\kappa,\nu,\epsilon,\m^{\epsilon}}^{j,d,l -\e,m,D,U} =
E -M_{\chi^j_{\Delta_j, X^0_j}} \mathtt R_{\Kappa^0_j +\kappa_j -\epsilon_j,
(\nu_{\Kappa^0 +\kappa -\epsilon}^{d, m,D,U}
(\n_{\epsilon}(\nu,\m^{\epsilon})))_j}^{1,l_j -1},
\end{multline*}
благодаря п. 2) леммы 1.3.2 и (1.3.4), а также в силу (2.2.10), (2.2.2),
(2.2.24), (2.2.9) используя (2.1.17), выводим
\begin{multline*} \tag{2.2.34}
\biggl\| (\prod_{j \in \s(\kappa)}
V_j(E -M_{\chi^j_{\Delta_j, X^0_j}} \mathtt R_{\Kappa^0_j +\kappa_j -\epsilon_j,
(\nu_{\Kappa^0 +\kappa -\epsilon}^{d, m,D,U}
(\n_{\epsilon}(\nu,\m^{\epsilon})))_j}^{1,l_j -1})) F
\biggr\|_{L_\infty(\mathcal D_{\Kappa^0 +\kappa,\nu,\epsilon,\m^{\epsilon}}^{d,m,D,U})} = \\
\biggl\| (\prod_{j \in \s(\kappa)}
V_j((E -M_{\chi^j_{\Delta_j, X^0_j}} \mathtt R_{\Kappa^0_j +\kappa_j -\epsilon_j,
(\nu_{\Kappa^0 +\kappa -\epsilon}^{d, m,D,U}
(\n_{\epsilon}(\nu,\m^{\epsilon})))_j}^{1,l_j -1}) \\
(E -M_{\chi^j_{\Delta_j, X^0_j}}
\mathcal S_{\Kappa^0 +\kappa,\nu,\epsilon,\m^{\epsilon}}^{j,d,l -\e,m,D,U}))) F
\biggr\|_{L_\infty(\mathcal D_{\Kappa^0 +\kappa,\nu,\epsilon,\m^{\epsilon}}^{d,m,D,U})} = \\
\biggl\| (\prod_{j \in \s(\kappa)}
V_j(E -M_{\chi^j_{\Delta_j, X^0_j}} \mathtt R_{\Kappa^0_j +\kappa_j -\epsilon_j,
(\nu_{\Kappa^0 +\kappa -\epsilon}^{d, m,D,U}
(\n_{\epsilon}(\nu,\m^{\epsilon})))_j}^{1,l_j -1})) \\
((\prod_{j \in \s(\kappa)}
V_j(E -M_{\chi^j_{\Delta_j, X^0_j}}
\mathcal S_{\Kappa^0 +\kappa,\nu,\epsilon,\m^{\epsilon}}^{j,d,l -\e,m,D,U})) F)
\biggr\|_{L_\infty(\mathcal D_{\Kappa^0 +\kappa,\nu,\epsilon,\m^{\epsilon}}^{d,m,D,U})} \le \\
c_9 \biggl\| (\prod_{j \in \s(\kappa)} V_j(E -M_{\chi^j_{\Delta_j, X^0_j}}
\mathcal S_{\Kappa^0 +\kappa,\nu,\epsilon,\m^{\epsilon}}^{j,d,l -\e,m,D,U})) F
\biggr\|_{L_\infty(\mathcal D_{\Kappa^0 +\kappa,\nu,\epsilon,\m^{\epsilon}}^{d,m,D,U})}, \\
\kappa \in \Z_+^d, \nu \in N_{\Kappa^0 +\kappa}^{d,m,U}, \epsilon \in
\Upsilon^d: \s(\epsilon) \subset \s(\kappa), \m^{\epsilon} \in
\M_{\epsilon}^m(\nu).
\end{multline*}

Соединяя (2.2.32), (2.2.33), (2.2.34), находим, что при $ \kappa \in \Z_+^d,
\nu \in N_{\Kappa^0 +\kappa}^{d,m,U} $ выполняется неравенство
\begin{multline*} \tag{2.2.35}
\biggl\| V_{\Kappa^0, \kappa,\nu}^{d,l -\e,m,D,U,\Nu} f
\biggr\|_{L_\infty(Q_{\Kappa^0 +\kappa, n_{\Kappa^0 +\kappa}^{d,m,D,U}(\nu)}^d)} \le \\
\sum_{ \gamma \in \Upsilon^d: \s(\kappa) \subset \s(\gamma)}
\sum_{\epsilon \in \Upsilon^d: \s(\epsilon) \subset \s(\kappa)}
\sum_{\m^{\epsilon} \in \M_{\epsilon}^m(\nu)} A_{\m^{\epsilon}}^m \\
c_8 c_9 \biggl\| (\prod_{j \in \s(\kappa)} V_j(E -M_{\chi^j_{\Delta_j, X^0_j}}
\mathcal S_{\Kappa^0 +\kappa,\nu,\epsilon,\m^{\epsilon}}^{j,d,l -\e,m,D,U})) F
\biggr\|_{L_\infty(\mathcal D_{\Kappa^0 +\kappa,\nu,\epsilon,\m^{\epsilon}}^{d,m,D,U})} \le \\
c_{10} \sum_{\epsilon \in \Upsilon^d: \s(\epsilon) \subset \s(\kappa)}
\sum_{\m^{\epsilon} \in \M_{\epsilon}^m(\nu)} A_{\m^{\epsilon}}^m \\
\biggl\| (\prod_{j \in \s(\kappa)} V_j(E -M_{\chi^j_{\Delta_j, X^0_j}}
\mathcal S_{\Kappa^0 +\kappa,\nu,\epsilon,\m^{\epsilon}}^{j,d,l -\e,m,D,U})) F
\biggr\|_{L_\infty(\mathcal D_{\Kappa^0 +\kappa,\nu,\epsilon,\m^{\epsilon}}^{d,m,D,U})}.
\end{multline*}

Для оценки слагаемых в правой части (2.2.35), применяя с учётом (2.2.2),
(2.2.24), (1.3.5), (2.1.46) неравенство (2.1.8) при $ \lambda = 0, q = \infty, $ а затем
используя (2.2.9), (1.3.5), (2.1.46) (см. (2.2.2)), приходим к неравенству
\begin{multline*} \tag{2.2.36}
\biggl\| (\prod_{j \in \s(\kappa)} V_j(E -M_{\chi^j_{\Delta_j, X^0_j}}
\mathcal S_{\Kappa^0 +\kappa,\nu,\epsilon,\m^{\epsilon}}^{j,d,l -\e,m,D,U})) F
\biggr\|_{L_\infty(\mathcal D_{\Kappa^0 +\kappa,\nu,\epsilon,\m^{\epsilon}}^{d,m,D,U})} \le \\
c_{11} (\bm \delta_{\Kappa^0 +\kappa,\nu,\epsilon,
\m^{\epsilon}}^{d,m,D,U})^{-p^{-1} \e}
\biggl(\| (\prod_{j \in \s(\kappa)} V_j(E -M_{\chi^j_{\Delta_j, X^0_j}}
\mathcal S_{\Kappa^0 +\kappa,\nu,\epsilon,\m^{\epsilon}}^{j,d,l -\e,m,D,U})) F
\|_{L_p(\mathcal D_{\Kappa^0 +\kappa,\nu,\epsilon,\m^{\epsilon}}^{d,m,D,U})} +\\
\sum_{ J \subset \Nu_{1,d}^1: J \ne \emptyset}
(\prod_{j \in J} (\bm \delta_{\Kappa^0 +\kappa,\nu,\epsilon,
\m^{\epsilon}}^{d,m,D,U})_j^{p^{-1}})
\int_{ (c_{12} \bm \delta_{\Kappa^0 +\kappa,\nu,\epsilon,
\m^{\epsilon}}^{d,m,D,U} I^d)^J} (\prod_{j \in J} t_j^{-2 p^{-1} -1})\times\\
\biggl(\int\limits_{(t B^d)^J}
\int\limits_{ (\mathcal D_{\Kappa^0 +\kappa,\nu,\epsilon,\m^{\epsilon}}^{d,m,D,U})_\xi^{l \chi_J}}
|\Delta_\xi^{l \chi_J} ((\prod_{j \in \s(\kappa)} V_j(E -M_{\chi^j_{\Delta_j, X^0_j}}
\mathcal S_{\Kappa^0 +\kappa,\nu,\epsilon,\m^{\epsilon}}^{j,d,l -\e,m,D,U}))
F)(x)|^p dx d\xi^J \biggr)^{1/p} dt^J \biggr) \le \\
c_{13} (2^{-\kappa})^{-p^{-1} \e}
\biggl(c_{14} \biggl(\prod_{j \in \s(\kappa)} (\bm \delta_{\Kappa^0 +\kappa,
\nu,\epsilon,\m^{\epsilon}}^{d,m,D,U})_j^{-1 /p} \biggr)\times\\
\biggl(\int_{(\bm \delta_{\Kappa^0 +\kappa,\nu,\epsilon,
\m^{\epsilon}}^{d,m,D,U} B^d)^{\s(\kappa)}}
\int_{(\mathcal D_{\Kappa^0 +\kappa,\nu,\epsilon,\m^{\epsilon}}^{d,m,D,U})_\xi^{l
\chi_{\s(\kappa)}}} | (\Delta_\xi^{l \chi_{\s(\kappa)}} F)(x)|^p dx
d\xi^{\s(\kappa)} \biggr)^{1 /p} +\\
c_{15} \sum_{ J \subset \Nu_{1,d}^1: J \ne \emptyset}
(\prod_{j \in J} 2^{-\kappa_j /p})
\int_{ (c_{16} 2^{-\kappa} I^d)^J} (\prod_{j \in J} t_j^{-2 p^{-1} -1})
\biggl(\int_{(t B^d)^J}
c_{17} \biggl(\prod_{j \in \s(\kappa) \setminus J}
(\bm \delta_{\Kappa^0 +\kappa,\nu,\epsilon,\m^{\epsilon}}^{d,m,D,U})_j^{-1} \biggr) \times\\
\int_{(\bm \delta_{\Kappa^0 +\kappa,\nu,\epsilon,
\m^{\epsilon}}^{d,m,D,U} B^d)^{\s(\kappa) \setminus J}}
\int_{ (\mathcal D_{\Kappa^0 +\kappa,\nu,\epsilon,
\m^{\epsilon}}^{d,m,D,U})_\xi^{l \chi_{J \cup \s(\kappa)}}}
|(\Delta_\xi^{l \chi_{J \cup \s(\kappa)}} F)(x)|^p dx d\xi^{\s(\kappa) \setminus J}
d\xi^J \biggr)^{1/p} dt^J \biggr) \le \\
c_{18} 2^{(\kappa, p^{-1} \e)}
\biggl( \biggl(\prod_{j \in \s(\kappa)} 2^{\kappa_j /p} \biggr)
\biggl(\int_{(c_{2} 2^{-\kappa} B^d)^{\s(\kappa)}}
\int_{(\mathcal D_{\Kappa^0 +\kappa,\nu,\epsilon,\m^{\epsilon}}^{d,m,D,U})_\xi^{l
\chi_{\s(\kappa)}}} | (\Delta_\xi^{l \chi_{\s(\kappa)}} f)(x)|^p dx
d\xi^{\s(\kappa)} \biggr)^{1 /p} +\\
\sum_{ J \subset \Nu_{1,d}^1: J \ne \emptyset} (\prod_{j \in J} 2^{-\kappa_j /p})
\biggl(\prod_{j \in \s(\kappa) \setminus J} 2^{\kappa_j /p} \biggr)
\int_{ (c_{3} 2^{-\kappa} I^d)^J} (\prod_{j \in J} t_j^{-2 p^{-1} -1})\times\\
\biggl( \int_{(t B^d)^J}
\int_{(c_{2} 2^{-\kappa} B^d)^{\s(\kappa) \setminus J}}
\int_{ (\mathcal D_{\Kappa^0 +\kappa,\nu,\epsilon,
\m^{\epsilon}}^{d,m,D,U})_\xi^{l \chi_{J \cup \s(\kappa)}}}
| (\Delta_\xi^{l \chi_{J \cup \s(\kappa)}} f)(x)|^p dx d\xi^{\s(\kappa) \setminus J}
d\xi^J \biggr)^{1/p} dt^J \biggr), \\
\kappa \in \Z_+^d, \nu \in N_{\Kappa^0 +\kappa}^{d,m,U}, \epsilon \in
\Upsilon^d: \s(\epsilon) \subset \s(\kappa), \m^{\epsilon} \in
\M_{\epsilon}^m(\nu).
\end{multline*}

Теперь при $ n \in \Z^d: Q_{\Kappa^0 +\kappa,n}^d \cap
G_{\Kappa^0 +\kappa}^{d,m,U} \ne \emptyset,
\nu \in N_{\Kappa^0 +\kappa}^{d,m,U}:  Q_{\Kappa^0 +\kappa,n}^d \cap
\supp g_{\Kappa^0 +\kappa, \nu}^{d,m} \ne \emptyset, $ из (2.2.35), (2.2.36)
ввиду (2.2.20), (1.4.21) получаем
\begin{multline*} \tag{2.2.37}
\biggl\| V_{\Kappa^0, \kappa,\nu}^{d,l -\e,m,D,U,\Nu} f
\biggr\|_{L_\infty(Q_{\Kappa^0 +\kappa, n_{\Kappa^0 +\kappa}^{d,m,D,U}(\nu)}^d)} \le
c_{10} \sum_{\epsilon \in \Upsilon^d: \s(\epsilon) \subset \s(\kappa)}
\sum_{\m^{\epsilon} \in \M_{\epsilon}^m(\nu)} A_{\m^{\epsilon}}^m
c_{18} 2^{(\kappa, p^{-1} \e)} \times\\
\biggl( \biggl(\prod_{j \in \s(\kappa)} 2^{\kappa_j /p} \biggr)
\biggl(\int_{(c_{2} 2^{-\kappa} B^d)^{\s(\kappa)}}
\int_{(\mathcal D_{\Kappa^0 +\kappa,\nu,\epsilon,\m^{\epsilon}}^{d,m,D,U})_\xi^{l
\chi_{\s(\kappa)}}} | (\Delta_\xi^{l \chi_{\s(\kappa)}} f)(x)|^p dx
d\xi^{\s(\kappa)} \biggr)^{1 /p} +\\
\sum_{ J \subset \Nu_{1,d}^1: J \ne \emptyset}
(\prod_{j \in J} 2^{-\kappa_j /p})
\biggl(\prod_{j \in \s(\kappa) \setminus J} 2^{\kappa_j /p} \biggr)
\int_{ (c_{3} 2^{-\kappa} I^d)^J} (\prod_{j \in J} t_j^{-2 p^{-1} -1})\times\\
\biggl( \int_{(t B^d)^J}
\int_{(c_{2} 2^{-\kappa} B^d)^{\s(\kappa) \setminus J}}
\int_{ (\mathcal D_{\Kappa^0 +\kappa,\nu,\epsilon,
\m^{\epsilon}}^{d,m,D,U})_\xi^{l \chi_{J \cup \s(\kappa)}}}
| (\Delta_\xi^{l \chi_{J \cup \s(\kappa)}} f)(x)|^p dx d\xi^{\s(\kappa) \setminus J}
d\xi^J \biggr)^{1/p} dt^J \biggr) \le \\
c_{19} \sum_{\epsilon \in \Upsilon^d: \s(\epsilon) \subset \s(\kappa)}
\sum_{\m^{\epsilon} \in \M_{\epsilon}^m(\nu)} A_{\m^{\epsilon}}^m
2^{(\kappa, p^{-1} \e)}\times\\
\biggl( \biggl(\prod_{j \in \s(\kappa)} 2^{\kappa_j /p} \biggr)
\biggl(\int_{(c_{2} 2^{-\kappa} B^d)^{\s(\kappa)}}
\int_{(D \cap D_{\Kappa^0 +\kappa,n}^{\prime d,m,D,U})_\xi^{l \chi_{\s(\kappa)}}}
| (\Delta_\xi^{l \chi_{\s(\kappa)}} f)(x)|^p dx
d\xi^{\s(\kappa)} \biggr)^{1 /p}+\\
\sum_{ J \subset \Nu_{1,d}^1: J \ne \emptyset}
(\prod_{j \in J} 2^{-\kappa_j /p})
\biggl(\prod_{j \in \s(\kappa) \setminus J} 2^{\kappa_j /p} \biggr)
\int_{ (c_{3} 2^{-\kappa} I^d)^J} (\prod_{j \in J} t_j^{-2 p^{-1} -1})\times\\
\biggl( \int_{(t B^d)^J}
\int_{(c_{2} 2^{-\kappa} B^d)^{\s(\kappa) \setminus J}}
\int_{(D \cap D_{\Kappa^0 +\kappa,n}^{\prime d,m,D,U})_\xi^{l \chi_{J \cup \s(\kappa)}}}
| (\Delta_\xi^{l \chi_{J \cup \s(\kappa)}} f)(x)|^p dx d\xi^{\s(\kappa) \setminus J}
d\xi^J \biggr)^{1/p} dt^J \biggr) \le \\
c_{20} 2^{(\kappa, p^{-1} \e)}
\biggl( \biggl(\prod_{j \in \s(\kappa)} 2^{\kappa_j /p} \biggr)
\biggl(\int_{(c_{2} 2^{-\kappa} B^d)^{\s(\kappa)}}
\int_{(D \cap D_{\Kappa^0 +\kappa,n}^{\prime d,m,D,U})_\xi^{l \chi_{\s(\kappa)}}}
| (\Delta_\xi^{l \chi_{\s(\kappa)}} f)(x)|^p dx
d\xi^{\s(\kappa)} \biggr)^{1 /p}+\\
\sum_{ J \subset \Nu_{1,d}^1: J \ne \emptyset}
(\prod_{j \in J} 2^{-\kappa_j /p})
\biggl(\prod_{j \in \s(\kappa) \setminus J} 2^{\kappa_j /p} \biggr)
\int_{ (c_{3} 2^{-\kappa} I^d)^J} (\prod_{j \in J} t_j^{-2 p^{-1} -1})\times\\
\biggl( \int_{(t B^d)^J}
\int_{(c_{2} 2^{-\kappa} B^d)^{\s(\kappa) \setminus J}}
\int_{(D \cap D_{\Kappa^0 +\kappa,n}^{\prime d,m,D,U})_\xi^{l \chi_{J \cup \s(\kappa)}}}
| (\Delta_\xi^{l \chi_{J \cup \s(\kappa)}} f)(x)|^p dx
d\xi^{\s(\kappa) \setminus J} d\xi^J \biggr)^{1/p} dt^J \biggr) \le \\
c_{20} 2^{(\kappa, p^{-1} \e)}
\biggl( \biggl(\prod_{j \in \s(\kappa)} 2^{\kappa_j /p} \biggr)
\biggl(\int_{(c_{2} 2^{-\kappa} B^d)^{\s(\kappa)}}
\int_{D_\xi^{l \chi_{\s(\kappa)}} \cap D_{\Kappa^0 +\kappa,n}^{\prime d,m,D,U}}
|(\Delta_\xi^{l \chi_{\s(\kappa)}} f)(x)|^p dx d\xi^{\s(\kappa)} \biggr)^{1 /p}+\\
\sum_{ J \subset \Nu_{1,d}^1: J \ne \emptyset}
(\prod_{j \in J} 2^{-\kappa_j /p})
\biggl(\prod_{j \in \s(\kappa) \setminus J} 2^{\kappa_j /p} \biggr)
\int_{ (c_{3} 2^{-\kappa} I^d)^J} (\prod_{j \in J} t_j^{-2 p^{-1} -1})\times\\
\biggl( \int_{(t B^d)^J}
\int_{(c_{2} 2^{-\kappa} B^d)^{\s(\kappa) \setminus J}}
\int_{D_\xi^{l \chi_{J \cup \s(\kappa)}} \cap D_{\Kappa^0 +\kappa,n}^{\prime d,m,D,U}}
|(\Delta_\xi^{l \chi_{J \cup \s(\kappa)}} f)(x)|^p dx
d\xi^{\s(\kappa) \setminus J} d\xi^J \biggr)^{1/p} dt^J \biggr).
\end{multline*}

Фиксируя $ \bm \epsilon \in \R_+^d $ так, чтобы соблюдалось условие
$ \alpha -p^{-1} \e -\bm \epsilon >0, $ и применяя неравенство Гёльдера,
при $ \kappa \in \Z_+^d, $ для $ J \subset \Nu_{1,d}^1: J \ne \emptyset, n \in \Z^d:
Q_{\Kappa^0 +\kappa,n}^d \cap G_{\Kappa^0 +\kappa}^{d,m,U} \ne \emptyset, $
имеем
\begin{multline*} \tag{2.2.38}
\int_{ (c_{3} 2^{-\kappa} I^d)^J} (\prod_{j \in J} t_j^{-2p^{-1} -1}) \times\\
\biggl(\int_{ (t B^d)^J}
\int_{ (c_{2} 2^{-\kappa} B^d)^{\s(\kappa) \setminus J}}
\int_{D_\xi^{l \chi_{J \cup \s(\kappa)}} \cap D_{\Kappa^0 +\kappa,n}^{\prime d,m,D,U}}
|(\Delta_\xi^{l \chi_{J \cup \s(\kappa)}} f)(x)|^p
dx d\xi^{\s(\kappa) \setminus J} d\xi^J \biggr)^{1/p} dt^J =\\
\int_{ (c_{3} 2^{-\kappa} I^d)^J} (\prod_{j \in J}
t_j^{\alpha_j -p^{-1} -\bm \epsilon_j -1/p^\prime})
(\prod_{j \in J} t_j^{-(\alpha_j -\bm \epsilon_j +p^{-1}) -1/p}) \times\\
\biggl(\int_{ (t B^d)^J}
\int_{ ( c_{2} 2^{-\kappa} B^d)^{\s(\kappa) \setminus J}}
\int_{D_\xi^{l \chi_{J \cup \s(\kappa)}} \cap D_{\Kappa^0 +\kappa,n}^{\prime d,m,D,U}}
| (\Delta_\xi^{l \chi_{J \cup \s(\kappa)}} f)(x)|^p
dx d\xi^{\s(\kappa) \setminus J} d\xi^J \biggr)^{1/p} dt^J \le\\
\biggl(\int_{ (c_{3} 2^{-\kappa} I^d)^J}
(\prod_{j \in J} t_j^{p^\prime (\alpha_j -p^{-1} -\bm \epsilon_j) -1})
dt^J \biggr)^{1 /p^\prime} \cdot \biggl(\int_{ (c_{3} 2^{-\kappa} I^d)^J}
(\prod_{j \in J} t_j^{-p(\alpha_j -\bm \epsilon_j +p^{-1}) -1}) \times\\
\int_{ (t B^d)^J} \int_{ ( c_{2} 2^{-\kappa} B^d)^{\s(\kappa) \setminus J}}
\int_{D_\xi^{l \chi_{J \cup \s(\kappa)}} \cap D_{\Kappa^0 +\kappa,n}^{\prime d,m,D,U}}
| (\Delta_\xi^{l \chi_{J \cup \s(\kappa)}} f)(x)|^p
dx d\xi^{\s(\kappa) \setminus J} d\xi^J dt^J \biggr)^{1/p} \le\\
c_{21} (\prod_{j \in J} 2^{-\kappa_j (\alpha_j -p^{-1} -\bm \epsilon_j)})
\biggl(\int_{ (c_{3} 2^{-\kappa} I^d)^J}
(\prod_{j \in J} t_j^{-p(\alpha_j -\bm \epsilon_j +p^{-1}) -1}) \times\\
\int_{ (t B^d)^J} \int_{ ( c_{2} 2^{-\kappa} B^d)^{\s(\kappa) \setminus J}}
\int_{D_\xi^{l \chi_{J \cup \s(\kappa)}} \cap D_{\Kappa^0 +\kappa,n}^{\prime d,m,D,U}}
| (\Delta_\xi^{l \chi_{J \cup \s(\kappa)}} f)(x)|^p
dx d\xi^{\s(\kappa) \setminus J} d\xi^J dt^J \biggr)^{1/p}.
\end{multline*}

Соединяя (2.2.23), (2.2.37), (2.2.38), заключаем, что при
$ n \in \Z^d: Q_{\Kappa^0 +\kappa,n}^d \cap G_{\Kappa^0 +\kappa}^{d,m,U} \ne \emptyset,
\nu \in N_{\Kappa^0 +\kappa}^{d,m,U}: Q_{\Kappa^0 +\kappa,n}^d \cap
\supp g_{\Kappa^0 +\kappa, \nu}^{d,m} \ne \emptyset, \mu \in \Z_+^d(\lambda) $
справедливо неравенство
\begin{multline*} \tag{2.2.39}
\biggl\| \D^\mu (V_{\Kappa^0, \kappa,\nu}^{d,l -\e,m,D,U,\Nu} f)
\D^{\lambda -\mu} g_{\Kappa^0 +\kappa, \nu}^{d,m} \biggr\|_{L_q(Q_{\Kappa^0 +\kappa,n}^d)} \le \\
c_7 2^{(\kappa, \lambda -q^{-1} \e)}
c_{20} 2^{(\kappa, p^{-1} \e)}
\biggl( \biggl(\prod_{j \in \s(\kappa)} 2^{\kappa_j /p} \biggr)\times\\
\biggl(\int_{(c_{2} 2^{-\kappa} B^d)^{\s(\kappa)}}
\int_{D_\xi^{l \chi_{\s(\kappa)}} \cap D_{\Kappa^0 +\kappa,n}^{\prime d,m,D,U}}
| (\Delta_\xi^{l \chi_{\s(\kappa)}} f)(x)|^p dx
d\xi^{\s(\kappa)} \biggr)^{1 /p}+\\
\sum_{ J \subset \Nu_{1,d}^1: J \ne \emptyset}
(\prod_{j \in J} 2^{-\kappa_j /p})
\biggl(\prod_{j \in \s(\kappa) \setminus J} 2^{\kappa_j /p} \biggr) \\ \times
c_{21} (\prod_{j \in J} 2^{-\kappa_j (\alpha_j -p^{-1} -\bm \epsilon_j)})
\biggl(\int_{ (c_{3} 2^{-\kappa} I^d)^J}
(\prod_{j \in J} t_j^{-p(\alpha_j -\bm \epsilon_j +p^{-1}) -1}) \times\\
\int_{ (t B^d)^J} \int_{ ( c_{2} 2^{-\kappa} B^d)^{\s(\kappa) \setminus J}}
\int_{D_\xi^{l \chi_{J \cup \s(\kappa)}} \cap D_{\Kappa^0 +\kappa,n}^{\prime d,m,D,U}}
| (\Delta_\xi^{l \chi_{J \cup \s(\kappa)}} f)(x)|^p
dx d\xi^{\s(\kappa) \setminus J} d\xi^J dt^J \biggr)^{1/p} \biggr) \le \\
c_{22} 2^{(\kappa, \lambda +p^{-1} \e -q^{-1} \e)} \times\\
\biggl( (\prod_{j \in \s(\kappa)} 2^{\kappa_j /p})
\biggl(\int_{ (c_{2} 2^{-\kappa} B^d)^{\s(\kappa)}}
\int_{D_\xi^{l \chi_{\s(\kappa)}} \cap D_{\Kappa^0 +\kappa,n}^{\prime d,m,D,U}}
| (\Delta_\xi^{l \chi_{\s(\kappa)}} f)(x)|^p
dx d\xi^{\s(\kappa)} \biggr)^{1/p} +\\
\sum_{J \subset \Nu_{1,d}^1: J \ne \emptyset}
(\prod_{j \in J} 2^{-\kappa_j (\alpha_j -\bm \epsilon_j)})
(\prod_{j \in \s(\kappa) \setminus J} 2^{\kappa_j /p})
\biggl(\int_{ (c_{3} 2^{-\kappa} I^d)^J}
(\prod_{j \in J} t_j^{-p(\alpha_j -\bm \epsilon_j +p^{-1}) -1}) \times\\
\int_{ (t B^d)^J} \int_{ ( c_{2} 2^{-\kappa} B^d)^{\s(\kappa) \setminus J}}
\int_{D_\xi^{l \chi_{J \cup \s(\kappa)}} \cap D_{\Kappa^0 +\kappa,n}^{\prime d,m,D,U}}
| (\Delta_\xi^{l \chi_{J \cup \s(\kappa)}} f)(x)|^p
dx d\xi^{\s(\kappa) \setminus J} d\xi^J dt^J \biggr)^{1/p} \biggr).
\end{multline*}

Объединяя (2.2.13) с (2.2.39) и используя (2.1.32), а затем применяя Неравенство
Гёльдера с показателем $ q $ и (1.1.2) при $ a = p /q \le 1, $ и, наконец,
пользуясь (2.2.22), при $ p \le q $
и $ \mu \in \Z_+^d(\lambda) $ приходим к
оценке
\begin{multline*}
\biggl\| \sum_{ \nu \in N_{\Kappa^0 +\kappa}^{d,m,U}}
\D^\mu (V_{\Kappa^0, \kappa,\nu}^{d,l -\e,m,D,U,\Nu} f)
\D^{\lambda -\mu} g_{\Kappa^0 +\kappa, \nu}^{d,m} \biggr\|_{L_q(\R^d)}^q \le \\
\sum_{\substack{n \in \Z^d: Q_{\Kappa^0 +\kappa,n}^d \cap\\ G_{\Kappa^0 +\kappa}^{d,m,U} \ne \emptyset}}
\biggl(c_{23} 2^{(\kappa, \lambda +p^{-1} \e -q^{-1} \e)} \times\\
\biggl((\prod_{j \in \s(\kappa)}
2^{\kappa_j /p})
\biggl(\int_{ (c_{2} 2^{-\kappa} B^d)^{\s(\kappa)}}
\int_{D_\xi^{l \chi_{\s(\kappa)}} \cap D_{\Kappa^0 +\kappa,n}^{\prime d,m,D,U}}
| (\Delta_\xi^{l \chi_{\s(\kappa)}} f)(x)|^p
dx d\xi^{\s(\kappa)} \biggr)^{1/p} +\\
\sum_{J \subset \Nu_{1,d}^1: J \ne \emptyset}
(\prod_{j \in J} 2^{-\kappa_j (\alpha_j -\bm \epsilon_j)})
(\prod_{j \in \s(\kappa) \setminus J} 2^{\kappa_j /p})
\biggl(\int_{ (c_{3} 2^{-\kappa} I^d)^J}
(\prod_{j \in J} t_j^{-p(\alpha_j -\bm \epsilon_j +p^{-1}) -1}) \times\\
\int\limits_{ (t B^d)^J} \int\limits_{ ( c_{2} 2^{-\kappa} B^d)^{\s(\kappa) \setminus J}}
\int\limits_{D_\xi^{l \chi_{J \cup \s(\kappa)}} \cap D_{\Kappa^0 +\kappa,n}^{\prime d,m,D,U}}
| (\Delta_\xi^{l \chi_{J \cup \s(\kappa)}} f)(x)|^p
dx d\xi^{\s(\kappa) \setminus J} d\xi^J dt^J \biggr)^{1/p} \biggr) \biggr)^q \le\\
\sum_{\substack{n \in \Z^d: Q_{\Kappa^0 +\kappa,n}^d \cap\\ G_{\Kappa^0 +\kappa}^{d,m,U} \ne \emptyset}}
(c_{24} 2^{(\kappa, \lambda +p^{-1} \e -q^{-1} \e)})^q
\biggl((\prod_{j \in \s(\kappa)} 2^{\kappa_j q /p}) \times\\
\biggl(\int_{ (c_{2} 2^{-\kappa} B^d)^{\s(\kappa)}}
\int_{D_\xi^{l \chi_{\s(\kappa)}} \cap D_{\Kappa^0 +\kappa,n}^{\prime d,m,D,U}}
| (\Delta_\xi^{l \chi_{\s(\kappa)}} f)(x)|^p
dx d\xi^{\s(\kappa)} \biggr)^{q /p} +\\
\sum_{J \subset \Nu_{1,d}^1: J \ne \emptyset}
(\prod_{j \in J} 2^{-\kappa_j q(\alpha_j -\bm \epsilon_j)})
(\prod_{j \in \s(\kappa) \setminus J} 2^{\kappa_j q /p})
\biggl(\int_{ (c_{3} 2^{-\kappa} I^d)^J}
(\prod_{j \in J} t_j^{-p(\alpha_j -\bm \epsilon_j +p^{-1}) -1}) \times\\
\int\limits_{ (t B^d)^J} \int\limits_{ ( c_{2} 2^{-\kappa} B^d)^{\s(\kappa) \setminus J}}
\int\limits_{D_\xi^{l \chi_{J \cup \s(\kappa)}} \cap D_{\Kappa^0 +\kappa,n}^{\prime d,m,D,U}}
| (\Delta_\xi^{l \chi_{J \cup \s(\kappa)}} f)(x)|^p
dx d\xi^{\s(\kappa) \setminus J} d\xi^J dt^J \biggr)^{q /p} \biggr) =\\
(c_{24} 2^{(\kappa, \lambda +p^{-1} \e -q^{-1} \e)})^q
\biggl((\prod_{j \in \s(\kappa)} 2^{\kappa_j q /p}) \times\\
\sum_{\substack{n \in \Z^d: Q_{\Kappa^0 +\kappa,n}^d \cap\\ G_{\Kappa^0 +\kappa}^{d,m,U} \ne \emptyset}}
\biggl(\int_{ (c_{2} 2^{-\kappa} B^d)^{\s(\kappa)}}
\int_{D_\xi^{l \chi_{\s(\kappa)}} \cap D_{\Kappa^0 +\kappa,n}^{\prime d,m,D,U}}
| (\Delta_\xi^{l \chi_{\s(\kappa)}} f)(x)|^p
dx d\xi^{\s(\kappa)} \biggr)^{q /p} +\\
\sum_{J \subset \Nu_{1,d}^1: J \ne \emptyset}
(\prod_{j \in J} 2^{-\kappa_j q(\alpha_j -\bm \epsilon_j)})
(\prod_{j \in \s(\kappa) \setminus J} 2^{\kappa_j q /p}) \times\\
\sum_{\substack{n \in \Z^d: Q_{\Kappa^0 +\kappa,n}^d \cap\\ G_{\Kappa^0 +\kappa}^{d,m,U} \ne \emptyset}}
\biggl(\int_{ (c_{3} 2^{-\kappa} I^d)^J}
(\prod_{j \in J} t_j^{-p(\alpha_j -\bm \epsilon_j +p^{-1}) -1})\times \\
\int\limits_{ (t B^d)^J} \int\limits_{ ( c_{2} 2^{-\kappa} B^d)^{\s(\kappa) \setminus J}}
\int\limits_{D_\xi^{l \chi_{J \cup \s(\kappa)}} \cap D_{\Kappa^0 +\kappa,n}^{\prime d,m,D,U}}
| (\Delta_\xi^{l \chi_{J \cup \s(\kappa)}} f)(x)|^p
dx d\xi^{\s(\kappa) \setminus J} d\xi^J dt^J \biggr)^{q /p} \biggr) \le \\
(c_{24} 2^{(\kappa, \lambda +p^{-1} \e -q^{-1} \e)})^q
\biggl((\prod_{j \in \s(\kappa)} 2^{\kappa_j q /p}) \times\\
\biggl(\sum_{\substack{n \in \Z^d: Q_{\Kappa^0 +\kappa,n}^d \cap\\ G_{\Kappa^0 +\kappa}^{d,m,U} \ne \emptyset}}
\int_{ (c_{2} 2^{-\kappa} B^d)^{\s(\kappa)}}
\int_{D_\xi^{l \chi_{\s(\kappa)}} \cap D_{\Kappa^0 +\kappa,n}^{\prime d,m,D,U}}
| (\Delta_\xi^{l \chi_{\s(\kappa)}} f)(x)|^p
dx d\xi^{\s(\kappa)} \biggr)^{q /p} +\\
\sum_{J \subset \Nu_{1,d}^1: J \ne \emptyset}
(\prod_{j \in J} 2^{-\kappa_j q(\alpha_j -\bm \epsilon_j)})
(\prod_{j \in \s(\kappa) \setminus J} 2^{\kappa_j q /p}) \times\\
\biggl(\sum_{\substack{n \in \Z^d: Q_{\Kappa^0 +\kappa,n}^d \cap\\ G_{\Kappa^0 +\kappa}^{d,m,U} \ne \emptyset}}
\int_{ (c_{3} 2^{-\kappa} I^d)^J}
(\prod_{j \in J} t_j^{-p(\alpha_j -\bm \epsilon_j +p^{-1}) -1})\times \\
\int\limits_{ (t B^d)^J} \int\limits_{ ( c_{2} 2^{-\kappa} B^d)^{\s(\kappa) \setminus J}}
\int\limits_{D_\xi^{l \chi_{J \cup \s(\kappa)}} \cap D_{\Kappa^0 +\kappa,n}^{\prime d,m,D,U}}
| (\Delta_\xi^{l \chi_{J \cup \s(\kappa)}} f)(x)|^p
dx d\xi^{\s(\kappa) \setminus J} d\xi^J dt^J \biggr)^{q /p} \biggr) =\\
(c_{24} 2^{(\kappa, \lambda +p^{-1} \e -q^{-1} \e)})^q
\biggl((\prod_{j \in \s(\kappa)} 2^{\kappa_j q /p}) \\
\times \biggl(\sum_{\substack{n \in \Z^d: Q_{\Kappa^0 +\kappa,n}^d \cap\\ G_{\Kappa^0 +\kappa}^{d,m,U} \ne \emptyset}}
\int_{ (c_{2} 2^{-\kappa} B^d)^{\s(\kappa)}}
\int_{ D_\xi^{l \chi_{\s(\kappa)}}}
\chi_{D_{\Kappa^0 +\kappa,n}^{\prime d,m,D,U}}(x)
|\Delta_\xi^{l \chi_{\s(\kappa)}} f(x)|^p
dx d\xi^{\s(\kappa)} \biggr)^{q /p} +\\
\sum_{J \subset \Nu_{1,d}^1: J \ne \emptyset}
(\prod_{j \in J} 2^{-\kappa_j q(\alpha_j -\bm \epsilon_j)})
(\prod_{j \in \s(\kappa) \setminus J} 2^{\kappa_j q /p})
\biggl(\sum_{\substack{n \in \Z^d: Q_{\Kappa^0 +\kappa,n}^d \cap\\ G_{\Kappa^0 +\kappa}^{d,m,U} \ne \emptyset}}
\int_{ (c_{3} 2^{-\kappa} I^d)^J}
(\prod_{j \in J} t_j^{-p(\alpha_j -\bm \epsilon_j +p^{-1}) -1}) \times\\
\int\limits_{ (t B^d)^J} \int\limits_{ ( c_{2} 2^{-\kappa} B^d)^{\s(\kappa) \setminus J}}
\int\limits_{ D_\xi^{l \chi_{J \cup \s(\kappa)}}}
\chi_{D_{\Kappa^0 +\kappa,n}^{\prime d,m,D,U}}(x)
|\Delta_\xi^{l \chi_{J \cup \s(\kappa)}} f(x)|^p
dx d\xi^{\s(\kappa) \setminus J} d\xi^J dt^J \biggr)^{q/p} \biggr) =\\
(c_{24} 2^{(\kappa, \lambda +p^{-1} \e -q^{-1} \e)})^q
\biggl((\prod_{j \in \s(\kappa)} 2^{\kappa_j q /p}) \\
\times \biggl(\int_{ (c_{2} 2^{-\kappa} B^d)^{\s(\kappa)}}
\int_{ D_\xi^{l \chi_{\s(\kappa)}}}
(\sum_{\substack{n \in \Z^d: Q_{\Kappa^0 +\kappa,n}^d \cap\\ G_{\Kappa^0 +\kappa}^{d,m,U} \ne \emptyset}}
\chi_{D_{\Kappa^0 +\kappa,n}^{\prime d,m,D,U}}(x))
| (\Delta_\xi^{l \chi_{\s(\kappa)}} f)(x)|^p
dx d\xi^{\s(\kappa)} \biggr)^{q /p} +\\
\sum_{J \subset \Nu_{1,d}^1: J \ne \emptyset}
(\prod_{j \in J} 2^{-\kappa_j q(\alpha_j -\bm \epsilon_j)})
(\prod_{j \in \s(\kappa) \setminus J} 2^{\kappa_j q /p})\times\\
\biggl(\int_{ (c_{3} 2^{-\kappa} I^d)^J}
(\prod_{j \in J} t_j^{-p(\alpha_j -\bm \epsilon_j +p^{-1}) -1})
\int\limits_{ (t B^d)^J} \int\limits_{ ( c_{2} 2^{-\kappa} B^d)^{\s(\kappa) \setminus J}}
\int\limits_{ D_\xi^{l \chi_{J \cup \s(\kappa)}}} \\
\biggl(\sum_{\substack{n \in \Z^d: Q_{\Kappa^0 +\kappa,n}^d \cap\\ G_{\Kappa^0 +\kappa}^{d,m,U} \ne \emptyset}}
\chi_{D_{\Kappa^0 +\kappa,n}^{\prime d,m,D,U}}(x)\biggr)
| (\Delta_\xi^{l \chi_{J \cup \s(\kappa)}} f)(x)|^p
dx d\xi^{\s(\kappa) \setminus J} d\xi^J dt^J \biggr)^{q /p} \biggr) \le\\
(c_{24} 2^{(\kappa, \lambda +p^{-1} \e -q^{-1} \e)})^q
\biggl((\prod_{j \in \s(\kappa)} 2^{\kappa_j q /p})
\biggl(\int_{ (c_{2} 2^{-\kappa} B^d)^{\s(\kappa)}}
\int_{ D_\xi^{l \chi_{\s(\kappa)}}} c_4
| (\Delta_\xi^{l \chi_{\s(\kappa)}} f)(x)|^p
dx d\xi^{\s(\kappa)} \biggr)^{q /p} +\\
\sum_{J \subset \Nu_{1,d}^1: J \ne \emptyset}
(\prod_{j \in J} 2^{-\kappa_j q(\alpha_j -\bm \epsilon_j)})
(\prod_{j \in \s(\kappa) \setminus J} 2^{\kappa_j q /p})
\biggl(\int_{ (c_{3} 2^{-\kappa} I^d)^J}
(\prod_{j \in J} t_j^{-p(\alpha_j -\bm \epsilon_j +p^{-1}) -1}) \times\\
\int\limits_{ (t B^d)^J} \int\limits_{ ( c_{2} 2^{-\kappa} B^d)^{\s(\kappa) \setminus J}}
\int\limits_{ D_\xi^{l \chi_{J \cup \s(\kappa)}}} c_{4}
| (\Delta_\xi^{l \chi_{J \cup \s(\kappa)}} f)(x)|^p
dx d\xi^{\s(\kappa) \setminus J} d\xi^J dt^J \biggr)^{q /p} \biggr),
\end{multline*}
откуда, применяя (1.1.2) при $ a = 1 /q $ и делая замену переменных, получаем
\begin{multline*} \tag{2.2.40}
\biggl\| \sum_{ \nu \in N_{\Kappa^0 +\kappa}^{d,m,U}}
\D^\mu (V_{\Kappa^0, \kappa,\nu}^{d,l -\e,m,D,U,\Nu} f)
\D^{\lambda -\mu} g_{\Kappa^0 +\kappa, \nu}^{d,m} \biggr\|_{L_q(\R^d)} \le\\
c_{25} 2^{(\kappa, \lambda +p^{-1} \e -q^{-1} \e)}
\biggl((\prod_{j \in \s(\kappa)} 2^{\kappa_j q /p})
\biggl(\int_{ (c_{2} 2^{-\kappa} B^d)^{\s(\kappa)}}
\int_{ D_\xi^{l \chi_{\s(\kappa)}}}
| (\Delta_\xi^{l \chi_{\s(\kappa)}} f)(x)|^p
dx d\xi^{\s(\kappa)} \biggr)^{q /p} +\\
\sum_{J \subset \Nu_{1,d}^1: J \ne \emptyset}
(\prod_{j \in J} 2^{-\kappa_j q(\alpha_j -\bm \epsilon_j)})
(\prod_{j \in \s(\kappa) \setminus J} 2^{\kappa_j q /p})
\biggl(\int_{ (c_{3} 2^{-\kappa} I^d)^J}
(\prod_{j \in J} t_j^{-p(\alpha_j -\bm \epsilon_j +p^{-1}) -1}) \times\\
\int_{ (t B^d)^J} \int_{ ( c_{2} 2^{-\kappa} B^d)^{\s(\kappa) \setminus J}}
\int_{ D_\xi^{l \chi_{J \cup \s(\kappa)}}}
| (\Delta_\xi^{l \chi_{J \cup \s(\kappa)}} f)(x)|^p
dx d\xi^{\s(\kappa) \setminus J} d\xi^J dt^J \biggr)^{q /p} \biggr)^{1/q} \le\\
c_{25} 2^{(\kappa, \lambda +p^{-1} \e -q^{-1} \e)}
\biggl((\prod_{j \in \s(\kappa)} 2^{\kappa_j /p})
\biggl(\int_{ (c_{2} 2^{-\kappa} B^d)^{\s(\kappa)}}
\int_{ D_\xi^{l \chi_{\s(\kappa)}}}
| (\Delta_\xi^{l \chi_{\s(\kappa)}} f)(x)|^p
dx d\xi^{\s(\kappa)} \biggr)^{1/p} +\\
\sum_{J \subset \Nu_{1,d}^1: J \ne \emptyset}
(\prod_{j \in J} 2^{-\kappa_j (\alpha_j -\bm \epsilon_j)})
(\prod_{j \in \s(\kappa) \setminus J} 2^{\kappa_j /p})
\biggl(\int_{ (c_{3} 2^{-\kappa} I^d)^J}
(\prod_{j \in J} t_j^{-p(\alpha_j -\bm \epsilon_j +p^{-1}) -1}) \times\\
\int_{ (t B^d)^J} \int_{ ( c_{2} 2^{-\kappa} B^d)^{\s(\kappa) \setminus J}}
\int_{ D_\xi^{l \chi_{J \cup \s(\kappa)}}}
| (\Delta_\xi^{l \chi_{J \cup \s(\kappa)}} f)(x)|^p
dx d\xi^{\s(\kappa) \setminus J} d\xi^J dt^J \biggr)^{1/p} \biggr) \le\\
c_{25} 2^{(\kappa, \lambda +p^{-1} \e -q^{-1} \e)}
\biggl( c_{26} \biggl((\prod_{j \in \s(\kappa)} (2 c_{2} 2^{-\kappa_j})^{-1}) \times\\
\int_{(c_{2} 2^{-\kappa} B^d)^{\s(\kappa)}}
\| \Delta_\xi^{l \chi_{\s(\kappa)}} f \|_{L_p(D_\xi^{l \chi_{\s(\kappa)}})}^p
d\xi^{\s(\kappa)} \biggr)^{1/p} +\\
\sum_{J \subset \Nu_{1,d}^1: J \ne \emptyset}
(\prod_{j \in J} 2^{-\kappa_j (\alpha_j -\bm \epsilon_j)})
\biggl(\int_{ (c_{3} 2^{-\kappa} I^d)^J}
(\prod_{j \in J} t_j^{-p(\alpha_j -\bm \epsilon_j) -1})\times\\
c_{27} \biggl((\prod_{j \in J} (2 t_j)^{-1})
(\prod_{j \in (\s(\kappa) \setminus J)} ( 2 c_{2} 2^{-\kappa_j})^{-1})\times\\
\int_{ (t B^d)^J \times ( c_{2} 2^{-\kappa} B^d)^{\s(\kappa) \setminus J}}
\| \Delta_\xi^{l \chi_{J \cup \s(\kappa)}} f \|_{L_p(D_\xi^{l \chi_{J \cup \s(\kappa)}})}^p
d\xi^J d\xi^{\s(\kappa) \setminus J} \biggr) dt^J \biggr)^{1/p} \biggr) \le\\
c_{28} 2^{(\kappa, \lambda +p^{-1} \e -q^{-1} \e)}
\biggl(\Omega^{\prime l \chi_{\s(\kappa)}}(f, (c_{2} 2^{-\kappa})^{\s(\kappa)})_{L_p(D)} +\\
\sum_{J \subset \Nu_{1,d}^1: J \ne \emptyset}
(\prod_{j \in J} 2^{-\kappa_j (\alpha_j -\bm \epsilon_j)})
\biggl(\int_{ (c_{3} 2^{-\kappa} I^d)^J}
(\prod_{j \in J} t_j^{-p(\alpha_j -\bm \epsilon_j) -1}) \times\\
(\Omega^{\prime l \chi_{J \cup \s(\kappa)}}(f,
(t \chi_J +c_{2} 2^{-\kappa}
\chi_{\s(\kappa) \setminus J})^{J \cup \s(\kappa)})_{L_p(D)})^p dt^J \biggr)^{1/p} \biggr)=\\
c_{28} 2^{(\kappa, \lambda +p^{-1} \e -q^{-1} \e)}
\biggl(\Omega^{\prime l \chi_{\s(\kappa)}}(f, (c_{2} 2^{-\kappa})^{\s(\kappa)})_{L_p(D)} +\\
\sum_{J \subset \Nu_{1,d}^1: J \ne \emptyset}
(\prod_{j \in J} 2^{-\kappa_j (\alpha_j -\bm \epsilon_j)})
\biggl(\int_{ (c_{3} I^d)^J}
(\prod_{j \in J} (2^{-\kappa_j} u_j)^{-p(\alpha_j -\bm \epsilon_j) -1}) \times\\
(\Omega^{\prime l \chi_{J \cup \s(\kappa)}}(f,
(2^{-\kappa} u \chi_J +c_{2} 2^{-\kappa}
\chi_{\s(\kappa) \setminus J})^{J \cup \s(\kappa)})_{L_p(D)})^p
(\prod_{j \in J} 2^{-\kappa_j}) du^J \biggr)^{1/p} \biggr) \\
= c_{28} 2^{(\kappa, \lambda +p^{-1} \e -q^{-1} \e)}
\biggl(\Omega^{\prime l \chi_{\s(\kappa)}}(f, (c_{2} 2^{-\kappa})^{\s(\kappa)})_{L_p(D)} +\\
\sum_{J \subset \Nu_{1,d}^1: J \ne \emptyset}
\biggl(\int_{ (c_{3} I^d)^J}
(\prod_{j \in J} u_j^{-p(\alpha_j -\bm \epsilon_j) -1})
(\Omega^{\prime l \chi_{J \cup \s(\kappa)}}(f, (u \chi_{J \setminus \s(\kappa)} \\
+2^{-\kappa} u \chi_{\s(\kappa) \cap J} +c_{2} 2^{-\kappa}
\chi_{\s(\kappa) \setminus J})^{J \cup \s(\kappa)})_{L_p(D)})^p du^J \biggr)^{1/p} \biggr).
\end{multline*}

Подставляя (2.2.40) в (2.2.12), находим, что при $ p \le q $ и $ \kappa
\in \Z_+^d \setminus \{0\} $ для $ f \in (S_p^\alpha H)^\prime(D) $
соблюдается неравенство
\begin{multline*}
\| \D^\lambda \mathcal R_{\Kappa^0,\kappa}^{d,l -\e,m,D,U,\Nu} f \|_{L_q(\R^d)} \le
c_{1} 2^{(\kappa, \lambda +p^{-1} \e -q^{-1} \e)}
\biggl(\Omega^{\prime l \chi_{\s(\kappa)}}(f, (c_{2} 2^{-\kappa})^{\s(\kappa)})_{L_p(D)} + \\
\sum_{J \subset \Nu_{1,d}^1: J \ne \emptyset}
\biggl(\int_{ (c_{3} I^d)^J}
(\prod_{j \in J} u_j^{-p(\alpha_j -\bm \epsilon_j) -1})
(\Omega^{\prime l \chi_{J \cup \s(\kappa)}}(f, (u \chi_{J \setminus \s(\kappa)} +\\
2^{-\kappa} u \chi_{\s(\kappa) \cap J} +c_{2} 2^{-\kappa}
\chi_{\s(\kappa) \setminus J})^{J \cup \s(\kappa)})_{L_p(D)})^p du^J \biggr)^{1/p} \biggr),
\end{multline*}
которое совпадает с (2.2.11) при $ p \le q. $
Неравенство (2.2.11) при $ q < p $ вытекает из неравенства (2.2.11) при
$ q = p $ и того факта, что в этом случае ввиду (2.2.24)
\begin{multline*}
\| \D^\lambda \mathcal R_{\Kappa^0, \kappa}^{d,l -\e,m,D,U,\Nu} f \|_{L_q(\R^d)} =
\| \D^\lambda \mathcal R_{\Kappa^0, \kappa}^{d,l -\e,m,D,U,\Nu} f
\|_{L_q(Q^{d,m,D})} \le \\
C(d,m,D,p,q) \| \D^\lambda \mathcal R_{\Kappa^0, \kappa}^{d,l -\e,m,D,U,\Nu} f \|_{L_p(Q^{d,m,D})} =\\
C(d,m,D,p,q) \| \D^\lambda \mathcal R_{\Kappa^0, \kappa}^{d,l -\e,m,D,U,\Nu} f \|_{L_p(\R^d)}. \square
\end{multline*}

Предложение 2.2.3

Пусть выполнены условия леммы 2.2.2 и соблюдается неравенство (2.1.6). Тогда
для любой функции $ f \in (S_p^\alpha H)^\prime(D) $ при $ l = l(\alpha) $

в $ L_q(U) $ имеет место равенство
\begin{equation*} \tag{2.2.41}
\D^\lambda (f \mid_U) = \sum_{\kappa \in \Z_+^d} (\D^\lambda (\mathcal R_{\Kappa^0, \kappa}^{d, l -\e,m,D,U,\Nu} f)) \mid_U.
\end{equation*}

Доказательство.

Прежде всего отметим, что в условиях предложения согласно предложению 2.1.7
в $ L_p(U) $ справедливо равенство (2.1.23) (с $ \Kappa^0 $ вместо $ \kappa^0 $).
Поэтому на основании леммы 1.2.1 и соотношений (1.4.16), (2.1.43)
заключаем, что в $ L_p(U) $ имеет место равенство
\begin{multline*} \tag{2.2.42}
f \mid_U = \sum_{\kappa \in \Z_+^d} (\sum_{\epsilon \in \Upsilon^d:
\s(\epsilon) \subset \s(\kappa)} (-\e)^\epsilon
(R_{\Kappa^0 +\kappa -\epsilon}^{d,l -\e,m,D,U,\nu_{\Kappa^0 +\kappa -\epsilon}} f) \mid_U) = \\
\sum_{\kappa \in \Z_+^d} (\sum_{\epsilon \in \Upsilon^d:
\s(\epsilon) \subset \s(\kappa)} (-\e)^\epsilon
(H_{\Kappa^0 +\kappa, \Kappa^0 +\kappa -\epsilon}^{d,l -\e,m,U}
R_{\Kappa^0 +\kappa -\epsilon}^{d,l -\e,m,D,U,\nu_{\Kappa^0 +\kappa -\epsilon}} f) \mid_U) = \\
\sum_{\kappa \in \Z_+^d} (\sum_{\epsilon \in \Upsilon^d:
\s(\epsilon) \subset \s(\kappa)} (-\e)^\epsilon
(H_{\Kappa^0 +\kappa, \Kappa^0 +\kappa -\epsilon}^{d,l -\e,m,U}
R_{\Kappa^0 +\kappa -\epsilon}^{d,l -\e,m,D,U,\nu_{\Kappa^0 +\kappa -\epsilon}} f)) \mid_U = \\
\sum_{ \kappa \in \Z_+^d} (\mathcal R_{\Kappa^0, \kappa}^{d,l -\e,m,D,U,\Nu} f) \mid_U.
\end{multline*}

Далее, для любой функции $ f \in (S_p^\alpha h)^\prime(D) $ и каждого множества
$ J^\prime \subset \Nu_{1,d}^1: J^\prime  \ne \emptyset, $ при $ \kappa \in \Z_+^d:
\s(\kappa) = J^\prime, $ ввиду (2.2.11) выполняется неравенство
\begin{multline*} \tag{2.2.43}
\| \D^\lambda ((\mathcal R_{\Kappa^0, \kappa}^{d, l -\e,m,D,U,\Nu} f) \mid_U) \|_{L_q(U)} =
\| (\D^\lambda \mathcal R_{\Kappa^0, \kappa}^{d, l -\e,m,D,U,\Nu} f) \mid_U \|_{L_q(U)} \le \\
\| \D^\lambda \mathcal R_{\Kappa^0,\kappa}^{d,l -\e,m,D,U,\Nu} f\|_{L_q(\R^d)} \le
c_{1} 2^{(\kappa, \lambda +(p^{-1} -q^{-1})_+ \e)}
\biggl(\Omega^{\prime l \chi_{J^\prime}}(f, (c_{2} 2^{-\kappa})^{J^\prime})_{L_p(D)} +\\
\sum_{J \subset \Nu_{1,d}^1: J \ne \emptyset}
\biggl(\int_{ (c_{3} I^d)^J} (\prod_{j \in J} u_j^{-p(\alpha_j -\bm \epsilon_j) -1})
(\Omega^{\prime l \chi_{J \cup J^\prime}}(f, (u \chi_{J \setminus J^\prime} +\\
2^{-\kappa} u \chi_{J^\prime \cap J} +c_{2} 2^{-\kappa}
\chi_{J^\prime \setminus J})^{J \cup J^\prime})_{L_p(D)})^p du^J \biggr)^{1/p} \biggr) \le \\
c_{1} 2^{(\kappa, \lambda +(p^{-1} -q^{-1})_+ \e)}
\biggl(\|f\|_{(S_p^\alpha H)^\prime(D)} ((c_{2} 2^{-\kappa})^{J^\prime})^{\alpha^{J^\prime}} +\\
\sum_{J \subset \Nu_{1,d}^1: J \ne \emptyset}
\biggl(\int_{ (c_{3} I^d)^J} (\prod_{j \in J} u_j^{-p(\alpha_j -\bm \epsilon_j) -1})
(\|f\|_{(S_p^\alpha H)^\prime(D)} (\prod_{j \in J \setminus J^\prime} u_j^{\alpha_j})\times\\
(\prod_{j \in J^\prime \cap J} (2^{-\kappa_j} u_j)^{\alpha_j})
(\prod_{j \in J^\prime \setminus J} (c_{2} 2^{-\kappa_j})^{\alpha_j}))^p du^J \biggr)^{1/p} \biggr) \le \\
c_{29} 2^{(\kappa, \lambda +(p^{-1} -q^{-1})_+ \e)}
\|f\|_{(S_p^\alpha H)^\prime(D)} \biggl((\prod_{j \in \s(\kappa)} 2^{-\kappa_j \alpha_j}) +\\
\sum_{J \subset \Nu_{1,d}^1: J \ne \emptyset}
\biggl(\int_{ (c_{3} I^d)^J} (\prod_{j \in J} u_j^{-p(\alpha_j -\bm \epsilon_j) -1})
(\prod_{j \in J} u_j^{p \alpha_j})
(\prod_{j \in J^\prime} 2^{-\kappa_j \alpha_j p})
du^J \biggr)^{1/p} \biggr) \le \\
c_{30} 2^{(\kappa, \lambda +(p^{-1} -q^{-1})_+ \e)}
\|f\|_{(S_p^\alpha H)^\prime(D)} \biggl(2^{-(\kappa, \alpha)} +\\
\sum_{J \subset \Nu_{1,d}^1: J \ne \emptyset}
(\prod_{j \in \s(\kappa)} 2^{-\kappa_j \alpha_j})
\biggl(\int_{ (c_{3} I^d)^J} (\prod_{j \in J} u_j^{p \bm \epsilon_j -1})
du^J \biggr)^{1/p} \biggr) \le \\
c_{31} \|f\|_{(S_p^\alpha H)^\prime(D)} 2^{-(\kappa, \alpha -\lambda -(p^{-1} -q^{-1})_+ \e)}.
\end{multline*}

Из (2.2.43) с учётом (2.1.6) следует, что ряд $ \sum_{ \kappa \in \Z_+^d}
\| \D^\lambda ((\mathcal R_{\Kappa^0,\kappa}^{d, l -\e,m,D,U,\Nu} f) \mid_U) \|_{L_q(U)} $
сходится, и, значит, ряд
$ \sum_{ \kappa \in \Z_+^d} \D^\lambda ((\mathcal R_{\Kappa^0,\kappa}^{d, l -\e,m,D,U,\Nu} f) \mid_U)
$
сходится в $ L_q(U). $

Принимая во внимание это обстоятельство и равенство (2.2.42), для
любой функции $ \phi \in C_0^\infty(U) $ имеем
\begin{multline*}
\langle \D^\lambda (f \mid_U), \phi \rangle
= \int_{U} (\sum_{ \kappa \in \Z_+^d} (\mathcal R_{\Kappa^0,\kappa}^{d,l -\e,m,D,U,\Nu} f) \mid_U)
(-1)^{(\lambda,\e)} \D^\lambda \phi dx =\\
\int_{U} (\sum_{ \kappa \in \Z_+^d}
(\D^\lambda (\mathcal R_{\Kappa^0,\kappa}^{d,l -\e,m,D,U,\Nu} f)) \mid_U) \phi dx.
\end{multline*}
А это значит, что в $ L_q(U) $ верно равенство (2.2.41). $ \square $

Напомним постановку задачи восстановления производных по значениям
функций в заданном числе точек.

Пусть $ T $ -- топологическое пространство, $ X $ -- банахово
пространство над $ \R, \mathcal V: D(\mathcal V) \mapsto X $ -- линейный оператор с областью
определения $ D(\mathcal V) \subset C(T), $ принимающий значения в $ X. $ Пусть ещё
$ K \subset D(\mathcal V) $ -- некоторый класс функций.
Для $ n \in \N $ через $ \Phi_n(C(T)) $ обозначим совокупнсть всех отображений
$ \phi: C(T) \mapsto \R^n, $ для каждого из которых существует набор точек
$ \{ t^j \in T, j =1,\ldots,n\} $ такой, что
$$
\phi(f) = (f(t^1), \ldots, f(t^n)), f \in C(T).
$$
При $ n \in \N $ обозначим также через $ \mathcal A^n(X)
(\overline {\mathcal A}^n(X)) $ множество всех отображений
(всех линейных отображений) $ A: \R^n \mapsto X. $

Тогда при $ n \in \N $ положим
$$
\sigma_n(\mathcal V,K,X) =
\inf_{ A \in \mathcal A^n(X), \phi \in \Phi_n(C(T))} \sup_{f \in K}
\| \mathcal V f -A \circ \phi (f)\|_X,
$$
а
$$
\overline \sigma_n(\mathcal V,K,X) =
\inf_{ A \in \overline{\mathcal A}^n(X), \phi \in \Phi_n(C(T))} \sup_{f \in K}
\| \mathcal V f -A \circ \phi (f)\|_X.
$$

Теперь установим следующее утверждение.

Предложение 2.2.4

Пусть выполнены условия предложения 2.2.3 (кроме условий $ m \in \N^d, 
\lambda \in \Z_+^d(m) $ ) и $ 1 \le \theta \le \infty. $
Тогда существуют константы $ c_{32}(d,\alpha,p,\theta,q,\lambda,D,U) >0 $ и
$ n_0(d,\alpha,p,\theta,q,\lambda,D,U) \in \N $ такие,что при $ n \ge n_0 $
можно построить отображения $ \phi \in \Phi_n(C(D)) $ и $ A \in
\overline {\mathcal A}^n(L_q(D)) $ такие, что для
любой функции $ f \in (\mathcal S_{p,\theta}^\alpha \mathcal B)^\prime(D) $
выполняется неравенство
\begin{equation*} \tag{2.2.44}
\| \D^\lambda (f \mid_U) -(A \circ \phi (f)) \mid_U\|_{L_q(U)} \le
c_{32} n^{-\mn} (\log n)^{(\mn +1 -1/\max(p,\theta))(\cmn -1)},
\end{equation*}
где
$$
\mn = \mn(\alpha -\lambda -(p^{-1} -q^{-1})_+ \e),
\cmn = \cmn(\alpha -\lambda -(p^{-1} -q^{-1})_+ \e) = \cmn(\alpha -\lambda).
$$

Доказательство.

В условиях предложения, полагая $ \mathcal J = \{j \in \Nu_{1,d}^1:
\alpha_j -\lambda_j -(p^{-1} -q^{-1})_+ = \mn \}, $
фиксируем вектор $ \beta \in \R_+^d, $ обладающий следующими свойствами:
$ \beta_j =1, j \in \mathcal J, \beta_j >1 $ и $ \beta_j^{-1}
(\alpha_j -\lambda_j -(p^{-1} -q^{-1})_+) > \mn,
j \in \Nu_{1,d}^1 \setminus \mathcal J. $

Далее, фиксируя $ m \in \N^d $ так, чтобы $ \lambda \in \Z_+^d(m), $ возьмём 
$ \Kappa^0 = \Kappa^0(d,m,D,U) \in \Z_+^d, 
\Nu = \{\nu_{\Kappa^0 +\kappa}^{d,m,D,U}, \kappa \in \Z_+^d\} $ 
из леммы 2.2.1, и при $ l = l(\alpha) $ определим семейство точек
\begin{multline*}
\mathtt x_{\Kappa^0 +\kappa, \nu}^{d, l -\e, \rho,m,D,U} = \xi_{\delta,
x^0}^{d,l-\e,\rho} = x^0 +\delta \xi_{\e,0}^{d,l-\e,\rho}
\text{ (см. п. 2.1.), при } x^0 = 2^{-\Kappa^0 -\kappa} \nu, \delta \\
= 2^{-\Kappa^0 -\kappa},
\rho \in \Z_+^d(l -\e), \nu \in \Z^d: Q_{\Kappa^0 +\kappa,\nu}^d \subset D,
\kappa \in \Z_+^d,
\end{multline*}
и для $ r \in \N $ рассмотрим множество точек
$$
\{ \mathtt x_{\Kappa^0 +\kappa, \nu}^{d, l -\e, \rho,m,D,U} \mid
\rho \in \Z_+^d(l -\e), \nu \in \Z^d: Q_{\Kappa^0 +\kappa,\nu}^d \subset D,
\kappa \in \Z_+^d: (\kappa, \beta) \le r \}.
$$
Число этих точек, учитывая (1.2.1), удовлетворяет неравенству
\begin{multline*}
\card \{ \mathtt x_{\Kappa^0 +\kappa, \nu}^{d, l -\e, \rho,m,D,U} \mid  \rho \in
\Z_+^d(l -\e), \nu \in \Z^d: Q_{\Kappa^0 +\kappa,\nu}^d \subset D,
\kappa \in \Z_+^d: (\kappa, \beta) \le r \} \le\\
\sum_{\kappa \in \Z_+^d: (\kappa,\beta) \le r} l^{\e} \card \{\nu \in
\Z^d: Q_{\Kappa^0 +\kappa,\nu}^d \subset D\} \le \\
\sum_{\kappa \in \Z_+^d: (\kappa,\beta) \le r} l^{\e} (\mes D) /
2^{-(\Kappa^0 +\kappa, \e)} \le
l^{\e} (\mes D) 2^{(\Kappa^0, \e)}
\sum_{\kappa \in \Z_+^d: (\kappa,\beta) \le r} 2^{(\kappa, \e)} \le \\
c_{33} 2^{\mx(\beta^{-1} \e) r} r^{\cmx(\beta^{-1} \e) -1} =
c_{33} 2^r r^{\card (\mathcal J) -1} = c_{33} 2^r r^{\cmn -1}.
\end{multline*}

Для $ n \ge n_0 = 2 c_{33} $ выберем $ r \in \N $ так, чтобы соблюдалось
соотношение
\begin{equation*} \tag{2.2.45}
c_{33} 2^r r^{\cmn -1} \le n < c_{33} 2^{r +1} (r +1)^{\cmn -1},
\end{equation*}
и построим соответствующую этому $ r \in \N $ систему точек
$$
\{ \mathtt x_{\Kappa^0 +\kappa, \nu}^{d, l -\e, \rho,m,D,U} \mid \rho \in
\Z_+^d(l -\e), \nu \in \Z^d: Q_{\Kappa^0 +\kappa,\nu}^d \subset D,
\kappa \in \Z_+^d: (\kappa, \beta) \le r \}.
$$
Тогда в условиях предложения (см. также (1.1.8)) нетрудно видеть, что существуют
отображения $ \phi \in \Phi_n(C(D)) $ и $ A \in \overline {\mathcal A}^n(L_q(D)) $
такие, что для $ f \in (\mathcal S_{p,\theta}^\alpha \mathcal B)^\prime(D) $
имеет место представление
$$
A \circ \phi (f) = \sum_{ \kappa \in \Z_+^d: (\kappa, \beta) \le r}
(\D^\lambda (\mathcal R_{\Kappa^0,\kappa}^{d, l -\e,m,D,U,\Nu} f)) \mid_D
$$
При этом в силу (2.2.41), (2.2.11) с учётом (1.1.8) для
$ f \in (\mathcal S_{p,\theta}^\alpha \mathcal B)^\prime(D) $ выполняется
неравенство
\begin{multline*} \tag{2.2.46}
\| \D^\lambda (f \mid_U) -(A \circ \phi(f)) \mid_U \|_{L_q(U)} =\\
\| \D^\lambda (f \mid_U) -(\sum_{ \kappa \in \Z_+^d: (\kappa, \beta) \le r}
(\D^\lambda (\mathcal R_{\Kappa^0,\kappa}^{d, l -\e,m,D,U,\Nu} f)) \mid_D) \mid_U \|_{L_q(U)} = \\
\| \D^\lambda (f \mid_U) -\sum_{ \kappa \in \Z_+^d: (\kappa, \beta) \le r}
(\D^\lambda (\mathcal R_{\Kappa^0,\kappa}^{d, l -\e,m,D,U,\Nu} f)) \mid_U \|_{L_q(U)} = \\
\| \sum_{ \kappa \in \Z_+^d: (\kappa, \beta) > r}
(\D^\lambda (\mathcal R_{\Kappa^0,\kappa}^{d,l -\e,m,D,U,\Nu} f)) \mid_U \|_{L_q(U)} \le \\
\sum_{ \kappa \in \Z_+^d: (\kappa,\beta) > r}
\| (\D^\lambda (\mathcal R_{\Kappa^0,\kappa}^{d, l -\e,m,D,U,\Nu} f)) \mid_U
\|_{L_q(U)} \le
\sum_{ \kappa \in \Z_+^d: (\kappa,\beta) > r}
\| \D^\lambda \mathcal R_{\Kappa^0,\kappa}^{d, l -\e,m,D,U,\Nu} f\|_{L_q(\R^d)} \le \\
\sum_{ \kappa \in \Z_+^d: (\kappa, \beta) > r}
c_{1} 2^{(\kappa, \lambda +(p^{-1} -q^{-1})_+ \e)} \biggl(\Omega^{\prime l
\chi_{\s(\kappa)}}(f, (c_{2} 2^{-\kappa})^{\s(\kappa)})_{L_p(D)} +\\
\sum_{J \subset \Nu_{1,d}^1: J \ne \emptyset} \biggl(\int_{ (c_{3} I^d)^J}
(\prod_{j \in J} u_j^{-p(\alpha_j -\bm \epsilon_j) -1})
(\Omega^{\prime l \chi_{J \cup \s(\kappa)}}(f, (u \chi_{J \setminus
\s(\kappa)} +\\
2^{-\kappa} u \chi_{\s(\kappa) \cap J} +c_{2}
2^{-\kappa} \chi_{\s(\kappa) \setminus J})^{J \cup \s(\kappa)})_{L_p(D)})^p
du^J \biggr)^{1/p} \biggr) = \\
\sum_{ \kappa \in \Z_+^d: (\kappa, \beta) > r} c_{1}
2^{(\kappa, \lambda +(p^{-1} -q^{-1})_+ \e)} \Omega^{\prime l \chi_{\s(\kappa)}}(f, (c_{2}
2^{-\kappa})^{\s(\kappa)})_{L_p(D)} +\\
\sum_{ \kappa \in \Z_+^d: (\kappa, \beta) > r}
\sum_{J \subset \Nu_{1,d}^1: J \ne \emptyset}
c_{1} 2^{(\kappa, \lambda +(p^{-1} -q^{-1})_+ \e)} \biggl(\int_{ (c_{3} I^d)^J}
(\prod_{j \in J} u_j^{-p(\alpha_j -\bm \epsilon_j) -1})\times \\
(\Omega^{\prime l \chi_{J \cup \s(\kappa)}}(f, (u \chi_{J \setminus
\s(\kappa)} + 2^{-\kappa} u \chi_{\s(\kappa) \cap J} +c_{2}
2^{-\kappa} \chi_{\s(\kappa) \setminus J})^{J \cup \s(\kappa)})_{L_p(D)})^p
du^J \biggr)^{1/p} = \\
c_{1} \biggl(\sum_{ \kappa \in \Z_+^d: (\kappa, \beta) > r}
2^{(\kappa, \lambda +(p^{-1} -q^{-1})_+ \e)} \Omega^{\prime l \chi_{\s(\kappa)}}(f,
(c_{2} 2^{-\kappa})^{\s(\kappa)})_{L_p(D)} +\\
\sum_{J \subset \Nu_{1,d}^1: J \ne \emptyset}
\sum_{ \kappa \in \Z_+^d: (\kappa,\beta) > r} 2^{(\kappa, \lambda +(p^{-1} -q^{-1})_+ \e)}
\biggl(\int_{(c_{3} I^d)^J} (\prod_{j \in J} u_j^{-p(\alpha_j -\bm \epsilon_j) -1}) \\
(\Omega^{\prime l \chi_{J \cup \s(\kappa)}}(f, (u \chi_{J \setminus
\s(\kappa)} + 2^{-\kappa} u \chi_{\s(\kappa) \cap J} +c_{2}
2^{-\kappa} \chi_{\s(\kappa) \setminus J})^{J \cup \s(\kappa)})_{L_p(D)})^p
du^J \biggr)^{1/p} \biggr) = \\
c_{1} \biggl(\sum_{ \kappa \in \Z_+^d: (\kappa, \beta) > r}
2^{(\kappa, \lambda +(p^{-1} -q^{-1})_+ \e)}
\Omega^{\prime l \chi_{\s(\kappa)}}(f, (c_{2} 2^{-\kappa})^{\s(\kappa)})_{L_p(D)} +\\
\sum_{J \subset \Nu_{1,d}^1: J \ne \emptyset}
\sum_{J^\prime \subset \Nu_{1,d}^1: J^\prime \ne \emptyset}
\sum_{ \kappa \in \Z_+^d: (\kappa, \beta) > r, \s(\kappa) = J^\prime}
2^{(\kappa, \lambda +(p^{-1} -q^{-1})_+ \e)} \times \\
\biggl(\int_{ (c_{3} I^d)^J} (\prod_{j \in J}
u_j^{-p(\alpha_j -\bm \epsilon_j) -1})
(\Omega^{\prime l \chi_{J \cup J^\prime}}(f, (u \chi_{J \setminus J^\prime} +\\
2^{-\kappa} u \chi_{J^\prime \cap J} +c_{2} 2^{-\kappa}
\chi_{J^\prime \setminus J})^{J \cup J^\prime})_{L_p(D)})^p du^J \biggr)^{1/p} \biggr).
\end{multline*}

Оценивая первую сумму в правой части (2.2.46), с помощью неравенства
Гёльдера для $ f \in (\mathcal S_{p,\theta}^\alpha \mathcal B)^\prime(D),
r \in \N $ получаем
\begin{multline*} \tag{2.2.47}
\sum_{ \kappa \in \Z_+^d: (\kappa, \beta) > r}
2^{(\kappa, \lambda +(p^{-1} -q^{-1})_+ \e)}
\Omega^{\prime l \chi_{\s(\kappa)}}(f, (c_{2} 2^{-\kappa})^{\s(\kappa)})_{L_p(D)} =\\
\sum_{ \kappa \in \Z_+^d: (\kappa, \beta) > r}
2^{-(\kappa, \alpha -\lambda -(p^{-1} -q^{-1})_+ \e)}
2^{(\kappa, \alpha)} \Omega^{\prime l \chi_{\s(\kappa)}}(f,
(c_{2} 2^{-\kappa})^{\s(\kappa)})_{L_p(D)} \le\\
\biggl(\sum_{ \kappa \in \Z_+^d: (\kappa, \beta) > r}
2^{-(\kappa, \alpha -\lambda -(p^{-1}  -q^{-1})_+ \e)
\theta^\prime} \biggr)^{1/\theta^\prime} \times\\
 \biggl(\sum_{ \kappa \in \Z_+^d: (\kappa, \beta) > r}
(2^{(\kappa, \alpha)} \Omega^{\prime l \chi_{\s(\kappa)}}(f,
(c_{2} 2^{-\kappa})^{\s(\kappa)})_{L_p(D)})^\theta \biggr)^{1/\theta},
\theta^\prime = \theta /(\theta -1).
\end{multline*}

Благодаря соблюдению (2.1.6) пользуясь (1.2.2), выводим
\begin{multline*} \tag{2.2.48}
\biggl(\sum_{ \kappa \in \Z_+^d: (\kappa, \beta) > r}
2^{-(\kappa, \alpha -\lambda -(p^{-1}  -q^{-1})_+ \e)
\theta^\prime} \biggr)^{1/\theta^\prime} \le\\
(c_{34} 2^{-\mn(\theta^\prime \beta^{-1} (\alpha -\lambda -(p^{-1}
-q^{-1})_+ \e)) r} r^{\cmn(\theta^\prime
\beta^{-1} (\alpha -\lambda -(p^{-1} -q^{-1})_+ \e)) -1})^{1/\theta^\prime} =\\
(c_{34} 2^{-\theta^\prime \mn(\alpha -\lambda -(p^{-1} -q^{-1})_+ \e) r}
r^{\cmn(\alpha -\lambda -(p^{-1} -q^{-1})_+ \e) -1})^{1/\theta^\prime} =\\
c_{35} 2^{-\mn(\alpha -\lambda -(p^{-1} -q^{-1})_+ \e) r}
r^{(1 -1/\theta)(\cmn(\alpha -\lambda -(p^{-1} -q^{-1})_+ \e) -1)}, r \in \N.
\end{multline*}

Как показано в [13, (см. (2.3.3))], существует константа
$ c_{36}(d,\alpha,p,\theta) >0 $ такая, что для $ f \in
(\mathcal S_{p,\theta}^\alpha \mathcal B)^\prime(D) $ при $ r \in \N $ имеет
место неравенство
\begin{equation*} \tag{2.2.49}
\biggl( \sum_{\kappa \in \Z_+^d: (\kappa, \beta) > r} (2^{(\kappa,\alpha)}
\Omega^{\prime l \chi_{\s(\kappa)}}(f,
(c_{2} 2^{-\kappa})^{\s(\kappa)})_{L_p(D)})^\theta \biggr)^{1/\theta} \le
c_{36}.
\end{equation*}

Подставляя (2.2.48) и (2.2.49) в (2.2.47), приходим к неравенству
\begin{multline*} \tag{2.2.50}
\sum_{ \kappa \in \Z_+^d: (\kappa, \beta) > r}
2^{(\kappa, \lambda +(p^{-1} -q^{-1})_+ \e)}
\Omega^{\prime l \chi_{\s(\kappa)}}(f, (c_{2} 2^{-\kappa})^{\s(\kappa)})_{L_p(D)} \le\\
c_{37} 2^{-\mn(\alpha -\lambda -(p^{-1} -q^{-1})_+ \e) r}
r^{(1 -1/\theta)(\cmn(\alpha -\lambda -(p^{-1} -q^{-1})_+ \e) -1)},
f \in (\mathcal S_{p,\theta}^\alpha \mathcal B)^\prime(D), r \in \N.
\end{multline*}

Продолжая оценку правой части (2.2.46), рассмотрим два случая.
В первом случае, когда $ \theta > p, $ для $ f \in
(\mathcal S_{p,\theta}^\alpha \mathcal B)^\prime(D), J \subset \Nu_{1,d}^1: J \ne \emptyset,
J^\prime \subset \Nu_{1,d}^1: J^\prime \ne \emptyset, r \in \N, $ благодаря
неравенству Гёльдера, с учётом (2.2.48) имеем
\begin{multline*} \tag{2.2.51}
\sum_{ \kappa \in \Z_+^d: (\kappa, \beta) > r, \s(\kappa) = J^\prime}
2^{(\kappa, \lambda +(p^{-1} -q^{-1})_+ \e)}
\biggl(\int_{ (c_{3} I^d)^J}
(\prod_{j \in J} u_j^{-p(\alpha_j -\bm \epsilon_j) -1}) \times \\
(\Omega^{\prime l \chi_{J \cup J^\prime}}(f, (u \chi_{J \setminus J^\prime} +
2^{-\kappa} u \chi_{J^\prime \cap J} +c_{2} 2^{-\kappa}
\chi_{J^\prime \setminus J})^{J \cup J^\prime})_{L_p(D)})^p du^J \biggr)^{1/p} = \\
\sum_{ \kappa \in \Z_+^d: (\kappa, \beta) > r, \s(\kappa) = J^\prime}
2^{-(\kappa, \alpha -\lambda -(p^{-1}  -q^{-1})_+ \e)}
\biggl(\int_{ (c_{3} I^d)^J}
(\prod_{j \in J} u_j^{-p(\alpha_j -\bm \epsilon_j) -1}) \times \\
(2^{(\kappa^{J^\prime}, \alpha^{J^\prime})}
\Omega^{\prime l \chi_{J \cup J^\prime}}(f, (u \chi_{J \setminus J^\prime} +
2^{-\kappa} u \chi_{J^\prime \cap J} +c_{2} 2^{-\kappa}
\chi_{J^\prime \setminus J})^{J \cup J^\prime})_{L_p(D)})^p du^J \biggr)^{1/p} \le \\
\biggl(\sum_{ \kappa \in \Z_+^d: (\kappa, \beta) > r, \s(\kappa) = J^\prime}
2^{-(\kappa, \alpha -\lambda -(p^{-1}  -q^{-1})_+ \e)
\theta^\prime} \biggr)^{1/\theta^\prime} \times \\
 \biggl(\sum_{ \kappa \in \Z_+^d: (\kappa, \beta) > r, \s(\kappa) = J^\prime}
(\int_{ (c_{3} I^d)^J} (\prod_{j \in J} u_j^{-p(\alpha_j -\bm \epsilon_j) -1}) \times \\
(2^{(\kappa^{J^\prime}, \alpha^{J^\prime})}
\Omega^{\prime l \chi_{J \cup J^\prime}}(f, (u \chi_{J \setminus J^\prime} +
2^{-\kappa} u \chi_{J^\prime \cap J} +c_{2} 2^{-\kappa}
\chi_{J^\prime \setminus J})^{J \cup J^\prime})_{L_p(D)})^p
du^J)^{\theta /p} \biggr)^{1 /\theta} \le \\
\biggl(\sum_{ \kappa \in \Z_+^d: (\kappa, \beta) > r}
2^{-(\kappa, \alpha -\lambda -(p^{-1}  -q^{-1})_+ \e)
\theta^\prime} \biggr)^{1 /\theta^\prime} \\
\biggl(\sum_{ \kappa \in \Z_+^d: \s(\kappa) = J^\prime}
(\int_{ (c_{3} I^d)^J} (\prod_{j \in J} u_j^{-p(\alpha_j -\bm \epsilon_j) -1}) \times \\
(2^{(\kappa^{J^\prime}, \alpha^{J^\prime})}
\Omega^{\prime l \chi_{J \cup J^\prime}}(f, (u \chi_{J \setminus J^\prime} +
2^{-\kappa} u \chi_{J^\prime \cap J} +c_{2} 2^{-\kappa}
\chi_{J^\prime \setminus J})^{J \cup J^\prime})_{L_p(D)})^p
du^J)^{\theta /p} \biggr)^{1 /\theta} \le \\
c_{35} 2^{-\mn(\alpha -\lambda -(p^{-1}  -q^{-1})_+ \e) r}
r^{(1 -1 /\theta)(\cmn(\alpha -\lambda -(p^{-1} -q^{-1})_+ \e) -1)} \times \\
\biggl(\sum_{ \kappa \in \Z_+^d: \s(\kappa) = J^\prime}
\biggl(\int_{ (c_{3} I^d)^J}
(\prod_{j \in J} u_j^{-p(\alpha_j -\bm \epsilon_j) -1})
(2^{(\kappa^{J^\prime}, \alpha^{J^\prime})}
\Omega^{\prime l \chi_{J \cup J^\prime}}(f, (u \chi_{J \setminus J^\prime} + \\
2^{-\kappa} u \chi_{J^\prime \cap J} +c_{2} 2^{-\kappa}
\chi_{J^\prime \setminus J})^{J \cup J^\prime})_{L_p(D)})^p
du^J \biggr)^{\theta /p} \biggr)^{1 /\theta}.
\end{multline*}

При оценке слагаемых в правой части (2.2.51), пользуясь неравенством
Гёльдера с показателем $ \theta /p >1, $ для $ f \in (\mathcal S_{p,\theta}^\alpha \mathcal B)^\prime(D) $
и $ J \subset \Nu_{1,d}^1: J \ne \emptyset,
J^\prime \subset \Nu_{1,d}^1: J^\prime \ne \emptyset,
\kappa \in \Z_+^d: \s(\kappa) = J^\prime, $ выводим
\begin{multline*} \tag{2.2.52}
\int_{ (c_{3} I^d)^J} (\prod_{j \in J} u_j^{-p(\alpha_j -\bm \epsilon_j) -1})
( 2^{(\kappa^{J^\prime}, \alpha^{J^\prime})} \times \\
\Omega^{\prime l \chi_{J \cup J^\prime}}(f, (u \chi_{J \setminus J^\prime} +
2^{-\kappa} u \chi_{J^\prime \cap J} +c_{2} 2^{-\kappa}
\chi_{J^\prime \setminus J})^{J \cup J^\prime})_{L_p(D)})^p du^J = \\
\int_{ (c_{3} I^d)^J}
(\prod_{j \in J} u_j^{\frac{1}{2} p \bm \epsilon_j -(\theta -p) /\theta})
(\prod_{j \in J} u_j^{-p(\alpha_j -\frac{1}{2} \bm \epsilon_j) -p /\theta})
( 2^{(\kappa^{J^\prime}, \alpha^{J^\prime})} \times \\
\Omega^{\prime l \chi_{J \cup J^\prime}}(f, (u \chi_{J \setminus J^\prime} +
2^{-\kappa} u \chi_{J^\prime \cap J} +c_{2} 2^{-\kappa}
\chi_{J^\prime \setminus J})^{J \cup J^\prime})_{L_p(D)})^p du^J \le \\
\biggl(\int_{ (c_{3} I^d)^J} (\prod_{j \in J}
u_j^{\frac{1}{2} p \bm \epsilon_j -(\theta -p) /\theta})^{\theta /(\theta -p)}
du^J \biggr)^{(\theta -p) /\theta} \times \\
\biggl(\int_{ (c_{3} I^d)^J} \biggl((\prod_{j \in J}
u_j^{-p(\alpha_j -\frac{1}{2} \bm \epsilon_j) -p /\theta})
( 2^{(\kappa^{J^\prime}, \alpha^{J^\prime})} \times \\
\Omega^{\prime l \chi_{J \cup J^\prime}}(f, (u \chi_{J \setminus J^\prime} +
2^{-\kappa} u \chi_{J^\prime \cap J} +c_{2} 2^{-\kappa}
\chi_{J^\prime \setminus J})^{J \cup J^\prime})_{L_p(D)})^p \biggr)^{\theta /p}
du^J \biggr)^{p /\theta} = \\
\biggl(\int_{ (c_{3} I^d)^J} (\prod_{j \in J}
u_j^{\frac{1}{2} p \bm \epsilon_j \theta /(\theta -p) -1})
du^J \biggr)^{(\theta -p) /\theta}
\biggl(\int_{ (c_{3} I^d)^J} (\prod_{j \in J}
u_j^{-\theta (\alpha_j -\frac{1}{2} \bm \epsilon_j) -1})
( 2^{(\kappa^{J^\prime}, \alpha^{J^\prime})} \times \\
\Omega^{\prime l \chi_{J \cup J^\prime}}(f, (u \chi_{J \setminus J^\prime} +
2^{-\kappa} u \chi_{J^\prime \cap J} +c_{2} 2^{-\kappa}
\chi_{J^\prime \setminus J})^{J \cup J^\prime})_{L_p(D)} )^\theta
du^J \biggr)^{p /\theta} \le \\
c_{38} \biggl(\int_{ (c_{3} I^d)^J}
(\prod_{j \in J} u_j^{-\theta (\alpha_j -\frac{1}{2} \bm \epsilon_j) -1})
( 2^{(\kappa^{J^\prime}, \alpha^{J^\prime})} \times \\
\Omega^{\prime l \chi_{J \cup J^\prime}}(f, (u \chi_{J \setminus J^\prime} +
2^{-\kappa} u \chi_{J^\prime \cap J} +c_{2} 2^{-\kappa}
\chi_{J^\prime \setminus J})^{J \cup J^\prime})_{L_p(D)} )^\theta
du^J \biggr)^{p /\theta}.
\end{multline*}

Из (2.2.52) с учётом замечания после леммы 1.2.1 и теоремы Лебега о
предельном переходе под знаком интеграла вытекает, что
при $ J \subset \Nu_{1,d}^1: J \ne \emptyset,
J^\prime \subset \Nu_{1,d}^1: J^\prime \ne \emptyset, $ для $ f \in
(\mathcal S_{p,\theta}^\alpha \mathcal B)^\prime(D) $ выполняется неравенство
\begin{multline*} \tag{2.2.53}
\biggl(\sum_{ \kappa \in \Z_+^d: \s(\kappa) = J^\prime}
\biggl(\int_{ (c_{3} I^d)^J}
(\prod_{j \in J} u_j^{-p(\alpha_j -\bm \epsilon_j) -1})
( 2^{(\kappa^{J^\prime}, \alpha^{J^\prime})} \times \\
\Omega^{\prime l \chi_{J \cup J^\prime}}(f, (u \chi_{J \setminus J^\prime} +
2^{-\kappa} u \chi_{J^\prime \cap J} +c_{2} 2^{-\kappa}
\chi_{J^\prime \setminus J})^{J \cup J^\prime})_{L_p(D)})^p
du^J \biggr)^{\theta /p} \biggr)^{1/\theta} \le \\
c_{39} \biggl(\sum_{ \kappa \in \Z_+^d: \s(\kappa) = J^\prime}
\int_{ (c_{3} I^d)^J}
(\prod_{j \in J} u_j^{-\theta (\alpha_j -\frac{1}{2} \bm \epsilon_j) -1})
( 2^{(\kappa^{J^\prime}, \alpha^{J^\prime})} \times \\
\Omega^{\prime l \chi_{J \cup J^\prime}}(f, (u \chi_{J \setminus J^\prime} +
2^{-\kappa} u \chi_{J^\prime \cap J} +c_{2} 2^{-\kappa}
\chi_{J^\prime \setminus J})^{J \cup J^\prime})_{L_p(D)} )^\theta
du^J \biggr)^{1 /\theta} = \\
c_{39} \biggl(\int_{ (c_{3} I^d)^J}
(\prod_{j \in J} u_j^{-\theta (\alpha_j -\frac{1}{2} \bm \epsilon_j) -1})
\biggl(\sum_{ \kappa \in \Z_+^d: \s(\kappa) = J^\prime}
( 2^{(\kappa^{J^\prime}, \alpha^{J^\prime})} \times \\
\Omega^{\prime l \chi_{J \cup J^\prime}}(f, (u \chi_{J \setminus J^\prime} +
2^{-\kappa} u \chi_{J^\prime \cap J} +c_{2} 2^{-\kappa}
\chi_{J^\prime \setminus J})^{J \cup J^\prime})_{L_p(D)} )^\theta \biggr)
du^J \biggr)^{1 /\theta}.
\end{multline*}

Далее, учитывая, что при $ J \subset \Nu_{1,d}^1: J \ne \emptyset,
J^\prime \subset \Nu_{1,d}^1: J^\prime \ne \emptyset,
\kappa \in \Z_+^d: \s(\kappa) = J^\prime,
u \in \R_+^d $ для $ f \in (\mathcal S_{p,\theta}^\alpha \mathcal B)^\prime(D) $
соблюдается неравенство
\begin{multline*}
( 2^{(\kappa^{J^\prime}, \alpha^{J^\prime})}
\Omega^{\prime l \chi_{J \cup J^\prime}}(f, (u \chi_{J \setminus J^\prime} +
2^{-\kappa} u \chi_{J^\prime \cap J} +c_{2} 2^{-\kappa}
\chi_{J^\prime \setminus J})^{J \cup J^\prime})_{L_p(D)} )^\theta = \\
(u^{J \setminus J^\prime})^{\theta \alpha^{J \setminus J^\prime}}
((u^{J \setminus J^\prime})^{-\alpha^{J \setminus J^\prime}}
(2^{-\kappa^{J^\prime}})^{-\alpha^{J^\prime}} \times \\
\Omega^{\prime l \chi_{J \cup J^\prime}}(f, (u \chi_{J \setminus J^\prime} +
2^{-\kappa} u \chi_{J^\prime \cap J} +c_{2} 2^{-\kappa}
\chi_{J^\prime \setminus J})^{J \cup J^\prime})_{L_p(D)} )^\theta = \\
(u^{J \setminus J^\prime})^{\theta \alpha^{J \setminus J^\prime}}
\int_{(u +u I^d)^{J \setminus J^\prime} \times (2^{-\kappa} +
2^{-\kappa} I^d)^{J^\prime}}
(u^{J \setminus J^\prime})^{-\e^{J \setminus J^\prime} -\theta
\alpha^{J \setminus J^\prime}}
(2^{-\kappa^{J^\prime}})^{-\e^{J^\prime} -\theta \alpha^{J^\prime}} \times \\
(\Omega^{\prime l \chi_{J \cup J^\prime}}(f, (u \chi_{J \setminus J^\prime} +
2^{-\kappa} u \chi_{J^\prime \cap J} +c_{2} 2^{-\kappa}
\chi_{J^\prime \setminus J})^{J \cup J^\prime})_{L_p(D)} )^\theta
d\tau^{J \cup J^\prime} \le \\
c_{40} (u^{J \setminus J^\prime})^{\theta \alpha^{J \setminus J^\prime}}
\int_{(u +u I^d)^{J \setminus J^\prime} \times (2^{-\kappa} +
2^{-\kappa} I^d)^{J^\prime}}
(\tau^{J \setminus J^\prime})^{-\e^{J \setminus J^\prime} -\theta
\alpha^{J \setminus J^\prime}}
(\tau^{J^\prime})^{-\e^{J^\prime} -\theta \alpha^{J^\prime}} \times \\
(\Omega^{\prime l \chi_{J \cup J^\prime}}(f, (\tau \chi_{J \setminus J^\prime} +
\tau u \chi_{J^\prime \cap J} +c_{2} \tau
\chi_{J^\prime \setminus J})^{J \cup J^\prime})_{L_p(D)} )^\theta
d\tau^{J \cup J^\prime} \le \\
c_{40} (u^{J \setminus J^\prime})^{\theta \alpha^{J \setminus J^\prime}}
\int_{(\R_+^d)^{J \setminus J^\prime} \times (2^{-\kappa} +
2^{-\kappa} I^d)^{J^\prime}}
(\tau^{J \cup J^\prime})^{-\e^{J \cup J^\prime} -\theta
\alpha^{J \cup J^\prime}} \times \\
(\Omega^{\prime l \chi_{J \cup J^\prime}}(f, (\tau \chi_{J \setminus J^\prime} +
\tau u \chi_{J^\prime \cap J} +c_{2} \tau
\chi_{J^\prime \setminus J})^{J \cup J^\prime})_{L_p(D)} )^\theta
d\tau^{J \cup J^\prime},
\end{multline*}
заключаем, что
\begin{multline*} \tag{2.2.54}
\sum_{ \kappa \in \Z_+^d: \s(\kappa) = J^\prime}
( 2^{(\kappa^{J^\prime}, \alpha^{J^\prime})}
\Omega^{\prime l \chi_{J \cup J^\prime}}(f, (u \chi_{J \setminus J^\prime} +
2^{-\kappa} u \chi_{J^\prime \cap J} +c_{2} 2^{-\kappa}
\chi_{J^\prime \setminus J})^{J \cup J^\prime})_{L_p(D)} )^\theta \le \\
c_{40} (u^{J \setminus J^\prime})^{\theta \alpha^{J \setminus J^\prime}}
\sum_{ \kappa \in \Z_+^d: \s(\kappa) = J^\prime}
\int_{(\R_+^d)^{J \setminus J^\prime} \times (2^{-\kappa} +
2^{-\kappa} I^d)^{J^\prime}}
(\tau^{J \cup J^\prime})^{-\e^{J \cup J^\prime} -\theta
\alpha^{J \cup J^\prime}} \times \\
(\Omega^{\prime l \chi_{J \cup J^\prime}}(f, (\tau \chi_{J \setminus J^\prime} +
\tau u \chi_{J^\prime \cap J} +c_{2} \tau
\chi_{J^\prime \setminus J})^{J \cup J^\prime})_{L_p(D)} )^\theta
d\tau^{J \cup J^\prime} = \\
c_{40} (u^{J \setminus J^\prime})^{\theta \alpha^{J \setminus J^\prime}}
\int_{\cup_{ \kappa \in \Z_+^d: \s(\kappa) = J^\prime}
(\R_+^d)^{J \setminus J^\prime} \times (2^{-\kappa} +2^{-\kappa} I^d)^{J^\prime}}
(\tau^{J \cup J^\prime})^{-\e^{J \cup J^\prime} -\theta
\alpha^{J \cup J^\prime}} \times \\
(\Omega^{\prime l \chi_{J \cup J^\prime}}(f, (\tau \chi_{J \setminus J^\prime} +
\tau u \chi_{J^\prime \cap J} +c_{2} \tau
\chi_{J^\prime \setminus J})^{J \cup J^\prime})_{L_p(D)} )^\theta
d\tau^{J \cup J^\prime} = \\
c_{40} (u^{J \setminus J^\prime})^{\theta \alpha^{J \setminus J^\prime}}
\int_{ (\R_+^d)^{J \setminus J^\prime} \times (I^d)^{J^\prime}}
(\tau^{J \cup J^\prime})^{-\e^{J \cup J^\prime} -\theta
\alpha^{J \cup J^\prime}} \times \\
(\Omega^{\prime l \chi_{J \cup J^\prime}}(f, (\tau \chi_{J \setminus J^\prime} +
\tau u \chi_{J^\prime \cap J} +c_{2} \tau
\chi_{J^\prime \setminus J})^{J \cup J^\prime})_{L_p(D)} )^\theta
d\tau^{J \cup J^\prime} \le \\
c_{40} (u^{J \setminus J^\prime})^{\theta \alpha^{J \setminus J^\prime}}
\int_{ (\R_+^d)^{J \cup J^\prime}}
(\tau^{J \cup J^\prime})^{-\e^{J \cup J^\prime} -\theta
\alpha^{J \cup J^\prime}} \times \\
(\Omega^{\prime l \chi_{J \cup J^\prime}}(f, (\tau \chi_{J \setminus J^\prime} +
\tau u \chi_{J^\prime \cap J} +c_{2} \tau
\chi_{J^\prime \setminus J})^{J \cup J^\prime})_{L_p(D)} )^\theta
d\tau^{J \cup J^\prime} = \\
c_{40} (u^{J \setminus J^\prime})^{\theta \alpha^{J \setminus J^\prime}}
\int_{ (\R_+^d)^{J \cup J^\prime}}
(\prod_{j \in (J \setminus J^\prime)}
t_j^{-1 -\theta \alpha_j})
(\prod_{j \in (J \cap J^\prime)} (t_j /u_j)^{-1 -\theta \alpha_j}) \times \\
(\prod_{j \in (J^\prime \setminus J)}
(t_j /c_{2})^{-1 -\theta \alpha_j})
(\Omega^{\prime l \chi_{J \cup J^\prime}}(f,
t^{J \cup J^\prime})_{L_p(D)} )^\theta
(\prod_{j \in (J \cap J^\prime)}
(1 /u_j))
(\prod_{j \in (J^\prime \setminus J)}
(1 /c_{2}))
dt^{J \cup J^\prime} \le \\
c_{41} (u^{J \setminus J^\prime})^{\theta \alpha^{J \setminus J^\prime}}
\int_{ (\R_+^d)^{J \cup J^\prime}}
(\prod_{j \in (J \cap J^\prime)} u_j^{\theta \alpha_j})
(t^{J \cup J^\prime})^{-\e^{J \cup J^\prime} -\theta \alpha^{J \cup J^\prime}}\times \\
(\Omega^{\prime l \chi_{J \cup J^\prime}}(f,
t^{J \cup J^\prime})_{L_p(D)} )^\theta dt^{J \cup J^\prime} =
c_{41} (u^J)^{\theta \alpha^J}
\int_{ (\R_+^d)^{J \cup J^\prime}}
(t^{J \cup J^\prime})^{-\e^{J \cup J^\prime} -\theta \alpha^{J \cup J^\prime}}\times \\
(\Omega^{\prime l \chi_{J \cup J^\prime}}(f,
t^{J \cup J^\prime})_{L_p(D)} )^\theta dt^{J \cup J^\prime} \le
c_{41} (u^J)^{\theta \alpha^J}.
\end{multline*}

Подставляя оценку (2.2.54) в (2.2.53), находим, что
\begin{multline*} \tag{2.2.55}
\biggl(\sum_{ \kappa \in \Z_+^d: \s(\kappa) = J^\prime}
\biggl(\int_{ (c_{3} I^d)^J}
(\prod_{j \in J} u_j^{-p(\alpha_j -\bm \epsilon_j) -1})
( 2^{(\kappa^{J^\prime}, \alpha^{J^\prime})} \times \\
\Omega^{\prime l \chi_{J \cup J^\prime}}(f, (u \chi_{J \setminus J^\prime} +
2^{-\kappa} u \chi_{J^\prime \cap J} +c_{2} 2^{-\kappa}
\chi_{J^\prime \setminus J})^{J \cup J^\prime})_{L_p(D)})^p
du^J \biggr)^{\theta /p} \biggr)^{1 /\theta} \le \\
c_{39} \biggl(\int_{ (c_{3} I^d)^J}
(\prod_{j \in J} u_j^{-\theta (\alpha_j -\frac{1}{2} \bm \epsilon_j) -1})
c_{41} (u^J)^{\theta \alpha^J} du^J \biggr)^{1 /\theta} = \\
c_{42} \biggl(\int_{ (c_{3} I^d)^J}
(u^J)^{\frac{1}{2} \theta \bm \epsilon^J -\e^J} du^J \biggr)^{1 /\theta} \le c_{43}, \\
J \subset \Nu_{1,d}^1: J \ne \emptyset,
J^\prime \subset \Nu_{1,d}^1: J^\prime \ne \emptyset, f \in
(\mathcal S_{p,\theta}^\alpha \mathcal B)^\prime(D).
\end{multline*}

Соединяя (2.2.51) с (2.2.55), приходим к неравенству
\begin{multline*} \tag{2.2.56}
\sum_{ \kappa \in \Z_+^d: (\kappa, \beta) > r, \s(\kappa) = J^\prime}
2^{(\kappa, \lambda +(p^{-1} -q^{-1})_+ \e)}
\biggl(\int_{ (c_{3} I^d)^J}
(\prod_{j \in J} u_j^{-p(\alpha_j -\bm \epsilon_j) -1}) \times \\
(\Omega^{\prime l \chi_{J \cup J^\prime}}(f, (u \chi_{J \setminus J^\prime} +
2^{-\kappa} u \chi_{J^\prime \cap J} +c_{2} 2^{-\kappa}
\chi_{J^\prime \setminus J})^{J \cup J^\prime})_{L_p(D)})^p du^J \biggr)^{1 /p} \le \\
c_{44} 2^{-\mn(\alpha -\lambda -(p^{-1} -q^{-1})_+ \e) r}
r^{(1 -1/\theta)(\cmn(\alpha -\lambda -(p^{-1} -q^{-1})_+ \e) -1)}, \\
f \in (\mathcal S_{p,\theta}^\alpha \mathcal B)^\prime(D), r \in \N,
J \subset \Nu_{1,d}^1: J \ne \emptyset,
J^\prime \subset \Nu_{1,d}^1: J^\prime \ne \emptyset.
\end{multline*}

После подстановки (2.2.50) и (2.2.56) в (2.2.46) получаем, что для
$ f \in (\mathcal S_{p,\theta}^\alpha \mathcal B)^\prime(D) $ и $ r \in \N, $
удовлетворяющего (2.2.45), при $ \theta > p $ выполняется неравенство
\begin{multline*} \tag{2.2.57}
\| \D^\lambda (f \mid_U) -(A \circ \phi(f)) \mid_U \|_{L_q(U)} \le \\
c_{45} 2^{-\mn(\alpha -\lambda -(p^{-1} -q^{-1})_+ \e) r}
r^{(1 -1 /\theta)(\cmn(\alpha -\lambda -(p^{-1} -q^{-1})_+ \e) -1)}.
\end{multline*}

Проводя оценку правой части (2.2.46) в случае, когда $ \theta \le p, $
для $ f \in (\mathcal S_{p,\theta}^\alpha \mathcal B)^\prime(D), J \subset
\Nu_{1,d}^1: J \ne \emptyset, J^\prime \subset \Nu_{1,d}^1: J^\prime \ne
\emptyset, r \in \N, $ благодаря неравенству Гёльдера, ввиду (2.2.48)
с заменой в нём $ \theta $ на $ p, $ а также используя с учётом замечания после
леммы 1.2.1 теорему Лебега о предельном переходе под знаком интеграла, выводим
\begin{multline*} \tag{2.2.58}
\sum_{ \kappa \in \Z_+^d: (\kappa, \beta) > r, \s(\kappa) = J^\prime}
2^{(\kappa, \lambda +(p^{-1} -q^{-1})_+ \e)}
\biggl(\int_{ (c_{3} I^d)^J}
(\prod_{j \in J} u_j^{-p(\alpha_j -\bm \epsilon_j) -1}) \times \\
(\Omega^{\prime l \chi_{J \cup J^\prime}}(f, (u \chi_{J \setminus J^\prime} +
2^{-\kappa} u \chi_{J^\prime \cap J} +c_{2} 2^{-\kappa}
\chi_{J^\prime \setminus J})^{J \cup J^\prime})_{L_p(D)})^p du^J \biggr)^{1 /p} = \\
\sum_{ \kappa \in \Z_+^d: (\kappa, \beta) > r, \s(\kappa) = J^\prime}
2^{-(\kappa, \alpha -\lambda -(p^{-1} -q^{-1})_+ \e)}
\biggl(\int_{ (c_{3} I^d)^J}
(\prod_{j \in J} u_j^{-p(\alpha_j -\bm \epsilon_j) -1})
( 2^{(\kappa^{J^\prime}, \alpha^{J^\prime})} \times \\
\Omega^{\prime l \chi_{J \cup J^\prime}}(f, (u \chi_{J \setminus J^\prime} +
2^{-\kappa} u \chi_{J^\prime \cap J} +c_{2} 2^{-\kappa}
\chi_{J^\prime \setminus J})^{J \cup J^\prime})_{L_p(D)})^p du^J \biggr)^{1 /p} \le \\
\biggl(\sum_{ \kappa \in \Z_+^d: (\kappa, \beta) > r, \s(\kappa) = J^\prime}
2^{-(\kappa, \alpha -\lambda -(p^{-1} -q^{-1})_+ \e)
p^\prime} \biggr)^{1 /p^\prime} \times \\
\biggl(\sum_{ \kappa \in \Z_+^d: (\kappa, \beta) > r, \s(\kappa) = J^\prime}
\int_{ (c_{3} I^d)^J}
(\prod_{j \in J} u_j^{-p(\alpha_j -\bm \epsilon_j) -1})
( 2^{(\kappa^{J^\prime}, \alpha^{J^\prime})} \times \\
\Omega^{\prime l \chi_{J \cup J^\prime}}(f, (u \chi_{J \setminus J^\prime} +
2^{-\kappa} u \chi_{J^\prime \cap J} +c_{2} 2^{-\kappa}
\chi_{J^\prime \setminus J})^{J \cup J^\prime})_{L_p(D)})^p
du^J \biggr)^{1 /p} \le \\
\biggl(\sum_{ \kappa \in \Z_+^d: (\kappa, \beta) > r}
2^{-(\kappa, \alpha -\lambda -(p^{-1} -q^{-1})_+ \e)
p^\prime} \biggr)^{1 /p^\prime} \times \\
\biggl(\sum_{ \kappa \in \Z_+^d: \s(\kappa) = J^\prime}
\int_{ (c_{3} I^d)^J} (\prod_{j \in J} u_j^{-p(\alpha_j -\bm \epsilon_j) -1})
( 2^{(\kappa^{J^\prime}, \alpha^{J^\prime})} \times \\
\Omega^{\prime l \chi_{J \cup J^\prime}}(f, (u \chi_{J \setminus J^\prime} +
2^{-\kappa} u \chi_{J^\prime \cap J} +c_{2} 2^{-\kappa}
\chi_{J^\prime \setminus J})^{J \cup J^\prime})_{L_p(D)})^p du^J \biggr)^{1 /p} \le \\
c_{46} 2^{-\mn(\alpha -\lambda -(p^{-1} -q^{-1})_+ \e) r}
r^{(1 -1 /p)(\cmn(\alpha -\lambda -(p^{-1} -q^{-1})_+ \e) -1)} \times \\
\biggl(\sum_{ \kappa \in \Z_+^d: \s(\kappa) = J^\prime}
\int_{ (c_{3} I^d)^J} (\prod_{j \in J} u_j^{-p(\alpha_j -\bm \epsilon_j) -1})
( 2^{(\kappa^{J^\prime}, \alpha^{J^\prime})} \times \\
\Omega^{\prime l \chi_{J \cup J^\prime}}(f, (u \chi_{J \setminus J^\prime} +
2^{-\kappa} u \chi_{J^\prime \cap J} +c_{2} 2^{-\kappa}
\chi_{J^\prime \setminus J})^{J \cup J^\prime})_{L_p(D)})^p du^J \biggr)^{1 /p} = \\
c_{46} 2^{-\mn(\alpha -\lambda -(p^{-1} -q^{-1})_+ \e) r}
r^{(1 -1 /p)(\cmn(\alpha -\lambda -(p^{-1} -q^{-1})_+ \e) -1)} \times \\
\biggl(\int_{ (c_{3} I^d)^J}
(\prod_{j \in J} u_j^{-p(\alpha_j -\bm \epsilon_j) -1})
(\sum_{ \kappa \in \Z_+^d: \s(\kappa) = J^\prime}
( 2^{(\kappa^{J^\prime}, \alpha^{J^\prime})} \times \\
\Omega^{\prime l \chi_{J \cup J^\prime}}(f, (u \chi_{J \setminus J^\prime} +
2^{-\kappa} u \chi_{J^\prime \cap J} +c_{2} 2^{-\kappa}
\chi_{J^\prime \setminus J})^{J \cup J^\prime})_{L_p(D)})^p )
du^J \biggr)^{1 /p}.
\end{multline*}

Применяя неравенство (1.1.2) при $ a = \theta /p \le 1, $
а также неравенство (2.2.54), для $ f \in (\mathcal S_{p,\theta}^\alpha \mathcal B)^\prime(D), J \subset
\Nu_{1,d}^1: J \ne \emptyset, J^\prime \subset \Nu_{1,d}^1: J^\prime \ne \emptyset,
u \in \R_+^d $ имеем
\begin{multline*} \tag{2.2.59}
\sum_{ \kappa \in \Z_+^d: \s(\kappa) = J^\prime}
( 2^{(\kappa^{J^\prime}, \alpha^{J^\prime})}
\Omega^{\prime l \chi_{J \cup J^\prime}}(f, (u \chi_{J \setminus J^\prime} +
2^{-\kappa} u \chi_{J^\prime \cap J} +c_{2} 2^{-\kappa}
\chi_{J^\prime \setminus J})^{J \cup J^\prime})_{L_p(D)})^p \le \\
(\sum_{ \kappa \in \Z_+^d: \s(\kappa) = J^\prime}
( 2^{(\kappa^{J^\prime}, \alpha^{J^\prime})}
\Omega^{\prime l \chi_{J \cup J^\prime}}(f, (u \chi_{J \setminus J^\prime} +
2^{-\kappa} u \chi_{J^\prime \cap J} +c_{2} 2^{-\kappa}
\chi_{J^\prime \setminus J})^{J \cup J^\prime})_{L_p(D)})^\theta )^{p /\theta} \le \\
(c_{41} (u^J)^{\theta \alpha^J} )^{p /\theta} = c_{47} (u^J)^{p \alpha^J}.
\end{multline*}

Из (2.2.59) вытекает, что при $ J \subset \Nu_{1,d}^1: J \ne \emptyset,
J^\prime \subset \Nu_{1,d}^1: J^\prime \ne \emptyset, \theta \le p, $
для $ f \in (\mathcal S_{p,\theta}^\alpha \mathcal B)^\prime(D) $ имеет
место неравенство
\begin{multline*} \tag{2.2.60}
\biggl(\int_{ (c_{3} I^d)^J}
(\prod_{j \in J} u_j^{-p(\alpha_j -\bm \epsilon_j) -1})
\biggl(\sum_{ \kappa \in \Z_+^d: \s(\kappa) = J^\prime}
( 2^{(\kappa^{J^\prime}, \alpha^{J^\prime})} \times \\
\Omega^{\prime l \chi_{J \cup J^\prime}}(f, (u \chi_{J \setminus J^\prime} +
2^{-\kappa} u \chi_{J^\prime \cap J} +c_{2} 2^{-\kappa}
\chi_{J^\prime \setminus J})^{J \cup J^\prime})_{L_p(D)})^p \biggr)
du^J \biggr)^{1/p} \le \\
\biggl(\int_{ (c_{3} I^d)^J}
(\prod_{j \in J} u_j^{-p(\alpha_j -\bm \epsilon_j) -1})
c_{47} (u^J)^{p \alpha^J} du^J \biggr)^{1 /p} = \\
c_{48} \biggl(\int_{ (c_{3} I^d)^J}
(u^J)^{p \bm \epsilon^J -\e^J} du^J \biggr)^{1 /p} \le c_{49}.
\end{multline*}

Подставляя оценку (2.2.60) в (2.2.58), находим, что для $ f \in (\mathcal S_{p,\theta}^\alpha \mathcal B)^\prime(D),
\theta \le p, $ при$ J \subset \Nu_{1,d}^1: J \ne \emptyset,
J^\prime \subset \Nu_{1,d}^1: J^\prime \ne \emptyset, r \in \N $
соблюдается неравенство
\begin{multline*} \tag{2.2.61}
\sum_{ \kappa \in \Z_+^d: (\kappa, \beta) > r, \s(\kappa) = J^\prime}
2^{(\kappa, \lambda +(p^{-1} -q^{-1})_+ \e)}
\biggl(\int_{ (c_{3} I^d)^J} (\prod_{j \in J} u_j^{-p(\alpha_j -\bm \epsilon_j) -1}) \times \\
(\Omega^{\prime l \chi_{J \cup J^\prime}}(f, (u \chi_{J \setminus J^\prime} +
2^{-\kappa} u \chi_{J^\prime \cap J} +c_{2} 2^{-\kappa}
\chi_{J^\prime \setminus J})^{J \cup J^\prime})_{L_p(D)})^p du^J \biggr)^{1 /p} \le \\
c_{50} 2^{-\mn(\alpha -\lambda -(p^{-1} -q^{-1})_+ \e) r}
r^{(1 -1 /p)(\cmn(\alpha -\lambda -(p^{-1} -q^{-1})_+ \e) -1)}.
\end{multline*}

Соединяя (2.2.50) и (2.2.61) с (2.2.46), получаем, что для
$ f \in (\mathcal S_{p,\theta}^\alpha \mathcal B)^\prime(D) $ и $ r \in \N, $
удовлетворяющего (2.2.45), при $ \theta \le p $ справедливо неравенство
\begin{multline*} \tag{2.2.62}
\| \D^\lambda (f \mid_U) -(A \circ \phi(f)) \mid_U \|_{L_q(U)} \le
c_{51} 2^{-\mn(\alpha -\lambda -(p^{-1} -q^{-1})_+ \e) r}
r^{(1 -1/ \theta)(\cmn(\alpha -\lambda -(p^{-1} -q^{-1})_+ \e) -1)} + \\
c_{52} 2^{-\mn(\alpha -\lambda -(p^{-1} -q^{-1})_+ \e) r}
r^{(1 -1 /p)(\cmn(\alpha -\lambda -(p^{-1} -q^{-1})_+ \e) -1)} \le \\
c_{53} 2^{-\mn(\alpha -\lambda -(p^{-1} -q^{-1})_+ \e) r}
r^{(1 -1 /p)(\cmn(\alpha -\lambda -(p^{-1} -q^{-1})_+ \e) -1)}.
\end{multline*}

Из (2.2.57), (2.2.62), (2.2.45)   следует (2.2.44). $ \square $
Отметим, что в проведенном доказательстве всюду при выводе
соотношений, содержащих ряды, неявно утверждается, что "если правая часть
выводимого соотношения определена, то левая его часть также определена и
справедливо выводимое соотношение."
Для завершения вывода оценки, объявленной в названии параграфа,
приведём ещё некоторые сведения из [13].

При $ d \in \N $ для $ \sigma \in \Sigma^d $ обозначим через $ \bm h_\sigma $
отображение, которое каждой функции $ f, $ заданной на некотором множестве
$ S \subset \R^d, $ ставит в соответствие функцию $ \bm h_\sigma f, $
определяемую на множестве $ \{ x \in \R^d: \sigma x \in S\} = \sigma^{-1} S =
\sigma S $ равенством $ (\bm h_\sigma f)(x) = f(\sigma x). $
Так как для $ \sigma \in \Sigma^d $ отображение
$ \R^d \ni x \mapsto \sigma x \in \R^d $ --- взаимно однозначно, то отображение
$ \bm h_\sigma $ является биекцией на себя множества всех функций с
областью определения в $ \R^d. $
При этом обратное отображение $ \bm h_\sigma^{-1} $ для $ f: S \mapsto \R $
задаётся равенством
\begin{equation*} \tag{2.2.63}
(\bm h_\sigma^{-1} f)(x) = f(\sigma^{-1} x) = f(\sigma x) =
(\bm h_\sigma f)(x), x \in \sigma S.
\end{equation*}

Отметим некоторые полезные для нас свойства отображений $ \bm h_\sigma. $
При $ d \in \N, \sigma \in \Sigma^d $ для любых множеств $ S \subset S^\prime
\subset \R^d $ и любой функции
$ f: S^\prime \mapsto \R $ верно равенство
\begin{equation*} \tag{2.2.64}
\bm h_\sigma (f \mid_S) = (\bm h_\sigma f) \mid_{\sigma^{-1} S}.
\end{equation*}

При $ d \in \N, \sigma \in \Sigma^d $ для открытого множества $ D \subset \R^d,
1 \le p \le \infty $ и $ f \in L_p(D) $ имеет место равенство
\begin{equation*} \tag{2.2.65}
\| \bm h_\sigma f \|_{L_p(\sigma^{-1} D)} = \| f \|_{L_p(D)},
\text{ а, значит, } \bm h_\sigma \in \mathcal B(L_p(D),L_p(\sigma^{-1} D)).
\end{equation*}

Лемма 2.2.5

Пусть $ d \in \N, \sigma \in \Sigma^d, D $ -- открытое множество в $ \R^d,
1 \le p, q \le \infty, \lambda \in \Z_+^d $ и $ f \in L_p(D),
\D^\lambda f \in L_q(D). $ Тогда соблюдается равенство
\begin{equation*} \tag{2.2.66}
\D^\lambda (\bm h_\sigma f) = \sigma^\lambda \bm h_\sigma (\D^\lambda f).
\end{equation*}

Лемма 2.2.6

Пусть $ d \in \N, D $ -- открытое множество в $ \R^d, 1 \le p < \infty,
\sigma \in \Sigma^d. $ Тогда для
$ f \in L_p(D), l \in \N^d, J \subset \Nu_{1,d}^1: J \ne \emptyset, $ при $ t^J
\in (\R_+^d)^J $ выполняется равенство
\begin{equation*}
\Omega^{\prime l \chi_{J}} (\bm h_\sigma f, t^{J})_{L_p(\sigma^{-1} D)} =
\Omega^{\prime l \chi_{J}} (f, t^{J})_{L_p(D)},
\end{equation*}
и
\begin{multline*} \tag{2.2.67}
\bm h_\sigma ((\mathcal S_{p,\theta}^\alpha \mathcal B)^\prime(D))
\subset (\mathcal S_{p,\theta}^\alpha \mathcal B)^\prime(\sigma^{-1} D),
\alpha \in \R_+^d, 1 \le p < \infty,  1 \le \theta \le \infty.
\end{multline*}

Теперь можно установить следующее утверждение.

Теорема 2.2.7

При $ d \in \N $ пусть $ D $ -- ограниченная область в $ \R^d, $ для 
которой существует система открытых подмножеств $ \{ U_i \subset D, i =1,\ldots,\mathcal I\} $
такая, что при $ i =1,\ldots,\mathcal I $ существуют $ \delta^i \in \R_+^d $ и
$ \sigma^i \in \Sigma^d, $ для которых имеет место включение
$ (U_i +\sigma^i \delta^i I^d) \subset D, $ и $ D = \cup_{i =1}^{\mathcal I} U_i. $
Пусть ещё $ \alpha \in \R_+^d, 1 \le p <
\infty, 1 \le \theta, q \le \infty, \lambda \in \Z_+^d $ и выполняются 
неравенства (2.1.9) и (2.1.6). Тогда существуют константы 
$ c_{54}(d,\alpha,p,\theta,q,\lambda,D) > 0 $ и $ N_0 \in \N $ 
такие, что для $ \mathcal V = \D^\lambda, D(\mathcal V) =
\{f \in C(D): \D^\lambda f \in L_q(D) \}, X = L_q(D), K = (\mathcal S_{p,\theta}^\alpha \mathcal B)^\prime(D) $
при $ N \ge N_0 $ имеет место неравенство
\begin{equation*} \tag{2.2.68}
\overline \sigma_N(\mathcal V,K,X) \le c_{54} N^{-\mn} (\log N)^{(\mn +1 -1/\max(p,\theta))(\cmn -1)}.
\end{equation*}

Доказательство.

Прежде всего заметим, что в условиях теоремы при $ i = 1,\ldots,\mathcal I, $
учитывая, что $ (U_i +\sigma^i \delta^i I^d) \subset D, $ а, значит,
$ ((\sigma^i)^{-1} U_i +\delta^i I^d) \subset (\sigma^i)^{-1} D, $ т.е. 
ограниченная область $ (\sigma^i)^{-1} D $ и её открытое подмножество 
$ (\sigma^i)^{-1} U_i $ удовлетворяют условиям леммы 2.2.1, а также 
выполняются и другие условия предложения 2.2.4.  
Полагая $ N_0 = \mathcal I \cdot \max_{i =1,\ldots.\mathcal I}
n_0(d,\alpha,p,\theta,q,\lambda,(\sigma^i)^{-1} D,(\sigma^i)^{-1} U_i) +\mathcal I $
(см. предложение 2.2.4), для $ N > N_0 $ выберем $ n \in \N, $ для которого
соблюдается неравенство
\begin{equation*} \tag{2.2.69}
\mathcal I n \le N < \mathcal I (n +1),
\end{equation*}
и при $ i =1,\ldots,\mathcal I $ на основании предложения 2.2.4 и (2.2.67)
построим отображения $ \phi_i \in \Phi_n(C((\sigma^i)^{-1} D)) $ и
$ A_i \in \overline {\mathcal A}^n(L_q((\sigma^i)^{-1} D)) $ такие, что для
любой функции $ f \in (\mathcal S_{p,\theta}^\alpha \mathcal B)^\prime(D) $
выполняется неравенство
\begin{multline*}
\| \D^\lambda ((\bm h_{\sigma^i} f) \mid_{(\sigma^i)^{-1} U_i}) -
(A_i \circ \phi_i (\bm h_{\sigma^i} f)) \mid_{(\sigma^i)^{-1} U_i}
\|_{L_q((\sigma^i)^{-1} U_i)} \le\\
c_{32}^i n^{-\mn} (\log n)^{(\mn +1 -1/\max(p,\theta))(\cmn -1)},
\end{multline*}
а, следовательно, в силу (2.2.63), (2.2.65) справедливо неравенство
\begin{multline*} \tag{2.2.70}
\| (\bm h_{\sigma^i})^{-1}
(\D^\lambda ((\bm h_{\sigma^i} f) \mid_{(\sigma^i)^{-1} U_i})) -
(\bm h_{\sigma^i})^{-1} ((A_i \circ \phi_i (\bm h_{\sigma^i} f)) \mid_{(\sigma^i)^{-1} U_i})
\|_{L_q(U_i)} \le\\
c_{32}^i n^{-\mn} (\log n)^{(\mn +1 -1/\max(p,\theta))(\cmn -1)}.
\end{multline*}
Принимая во внимание, что в силу (2.2.63), (2.2.64), (2.2.66) имеют место
равенства
\begin{multline*}
(\bm h_{\sigma^i})^{-1}
(\D^\lambda ((\bm h_{\sigma^i} f) \mid_{(\sigma^i)^{-1} U_i})) =
\bm h_{\sigma^i}
(\D^\lambda ((\bm h_{\sigma^i} f) \mid_{(\sigma^i)^{-1} U_i})) = \\
(\sigma^i)^{-\lambda} \D^\lambda (\bm h_{\sigma^i}
((\bm h_{\sigma^i} f) \mid_{(\sigma^i)^{-1} U_i})) =
(\sigma^i)^{\lambda} \D^\lambda (\bm h_{\sigma^i}
(\bm h_{\sigma^i} (f \mid_{U_i}))) = \\
(\sigma^i)^{\lambda} \D^\lambda ((\bm h_{\sigma^i})^{-1}
(\bm h_{\sigma^i} (f \mid_{U_i}))) =
(\sigma^i)^{\lambda} \D^\lambda (f \mid_{U_i}) =
(\sigma^i)^{\lambda} (\D^\lambda f) \mid_{U_i},
\end{multline*}
\begin{multline*}
(\bm h_{\sigma^i})^{-1} ((A_i \circ \phi_i (\bm h_{\sigma^i} f)) \mid_{(\sigma^i)^{-1} U_i}) =
(\bm h_{\sigma^i} (A_i \circ \phi_i (\bm h_{\sigma^i} f))) \mid_{U_i} = \\
(\bm h_{\sigma^i} (A_i (\phi_i (\bm h_{\sigma^i} f)))) \mid_{U_i} = 
((\bm h_{\sigma^i} \circ A_i) \circ (\phi_i \circ \bm h_{\sigma^i}) (f)) \mid_{U_i} = 
(\bm a_i \circ \bm \phi_i (f)) \mid_{U_i},
\end{multline*}
а в силу (2.2.65) справедливо включение
$$
\bm a_i = \bm h_{\sigma^i} \circ A_i \in \overline {\mathcal A}^n(L_q(D)),
$$
и
$$
\bm \phi_i = \phi_i \circ \bm h_{\sigma^i} \in \Phi_n(C(D)),
$$
из (2.2.70) получаем,что
\begin{equation*}
\| (\sigma^i)^{\lambda} (\D^\lambda f) \mid_{U_i} -
(\bm a_i \circ \bm \phi_i (f)) \mid_{U_i} \|_{L_q(U_i)} \le
c_{32}^i n^{-\mn} (\log n)^{(\mn +1 -1/\max(p,\theta))(\cmn -1)}
\end{equation*}
или
\begin{multline*} \tag{2.2.71}
\| (\D^\lambda f) \mid_{U_i} -(\sigma^i)^{\lambda}
(\bm a_i \circ \bm \phi_i (f)) \mid_{U_i} \|_{L_q(U_i)} \le
c_{32}^i n^{-\mn} (\log n)^{(\mn +1 -1/\max(p,\theta))(\cmn -1)},\\
f \in (\mathcal S_{p,\theta}^\alpha \mathcal B)^\prime(D), i =1,\ldots,\mathcal I.
\end{multline*}
Обозначая через $ \bm u^i = U_i \setminus (\cup_{j =1}^{i -1}
U_j), \chi^i = \chi_{\bm u^i}, i =1,\ldots,\mathcal I, $ и
учитывая, что $ \cup_{i =1}^{\mathcal I} \bm u^i = D, $
видим, что для $ x \in D $ справедливо равенство $ \sum_{i =1}^{\mathcal I} \chi^i(x) =1. $
Теперь, учитывая (2.2.69), возьмём отображения
$ A \in \overline {\mathcal A}^N(L_q(D)), \phi \in
\Phi_N(C(D)), $ для которых соблюдается равенство
$$
A \circ \phi(f) = \sum_{i =1}^{\mathcal I} \chi^i (\sigma^i)^{\lambda}
\bm a_i \circ \bm \phi_i (f), f \in C(D).
$$
Тогда, принимая во внимание сказанное и теорему 2.2.5 из [13] (см. также (1.1.8)), 
для $ f \in (\mathcal S_{p,\theta}^\alpha \mathcal B)^\prime(D) $ имеем
\begin{multline*}
\| \D^\lambda f -A \circ \phi(f) \|_{L_q(D)} =
\| (\sum_{i =1}^{\mathcal I} \chi^i) \D^\lambda f -
\sum_{i =1}^{\mathcal I} \chi^i (\sigma^i)^{\lambda}
\bm a_i \circ \bm \phi_i (f)\|_{L_q(D)} = \\
\| \sum_{i =1}^{\mathcal I} \chi^i (\D^\lambda f -
(\sigma^i)^{\lambda} \bm a_i \circ \bm \phi_i (f))\|_{L_q(D)} \le \\
\sum_{i =1}^{\mathcal I} \| \chi^i (\D^\lambda f -
(\sigma^i)^{\lambda} \bm a_i \circ \bm \phi_i (f))\|_{L_q(D)} = \\
\sum_{i =1}^{\mathcal I} \| (\chi^i (\D^\lambda f -(\sigma^i)^{\lambda}
\bm a_i \circ \bm \phi_i (f))) \mid_{U_i} \|_{L_q(U_i)} \le \\
\sum_{i =1}^{\mathcal I} \| (\D^\lambda f -(\sigma^i)^{\lambda}
\bm a_i \circ \bm \phi_i (f)) \mid_{U_i} \|_{L_q(U_i)} = \\
\sum_{i =1}^{\mathcal I} \| (\D^\lambda f) \mid_{U_i} -(\sigma^i)^{\lambda}
(\bm a_i \circ \bm \phi_i (f)) \mid_{U_i} \|_{L_q(U_i)}.
\end{multline*}
Из последнего неравенства в соединении с (2.2.71) и с учётом (2.2.69) приходим к
(2.2.68). $ \square $
\bigskip

\centerline{\S 3. Оценка снизу наилучшей точности восстановления
в $ L_q(D) $ производных $ \D^\lambda f$}
\centerline{по значениям в $n$ точках функций $ f$ из $ B((S_{p,\theta}^\alpha B)^\prime(D)) $}
\bigskip

3.1. В этом пункте напомним некоторые сведения, используемые при выводе нижней
оценки изучаемой величины.

При $ n \in \N $ и $ 1 \le p \le \infty $ через $ l_p^n $ обозначается
пространство $ \R^n $ с фиксированной в нём нормой
$$
\|x\|_{l_p^n} = \begin{cases} (\sum_{j =1}^n |x_j|^p)^{1/p} \text{ при } p < \infty; \\
\max_{j =1}^n |x_j| \text{ при } p = \infty, \ x \in \R^n.
\end{cases}
$$

При $ 1 \le p,q \le \infty, n \in \N $ справедливо равенство
\begin{equation*} \tag{3.1.1}
\max_{x \in B(l_p^n)} \|x\|_{l_q^n} = n^{(1 /q -1 /p)_+}.
\end{equation*}

При $ m,n \in \N, 1 \le p,q \le \infty $ через $ l_{p,q}^{m,n} $
обозначается пространство $ \R^{mn}, $ состоящее из векторов $ x =
\{ x_{i,j} \in \R: i =1,\ldots,m; j =1,\ldots,n\}, $ с нормой
$$
\| x \|_{l_{p,q}^{m,n}} = (\sum_{j =1}^n (\sum_{i =1}^m | x_{i,j} |^p )^{q /p} )^{1 /q}.
$$
При $ p = \infty $ и/или $ q = \infty $ в последнем соотношении следует сделать
соответствующие изменения.

Для $ d \in \N, r \in \Z_+ $ положим
$$
c(d,r) = \card \{ \kappa \in \Z_+^d: (\kappa, \e) = r \}.
$$

Отметим (см. [8]), что существуют константы $ c_{1}(d) >0, c_{2}(d) >0 $
такие, что при $ r \in \N $ выполняется неравенство
\begin{equation*} \tag{3.1.2}
c_{1} r^{d -1} \le c(d, r) \le c_{2} r^{d -1}.
\end{equation*}

Пусть $ C $ -- подмножество банахова пространства $X$ и
$ n \in \Z_+.$ Тогда $n$-поперечни\-ком по Колмогорову множества $C$
в пространстве $X$ называется величина
$$
d_n(C,X) = \inf_{M \in \mathcal M_n(X)} \sup_{x \in C} \inf_{y \in M} \|x -y\|_X,
$$
где $ \mathcal M_n(X) $ -- совокупность всех плоскостей $ M $ в $ X, $ у которых
$ \dim M \le n. $

$n$-поперечником по Гельфанду множества $ C $ в пространстве $ X $ называется
величина
$$
d^n(C,X) = \inf_{M \in \mathcal M^n(X)} \sup_{x \in C \cap M} \|x\|_X,
$$
где $ \mathcal M^n(X)$ -- множество всех замкнутых линейных подпространств $ M $
в $ X, $ у которых $ \codim M \le n. $

Будет полезно следующее предложение, являющееся частным случаем
соответствующего утверждения из [15].

Предложение 3.1.1

Пусть $ U: X \mapsto Y$ -- непрерывное линейное отображение
банахова пространства $ X $ в банахово пространство $ Y $ и $ L \subset X$ --
замкнутое линейное подпространство, для которого сужение $ U \mid_L $ является
гомеоморфизмом $ L $ на $ U(L), $ и пусть $ C \subset L$ -- некоторое
множество. Тогда при $ n \in \Z_+ $
имеет место неравенство
\begin{equation*} \tag{3.1.3}
d^n(U(C), Y) \le \| U \|_{\mathcal B(X,Y)} d^n(C,X).
\end{equation*}

Согласно [16], [17], [18] при $ 1 \le p,q \le \infty $ существует
константа $ c_3(p,q) >0 $ такая, что для $ n \in \N $ выполняется неравенство
\begin{multline*} \tag{3.1.4}
d_n(B(l_p^{2n}),l_q^{2n}) \ge c_3
\begin{cases}
n^{1 /q -1 /p}, & \text{ при $ q \le p $ или $ 2 \le p < q;$ } \\
n^{1 /q -1/2}, & \text{ при $ p < 2 \le q; $ } \\
1, & \text{ при $ p < q < 2.$ }
\end{cases}
\end{multline*}

Из результатов, содержащихся в [19], вытекает, что при
$ 1 \le p, q \le \infty, n, m \in \Z_+: n \le m, $ гельфандовский и
колмогоровский поперечники связаны равенством
\begin{equation*} \tag{3.1.5}
d^n(B(l_p^m), l_q^m) = d_n(B(l_{q^{\prime}}^m), l_{p^{\prime}}^m).
\end{equation*}

Из (3.1.4), (3.1.5) вытекает, что при $ 1 \le p, q \le \infty, n \in \N $ верно
неравенство
\begin{multline*} \tag{3.1.6}
d^n(B(l_p^{2n}), l_q^{2n}) =
d_n(B(l_{q^\prime}^{2n}),l_{p^\prime}^{2n}) \ge \\
c_4(p,q)
\begin{cases}
n^{1 /q -1 /p}, & \text{ при $ q \le p $ или $ p < q \le 2;$ } \\
n^{1/2 -1 /p}, & \text{ при $ p \le 2 < q; $ } \\
1, & \text{ при $ 2 < p < q. $}
\end{cases}
\end{multline*}

Также понадобится следующее утверждение.

Предложение 3.1.2

Пусть $ d, n \in \N, 1 < p,q < \infty. $ Тогда для любого набора
комплекснозначных функций $ \{f_j \in L_p(\R^d), j =1,\ldots,n\} $ имеет место
равенство
\begin{equation*} \tag{3.1.7}
(\sum_{j =1}^n \| f_j \|_{L_p(\R^d)}^q )^{1 /q} = 
\max_{\{ g_j \in L_{p^\prime}(\R^d), j =1,\ldots,n\}:
(\sum_{j =1}^n \| g_j \|_{L_{p^\prime}(\R^d)}^{q^\prime} )^{1 /q^\prime} \le 1}
\biggl| \int_{ \R^d} \sum_{j =1}^n f_j \overline {g_j} dx\biggr|.
\end{equation*}

Доказательство.

В условиях предложения для наборов функций $ \{f_j \in L_p(\R^d), j =1,\ldots,n\} $
и $ \{ g_j \in L_{p^\prime}(\R^d), j =1,\ldots,n\}:
(\sum_{j =1}^n \| g_j \|_{L_{p^\prime}(\R^d)}^{q^\prime} )^{1 /q^\prime} \le 1, $
в силу неравенства Гёльдера имеем
\begin{multline*}
| \int_{ \R^d} \sum_{j =1}^n f_j \overline {g_j} dx| =
| \sum_{j =1}^n \int_{\R^d} f_j \overline {g_j} dx| \le
\sum_{j =1}^n | \int_{\R^d} f_j \overline {g_j} dx| \le \\
\sum_{j =1}^n \int_{\R^d} | f_j | \cdot | g_j | dx \le
\sum_{j =1}^n \| f_j \|_{L_p(\R^d)} \cdot \| g_j \|_{L_{p^\prime}(\R^d)} \le \\
(\sum_{j =1}^n \| f_j \|_{L_p(\R^d)}^q )^{1 /q} \cdot
(\sum_{j =1}^n \| g_j \|_{L_{p^\prime}(\R^d)}^{q^\prime} )^{1 /q^\prime} \le
(\sum_{j =1}^n \| f_j \|_{L_p(\R^d)}^q )^{1 /q}.
\end{multline*}

С другой стороны, для набора функций $ \{f_j \in L_p(\R^d), j =1,\ldots,n\} $
выберем набор чисел $ \{ y_j \ge 0, j =1,\ldots,n\} $ такой, что
\begin{equation*} \tag{3.1.8}
(\sum_{j =1}^n y_j^{q^\prime} )^{1 /q^\prime} =1,
\sum_{j =1}^n \| f_j \|_{L_p(\R^d)} y_j =
(\sum_{j =1}^n \| f_j \|_{L_p(\R^d)}^q )^{1 /q}.
\end{equation*}
Кроме того, для каждого $ j =1, \ldots,n $ возьмём функцию $ \mathfrak g_j \in
L_{p^\prime}(\R^d) $ вида $ \mathfrak g_j = \lambda_j f_j | f_j|^{p -2},
\lambda_j \in \R, $ такую, что
\begin{equation*} \tag{3.1.9}
\| \mathfrak g_j \|_{L_{p^\prime}(\R^d)} =1, \\
\int_{\R^d} f_j \cdot \overline{\mathfrak g_j} dx = \| f_j \|_{L_p(\R^d)}.
\end{equation*}
Полагая $ g_j = y_j \mathfrak g_j, j =1,\ldots,n, $ и принимая во внимание
(3.1.8), (3.1.9), получаем
\begin{multline*}
(\sum_{j =1}^n \| g_j \|_{L_{p^\prime}(\R^d)}^{q^\prime} )^{1 /q^\prime} =
(\sum_{j =1}^n \| y_j \mathfrak g_j \|_{L_{p^\prime}(\R^d)}^{q^\prime} )^{1 /q^\prime} =
(\sum_{j =1}^n y_j^{q^\prime} )^{1 /q^\prime} = 1, \\
(\sum_{j =1}^n \| f_j \|_{L_p(\R^d)}^q )^{1 /q} =
\sum_{j =1}^n \| f_j \|_{L_p(\R^d)} y_j = \sum_{j =1}^n y_j \cdot
\int_{\R^d} f_j \cdot \overline{\mathfrak g_j} dx = \\
\sum_{j =1}^n \int_{\R^d} f_j \cdot \overline{(y_j \mathfrak g_j)} dx =
\sum_{j =1}^n \int_{\R^d} f_j \cdot \overline{g_j} dx =
| \sum_{j =1}^n \int_{\R^d} f_j \cdot \overline{g_j} dx | =
| \int_{\R^d} \sum_{j =1}^n f_j \cdot \overline{g_j} dx |.
\end{multline*}

Сопоставляя сказанное, приходим к (3.1.7). $ \square $

Лемма 3.1.3

При $ m,n \in \N, 1 \le p_0, p_1, q_0, q_1 \le \infty $ имеет место равенство
\begin{equation*} \tag{3.1.10}
\max_{x \in B(l_{p_0,q_0}^{m,n})} \|x\|_{l_{p_1,q_1}^{m,n}} =
m^{(1 /p_1 -1 /p_0)_+} n^{(1 /q_1 -1 /q_0)_+}.
\end{equation*}

Доказательство.

В условиях леммы, с одной стороны, для
$ x = \{ x_{i,j} \in \R: i =1,\ldots,m; j =1,\ldots,n\} $ в силу (3.1.1) имеет
место неравенство
\begin{multline*}
\| x \|_{l_{p_1,q_1}^{m,n}} = (\sum_{j =1}^n (\sum_{i =1}^m
| x_{i,j} |^p_1 )^{q_1 /p_1} )^{1 /q_1} \le \\
(\sum_{j =1}^n (m^{(1 /p_1 -1 /p_0)_+} (\sum_{i =1}^m
| x_{i,j} |^{p_0} )^{1 /p_0})^{q_1} )^{1 /q_1} = \\
m^{(1 /p_1 -1 /p_0)_+} (\sum_{j =1}^n
((\sum_{i =1}^m | x_{i,j} |^{p_0} )^{1 /p_0})^{q_1} )^{1 /q_1} \le \\
m^{(1 /p_1 -1 /p_0)_+} n^{(1 /q_1 -1 /q_0)_+} (\sum_{j =1}^n ((\sum_{i =1}^m
| x_{i,j} |^{p_0} )^{1 /p_0})^{q_0} )^{1 /q_0} = \\
m^{(1 /p_1 -1 /p_0)_+} n^{(1 /q_1 -1 /q_0)_+} \| x \|_{l_{p_0,q_0}^{m,n}}.
\end{multline*}
Откуда
\begin{equation*} \tag{3.1.11}
\sup_{x \in B(l_{p_0,q_0}^{m,n})} \|x\|_{l_{p_1,q_1}^{m,n}} \le
m^{(1 /p_1 -1 /p_0)_+} n^{(1 /q_1 -1 /q_0)_+}.
\end{equation*}

С другой стороны, для $ \xi = (\xi_1,\ldots,\xi_m) \in \R^m, \eta =
(\eta_1,\ldots,\eta_n) \in \R^n, x = \{ x_{i,j} = \xi_i \eta_j, i =1,\ldots,m;
j =1,\ldots,n\} \in \R^{mn} $ при $ 1 \le p,q \le \infty $ соблюдается
равенство
\begin{multline*} \tag{3.1.12}
\| x \|_{l_{p,q}^{m,n}} = (\sum_{j =1}^n (\sum_{i =1}^m
| x_{i,j} |^p )^{q /p} )^{1 /q} =
(\sum_{j =1}^n (\sum_{i =1}^m | \xi_i \eta_j |^p )^{q /p} )^{1 /q} = \\
(\sum_{j =1}^n (\sum_{i =1}^m | \xi_i|^p | \eta_j |^p )^{q /p} )^{1 /q} =
(\sum_{j =1}^n (| \eta_j |^p \sum_{i =1}^m | \xi_i|^p )^{q /p} )^{1 /q} = \\
(\sum_{j =1}^n | \eta_j |^q (\sum_{i =1}^m | \xi_i|^p )^{q /p} )^{1 /q} =
((\sum_{i =1}^m | \xi_i|^p )^{q /p} \sum_{j =1}^n | \eta_j |^q )^{1 /q} = \\
(\sum_{i =1}^m | \xi_i|^p )^{1 /p} (\sum_{j =1}^n | \eta_j |^q )^{1 /q} =
\| \xi \|_{l_p^m} \| \eta\|_{l_q^n}.
\end{multline*}

Теперь в условиях леммы, опираясь на (3.1.1), выберем вектор $ \xi \in \R^m $
такой, что $ \| \xi \|_{l_{p_0}^m} =1,
\| \xi \|_{l_{p_1}^m} = m^{(1 /p_1 -1 /p_0)_+}, $ и вектор $ \eta \in \R^n $
такой, что $ \| \eta \|_{l_{q_0}^n} =1, \| \eta\|_{l_{q_1}^n} = n^{(1 /q_1 -1 /q_0)_+}. $
Тогда для $ x = \{ x_{i,j} = \xi_i \eta_j, i =1,\ldots,m; j =1,\ldots,n\}
\in \R^{mn} $ в силу (3.1.12) получаем, что
\begin{equation*}
\| x \|_{l_{p_0,q_0}^{m,n}} = \| \xi \|_{l_{p_0}^m} \| \eta\|_{l_{q_0}^n} =1, \\
\| x \|_{l_{p_1,q_1}^{m,n}} = \| \xi \|_{l_{p_1}^m} \| \eta\|_{l_{q_1}^n} =
m^{(1 /p_1 -1 /p_0)_+} n^{(1 /q_1 -1 /q_0)_+}.
\end{equation*}
Сопоставляя эти равенства с (3.1.11), приходим к (3.1.10). $ \square $
\bigskip

3.2. В этом пункте приведём используемые ниже соотношения между нормами образов
и прообразов при некоторых отображениях рассматриваемых пространств.

При $ d \in \N $ для $ \delta \in \R_+^d $ и $ x^0 \in \R^d $ обозначим через
$ h_{\delta, x^0} $ отображение, которое каждой функции $ f, $ заданной на
некотором множестве $ S \subset \R^d, $  ставит в соответствие функцию
$ h_{\delta, x^0} f, $ определяемую на множестве $ \{ x \in \R^d: x^0 +\delta
x \in S\} = \delta^{-1} (S -x^0) $ равенством $ (h_{\delta, x^0} f)(x) =
f(x^0 +\delta x). $ Так как для $ \delta \in \R_+^d, x^0 \in \R^d $
отображение $ \R^d \ni x \mapsto x^0 +\delta x \in \R^d $ ---
взаимно однозначно, то отображение $ h_{\delta, x^0} $ является
биекцией на себя множества всех функций с областью определения в $ \R^d. $
При этом обратное отображение $ h_{\delta, x^0}^{-1} $
задаётся равенством
\begin{equation*} \tag{3.2.1}
(h_{\delta, x^0}^{-1} f)(x) = f(\delta^{-1} (x -x^0)) =
(h_{\delta^\prime, x^{\prime 0}} f)(x)  \text{ с } \delta^\prime = \delta^{-1},
x^{\prime 0} =-\delta^{-1} x^0.
\end{equation*}

Отметим, что при $ 1 \le p \le \infty $ для $ f \in L_p(x^0 +\delta D), $ где
$ D $ -- область в $ \R^d, \delta \in \R_+^d, x^0 \in \R^d, $ имеет место
равенство
\begin{equation*} \tag{3.2.2}
\| h_{\delta,x^0} f\|_{L_p(D)} = \delta^{-p^{-1} \e} \|f\|_{L_p(x^0 +\delta D)},
\end{equation*}
а, следовательно, для $ f \in L_p(D) $ выполняется равенство
\begin{equation*}
\| h_{\delta,x^0}^{-1} f\|_{L_p(x^0 +\delta D)} = \delta^{p^{-1} \e} \|f\|_{L_p(D)}.
\end{equation*}

Лемма  3.2.1

Пусть $ d \in \N, l \in \N^d, D $ -- область в $ \R^d, 1 \le p < \infty,
\delta \in \R_+^d, x^0 \in \R^d. $ Тогда при $ J \subset \{1,\ldots,d\}: J \ne
\emptyset, t \in \R_+^d $ для $ f \in L_p(x^0 +\delta D) $ справедливо
равенство
\begin{equation*} \tag{3.2.3}
\Omega^{\prime l \chi_J}((h_{\delta, x^0} f),t^J)_{L_p(D)} =
\delta^{-p^{-1} \e} \Omega^{\prime l \chi_J}(f, (\delta t)^J)_{L_p(x^0 +\delta D)}.
\end{equation*}

Доказательство.

В условиях леммы при $ J \subset \{1,\ldots,d\}: J \ne \emptyset $ для $ f \in
L_p(x^0 +\delta D), t \in \R_+^d, $ делая замену переменных
$ \eta = \delta^{-1} \xi, y = \delta^{-1}(x -x^0), $ находим, что соблюдается
равенство
\begin{multline*}
\int_{ (t B^d)^J} \int_{ D_\eta^{l \chi_J}} | (\Delta_{\eta}^{l \chi_J}
(h_{\delta, x^0} f))(y)|^p dy d\eta^J = \\
\int\limits_{\{(\eta^J,y): \eta^J \in (t B^d)^J, y \in D_\eta^{l \chi_J}\}}
\biggl| \sum_{k^J \in (\Z_+^d(l))^J} C_{l^J}^{k^J} (-\e^J)^{l^J -k^J}
(h_{\delta, x^0} f)(y +k \eta \chi_J)\biggr|^p d\eta^J dy = \\
\int\limits_{\{(\eta^J,y): \eta^J \in (t B^d)^J, y \in D_\eta^{l \chi_J}\}}
\biggl| \sum_{k^J \in (\Z_+^d(l))^J} C_{l^J}^{k^J} (-\e^J)^{l^J -k^J}
f(x^0 +\delta y +k \delta \eta \chi_J)\biggr|^p d\eta^J dy = \\
\int\limits_{\substack{\{ (\xi^J,x): \xi^J \in \delta^J (t B^d)^J, \\ x \in (x^0 +\delta D)_\xi^{l \chi_J} \}}}
\biggl| \sum_{k^J \in (\Z_+^d(l))^J} C_{l^J}^{k^J} (-\e^J)^{l^J -k^J}
f(x +k \xi \chi_J)\biggr|^p \delta^{-\e} (\delta^J)^{-\e^J} d\xi^J dx = \\
\delta^{-\e} (\delta^J)^{-\e^J} \int_{(\delta t B^d)^J} \int_{ (x^0 +\delta D)_\xi^{l \chi_J}}
| (\Delta_\xi^{l \chi_J} f)(x)|^p dx d\xi^J.
\end{multline*}
Откуда выводим
\begin{multline*}
\Omega^{\prime l \chi_J}((h_{\delta, x^0} f),t^J)_{L_p(D)} =
((2 t^J)^{-\e^J} \int_{ (t B^d)^J} \int_{D_\eta^{l \chi_J}}
| (\Delta_\eta^{l \chi_J} (h_{\delta, x^0} f))(y)|^p dy d\eta^J)^{1 /p} = \\
((2 t^J)^{-\e^J} \delta^{-\e} (\delta^J)^{-\e^J} \int_{(\delta t B^d)^J}
\int_{ (x^0 +\delta D)_\xi^{l \chi_J}}
| (\Delta_\xi^{l \chi_J} f)(x)|^p dx d\xi^J)^{1 /p} = \\
\delta^{-p^{-1} \e} ((2 \delta^J t^J)^{-\e^J} \int_{\delta^J t^J (B^d)^J}
\int_{ (x^0 +\delta D)_\xi^{l \chi_J}}
| (\Delta_\xi^{l \chi_J} f)(x)|^p dx d\xi^J)^{1 /p} = \\
\delta^{-p^{-1} \e} \Omega^{\prime l \chi_J}(f, \delta^J t^J)_{L_p(x^0 +\delta D)},
t \in \R_+^d, f \in L_p(x^0 +\delta D), J \subset \{1,\ldots,d\}: J \ne \emptyset,
\end{multline*}
что совпадает с (3.2.3). $ \square $

Лемма 3.2.2

Пусть $ d \in \N, D $ -- область в $ \R^d, \alpha \in \R_+^d, 1 \le p < \infty,
1 \le \theta \le \infty, \delta \in \R_+^d, x^0 \in \R^d. $ Тогда существуют
константы $ c_1(d,\alpha,p,\delta) > 0, c_2(d,\alpha,p,\delta) > 0 $
такие, что для любой функции $ f \in (S_{p,\theta}^\alpha B)^\prime(x^0 +\delta D) $
соблюдается неравенство
\begin{equation*} \tag{3.2.4}
\| h_{\delta, x^0} f \|_{(S_{p,\theta}^\alpha B)^\prime(D)} \le
c_1 \| f \|_{(S_{p,\theta}^\alpha B)^\prime(x^0 +\delta D)},
\end{equation*}
а для $ f \in (S_{p,\theta}^\alpha B)^\prime(D) $ выполняется неравенство
\begin{equation*} \tag{3.2.5}
\| h_{\delta, x^0}^{-1} f \|_{(S_{p,\theta}^\alpha B)^\prime(x^0 +\delta D)} \le
c_2 \| f \|_{(S_{p,\theta}^\alpha B)^\prime(D)}.
\end{equation*}

Доказательство.

В условиях леммы, полагая $ l = l(\alpha), $ для $ f \in
(S_{p,\theta}^\alpha B)^\prime(x^0 +\delta D) $ при $ J \subset
\{1,\ldots,d\}: J \ne \emptyset, $ в силу (3.2.3) функция $
(t^J)^{-\e^J -\theta \alpha^J} (\Omega^{\prime l
\chi_J}((h_{\delta, x^0} f),t^J)_{L_p(D)})^\theta $ суммируема на
$ (\R_+^d)^J, $ и справедливо соотношение
\begin{multline*} \tag{3.2.6}
\biggl(\int_{(\R_+^d)^J} (t^J)^{-\e^J -\theta \alpha^J}
(\Omega^{\prime l \chi_J}((h_{\delta, x^0} f),t^J)_{L_p(D)})^\theta dt^J \biggr)^{1 /\theta} = \\
\biggl(\int_{(\R_+^d)^J} (t^J)^{-\e^J -\theta \alpha^J} (\delta^{-p^{-1} \e}
\Omega^{\prime l \chi_J}(f, \delta^J t^J)_{L_p(x^0 +\delta D)})^\theta dt^J \biggr)^{1 /\theta} = \\
\delta^{-p^{-1} \e} \biggl(\int_{(\R_+^d)^J} (t^J)^{-\e^J -\theta \alpha^J}
(\Omega^{\prime l \chi_J}(f, \delta^J t^J)_{L_p(x^0 +\delta D)})^\theta dt^J \biggr)^{1 /\theta} = \\
\delta^{-p^{-1} \e} \biggl(\int_{(\R_+^d)^J} ((\delta^J)^{-1} \tau^J)^{-\e^J -\theta \alpha^J}
(\Omega^{\prime l \chi_J}(f, \delta^J (\delta^J)^{-1} \tau^J)_{L_p(x^0 +\delta D)})^\theta
(\delta^J)^{-\e^J} d\tau^J \biggr)^{1 /\theta} = \\
\delta^{-p^{-1} \e} \biggl(\int_{(\R_+^d)^J} (\delta^J)^{\theta \alpha^J} (\tau^J)^{-\e^J -\theta \alpha^J}
(\Omega^{\prime l \chi_J}(f, \tau^J)_{L_p(x^0 +\delta D)})^\theta d\tau^J \biggr)^{1 /\theta} = \\
\delta^{-p^{-1} \e} (\delta^J)^{\alpha^J} \biggl(\int_{(\R_+^d)^J} (\tau^J)^{-\e^J -\theta \alpha^J}
(\Omega^{\prime l \chi_J}(f, \tau^J)_{L_p(x^0 +\delta D)})^\theta d\tau^J \biggr)^{1 /\theta}.
\end{multline*}
Объединяя неравенства (3.2.6) и (3.2.2), приходим к (3.2.4).
Наконец, для $ f \in (S_{p,\theta}^\alpha B)^\prime(D), $ применяя (3.2.1),
(3.2.4), получаем, что $ h_{\delta,x^0}^{-1} f \in
(S_{p,\theta}^\alpha B)^\prime(x^0 +\delta D) $ и
соблюдается (3.2.5). $ \square $
\bigskip

3.3. В этом пункте устанавливаются некоторые метрические соотношения,
используемые для получения интересующих нас оценок.

В [20] доказана следующая теорема (см. там теорему 4).

Теорема 3.3.1

Пусть $ d \in \N, 1 < p < \infty $ и наборы функций
$$
\bm \phi = \{\phi_1, \ldots, \phi_d\},
\tilde{\bm \phi} = \{\tilde \phi_1, \ldots, \tilde \phi_d \},
$$
таковы, что соблюдаются условия предложения 10 из [20].
Тогда существуют константы $ c_{1}(d, \bm \phi, \tilde{\bm \phi}, p) >0, \\
c_{2}(d,\bm \phi, \tilde{\bm \phi}, p) >0 $ такие, что для любой
функции $ f \in L_p(\R^d) $ выполняются неравенства
\begin{equation*} \tag{3.3.1}
c_{1} \| f \|_{L_p(\R^d)} \le
\biggl( \int_{\R^d} \biggl(\sum_{\kappa \in \Z_+^d}
| (\mathcal E_\kappa^{\bm \phi,\tilde{\bm \phi},p} f)(x)|^2 \biggr)^{p/2}
dx \biggr)^{1/p} \le c_{2} \| f \|_{L_p(\R^d)},
\end{equation*}
где операторы $ \mathcal E_\kappa^{\bm \phi,\tilde{\bm \phi},p}:
L_p(\R^d) \mapsto L_p(\R^d), $ определяются следующим образом.

При $ d \in \N, 1 < p < \infty, \kappa \in \Z_+^d $ в условиях теоремы 3.3.1
рассмотрим линейные непрерывные операторы
$$
E_{\kappa_j}^{j, p} = E_{\kappa_j}^{\phi_j, \tilde \phi_j, p}, j =1, \ldots, d,
$$
значения которых для $ f \in L_p(\R) $ задаются равенствами
\begin{multline*}
(E_{\kappa_j}^{\phi_j, \tilde \phi_j, p} f)(x) = \sum_{\nu \in \Z} 2^{\kappa_j}
\biggl(\int_{\R} f(y) \overline {\tilde \phi_j(2^{\kappa_j} y -\nu)} dy \biggr)
\phi_j(2^{\kappa_j} x -\nu) \\
\text{ почти для всех } x \in \R, \kappa_j \in \Z_+, 
(E_{-1}^{\phi_j, \tilde \phi_j, p} f)(x) = 0 \text{ почти для всех } x \in \R, \\
(\text{см. предложение 2 из [20] и (2.1.8) из [20]}).
\end{multline*}
Определим операторы
\begin{equation*}
E_\kappa^p = E_\kappa^{\bm \phi, \tilde{\bm \phi}, p} =
\prod_{j =1}^d V_j^{L_p} (E_{\kappa_j}^{\phi_j,\tilde \phi_j, p}) =
\prod_{j =1}^d V_j (E_{\kappa_j}^{j, p}), \kappa \in \Z_+^d.
\end{equation*}

Далее, определим линейные операторы
\begin{equation*}
\begin{gathered}
\mathcal E_\kappa^{p}  =  \mathcal E_\kappa^{\bm \phi,\tilde{\bm \phi},p}:
L_p(\R^d) \mapsto L_p(\R^d), \\
\mathcal E_{\kappa_j}^{j, p}: L_p(\R) \mapsto L_p(\R), j =1,\ldots,d, \kappa \in \Z_+^d,
\end{gathered}
\end{equation*}
следующими соотношениями
\begin{equation*}
\begin{gathered}
\mathcal E_\kappa^{p} = \mathcal E_\kappa^{\bm \phi,\tilde{\bm \phi},p} =
\sum_{\epsilon \in \Upsilon^d: \s(\epsilon) \subset \s(\kappa)}
(-\e)^\epsilon E_{\kappa -\epsilon}^{\bm \phi,\tilde{\bm \phi},p} =
\sum_{\epsilon \in \Upsilon^d: \s(\epsilon) \subset \s(\kappa)}
(-\e)^\epsilon E_{\kappa -\epsilon}^{p} \\
\mathcal E_{\kappa_j}^{j, p} = E_{\kappa_j}^{j, p} -
E_{\kappa_j -1}^{j, p}, j =1,\ldots,d, \kappa \in \Z_+^d.
\end{gathered}
\end{equation*}

При этом в условиях теоремы 3.3.1 имеет место равенство (см. [20]):
\begin{equation*} \tag{3.3.2}
\mathcal E_\kappa^{p} = \prod_{j =1}^d V_j(\mathcal E_{\kappa_j}^{j, p}).
\end{equation*}

Из (3.3.1) следует, что в условиях теоремы 3.3.1 соблюдается неравенство
\begin{equation*} \tag{3.3.3}
\| \mathcal E_\kappa^{p} \|_{\mathcal B(L_p(\R^d),L_p(\R^d)} \le c_2,
\kappa \in \Z_+^d.
\end{equation*}

Следствие

Пусть выполнены условия теоремы 3.3.1. Тогда существуют константы
$ c_{3}(d,\bm \phi, \tilde{\bm \phi},p) >0,
c_{4}(d,\bm \phi, \tilde{\bm \phi},p) >0 $ такие, что
для $ f \in L_p(\R^d) $ при $ \p = \min(2,p), \P = \max(2,p) $ имеют место
неравенства
\begin{equation*} \tag{3.3.4}
c_{3} (\sum_{\kappa \in \Z_+^d }
\| \mathcal E_\kappa^{\bm \phi,\tilde{\bm \phi},p} f \|_{L_p(\R^d)}^{\P} )^{1 /\P} \le \\
\| f \|_{L_p(\R^d)} \le c_{4} (\sum_{\kappa \in \Z_+^d }
\| \mathcal E_\kappa^{\bm \phi,\tilde{\bm \phi},p} f \|_{L_p(\R^d)}^{\p} )^{1 /\p}.
\end{equation*}

Доказательство.

Сначала установим второе неравенство в (3.3.4).
В условиях теоремы 3.3.1 с помощью неравенства треугольника для нормы в
пространстве $ L_{p/2}(\R^d) $ при $ p > 2, $ и неравенства (1.1.2) с $ a = p/2 $
при $ 1 < p \le 2, $ из (2.4.10) в [20] выводим
\begin{multline*} \tag{3.3.5}
\| E_k f \|_{L_p(\R^d)} \le c_{4}  \biggl( \| \sum_{\kappa \in \Z_+^d(k)}
| \mathcal E_\kappa f |^2 \|_{L_{p/2}(\R^d)} \biggr)^{1/2} \le 
c_{4}  \biggl( \sum_{\kappa \in \Z_+^d(k)} \| | \mathcal E_\kappa f |^2
\|_{L_{p/2}(\R^d)} \biggr)^{1/2} = \\
c_{4}  \biggl(\sum_{\kappa \in
\Z_+^d(k)} \| \mathcal E_\kappa f \|_{L_p(\R^d)}^2 \biggr)^{1/2}, 
f \in L_p(\R^d), k \in \Z_+^d, p > 2,
\end{multline*}
и
\begin{multline*} \tag{3.3.6}
\| E_k f \|_{L_p(\R^d)} \le c_{4} \biggl(\int_{\R^d} \biggl(\sum_{\kappa
\in \Z_+^d(k)} | (\mathcal E_\kappa f)(x)|^2 \biggr)^{p/2} dx \biggr)^{1/p} \le \\
c_{4} \biggl(\int_{\R^d} \sum_{\kappa \in \Z_+^d(k)} | (\mathcal E_\kappa f)(x)|^p
dx\biggr)^{1/p} = c_{4} \biggl( \sum_{\kappa \in \Z_+^d(k)}
\| \mathcal E_\kappa f \|_{L_p(\R^d)}^p \biggr)^{1/p}, \\
f \in L_p(\R^d), k \in \Z_+^d, 1 < p \le 2.
\end{multline*}

Отметим ещё, что для любого семейства неотрицательных чисел
$ \{a_\kappa \in \R: a_\kappa \ge 0, \kappa \in \Z_+^d \} $
существует предел (конечный или бесконечный)
$$
\lim_{ \mn(k) \to \infty} \sum_{ \kappa \in \Z_+^d(k)} a_\kappa =
\sup_{k \in \Z_+^d} \sum_{ \kappa \in \Z_+^d(k)} a_\kappa,
$$
поскольку для любых $ k, k^\prime \in \Z_+^d: \mn(k) \ge \mx(k^\prime), $ соблюдается неравенство
$$
\sum_{ \kappa \in \Z_+^d(k)} a_\kappa \ge
\sum_{ \kappa \in \Z_+^d(k^\prime)} a_\kappa.
$$

Пусть теперь $ f \in L_p(\R^d) $ и ряд $ \sum_{\kappa \in \Z_+^d }
\| \mathcal E_\kappa f \|_{L_p(\R^d)}^{\p} $ сходится. Тогда,
с учётом замечаний, на основании (3.3.5), (3.3.6) заключаем, что
при $ k \in \Z_+^d $ справедливо неравенство
\begin{equation*} \tag{3.3.7}
\| E_k f \|_{L_p(\R^d)} \le c_{4} (\sum_{\kappa \in \Z_+^d}
\| \mathcal E_\kappa f \|_{L_p(\R^d)}^{\p} )^{1 /\p} < \infty.
\end{equation*}
Причём в силу (2.3.20) из [20] в $ L_p(\R^d) $ соблюдается равенство
\begin{equation*}
f = \lim_{ \mn(k) \to \infty} E_k f,
\end{equation*}
а, следовательно,
\begin{equation*}
\| f \|_{L_p(\R^d)} = \lim_{ \mn(k) \to \infty} \| E_k f \|_{L_p(\R^d)}.
\end{equation*}
Соединяя последнее равенство с (3.3.7), приходим ко второму неравенству в  (3.3.4).

Проверим соблюдение первого неравенства в (3.3.4). При $ 1 < p \le 2, $ принимая
во внимание, что в силу второго неравенства в (3.3.4), а также
вследствие (2.3.6) из [20] и (3.3.3) верно неравенство
\begin{multline*}
\| \sum_{\kappa \in \Z_+^d(k) }
\mathcal E_\kappa g_\kappa \|_{L_{p^\prime}(\R^d)} \le
c_{5} ( \sum_{\kappa^\prime \in \Z_+^d }
\| \mathcal E_{\kappa^\prime}( \sum_{\kappa \in \Z_+^d(k) }
\mathcal E_\kappa g_\kappa)\|_{L_{p^\prime}(\R^d)}^2)^{1 /2} = \\
c_{5} ( \sum_{\kappa^\prime \in \Z_+^d }
\| \sum_{\kappa \in \Z_+^d(k) } \mathcal E_{\kappa^\prime}
\mathcal E_\kappa g_\kappa\|_{L_{p^\prime}(\R^d)}^2)^{1 /2} = \\
c_{5} ( \sum_{\kappa \in \Z_+^d(k) }
\| \mathcal E_\kappa g_\kappa\|_{L_{p^\prime}(\R^d)}^2)^{1 /2} \le
c_{5} ( \sum_{\kappa \in \Z_+^d(k) }
(c_{6} \| g_\kappa \|_{L_{p^\prime}(\R^d)})^2)^{1 /2} = \\
c_{7} ( \sum_{\kappa \in \Z_+^d(k) }
\| g_\kappa \|_{L_{p^\prime}(\R^d)}^2)^{1 /2},
g_\kappa \in L_{p^\prime}(\R^d), \kappa \in \Z_+^d(k), k \in \Z_+^d,
\end{multline*}
а также учитывая неравенство Гёльдера, на основании (3.1.7), (2.3.5) из [20]
заключаем, что для $ f \in L_p(\R^d), k \in \Z_+^d $ имеет место соотношение
\begin{multline*}
(\sum_{\kappa \in \Z_+^d(k)} \| \mathcal E_\kappa f \|_{L_p(\R^d)}^2 )^{1/2} = \\
\max_{\{g_\kappa \in L_{p^\prime}(\R^d), \kappa \in \Z_+^d(k)\}:
(\sum_{\kappa \in \Z_+^d(k)} \| g_\kappa \|_{L_{p^\prime}(\R^d)}^2 )^{1/2} \le 1}
| \int_{\R^d} \sum_{\kappa \in \Z_+^d(k) }
( \mathcal E_\kappa f ) \cdot \overline{g_\kappa} dx | = \\
\max_{\{g_\kappa \in L_{p^\prime}(\R^d), \kappa \in \Z_+^d(k)\}:
(\sum_{\kappa \in \Z_+^d(k)} \| g_\kappa \|_{L_{p^\prime}(\R^d)}^2 )^{1/2} \le 1}
| \sum_{\kappa \in \Z_+^d(k) }
\int_{\R^d} ( \mathcal E_\kappa f ) \cdot \overline{g_\kappa} dx | = \\
\max_{\{g_\kappa \in L_{p^\prime}(\R^d), \kappa \in \Z_+^d(k)\}:
(\sum_{\kappa \in \Z_+^d(k)} \| g_\kappa \|_{L_{p^\prime}(\R^d)}^2)^{1/2} \le 1}
| \sum_{\kappa \in \Z_+^d(k) }
\int_{\R^d} f \cdot \overline{(\mathcal E_\kappa g_\kappa )} dx | = \\
\max_{\{g_\kappa \in L_{p^\prime}(\R^d), \kappa \in \Z_+^d(k)\}:
(\sum_{\kappa \in \Z_+^d(k)} \| g_\kappa \|_{L_{p^\prime}(\R^d)}^2)^{1/2} \le 1}
| \int_{\R^d} f \cdot
\overline{( \sum_{\kappa \in \Z_+^d(k) } \mathcal E_\kappa g_\kappa )} dx | \le \\
\max_{g \in L_{p^\prime}(\R^d): \| g \|_{L_{p^\prime}(\R^d)} \le c_{7}}
| \int_{\R^d} f \cdot \overline g dx | =
c_{7} \max_{g \in L_{p^\prime}(\R^d): \| g \|_{L_{p^\prime}(\R^d)} \le 1}
| \int_{\R^d} f \cdot \overline g dx | = c_{7} \| f \|_{L_{p}(\R^d)}.
\end{multline*}
Отсюда, переходя к пределу при $ \mn(k) \to \infty, $ получаем первое
неравенство в (3.3.4) при $ 1 < p \le 2. $
При $ 2 < p < \infty $ в силу тех же соображений, что и выше, пользуясь тем,
что
\begin{multline*}
\| \sum_{\kappa \in \Z_+^d(k) }
\mathcal E_\kappa g_\kappa \|_{L_{p^\prime}(\R^d)} \le
c_{5} ( \sum_{\kappa^\prime \in \Z_+^d }
\| \mathcal E_{\kappa^\prime}( \sum_{\kappa \in \Z_+^d(k) }
\mathcal E_\kappa g_\kappa)\|_{L_{p^\prime}(\R^d)}^{p^\prime})^{1 /p^\prime} = \\
c_{5} ( \sum_{\kappa^\prime \in \Z_+^d }
\| \sum_{\kappa \in \Z_+^d(k) } \mathcal E_{\kappa^\prime}
\mathcal E_\kappa g_\kappa\|_{L_{p^\prime}(\R^d)}^{p^\prime})^{1 /p^\prime} = \\
c_{5} ( \sum_{\kappa \in \Z_+^d(k) }
\| \mathcal E_\kappa g_\kappa\|_{L_{p^\prime}(\R^d)}^{p^\prime})^{1 /p^\prime} \le
c_{5} ( \sum_{\kappa \in \Z_+^d(k) }
(c_{6} \| g_\kappa \|_{L_{p^\prime}(\R^d)})^{p^\prime})^{1 /p^\prime} = \\
c_{7} ( \sum_{\kappa \in \Z_+^d(k) }
\| g_\kappa \|_{L_{p^\prime}(\R^d)}^{p^\prime})^{1 /p^\prime},
g_\kappa \in L_{p^\prime}(\R^d), \kappa \in \Z_+^d(k), k \in \Z_+^d,
\end{multline*}
получаем, что для $ f \in L_p(\R^d), k \in \Z_+^d $ выполняется неравенство
\begin{multline*}
\biggl(\sum_{\kappa \in \Z_+^d(k)} \| \mathcal E_\kappa f \|_{L_p(\R^d)}^p \biggr)^{1/p} = \\
\max_{\{g_\kappa \in L_{p^\prime}(\R^d), \kappa \in \Z_+^d(k)\}:
(\sum_{\kappa \in \Z_+^d(k)} \| g_\kappa \|_{L_{p^\prime}(\R^d)}^{p^\prime} )^{1/p^\prime} \le 1}
\biggl| \int_{\R^d} \sum_{\kappa \in \Z_+^d(k) }
( \mathcal E_\kappa f ) \cdot \overline{g_\kappa} dx \biggr| = \\
\max_{\{g_\kappa \in L_{p^\prime}(\R^d), \kappa \in \Z_+^d(k)\}:
(\sum_{\kappa \in \Z_+^d(k)} \| g_\kappa \|_{L_{p^\prime}(\R^d)}^{p^\prime} )^{1/p^\prime} \le 1}
\biggl| \sum_{\kappa \in \Z_+^d(k) }
\int_{\R^d} ( \mathcal E_\kappa f ) \cdot \overline{g_\kappa} dx \biggr| = \\
\max_{\{g_\kappa \in L_{p^\prime}(\R^d), \kappa \in \Z_+^d(k)\}:
(\sum_{\kappa \in \Z_+^d(k)} \| g_\kappa \|_{L_{p^\prime}(\R^d)}^{p^\prime})^{1/p^\prime} \le 1}
\biggl| \sum_{\kappa \in \Z_+^d(k) }
\int_{\R^d} f \cdot \overline{(\mathcal E_\kappa g_\kappa )} dx \biggr| = \\
\max_{\{g_\kappa \in L_{p^\prime}(\R^d), \kappa \in \Z_+^d(k)\}:
(\sum_{\kappa \in \Z_+^d(k)} \| g_\kappa \|_{L_{p^\prime}(\R^d)}^{p^\prime})^{1/p^\prime} \le 1}
\biggl| \int_{\R^d} f \cdot
\overline{( \sum_{\kappa \in \Z_+^d(k) } \mathcal E_\kappa g_\kappa )} dx \biggr| \le \\
\max_{g \in L_{p^\prime}(\R^d): \| g \|_{L_{p^\prime}(\R^d)} \le c_{7}}
\biggl| \int_{\R^d} f \cdot \overline g dx \biggr| =
c_{7} \max_{g \in L_{p^\prime}(\R^d):
\| g \|_{L_{p^\prime}(\R^d)} \le 1} \biggl| \int_{\R^d} f \cdot \overline g dx \biggr| =
c_{7} \| f \|_{L_{p}(\R^d)}.
\end{multline*}
Переходя в этом неравенстве к пределу при $ \mn(k) \to \infty, $ приходим к
первому неравенству в (3.3.4) при $ 2 < p < \infty. \square $

Приведём ещё некоторые полезные для нас сведения.

Для $ d \in \N, m \in \Z_+^d $ и области $ D \subset \R^d $ обозначим через
$ C^m(D) $ пространство всех функций $ f: D \mapsto \R, $ у которых для каждого
$ \lambda \in \Z_+^d(m) $ существует непрерывная в области $ D $ частная
производная $ \D^\lambda f $ порядка $ \lambda, $ а через $ C_0^m(d) $ обозначим
пространство всех функций $ f \in C^m(\R^d), $ у которых носитель $ \supp f
\subset D. $

Пусть $ d \in \N $ и функция $ \psi \in C_0^0(I^d), $ т.е. $ f \in C(\R^d) $
такоова, что её носитель $ \supp \psi \subset I^d. $ При $ \kappa \in \Z_+^d $
через $ \L_\kappa^{d,\psi} $ обозначим линейную оболочку системы функций
$ \{ \psi(2^\kappa \cdot -\nu), \nu \in \Nu_{0, 2^\kappa -\e}^d \} $
и определим линейное изоморфное отображение $ V_\kappa^{d,\psi}:
\R^{2^{(\kappa, \e)}} \mapsto \L_\kappa^{d,\psi}, $ полагая для
$ Y_\kappa = \{y_{\kappa,\nu^\kappa}, \nu^\kappa \in \Nu_{0, 2^\kappa -\e}^d \}
\in \R^{2^{(\kappa, \e)}} $ значение
$$
(V_\kappa^{d,\psi} Y_\kappa)(x) = \sum_{ \nu^\kappa \in \Nu_{0, 2^\kappa -\e}^d}
y_{\kappa, \nu_\kappa} \psi(2^\kappa x -\nu^\kappa), x \in \R^d.
$$
Понятно, что если при $ m \in \Z_+^d $ функция $ \psi \in C_0^m(I^d), $ то
\begin{equation*} \tag{3.3.8}
\L_\kappa^{d,\psi} \subset C_0^m(I^d) \subset L_q(\R^d), 1 \le q \le \infty,
\kappa \in \Z_+^d.
\end{equation*}

Лемма 3.3.2

Пусть $ d \in \N, m \in \Z_+^d $ и функция $ \psi \in C_0^m(I^d). $
Пусть ещё $ \lambda \in \Z_+^d(m), 1 \le q \le \infty. $
Тогда существует константа $ c_8(d,\psi,\lambda,q) >0 $ такая, что при $ \kappa
\in \Z_+^d $ для $ Y_\kappa \in \R^{2^{(\kappa,\e)}} $ выполняются равенства
\begin{multline*} \tag{3.3.9}
\| \D^\lambda (V_\kappa^{d,\psi} Y_\kappa) \|_{L_q(\R^d)} =
\| \D^\lambda (V_\kappa^{d,\psi} Y_\kappa) \|_{L_q(I^d)} = \\
\| \D^\lambda \psi \|_{L_q(I^d)} 2^{(\kappa, \lambda -q^{-1} \e)}
\| Y_\kappa \|_{l_q^{2^{(\kappa, \e)}}} = c_8 2^{(\kappa, \lambda)} \| V_\kappa^{d,\psi} Y_\kappa \|_{L_q(\R^d)}.
\end{multline*}

Лемма 3.3.3

Пусть $ d \in \N, \alpha \in \R_+^d, l = l(\alpha), m \in \N^d: l \in \Z_+^d(m), $
функция $ \psi \in C_0^m(I^d), 1 \le p \le \infty, 1 \le \theta \le \infty. $
Тогда существует константа $ c_9(d,\alpha,p,\theta,m,\psi) >0 $ такая,
что для любого семейства функций $ \{ g_\kappa \in \L_\kappa^{d,\psi},
\kappa \in \Z_+^d \}, $ удовлетворяющего условию
\begin{equation*} \tag{3.3.10}
( \sum_{ \kappa \in \Z_+^d} (2^{(\kappa, \alpha)}
\| g_\kappa \|_{L_p(\R^d)})^\theta)^{1 / \theta} \le c_9,
\end{equation*}
ряд $ \sum_{\kappa \in \Z_+^d} g_\kappa $ сходится в $ L_p(\R^d) $
и справедливо включение
\begin{equation*} \tag{3.3.11}
( \sum_{ \kappa \in \Z_+^d} g_\kappa) \in B((S_{p, \theta}^\alpha B)^0(\R^d)).
\end{equation*}

Доказательство лемм 3.3.2 и 3.3.3 по существу проведено в [8] (см. лемма
2.1.3 и лемма 2.1.4).

Лемма 3.3.4

Пусть $ d \in \N, m \in \N^d, \lambda \in \Z_+^d. $ Тогда можно построить
функцию $ g^{d,m,\lambda} \in C_0^{m +\lambda}(I^d), $ для которой
при $ 1 < q < \infty $ существует константа $ c_{10}(g^{d,m,\lambda},q) >0 $
такая, что при любом $ r \in \N $ для любого набора векторов $ \{Y_\kappa \in
\R^{2^{(\kappa, \e)}} \mid \kappa \in \Z_+^d: (\kappa,\e) = r \} $ соблюдается
неравенство
\begin{multline*} \tag{3.3.12}
\| \D^\lambda (\sum_{\kappa \in \Z_+^d: (\kappa, \e) = r}
V_\kappa^{d,g^{d,m,\lambda}} Y_\kappa) \|_{L_q(I^d)} \ge \\
c_{10} \biggl(\sum_{ \kappa \in \Z_+^d: (\kappa,\e) = r}
\| 2^{(\kappa,\lambda)} V_\kappa^{d,\D^\lambda g^{d,m,\lambda}} Y_\kappa
\|_{L_q(I^d)}^{\Q}\biggr)^{1 / \Q}
\end{multline*}
или
\begin{multline*}
\| \D^\lambda (\sum_{\kappa \in \Z_+^d: (\kappa, \e) = r} \sum_{\nu^\kappa \in
\Nu_{0,2^\kappa -\e}^d} y_{\kappa, \nu^\kappa} g^{d,m,\lambda}(2^\kappa \cdot -
\nu^\kappa)) \|_{L_q(I^d)} \ge\\
 c_{10} \biggl(\sum_{\kappa \in \Z_+^d: (\kappa,\e) = r}
\| 2^{(\kappa,\lambda)} \sum_{\nu^\kappa \in \Nu_{0,2^\kappa -\e}^d} y_{\kappa,\nu^\kappa}
\D^\lambda g^{d,m,\lambda}(2^\kappa \cdot -\nu^\kappa) \|_{L_q(I^d)}^{\Q}\biggr)^{1 /\Q},
\end{multline*}
где $ \Q = \max(2,q). $

Доказательство.

Для получения (3.3.12) сначала проведём следующие построения. При $ m \in \N, $
опираясь на [21, гл. 4], построим вещественную финитную масштабирующую
функцию $ \phi^m \in C^m(\R) $ для кратномасштабного анализа
$ \{X_k \subset L_2(\R), k \in \Z\} $ такую, что система функций
$ \{ \phi^m(\cdot -\nu), \nu \in \Z\} $ -- ортонормированный базис в $ X_0, $ а
$ \psi^m \in C^m(\R) $ -- соответствующая вещественная финитная
всплеск-функция, для которой система функций $ \{\psi^m(\cdot -\nu),
\nu \in \Z\} $ -- ортонормированная. При этом
\begin{equation*} \tag{3.3.13}
\psi^m(\cdot -\nu) \in X_1, \int_\R f(x) \psi^m(x -\nu) dx =0,
f \in X_0, \nu \in \Z.
\end{equation*}
Кроме того, сдвигая, если нужно, на целочисленный шаг, можно считать, что
носитель $ \supp \psi^m \subset n^{m,0} I, $ где $ n^{m,0} \in \N. $

Далее, полагая $ \psi^{m,0} = \psi^m, $ при $ l \in \Z_+ $ определим число
$ n^{m,l +1} = 2 n^{m,l} \in \N $ и функцию
$$
\psi^{m,l +1}(x) = \int_{-\infty}^x \psi^{m,l}(t) dt -
\int_{-\infty}^x \psi^{m,l}(t -n^{m,l}) dt, x \in \R.
$$
Принимая во внимание эти обозначения, по индукции легко проверить, что
\begin{equation*} \tag{3.3.14}
\supp \psi^{m,l} \subset n^{m,l} I, l \in \Z_+,
\end{equation*}
а также, учитывая соотношение
\begin{multline*}
\frac{d^{l +1} \psi^{m,l +1}}{dx^{l +1}}(x) =
\frac{d^{l +1}}{dx^{l +1}} (\psi^{m,l +1}(x)) = \\
\frac{d^l}{dx^l} (\frac{d}{dx}
(\int_{-\infty}^x \psi^{m,l}(t) dt)) -\frac{d^l}{dx^l} (\frac{d}{dx}
(\int_{-\infty}^x \psi^{m,l}(t -n^{m,l}) dt)) =\\
 \frac{d^l}{dx^l} (\psi^{m,l}(x))
-\frac{d^l}{dx^l} (\psi^{m,l}(x -n^{m,l})) =
\frac{d^l \psi^{m,l}}{dx^l}(x) -\frac{d^l \psi^{m,l}}{dx^l}(x -n^{m,l}),
\end{multline*}
несложно убедиться в том, что при $ l \in \Z_+, n \in \Z $ имеет место включение
\begin{equation*} \tag{3.3.15}
\frac{d^l \psi^{m,l}}{dx^l}(\cdot -n) \in \span \{\psi^m(\cdot -\nu), \nu \in \Z\}.
\end{equation*}

При этом при $ k \in \Z_+, \kappa \in \Z_+, $ принимая во внимание (3.3.13), а
также (1.2.14), (1.2.16), (1.2.17), пп. 4), 3) предложения 4 из [20],
для $ \psi \in \span \{ \psi^m(2^k \cdot -\nu), \nu \in \Z\} \subset
\Im \mathcal U_{k +1}, $ где $ \mathcal U_k = \mathcal E_k^{\phi^m,\phi^m,2},
U_k = E_k^{\phi^m,\phi^m,2}, k \in \Z_+, $ (см. (1.2.26) из [20]) имеем
равенства
\begin{multline*} \tag{3.3.16}
\mathcal E_\kappa^p \psi = \mathcal E_\kappa^{\phi^m,\phi^m,p} \psi =
E_\kappa^p \psi -E_{\kappa -1}^p \psi = U_\kappa \psi -U_{\kappa -1} \psi =\\
\mathcal U_\kappa \psi = \mathcal U_\kappa \mathcal U_{k +1} \psi =
\begin{cases}
\psi, \text{ при } \kappa = k +1; \\
0, \text{ при } \kappa \ne k +1.
\end{cases}
\end{multline*}
Выбирая $ k^{m,l} \in \Z_+ $ так, что $ 2^{-k^{m,l}} n^{m,l} \le 1, l \in \Z_+, $
и определяя функцию $ g^{m,l} $ равенством $ g^{m,l}(\cdot) =
\psi^{m,l}(2^{k^{m,l}} \cdot) \in C^{m +l}(\R), l \in \Z_+, $ ввиду (3.3.14)
получаем, что
\begin{equation*} \tag{3.3.17}
\supp g^{m,l} = 2^{-k^{m,l}} \supp \psi^{m,l} \subset 2^{-k^{m,l}} n^{m,l} I
\subset I, l \in \Z_+.
\end{equation*}
Отметим ещё, что при $ l \in \Z_+, k \in \Z_+, n \in \Z $ вследствие (3.3.15)
справедливо соотношение
\begin{multline*} \tag{3.3.18}
\frac{d^l g^{m,l}}{dx^l}(2^k \cdot -n) =
\frac{d^l}{dx^l}(\psi^{m,l}(2^{k^{m,l}} x)) \mid_{x = 2^k \cdot -n} = \\
2^{l k^{m,l}} \frac{d^l \psi^{m,l}}{dx^l}(2^{k^{m,l}} (2^k \cdot -n)) =
2^{l k^{m,l}} \frac{d^l \psi^{m,l}}{dx^l}(2^{k^{m,l} +k} \cdot -2^{k^{m,l}} n) \\
\in \span \{\psi^m(2^{k^{m,l} +k} \cdot -\nu), \nu \in \Z\}.
\end{multline*}

Теперь, имея в виду проведенные построения, в условиях леммы задавая
наборы функций
$$
\bm \phi = \{\phi_1, \ldots, \phi_d\},
\tilde{\bm \phi} = \{ \tilde \phi_1, \ldots, \tilde \phi_d \},
$$
так, что при $ j =1, \ldots, d $ соблюдаются равенства
$ \phi_j = \tilde \phi_j = \phi^{m_j} $ (см. выше), определим функцию
$ g^{d,m,\lambda} \in C^{m +\lambda}(\R^d) $ равенством
$$
g^{d,m,\lambda}(x) = \prod_{j =1}^d g^{m_j,\lambda_j}(x_j), x \in \R^d.
$$
Тогда из (3.3.17) следует, что
\begin{equation*} \tag{3.3.19}
\supp g^{d,m,\lambda} \subset I^d, \\
\supp \D^\lambda g^{d,m,\lambda} \subset I^d, \\
\D^\lambda g^{d,m,\lambda} \in C^m(\R^d).
\end{equation*}

Поскольку построенные наборы функций $ \bm \phi = \tilde{\bm \phi} $
удовлетворяют условиям следствия из теоремы 3.3.1, то применяя первое
неравенство в (3.3.4), с учётом (3.3.19), (3.3.8) при $ r \in \N $
получаем
\begin{multline*} \tag{3.3.20}
\| \D^\lambda (\sum_{\kappa \in \Z_+^d: (\kappa, \e) = r}
V_\kappa^{d,g^{d,m,\lambda}} Y_\kappa) \|_{L_q(I^d)} = \\
\| \D^\lambda (\sum_{\kappa \in \Z_+^d: (\kappa, \e) = r} \sum_{\nu^\kappa \in
\Nu_{0,2^\kappa -\e}^d} y_{\kappa, \nu^\kappa} g^{d,m,\lambda}(2^\kappa \cdot -
\nu^\kappa)) \|_{L_q(I^d)} = \\
\| \sum_{\kappa \in \Z_+^d: (\kappa, \e) = r} \sum_{\nu^\kappa \in
\Nu_{0,2^\kappa -\e}^d} y_{\kappa, \nu^\kappa}
\D^\lambda (g^{d,m,\lambda}(2^\kappa \cdot -\nu^\kappa)) \|_{L_q(I^d)} = \\
\| \sum_{\kappa \in \Z_+^d: (\kappa, \e) = r} \sum_{\nu^\kappa \in
\Nu_{0,2^\kappa -\e}^d} y_{\kappa, \nu^\kappa} 2^{(\kappa,\lambda)}
(\D^\lambda g^{d,m,\lambda})(2^\kappa \cdot -\nu^\kappa) \|_{L_q(I^d)} = \\
\| \sum_{\kappa \in \Z_+^d: (\kappa, \e) = r} 2^{(\kappa,\lambda)}
\sum_{\nu^\kappa \in \Nu_{0,2^\kappa -\e}^d} y_{\kappa, \nu^\kappa}
(\D^\lambda g^{d,m,\lambda})(2^\kappa \cdot -\nu^\kappa) \|_{L_q(\R^d)} \ge\\
c_{10} \biggl(\sum_{\kappa \in \Z_+^d }
\| \mathcal E_\kappa(\sum_{k \in \Z_+^d: (k, \e) = r} 2^{(k,\lambda)}
\sum_{\nu^k \in \Nu_{0,2^k -\e}^d} y_{k, \nu^k}
(\D^\lambda g^{d,m,\lambda})(2^k \cdot -\nu^k))
\|_{L_q(\R^d)}^{\Q} \biggr)^{1 /\Q} = \\
c_{10} \biggl(\sum_{\kappa \in \Z_+^d }
\| \sum_{k \in \Z_+^d: (k, \e) = r} 2^{(k,\lambda)}
\sum_{\nu^k \in \Nu_{0,2^k -\e}^d} y_{k, \nu^k}
\mathcal E_\kappa((\D^\lambda g^{d,m,\lambda})(2^k \cdot -\nu^k))
\|_{L_q(\R^d)}^{\Q} \biggr)^{1 /\Q}.
\end{multline*}

Теперь заметим, что при $ \kappa \in \Z_+^d, k \in \Z_+^d: (k, \e) = r,
\nu^k \in \Nu_{0,2^k -\e}^d, $ благодаря (3.3.2), (1.3.1), (3.3.18), (3.3.16)
имеют место равенства
\begin{multline*} \tag{3.3.21}
\mathcal E_\kappa^{q}((\D^\lambda g^{d,m,\lambda})(2^k \cdot -\nu^k)) =
(\prod_{j =1}^d V_j(\mathcal E_{\kappa_j}^{j,q}))
(\prod_{j =1}^d (\D^{\lambda_j} g^{m_j,\lambda_j})(2^{k_j} \cdot -(\nu^k)_j)) = \\
\prod_{j =1}^d \mathcal E_{\kappa_j}^{j,q}
((\D^{\lambda_j} g^{m_j,\lambda_j})(2^{k_j} \cdot -(\nu^k)_j)) = \\
\prod_{j =1}^d \mathcal E_{\kappa_j}^{\phi^{m_j},\phi^{m_j},q}
((\D^{\lambda_j} g^{m_j,\lambda_j})(2^{k_j} \cdot -(\nu^k)_j)) = \\
\begin{cases}
\prod_{j =1}^d ((\D^{\lambda_j} g^{m_j,\lambda_j})(2^{k_j} \cdot -(\nu^k)_j)),
\text{ если } \kappa_j = k^{m_j,\lambda_j} +k_j +1, \text{ при } j =1,\ldots,d; \\
0, \text{ если } \exists j =1,\ldots,d: \kappa_j \ne k^{m_j,\lambda_j} +k_j +1, \\
\end{cases}
=\\
\begin{cases}
(\D^\lambda g^{d,m,\lambda})(2^k \cdot -\nu^k),
\text{ при } \kappa = k^{d,m,\lambda} +k +\e; \\
0, \text{ при } \kappa \ne k^{d,m,\lambda} +k +\e,
\end{cases}
\end{multline*}
где $ k^{d,m,\lambda} \in \Z_+^d: (k^{d,m,\lambda})_j = k^{m_j,\lambda_j},
j =1,\ldots,d. $

Используя для преобразования правой части (3.3.20) соотношение (3.3.21),
с учётом (3.3.19) выводим
\begin{multline*} \tag{3.3.22}
\biggl(\sum_{\kappa \in \Z_+^d }
\| \sum_{k \in \Z_+^d: (k, \e) = r} 2^{(k,\lambda)}
\sum_{\nu^k \in \Nu_{0,2^k -\e}^d} y_{k, \nu^k}
\mathcal E_\kappa((\D^\lambda g^{d,m,\lambda})(2^k \cdot -\nu^k))
\|_{L_q(\R^d)}^{\Q} \biggr)^{1 /\Q} = \\
\biggl(\sum_{\substack{\kappa \in \Z_+^d \mid \kappa = k^{d,m,\lambda} +k^\prime +\e,\\
k^\prime \in \Z_+^d: (k^\prime, \e) = r}}
\| \sum_{k \in \Z_+^d: (k, \e) = r} 2^{(k,\lambda)}\times\\
\sum_{\nu^k \in \Nu_{0,2^k -\e}^d} y_{k, \nu^k}
\mathcal E_\kappa((\D^\lambda g^{d,m,\lambda})(2^k \cdot -\nu^k))
\|_{L_q(\R^d)}^{\Q} \biggr)^{1 /\Q} = \\
\biggl(\sum_{k^\prime \in \Z_+^d: (k^\prime, \e) = r}
\| \sum_{k \in \Z_+^d: (k, \e) = r} 2^{(k,\lambda)}\times\\
\sum_{\nu^k \in \Nu_{0,2^k -\e}^d} y_{k, \nu^k}
\mathcal E_{k^{d,m,\lambda} +k^\prime +\e}((\D^\lambda g^{d,m,\lambda})(2^k \cdot -\nu^k))
\|_{L_q(\R^d)}^{\Q} \biggr)^{1 /\Q} = \\
\biggl(\sum_{k \in \Z_+^d: (k, \e) = r}
\| 2^{(k,\lambda)} \sum_{\nu^k \in \Nu_{0,2^k -\e}^d} y_{k, \nu^k}
(\D^\lambda g^{d,m,\lambda})(2^k \cdot -\nu^k)\|_{L_q(\R^d)}^{\Q} \biggr)^{1 /\Q} = \\
\biggl(\sum_{k \in \Z_+^d: (k, \e) = r}
\| 2^{(k,\lambda)} \sum_{\nu^k \in \Nu_{0,2^k -\e}^d} y_{k, \nu^k}
(\D^\lambda g^{d,m,\lambda})(2^k \cdot -\nu^k)\|_{L_q(I^d)}^{\Q} \biggr)^{1 /\Q} = \\
\biggl(\sum_{ \kappa \in \Z_+^d: (\kappa,\e) = r}
\| 2^{(\kappa,\lambda)} V_\kappa^{d,\D^\lambda g^{d,m,\lambda}} Y_\kappa
\|_{L_q(I^d)}^{\Q}\biggr)^{1 /\Q}.
\end{multline*}
Подставляя (3.3.22) в (3.3.20), приходим к (3.3.12). $ \square $
\bigskip

3.4. В этом пункте устанавливается оценка снизу величины
$ \sigma_n(\D^\lambda, B((S_{p,\theta}^\alpha B)^\prime(D)), L_q(D)). $

Лемма 3.4.1

Пусть $ d \in \N, \alpha \in \R_+^d, 1 \le p < \infty, 1 \le \theta \le \infty,
1 < q < \infty, \lambda \in \Z_+^d $ удовлетворяют условиям (2.1.6), (2.1.9)
и $ J = \{j \in \Nu_{1,d}^1: \alpha_j - \lambda_j = \mn(\alpha -\lambda)\}. $
Тогда существует константа $ c_{1}(d, \alpha, p,\theta,q,\lambda) > 0 $ такая,
что при $ n \in \N $ имеет место
неравенство
\begin{multline*} \tag{3.4.1}
\inf_{ \phi \in \Phi_n(C(I^d))} \sup_{f \in
B((S_{p,\theta}^\alpha B)^0(\R^d)) \cap C_0^0(I^d): \phi f =0}
\| \D^\lambda f \|_{L_q(I^d)} \ge \\c_{1}
\begin{cases}
n^{-\mu} (\log n)^{(\gamma -1)(\mu -(1 /p -1 /q) -(1 /\theta -1 /p)_+ -(1 /q -1 /\Q)_+)}, &
\text{ при $ q \le p $ или $ p < q \le 2;$ } \\
n^{-(\mu -(1/2 -1/q))} (\log n)^{(\gamma -1)(\mu -(1/p -1/q) -(1 /\theta -1 /p)_+)}, &
\text{ при $ p \le 2 < q; $ } \\
n^{-(\mu -(1/p -1/q))}
(\log n)^{(\gamma -1)(\mu -(1/p -1/q) -(1 /\theta -1 /p)_+)}, &
\text{ при $ 2 < p < q, $ }
\end{cases}
\end{multline*}
где $ \mu = \mn(\alpha -\lambda), \gamma = \cmn(\alpha -\lambda) = \card J. $

Доказательство.

Сначала в условиях леммы, фиксируя $ m \in \N^d $ так, чтобы $ l(\alpha) \in \Z_+^d(m), $
в соответствии с леммой 3.3.4 построим функцию $ g^{d,m, \lambda} \in C_0^{m +\lambda}(I^d). $
Теперь для $ n \in \N $ выберем $ r \in \N, $ для которого соблюдается
соотношение
\begin{equation*} \tag{3.4.2}
c(\gamma,r -1) 2^{r -1} \le 2n < c(\gamma,r) 2^r,
\end{equation*}
где $ c(\gamma,r) = \card \{ \kappa \in \Z_+^\gamma: (\kappa, \e) = r\}. $
Тогда применяя лемму 3.3.3 при $ \psi = g^{d,m,\lambda}, $ согласно (3.3.10),
(3.3.11) для $ n, r \in \N, $ связанных соотношением (3.4.2), обозначая
$ \Z_{+ r}^{d,J} = \{ \kappa \in \Z_+^d: (\kappa,\e) = r, \kappa = \kappa \chi_J\}, $
имеем
\begin{multline*} \tag{3.4.3}
\inf_{ \phi \in \Phi_n(C(I^d))} \sup_{f \in B((S_{p,\theta}^\alpha B)^0(\R^d)) \cap
C_0^0(I^d): \phi f =0}
\| \D^\lambda f \|_{L_q(I^d)} \ge \\
\inf_{ \phi \in \Phi_n(C(I^d))} \sup_{\substack{ \{Y_\kappa \in \R^{2^{(\kappa, \e)}},
\kappa \in \Z_{+ r}^{d,J}\}:\\
(\sum_{\kappa \in \Z_{+ r}^{d,J}} (2^{(\kappa,  \alpha)}
\| V_\kappa^{d,g^{d,m, \lambda}} Y_\kappa\|_{L_p(\R^d)})^\theta)^{1 /\theta}\le c_2, \\
\phi (\sum_{\kappa \in \Z_{+ r}^{d,J}}
V_\kappa^{d,g^{d,m,\lambda}} Y_\kappa) =0}}
\| \D^\lambda (\sum_{\kappa \in \Z_{+ r}^{d,J}}
V_\kappa^{d,g^{d,m,\lambda}} Y_\kappa) \|_{L_q(I^d)} = V_1.
\end{multline*}

Для оценки $ V_1, $ используя (3.3.12) и делая замену переменных $ Y_\kappa =
2^{-(\kappa, \alpha)} X_\kappa, \kappa \in \Z_{+ r}^{d,J}, $ а затем, учитывая,
что $ (\alpha -\lambda) \chi_J = \mu \ch_J, $
а, следовательно, при $ \kappa \in \Z_{+ r}^{d,J}, $
соблюдается соотношение
\begin{multline*}
(\kappa, \alpha -\lambda) = (\kappa \chi_J, \alpha -\lambda) =
(\kappa \chi_J, (\alpha -\lambda) \chi_J) =\\
(\kappa \chi_J, \mu \chi_J) = \mu (\kappa \chi_J, \chi_J) =
\mu (\kappa \chi_J, \e) = \mu (\kappa, \e) = \mu r,
\end{multline*}
выводим
\begin{multline*} \tag{3.4.4}
V_1 \ge \\
\inf_{ \phi \in \Phi_n(C(I^d))} \sup_{\substack{ \{Y_\kappa \in
\R^{2^{(\kappa,\e)}}, \kappa \in \Z_{+ r}^{d,J}\}:\\
(\sum_{\kappa \in \Z_{+ r}^{d,J}} (2^{(\kappa,  \alpha)}
\| V_\kappa^{d,g^{d,m,\lambda}} Y_\kappa\|_{L_p(\R^d)})^\theta)^{1 /\theta} \le c_2,\\
\phi (\sum_{\kappa \in \Z_{+ r}^{d,J}}
V_\kappa^{d,g^{d,m,\lambda}} Y_\kappa) =0}}
c_3 \biggl(\sum_{ \kappa \in \Z_{+ r}^{d,J}}
\| 2^{(\kappa, \lambda)} V_\kappa^{d,\D^\lambda g^{d,m,\lambda}}
Y_\kappa \|_{L_q(I^d)}^{\Q}\biggr)^{1 /\Q} = \\
\inf_{ \phi \in \Phi_n(C(I^d))} \sup_{\substack{ \{X_\kappa \in \R^{2^{(\kappa, \e)}},
\kappa \in \Z_{+ r}^{d,J}\}:\\
(\sum_{\kappa \in \Z_{+ r}^{d,J}}
(\| V_\kappa^{d,g^{d,m,\lambda}} X_\kappa\|_{L_p(\R^d)})^\theta)^{1 /\theta}\le c_2,\\
\phi (\sum_{\kappa \in \Z_{+ r}^{d,J}}
2^{-(\kappa, \alpha)} V_\kappa^{d,g^{d,m,\lambda}} X_\kappa) =0}}
c_3 \biggl(\sum_{ \kappa \in \Z_{+ r}^{d,J}}
\| 2^{-(\kappa, \alpha -\lambda)}
V_\kappa^{d,\D^\lambda g^{d,m,\lambda}}
X_\kappa \|_{L_q(I^d)}^{\Q}\biggr)^{1 /\Q} = \\
c_3 2^{-\mu r} \inf_{ \phi \in \Phi_n(C(I^d))} \sup_{\substack{ \{X_\kappa \in
\R^{2^{(\kappa, \e)}}, \kappa \in \Z_{+ r}^{d,J}\}:\\
(\sum_{\kappa \in \Z_{+ r}^{d,J}}
(\| V_\kappa^{d,g^{d,m,\lambda}} X_\kappa\|_{L_p(\R^d)})^\theta)^{1 /\theta}\le c_2,\\
\phi (\sum_{\kappa \in \Z_{+ r}^{d,J}} 2^{-(\kappa, \alpha)}
V_\kappa^{d,g^{d,m,\lambda}} X_\kappa) =0}}
\biggl(\sum_{\kappa \in \Z_{+ r}^{d,J}}
\| V_\kappa^{d,\D^\lambda g^{d,m,\lambda}}
X_\kappa \|_{L_q(I^d)}^{\Q}\biggr)^{1 /\Q} = \\
c_3 2^{-\mu r} V_2.
\end{multline*}

Применяя (3.3.9) и делая замену переменных $ X_\kappa = (c_2 /c_4)
2^{p^{-1} (\kappa,\e)} U_\kappa, \kappa \in \Z_{+ r}^{d,J}, $
получаем
\begin{multline*} \tag{3.4.5}
V_2 = \inf_{ \phi \in \Phi_n(C(I^d))} \sup_{\substack{ \{X_\kappa \in \R^{2^{(\kappa, \e)}},
\kappa \in \Z_{+ r}^{d,J}\}:\\
(\sum_{\kappa \in \Z_{+ r}^{d,J}}
(c_4 2^{-p^{-1} (\kappa, \e)} \| X_\kappa\|_{l_p^{2^{(\kappa,\e)}}})^\theta)^{1 /\theta} \le c_2,\\
\phi (\sum_{\kappa \in \Z_{+ r}^{d,J}}
2^{-(\kappa, \alpha)} V_\kappa^{d,g^{d,m,\lambda}} X_\kappa) =0}}
\biggl(\sum_{\kappa \in \Z_{+ r}^{d,J}}
(c_5 2^{-q^{-1} (\kappa,\e)} \| X_\kappa \|_{l_q^{2^{(\kappa,\e)}}})^{\Q}\biggr)^{1 /\Q} = \\
\inf_{ \phi \in \Phi_n(C(I^d))} \sup_{\substack{ \{U_\kappa \in \R^{2^{(\kappa, \e)}},
\kappa \in \Z_{+ r}^{d,J}\}:\\
(\sum_{\kappa \in \Z_{+ r}^{d,J}}
(\| U_\kappa\|_{l_p^{2^{(\kappa,\e)}}})^\theta)^{1 /\theta} \le 1,\\
\phi (\sum_{\kappa \in \Z_{+ r}^{d,J}} c_6
2^{-(\kappa, \alpha) +p^{-1} (\kappa,\e)}
V_\kappa^{d,g^{d,m,\lambda}} U_\kappa) =0}}
\biggl(\sum_{\kappa \in \Z_{+ r}^{d,J}}
(c_7 2^{-q^{-1} (\kappa,\e)} 2^{p^{-1} (\kappa,\e)}
\| U_\kappa \|_{l_q^{2^{(\kappa,\e)}}})^{\Q}\biggr)^{1 /\Q} = \\
c_7 2^{(p^{-1} -q^{-1}) r} \inf_{ \phi \in \Phi_n(C(I^d))}
\sup_{\substack{ \{U_\kappa \in \R^{2^{(\kappa, \e)}},
\kappa \in \Z_{+ r}^{d,J}\}:\\
(\sum_{\kappa \in \Z_{+ r}^{d,J}}
(\| U_\kappa\|_{l_p^{2^{(\kappa,\e)}}})^\theta)^{1 /\theta} \le 1,\\
\phi (\sum_{\kappa \in \Z_{+ r}^{d,J}} c_6
2^{-(\kappa, \alpha) +p^{-1} (\kappa,\e)}
V_\kappa^{d,g^{d,m,\lambda}} U_\kappa) =0}}
\biggl(\sum_{\kappa \in \Z_{+ r}^{d,J}}
(\| U_\kappa \|_{l_q^{2^{(\kappa,\e)}}})^{\Q}\biggr)^{1 /\Q}.
\end{multline*}
Соединяя (3.4.3), (3.4.4), (3.4.5) и обозначая для $ \phi \in \Phi_n(C(I^d)) $
через $ F_\phi $ линейное отображение пространства
$ \prod_{\kappa \in \Z_{+ r}^{d,J}} \R^{2^{(\kappa, \e)}} $
в $ \R^n, $ задаваемое равенством
$$
F_\phi (U) = \phi (\sum_{\kappa \in \Z_{+ r}^{d,J}}
c_6 2^{-(\kappa, \alpha) +p^{-1} (\kappa,\e)}
V_\kappa^{d,g^{d,m,\lambda}} U_\kappa), U = \{U_\kappa \in \R^{2^{(\kappa,\e)}},
\kappa \in \Z_{+ r}^{d,J}\},
$$
приходим к неравенству
\begin{multline*} \tag{3.4.6}
\inf_{ \phi \in \Phi_n(C(I^d))}
\sup_{f \in B((S_{p,\theta}^\alpha B)^0(\R^d)) \cap C_0^0(I^d): \phi f =0}
\| \D^\lambda f \|_{L_q(I^d)} \ge \\
c_8 2^{-\mu r +(p^{-1} -q^{-1}) r}
\inf_{ \phi \in \Phi_n(C(I^d))}
\sup_{\substack{ \{U_\kappa \in \R^{2^r}, \kappa \in \Z_{+ r}^{d,J}\}:\\
(\sum_{\kappa \in \Z_{+ r}^{d,J}}
(\| U_\kappa\|_{l_p^{2^r}})^\theta)^{1 /\theta} \le 1,\\
\phi (\sum_{\kappa \in \Z_{+ r}^{d,J}} c_6
2^{-(\kappa, \alpha) +p^{-1} (\kappa,\e)}
V_\kappa^{d,g^{d,m,\lambda}} U_\kappa) =0}}
\biggl(\sum_{\kappa \in \Z_{+ r}^{d,J}}
(\| U_\kappa \|_{l_q^{2^r}})^{\Q}\biggr)^{1 /\Q} = \\
c_8 2^{-\mu r +(p^{-1} -q^{-1}) r} \inf_{ \phi \in \Phi_n(C(I^d))}
\sup_{ U \in B(l_{p,\theta}^{2^r, c(\gamma,r)}): F_\phi U = 0}
\|U\|_{l_{q,\Q}^{2^r,c(\gamma,r)}} =\\
c_8 2^{-(\mu -(p^{-1} -q^{-1})) r}
\inf_{ \phi \in \Phi_n(C(I^d))}
\sup_{ U \in (B(l_{p,\theta}^{2^r, c(\gamma,r)}) \cap \Ker F_\phi)}
\|U\|_{l_{q,\Q}^{2^r,c(\gamma,r)}},
\end{multline*}
ибо
\begin{equation*} \tag{3.4.7}
\card \{\kappa \in \Z_{+ r}^{d,J}\} =
\card \{\kappa^J \in (\Z_+^d)^J: (\kappa^J,\e^J) = r\} =
c(\card J,r) = c(\gamma,r).
\end{equation*}
Из (3.4.6), учитывая, что для $ \phi \in \Phi_n(C(I^d)) $ ядро $ \Ker F_\phi $
является пересечением ядер $ n $ линейных функционалов на
$ \R^{2^r c(\gamma,r)}, $ и, следовательно, $ \codim \Ker F_\phi \le n, $
получаем, что
\begin{multline*} \tag{3.4.8}
\inf_{ \phi \in \Phi_n(C(I^d))} \sup_{f \in B((S_{p,\theta}^\alpha B)^0(\R^d)) \cap
C_0^0(I^d): \phi f =0}
\| \D^\lambda f \|_{L_q(I^d)} \ge\\
c_8 2^{-(\mu -(p^{-1} -q^{-1})) r}
\inf_{M \in \mathcal M^n(l_{q,\Q}^{2^r,c(\gamma,r)})}
\sup_{ U \in B(l_{p,\theta}^{2^r, c(\gamma,r)}) \cap M}
\|U\|_{l_{q,\Q}^{2^r,c(\gamma,r)}} =\\
c_8 2^{-(\mu -(p^{-1} -q^{-1})) r}
d^n(B(l_{p,\theta}^{2^r, c(\gamma,r)}), l_{q,\Q}^{2^r,c(\gamma,r)}).
\end{multline*}

Для оценки правой части (3.4.8), учитывая, что в силу (3.1.1) справедливы
соотношения
\begin{equation*}
\|x\|_{l_{q, \Q}^{2^r, c(\gamma,r)}} \ge
(c(\gamma,r))^{-(1 /q -1 /\Q)_+} \cdot
\|x\|_{l_{q, q}^{2^r, c(\gamma,r)}} =
(c(\gamma,r))^{-(1 /q -1 /\Q)_+} \|x\|_{l_{q}^{2^r c(\gamma,r)}},
\end{equation*}
\begin{equation*}
\|x\|_{l_{p, \theta}^{2^r, c(\gamma,r)}} \le
(c(\gamma,r))^{(1 /\theta -1 /p)_+} \cdot
\|x\|_{l_{p, p}^{2^r, c(\gamma,r)}} =
(c(\gamma,r))^{(1 /\theta -1 /p)_+} \cdot
\|x\|_{l_{p}^{2^r c(\gamma,r)}},
\end{equation*}
или
\begin{equation*}
(c(\gamma,r))^{-(1 /\theta -1 /p)_+} \cdot
B(l_{p}^{2^r c(\gamma,r)}) \subset
B(l_{p, \theta}^{2^r, c(\gamma,r)}),
\end{equation*}
благодаря (3.1.3) выводим
\begin{multline*} \tag{3.4.9}
d^n(B(l_{p,\theta}^{2^r, c(\gamma,r)}), l_{q,\Q}^{2^r,c(\gamma,r)}) \ge\\
(c(\gamma,r))^{-(1 /q -1 /\Q)_+} d^n((c(\gamma,r))^{-(1 /\theta -1 /p)_+} \cdot
B(l_{p}^{2^r c(\gamma,r)}), l_{q}^{2^r c(\gamma,r)}) = \\
(c(\gamma,r))^{-(1 /q -1 /\Q)_+} (c(\gamma,r))^{-(1 /\theta -1 /p)_+} \cdot
d^n(B(l_{p}^{2^r c(\gamma,r)}), l_{q}^{2^r c(\gamma,r)}).
\end{multline*}
Учитывая (3.4.2) и используя (3.1.6), приходим к оценке
\begin{multline*} \tag{3.4.10}
d^n(B(l_{p}^{2^r c(\gamma,r)}), l_{q}^{2^r c(\gamma,r)}) \ge
d^m(B(l_p^{2m}), l_q^{2m}) \mid_{m = \frac|{1}{2} c(\gamma,r) 2^r} \ge\\
c_9
\begin{cases}
m^{1 /q -1 /p}, & \text{ при $ q \le p$ или $ p < q \le 2;$ } \\
m^{1/2 -1 /p}, & \text{ при $ p \le 2 < q; $ } \\
1, & \text{ при $ 2 < p < q $,} \\
\end{cases}
\ge\\
c_{10}
\begin{cases}
(2^r c(\gamma,r))^{1 /q -1 /p}, & \text{ при $ q \le p$ или $ p < q \le 2;$ } \\
(2^r c(\gamma,r))^{1/2 -1 /p}, & \text{ при $ p \le 2 < q; $ } \\
1, & \text{ при $ 2 < p < q $ }.
\end{cases}
\end{multline*}
Соединяя (3.4.9) и (3.4.10) и принимая во внимание (3.1.2), выводим
\begin{multline*} \tag{3.4.11}
d^n(B(l_{p,\theta}^{2^r, c(\gamma,r)}), l_{q,\Q}^{2^r,c(\gamma,r)}) \ge\\
c_{11}
\begin{cases}
2^{r (1/q -1/p)} r^{(\gamma -1)(1 /q -1 /p -(1 /\theta -1 /p)_+ -(1 /q -1 /\Q)_+)}, &
\text{ при $ q \le p$ или $ p < q \le 2;$ } \\
2^{r (1/2 -1/p)} r^{(\gamma -1)(1/2 -1 /p -(1 /\theta -1 /p)_+ -(1 /q -1 /\Q)_+)}, &
\text{ при $ p \le 2 < q; $ } \\
r^{-(\gamma -1)((1 /\theta -1 /p)_+ +(1 /q -1 /\Q)_+)}, & \text{ при $ 2 < p < q $}.
\end{cases}
\end{multline*}
Подставляя (3.4.11) в (3.4.8) и пользуясь (3.4.2), (3.1.2), получаем
соотношение
\begin{multline*}
\inf_{ \phi \in \Phi_n(C(I^d))} \sup_{f \in B((S_{p,\theta}^\alpha B)^0(\R^d)) \cap
C_0^0(I^d): \phi f =0} \| \D^\lambda f \|_{L_q(I^d)} \ge \\
c_{12}
\begin{cases}
2^{-r \mu} r^{(\gamma -1)(1 /q -1 /p -(1 /\theta -1 /p)_+ -(1 /q -1 /\Q)_+)},
\text{ при $ q \le p$ или $ p < q \le 2;$ } \\
2^{-r(\mu -(1/2 -1/q))} r^{(\gamma -1)(1/2 -1 /p -(1 /\theta -1 /p)_+ -(1 /q -1 /\Q)_+)},
\text{ при $ p \le 2 < q; $ } \\
2^{-r(\mu -(1/p -1/q))} r^{-(\gamma -1)((1 /\theta -1 /p)_+ +(1 /q -1 /\Q)_+)},
\text{ при $ 2 < p < q $,}
\end{cases}=\\
c_{12}
\begin{cases}
(2^r r^{\gamma -1})^{-\mu} r^{(\gamma -1)(\mu +1 /q -1 /p -(1 /\theta -
1 /p)_+ -(1 /q -1 /\Q)_+)},
\text{ при $ q \le p$ или $ p < q \le 2;$ } \\
(2^r r^{\gamma -1})^{-(\mu -(1/2 -1/q))} r^{(\gamma -1)(\mu -(1/p -1/q) -(1 /\theta -1 /p)_+)},
\text{ при $ p \le 2 < q; $ } \\
(2^r r^{\gamma -1})^{-(\mu -(1/p -1/q))}
r^{(\gamma -1)(\mu -(1/p -1/q) -(1 /\theta -1 /p)_+)},
\text{ при $ 2 < p < q, $}
\end{cases}\ge\\
c_1
\begin{cases}
n^{-\mu} (\log n)^{(\gamma -1)(\mu -(1 /p -1 /q) -(1 /\theta -1 /p)_+ -(1 /q -1 /\Q)_+)},
\text{ при $ q \le p$ или $ p < q \le 2;$ } \\
n^{-(\mu -(1/2 -1/q))} (\log n)^{(\gamma -1)(\mu -(1/p -1/q) -(1 /\theta -1 /p)_+)},
\text{ при $ p \le 2 < q; $ } \\
n^{-(\mu -(1/p -1/q))}
(\log n)^{(\gamma -1)(\mu -(1/p -1/q) -(1 /\theta -1 /p)_+)},
\text{ при $ 2 < p < q, $}
\end{cases}
\end{multline*}
что совпадает с (3.4.1). $ \square $

Лемма 3.4.2

При соблюдении условий леммы 3.4.1 существует константа
$ c_{13}(d, \alpha, p,\theta,q,\lambda) > 0 $ такая, что при $ n \in \N $ выполняется неравенство
\begin{multline*} \tag{3.4.12}
\inf_{ \phi \in \Phi_n(C(I^d))} \sup_{f \in B((S_{p,\theta}^\alpha B)^0(\R^d))
\cap C_0^0(I^d): \phi f =0} \| \D^\lambda f \|_{L_q(I^d)} \ge\\
c_{13}
n^{-\mu +(1 /p -1 /q)_+} (\log n)^{(\gamma -1)(1 /\Q -1 /\theta)_+},
\end{multline*}
где $ \mu = \mn(\alpha -\lambda), \gamma = \cmn(\alpha -\lambda), \Q =
\max(2,q).$

Доказательство.

При доказательстве леммы 3.4.2 будем использовать те же обозначения, что при
доказательстве леммы 3.4.1. Сначала так же, как при доказательстве леммы 3.4.1,
фиксируя $ m \in \N^d $ так, чтобы $ l(\alpha) \in \Z_+^d(m), $
в соответствии с леммой 3.3.4 построим функцию $ g^{d,m, \lambda} \in C_0^{m +\lambda}(I^d). $
Затем для $ n \in \N $ выберем $ r \in \N, $ для которого соблюдается условие
\begin{equation*} \tag{3.4.13}
2^{r -1} < 2n \le 2^r.
\end{equation*}
Тогда, повторяя соответствующую часть доказательства леммы 3.4.1, для
$ n \in \N, r \in \N, $ удовлетворяющих условию (3.4.13), получим соотношения
(3.4.3), (3.4.4), (3.4.5), соединение которых приводит к оценке
\begin{multline*} \tag{3.4.14}
\inf_{ \phi \in \Phi_n(C(I^d))}
\sup_{f \in B((S_{p,\theta}^\alpha B)^0(\R^d)) \cap C_0^0(I^d): \phi f =0}
\| \D^\lambda f \|_{L_q(I^d)} \ge\\
c_8 2^{-\mu r +(p^{-1} -q^{-1}) r}
\inf_{ \phi \in \Phi_n(C(I^d))}
\sup_{\substack{ \{U_\kappa \in \R^{2^r}, \kappa \in \Z_{+ r}^{d,J}\}:\\
(\sum_{\kappa \in \Z_{+ r}^{d,J}}
(\| U_\kappa\|_{l_p^{2^r}})^\theta)^{1 /\theta} \le 1,\\
\phi (\sum_{\kappa \in \Z_{+ r}^{d,J}} c_6 2^{-(\kappa, \alpha) +p^{-1} (\kappa,\e)}
V_\kappa^{d,g^{d,m,\lambda}} U_\kappa) =0}}
\biggl(\sum_{\kappa \in \Z_{+ r}^{d,J}} (\| U_\kappa \|_{l_q^{2^r}})^{\Q}\biggr)^{1 /\Q}.
\end{multline*}

Чтобы оценить правую часть (3.4.14), для каждого $ \phi \in \Phi_n(C(I^d)) $
рассмотрим систему точек $ \{x^i \in I^d, i =1,\ldots,n\} $ такую, что
$ \phi f = (f(x^1),\ldots,f(x^n)) \in \R^n $ для $ f \in C(I^d). $ Тогда имеем
\begin{multline*} \tag{3.4.15}
\{ \{U_\kappa \in \R^{2^r}, \kappa \in \Z_{+ r}^{d,J}\}:
\phi (\sum_{\kappa \in \Z_{+ r}^{d,J}} c_6 2^{-(\kappa, \alpha) +p^{-1} (\kappa,\e)}
V_\kappa^{d,g^{d,m,\lambda}} U_\kappa) =0 \} = \\
\{ \{U_\kappa = \{u_{\kappa,\nu^\kappa} \in \R, \nu^\kappa \in \Nu_{0,2^\kappa -\e}^d\},
\kappa \in \Z_{+ r}^{d,J}\}:
\sum_{\kappa \in \Z_{+ r}^{d,J}} c_6 2^{-(\kappa, \alpha) +p^{-1} (\kappa,\e)} \\
\sum_{\nu^\kappa \in \Nu_{0,2^\kappa -\e}^d}
u_{\kappa,\nu^\kappa} g^{d,m,\lambda}(2^\kappa x^i -\nu^\kappa) =0, i =1,\ldots,n\}.
\end{multline*}

Для $ \kappa \in \Z_{+ r}^{d,J} $ обозначая через $ M_\kappa = M_\kappa(\phi) $
множество тех $ \nu^\kappa \in \Nu_{0,2^\kappa -\e}^d, $ для которых
существует $ i \in \{1,\ldots,n\} $ такое, что $ x^i \in
(2^{-\kappa} \nu^\kappa +2^{-\kappa} I^d), $ и сопоставляя каждому
$ \nu^\kappa \in M_\kappa $ некоторое $ i_{\nu^\kappa} \in \{1,\ldots,n\} $
такое, что $ x^{i_{\nu^\kappa}} \in (2^{-\kappa} \nu^\kappa +2^{-\kappa} I^d), $
заметим, что отображение $ M_\kappa \ni \nu^\kappa \mapsto i_{\nu^\kappa} \in
\{1,\ldots,n\} $ инъективно, ибо если для $ \nu^\kappa, \nu^{\prime \kappa}
\in M_\kappa: \nu^\kappa \ne \nu^{\prime \kappa}, $
выполняется равенство $ i_{\nu^\kappa} =
i_{\nu^{\prime \kappa}}, $ то $ x^{i_{\nu^\kappa}} = x^{i_{\nu^{\prime \kappa}}}, $
а, следовательно, $ (2^{-\kappa} \nu^\kappa +2^{-\kappa} I^d) \cap
(2^{-\kappa} \nu^{\prime \kappa} +2^{-\kappa} I^d) \ne \emptyset, $ что не верно.
Отсюда заключаем, что $ \card M_\kappa \le n, $ и в силу (3.4.13) видим, что
$ \card(\Nu_{0,2^\kappa -\e}^d \setminus M_\kappa) = \card \Nu_{0,2^\kappa -\e}^d
-\card M_\kappa \ge 2^{(\kappa,\e)} -n = 2^r -n \ge 2n -n = n,
\kappa \in \Z_{+ r}^{d,J}. $

Для каждого $ \kappa \in \Z_{+ r}^{d,J} $ возьмём некоторое множество $ N_\kappa
= N_\kappa(\phi) \subset (\Nu_{0,2^\kappa -\e}^d \setminus M_\kappa), $ для
которого $ \card N_\kappa = n. $
Тогда при $ \kappa \in \Z_{+ r}^{d,J}, \\nu^\kappa \in N_\kappa $ для
$ i =1,\ldots,n $ справедливо соотношение
$ x^i \notin (2^{-\kappa} \nu^\kappa +2^{-\kappa} I^d), $
а поскольку носитель $ \supp g^{d,m, \lambda}(2^{\kappa} \cdot -\nu^\kappa)
\subset (2^{-\kappa} \nu^\kappa +2^{-\kappa} I^d), $
то $ g^{d,m, \lambda}(2^{\kappa} x^i -\nu^\kappa) =0 $ при $ \nu^\kappa \in
N_\kappa, \kappa \in \Z_{+ r}^{d,J}, i =1,\ldots,n. $
Принимая во внимание это обстоятельство, получаем, что
\begin{multline*} \tag{3.4.16}
\{ \{U_\kappa = \{u_{\kappa,\nu^\kappa} \in \R, \nu^\kappa \in \Nu_{0,2^\kappa -\e}^d\},
\kappa \in \Z_{+ r}^{d,J}\}:\\
\sum_{\kappa \in \Z_{+ r}^{d,J}} c_6 2^{-(\kappa, \alpha) +p^{-1} (\kappa,\e)}
\sum_{\nu^\kappa \in \Nu_{0,2^\kappa -\e}^d}
u_{\kappa,\nu^\kappa} g^{d,m,\lambda}(2^\kappa x^i -\nu^\kappa) =0, i =1,\ldots,n\} \supset \\
\{ \{U_\kappa = \{u_{\kappa,\nu^\kappa} \in \R, \nu^\kappa \in \Nu_{0,2^\kappa -\e}^d\},
\kappa \in \Z_{+ r}^{d,J}\}:
u_{\kappa,\nu^\kappa} \in \R, \nu^\kappa \in N_\kappa; \\
u_{\kappa,\nu^\kappa} =0,
\nu^\kappa \in \Nu_{0,2^\kappa -\e}^d \setminus N_\kappa,
\kappa \in \Z_{+ r}^{d,J}\}.
\end{multline*}
Из (3.4.15), (3.4.16) вытекает, что
\begin{multline*}
\{ \{ U_\kappa \in \R^{2^r}, \kappa \in \Z_{+ r}^{d,J}\}:\\
(\sum_{\kappa \in \Z_{+ r}^{d,J}}
(\| U_\kappa\|_{l_p^{2^r}})^\theta)^{1 /\theta} \le 1,
\phi (\sum_{\kappa \in \Z_{+ r}^{d,J}} c_6 2^{-(\kappa, \alpha) +p^{-1} (\kappa,\e)}
V_\kappa^{d,g^{d,m,\lambda}} U_\kappa) =0 \} \supset \\
\{ \{U_\kappa = \{u_{\kappa,\nu^\kappa} \in \R, \nu^\kappa \in \Nu_{0,2^\kappa -\e}^d\},
\kappa \in \Z_{+ r}^{d,J}\}:
u_{\kappa,\nu^\kappa} \in \R, \nu^\kappa \in N_\kappa; \\
u_{\kappa,\nu^\kappa} =0,
\nu^\kappa \in \Nu_{0,2^\kappa -\e}^d \setminus N_\kappa,
\kappa \in \Z_{+ r}^{d,J};
(\sum_{\kappa \in \Z_{+ r}^{d,J}}
(\| U_\kappa\|_{l_p^{2^r}})^\theta)^{1 /\theta} \le 1\}.
\end{multline*}
Отсюда, используя (3.1.10), (3.4.7), (3.1.2), (3.4.13), получаем, что при
$ n \in \N $ для каждого $ \phi \in \Phi_n(C(I^d)) $ справедливо неравенство
\begin{multline*}
\sup_{\substack{ \{U_\kappa \in \R^{2^r}, \kappa \in \Z_{+ r}^{d,J}\}:\\
(\sum_{\kappa \in \Z_{+ r}^{d,J}}
(\| U_\kappa\|_{l_p^{2^r}})^\theta)^{1 /\theta} \le 1,\\
\phi (\sum_{\kappa \in \Z_{+ r}^{d,J}} c_6 2^{-(\kappa, \alpha) +p^{-1} (\kappa,\e)}
V_\kappa^{d,g^{d,m,\lambda}} U_\kappa) =0}}
\biggl(\sum_{\kappa \in \Z_{+ r}^{d,J}}
(\| U_\kappa \|_{l_q^{2^r}})^{\Q}\biggr)^{1 /\Q} \ge \\
\sup_{\substack{ \{U_\kappa = \{u_{\kappa,\nu^\kappa} \in \R, \nu^\kappa \in \Nu_{0,2^\kappa -\e}^d\},
\kappa \in \Z_{+ r}^{d,J}\}:\\
u_{\kappa,\nu^\kappa} \in \R, \nu^\kappa \in N_\kappa; u_{\kappa,\nu^\kappa} =0,
\nu^\kappa \in \Nu_{0,2^\kappa -\e}^d \setminus N_\kappa,
\kappa \in \Z_{+ r}^{d,J};\\
(\sum_{\kappa \in \Z_{+ r}^{d,J}}
(\| U_\kappa\|_{l_p^{2^r}})^\theta)^{1 /\theta} \le 1}}
\biggl(\sum_{\kappa \in \Z_{+ r}^{d,J}}
(\| U_\kappa \|_{l_q^{2^r}})^{\Q}\biggr)^{1 /\Q} = \\
\sup_{\substack{ \{u_{\kappa,\nu^\kappa} \in \R, \kappa \in \Z_{+ r}^{d,J}, \nu^\kappa \in N_\kappa\}:\\
(\sum_{\kappa \in \Z_{+ r}^{d,J}}
((\sum_{\nu^\kappa \in \N_\kappa} | u_{\kappa,\nu^\kappa}|^p)^{1 /p})^\theta)^{1 /\theta} \le 1}}
\biggl(\sum_{\kappa \in \Z_{+ r}^{d,J}}
((\sum_{\nu^\kappa \in \N_\kappa} | u_{\kappa,\nu^\kappa}|^q)^{1 /q})^{\Q}\biggr)^{1 /\Q} = \\
n^{(1 /q -1 /p)_+} (\card \Z_{+ r}^{d,J})^{(1 /\Q -1 /\theta)_+} = \\
n^{(1 /q -1 /p)_+} (c(\gamma,r))^{(1 /\Q -1 /\theta)_+} \ge \\
c_{14} n^{(1 /q -1 /p)_+} r^{(\gamma -1)(1 /\Q -1 /\theta)_+} \ge \\
c_{15} n^{(1 /q -1 /p)_+} (\log n)^{(\gamma -1)(1 /\Q -1 /\theta)_+},
\end{multline*}
а, значит,
\begin{multline*}
\inf_{ \phi \in \Phi_n(C(I^d))}
\sup_{\substack{ \{U_\kappa \in \R^{2^r}, \kappa \in \Z_{+ r}^{d,J}\}:\\
(\sum_{\kappa \in \Z_{+ r}^{d,J}}
(\| U_\kappa\|_{l_p^{2^r}})^\theta)^{1 /\theta} \le 1,\\
\phi (\sum_{\kappa \in \Z_{+ r}^{d,J}} c_6 2^{-(\kappa, \alpha) +p^{-1} (\kappa,\e)}
V_\kappa^{d,g^{d,m,\lambda}} U_\kappa) =0}}
\biggl(\sum_{\kappa \in \Z_{+ r}^{d,J}}
(\| U_\kappa \|_{l_q^{2^r}})^{\Q}\biggr)^{1 /\Q} \ge \\
c_{15} n^{(1 /q -1 /p)_+} (\log n)^{(\gamma -1)(1 /\Q -1 /\theta)_+}.
\end{multline*}
Подставляя эту оценку в (3.4.14) и используя (3.4.13), приходим к
неравенству
\begin{multline*}
\inf_{ \phi \in \Phi_n(C(I^d))}
\sup_{f \in B((S_{p,\theta}^\alpha B)^0(\R^d)) \cap C_0^0(I^d): \phi f =0}
\| \D^\lambda f \|_{L_q(I^d)} \ge\\
c_8 2^{-\mu r +(p^{-1} -q^{-1}) r}
c_{15} n^{(1 /q -1 /p)_+} (\log n)^{(\gamma -1)(1 /\Q -1 /\theta)_+} \ge \\
c_{13} n^{-\mu +(1 /p -1 /q)_+} (\log n)^{(\gamma -1)(1 /\Q -1 /\theta)_+},
\end{multline*}
что совпадает с (3.4.12). $ \square $

Теорема 3.4.3

При $ d \in \N $ пусть $ D \subset \R^d $ -- ограниченная область, а для
$ \alpha \in \R_+^d, 1 \le p < \infty,  1 < q < \infty, \lambda \in \Z_+^d $
соблюдаются соотношения (2.1.6), (2.1.9) и $ 1 \le \theta \le \infty. $
Пусть ещё $ \mathcal V = \D^\lambda, D(\mathcal V) = \{f \in C(D): \D^\lambda f \in L_q(D) \},
X = L_q(D), K = B((S_{p,\theta}^\alpha B)^\prime(D)). $ Тогда существует
константа $ c_{16}(\mathcal V,K,X) >0 $ такая, что при $ n \in \N $ справедливо
неравенство
\begin{multline*} \tag{3.4.17}
\sigma_n(\mathcal V,K,X) \ge\\
c_{16}
\begin{cases}
\begin{aligned}
&n^{-(\mu -(1/p -1/q)_+)}\times\\
&(\log n)^{(\gamma -1)(\mu -(1/p -1/q) -(1 /\theta -1 /p)_+) -(1 /q -1 /\Q)_+},
\end{aligned}&
\text{ при $ q \le p $ или $ 2 \le p < q, $ } \\
n^{-(\mu -(1 /p -1 /q)_+)} (\log n)^{(\gamma -1)(1 /\Q -1 /\theta)_+}, &
\text{ при $ p < \min(2,q). $ }
\end{cases}
\end{multline*}

Доказательство.

В условиях теоремы для $ N \in \N $ полагая
\begin{equation*} \tag{3.4.18}
V_0^N = \inf_{ \bm \phi \in \Phi_N(C(I^d))} \sup_{f \in
B((S_{p,\theta}^\alpha B)^0(\R^d)) \cap C_0^0(I^d): \bm \phi f =0}
\| \D^\lambda f \|_{L_q(I^d)},
\end{equation*}
фиксируем точку $ x^0 \in \R^d $ и вектор $ \delta \in \R_+^d $
такие, что $ Q = (x^0 +\delta I^d) \subset D. $ Пусть $ n \in \N,
A \in \mathcal A^n(L_q(D)) $ и $ \phi \in \Phi_n(C(D)) $ определяется системой
точек $ y^1,\ldots,y^n \in D. $ Тогда, ввиду того, что для $ \bm \phi \in
\Phi_n(C(I^d)) $ верно неравенство

$$
\sup_{\bm f \in B((S_{p,\theta}^\alpha B)^0(\R^d)) \cap C_0^0(I^d): \bm \phi \bm f =0}
\| \D^\lambda \bm f \|_{L_q(I^d)} \ge V_0^n > 0,
$$
возьмём функцию $ f \in B((S_{p,\theta}^\alpha B)^0(\R^d)) \cap
C_0^0(I^d), $ обладающую следующими свойствами:
\begin{equation*} \tag{3.4.19}
\| \D^\lambda f \|_{L_q(I^d)} > \frac{1}{2} V_0^n,
\end{equation*}
и
\begin{equation*} \tag{3.4.20}
f(x^i) =0, i =1,\ldots,n,
\end{equation*}
где
$ x^i = \delta^{-1}(y^i -x^0), i =1,\ldots,n. $ Рассмотрим функцию
$ F = h_{\delta,x^0}^{-1} f $ (см. п. 3.2.). Тогда, поскольку функция
$ f \in B((S_{p,\theta}^\alpha B)^0(\R^d) \cap C_0^0(I^d), $ то в силу
(1.1.9) и (3.2.5) функция
$ F \in (c_{17} B((S_{p,\theta}^\alpha B)^\prime(\R^d)) \cap C_0^0(Q). $
Полагая $ \mathcal F = (1 / c_{17}) F \mid_D, $ получаем, что
$ (\supp \mathcal F) \subset Q \subset D, \mathcal F \in B((S_{p,\theta}^\alpha B)^\prime(D)), $
благодаря (3.4.20) -- $ \mathcal F(y^i) = (1 / c_{17}) f(\delta^{-1}(y^i -x^0))
= (1 / c_{17}) f(x^i) =0, i =1,\ldots,n. $
При этом ввиду (3.4.19) соблюдается неравенство
\begin{multline*}
\frac{1}{2} V_0^n < \| \D^\lambda f\|_{L_q(I^d)}
= \| \D^\lambda (h_{\delta,x^0} F)\|_{L_q(I^d)} =\\
c_{18}(\delta,\lambda,q) \| \D^\lambda F\|_{L_q(Q)} =
c_{19} \| \D^\lambda \mathcal F\|_{L_q(Q)} = c_{19} \| \D^\lambda \mathcal F\|_{L_q(D)} = \\
c_{19} \| \D^\lambda \mathcal F -A \circ \phi(\mathcal F) +A \circ \phi(0) -\D^\lambda 0\|_{L_q(D)}\le\\
2 c_{19} \sup_{g \in K} \| \D^\lambda g -A \circ \phi(g)\|_X.
\end{multline*}
Отсюда в силу произвольности $ A \in \mathcal A^n(L_q(D)) $ и $ \phi \in
\Phi_n(C(D)), $ заключаем, что имеет место оценка
\begin{equation*}
\sigma_n(\mathcal V,K,X) \ge c_{20} V_0^n.
\end{equation*}
Объединяя эту оценку с (3.4.18), (3.4.1), (3.4.12), приходим к (3.4.17). $ \square $

\end{document}